\documentclass{article}
\usepackage{amsmath}
\usepackage{amsthm}
\usepackage{amssymb}
\usepackage{enumerate, fullpage}
\usepackage{color, xcolor}
\usepackage{graphicx}
\usepackage{verbatim}
\usepackage{epstopdf}
\usepackage{epsfig} 
\usepackage{tablefootnote}
\usepackage{tabularx,booktabs}
\usepackage{multirow}
\usepackage{hyperref}
\usepackage{algorithm2e}
\usepackage{algpseudocode}
\usepackage{soul}
\RestyleAlgo{ruled}
\usepackage{graphbox}
\newtheorem{theorem}{Theorem}

\newcommand{\commentout}[1]{}
\newcommand{\bW}{\boldsymbol{W}}
\newcommand{\bD}{\boldsymbol{D}}
\newcommand{\bc}{\boldsymbol{c}}
\newcommand{\bb}{\boldsymbol{b}}
\newcommand{\bw}{\boldsymbol{w}}
\DeclareMathOperator*{\argmin}{arg\,min}

\usepackage{multirow, tabularx}
\newcolumntype{C}{>{\centering\arraybackslash}X}
\newcolumntype{b}{X}
\newcolumntype{s}{>{\hsize=.4\hsize}X}
\newcolumntype{o}{>{\hsize=.3\hsize}X}

\title{WeakIdent: Weak formulation for  Identifying Differential Equations using Narrow-fit and Trimming}
\author{Mengyi Tang, Wenjing Liao, Rachel Kuske  and Sung Ha Kang \thanks{School of Mathematics, Georgia Institute of Technology, Atlanta, GA 30332-0160. Mengyi Tang, the corresponding author  (tangmengyi@gatech.edu), Liao (wliao60@gatech.edu),  Kuske (rachel@math.gatech.edu) and Kang (kang@math.gatech.edu). This work was partially funded by Simons Foundation grant 584960 and NSF 2145167. We provide the code of WeakIdent at \url{http://github.com/sunghakang/WeakIdent}.}}

\begin{document}
\maketitle
\begin{abstract}
Data-driven identification of differential equations is an interesting but challenging problem, especially when the given data are corrupted by noise. 
When the governing differential equation is a linear combination of various differential terms, the identification problem can be formulated as solving a linear system, with the feature matrix consisting of linear and nonlinear terms multiplied by a coefficient vector. This product is equal to the time derivative term, and thus generates dynamical behaviors.
The goal is to identify the correct terms that form the equation to capture the dynamics of the given data.
We propose a general and robust framework to recover differential equations using a weak formulation, for both ordinary and partial differential equations (ODEs and PDEs). The weak formulation facilitates an efficient and robust way to handle noise. For a robust recovery against noise and the choice of hyper-parameters, we introduce two new mechanisms, narrow-fit and trimming, for the coefficient support and value recovery, respectively.   
For each sparsity level, Subspace Pursuit is utilized to find an initial set of support from the large dictionary. Then, we focus on highly dynamic regions (rows of the feature matrix), and error normalize the feature matrix in the narrow-fit step.   The support is further updated via trimming of the  terms that contribute the least. Finally, the support set of features with the smallest Cross-Validation error is chosen as the result.
A comprehensive set of numerical experiments are presented for both systems of ODEs and PDEs with various noise levels. 
The proposed method gives a robust recovery of the coefficients, and a significant denoising effect which can handle up to $100\%$ noise-to-signal ratio for some equations.  We compare the proposed method with several state-of-the-art algorithms for the recovery of differential equations.
\end{abstract}

\section{Introduction}

In recent years, there has been an increasing interest in discovering physical or biological dynamics from complex data. The discovery of  differential equations can offer  important insights into contemporary neuroscience \cite{favela2021dynamical}, fluid mechanics, physical systems\cite{bongard2007automated,schmidt2009distilling}, and biology \cite{nardini2020learning}.

In this paper, we focus on the inverse problem of identifying a  differential equation corresponding to given data corrupted by noise. 
Given a time-dependent discrete data set, we aim to discover the underlying equation of the form
\begin{equation}
\partial_t u = f(u, \partial_{\boldsymbol x} u, \ldots, \partial ^k_{\boldsymbol x} u,\ldots,u^2,\partial_{\boldsymbol x} u^2, ..., \partial ^k_{\boldsymbol x} u^2,\ldots,u^3,\partial_{\boldsymbol x} u^3, ..., \partial ^k_{\boldsymbol x} u^3,\ldots)
\label{eqdifferential1}
\end{equation}
where each differential term in the right hand side of \eqref{eqdifferential1} is called a feature in this prescribed dictionary. In particular, $f$ is called the governing equation of \eqref{eqdifferential1}. 
We assume that $f$ in \eqref{eqdifferential1} is a linear combination of the features, so that 
this inverse problem becomes identification of a sparse coefficient vector with both the support and the values unknown. 
Since feature terms include linear and nonlinear terms, this $f$ in \eqref{eqdifferential1} includes nonlinear differential equations.   This  model identification problem is very challenging when the given data are corrupted by noise.

Parameter identification in differential equations and dynamical systems has been studied by scientists in various fields. Earlier works include  \cite{baake1992fitting,bar1999fitting,bock1983recent,bongard2007automated,muller2004parameter,schmidt2009distilling}, where the differential equation \eqref{eqdifferential1} is considered in \cite{bar1999fitting,muller2004parameter}, symbolic regression is used in  \cite{bongard2007automated,schmidt2009distilling}, and an optimization approach is taken in \cite{baake1992fitting,bar1999fitting,muller2004parameter}.
In recent years, sparse regression is incorporated into the  model  identification problem to promote sparsity in coefficient recovery \cite{brunton2016discovering,Zhang2019OnTC,Kaheman2020SINDyPIAR,Rudy2019DataDrivenIO,loiseau2018sparse, Guan2021SparseNM,champion2019discovery,kang2021ident, he2020robust,tran2017exact,schaeffer2020extracting,rudy2017data, wu2019numerical,Gurevich2019RobustAO,messenger2021weakPDE,messenger2021weak}. Representative works include Sparse Identification of Nonlinear Dynamics (SINDy)\cite{brunton2016discovering,Zhang2019OnTC,Kaheman2020SINDyPIAR,Rudy2019DataDrivenIO}, Identifying Differential Equations with Numerical Time evolution (IDENT)\cite{kang2021ident, he2020robust}, Weak SINDy \cite{messenger2021weakPDE,messenger2021weak}, RGG \cite{reinbold2020using} and many others\cite{tran2017exact,schaeffer2017learning,schaeffer2020extracting,zhang2018robust}. The PDE and dynamics identification problem is also addressed by deep learning approaches \cite{Chen2022DeepNN,qin2019data,lusch2018deep,raissi2018hidden,Wu2020DataDrivenDL,xu2019dl, Xu2020DLGAPDEDO}.

The majority of existing works apply sparse regression on a linear system formed from \eqref{eqdifferential1} with differential features \cite{brunton2016discovering,Zhang2019OnTC,Kaheman2020SINDyPIAR,Rudy2019DataDrivenIO,kang2021ident, he2020robust,tran2017exact,rudy2017data}. From the given data, differential features are approximated via numerical differentiation. When the given data contain noise, a denoising step is applied before  numerical differentiation.  Least-squares moving average is applied in \cite{kang2021ident}, successively denoised differentiation is proposed in \cite{he2020robust} and regularization is used in \cite{he2022numerical}.  
In terms of sparse regression, $L_1$ or regularized $L_1$ minimization has been widely used \cite{brunton2016discovering,kang2021ident,tran2017exact,schaeffer2017learning}; Sequentially thresholded least-squares is used in \cite{Kaheman2020SINDyPIAR,loiseau2018sparse, Guan2021SparseNM,champion2019discovery}; Greedy algorithms are used in \cite{he2020robust}. More generally, the coefficients are allowed to be spatially dependent in \cite{rudy2019data,kang2021ident}, and the Group Lasso is used to promote group sparsity where each group represents a feature, which is also used for varying coefficient case in \cite{kang2021ident}.  While these methods using differential features give good  results, numerical differentiation can be unstable for high-order features, and the coefficient recovery may not be robust when the given data is corrupted by noise. 

Recent progress using a weak/integral formulation \cite{Gurevich2019RobustAO,reinbold2020using,messenger2021weakPDE,messenger2021weak} shows improvements in the robustness of the sparse coefficient identification. Using a weak form for \eqref{eqdifferential1} with a set of test functions gives a linear system with integral features instead of differential features. Noise is tackled through the weak form, since a proper test function gives a denoising effect. The test functions are chosen to be  localized smooth functions vanishing on the boundaries, thus resembling kernel functions commonly used in kernel denoising methods. It is shown in \cite{messenger2021weakPDE,messenger2021weak} that using the weak form with the standard sequentially threshold least-squares algorithm gives rise to superior numerical performance. Differential equations with high-order derivatives, including the Korteweg–De Vries (KdV) equation, the Kuramoto–Sivashinsky (KS) equation, and 2D reaction-diffusion equations can be recovered even with a significant amount of noise. 
A related work \cite{Chen2020MethodsTR} focuses on the  identification of advection-diffusion equations, and shows that a Galerkin-type algorithm using the weak form out-performs the collocation-type algorithm based on  using a  differential form.

In this paper, we propose a Weak formulation for Identification of Differential Equations  with Narrow-fit and Trimming (WeakIDENT).  
To recover \eqref{eqdifferential1} where $f$ is  a linear combination of various differential terms, we construct a linear system: the feature matrix consisting of linear and nonlinear terms, called features,  multiplied by a coefficient vector, where this product is set equal to the time derivative.  We use the term coefficient support  to refer to a collection of nonzero components in the coefficient vector, i.e., yielding the linear system composed of the collection of features that contribute to the dynamics represented by the data. 
For our sparse coefficient recovery, we perform an iterative greedy support identification scheme as in \cite{he2020robust, kang2021ident} to find the support which gives the collection of linear and nonlinear differential terms.   
For each sparsity level, we use the Subspace Pursuit (SP) algorithm \cite{dai2009subspace} to first find the initial guess of the coefficient support.    We propose new narrow-fit and trimming  steps which improve the support selection as well as coefficient value recovery.  Among different sparsity results, we choose the one with the minimum Cross-Validation (CV) error as the final result.  For  Cross-Validation, we randomly separate the given data in half, use one set to find the coefficients, then use this coefficient vector with the other set of data  to compute the error.  We provide an error analysis to show the improvement in denoising when using the weak form. 
Our contributions can be summarized as follows:
\begin{enumerate}
    \item Proposing WeakIDENT to robustly identify differential equations in \eqref{eqdifferential1} from highly corrupted noisy data. Using a weak form facilitates handling noise efficiently by moving the derivative to the test function. 
    \item Proposing two  novel mechanisms, narrow-fit and trimming, to improve coefficient support recovery and coefficient value recovery, respectively. These mechanisms utilize column-wise error normalization for improved coefficient value recovery.  Narrow-fit is focused on highly dynamic regions to reduce the size of the feature matrix, and trimming the features with small contributions to the result further contributes to this improvement. 
    \item We provide comprehensive numerical experiments for ordinary differential systems (ODEs) and partial differential equations (PDEs), and compare with existing methods such as  \cite{kang2021ident, he2020robust,Gurevich2019RobustAO,messenger2021weakPDE,messenger2021weak}.
\end{enumerate}

We organize this paper as follows. In Section \ref{sec:model}, we state the identification problem of differential equations, and give details about the feature formulation in the weak form, the choice of test functions, and provide an error analysis of the weak form.   We present our WeakIDENT Algorithm in Section \ref{sec: weakIdent Algorithm} with the details of error normalization,  selection of highly dynamic regions with dominant contributions to the identification,  and trimming of the features with the least contribution to the support. 
A comprehensive set of numerical experiments is provided in Section \ref{sec: numerical experiments}, including various comparisons against state-of-the-art algorithms. We conclude the paper in Section \ref{sec: conclusion}, and provide additional experiments in Appendix \ref{sec: appendix}.

\section{Problem set-up, Weak formulation and Error analysis}
\label{sec:model}

In this section, we state the identification problem for differential equations and formulate a linear system in a weak form. We also discuss the choice of test functions and provide an error analysis of the weak formulation.

\subsection{Problem set-up}
We  present the identification problem with one spatial variable for simplicity. It can be easily extended to multi-variables, and numerical results are provided for the multi-variable case. 
We consider a spatial-temporal domain  $\Omega = [X_1,X_2]  \times [0,T]$ with $X_1<X_2$ and $T>0$.
We assume a set of discrete time-dependent noisy data is given:
\begin{equation} 
\mathcal{D} =  \{ {\hat{U}}_{{i}}^n |i = 1,2,...,\mathbb{N}_x; n = 1,...,\mathbb{N}_t \} \in \mathbb{R}^{\mathbb{N}_x \times \mathbb{N}_t}, 
\label{e: given data} \end{equation}
where $\mathbb{N}_x$, and $\mathbb{N}_t \in \mathbb{N}$ are the size of discretization in spatial and temporal dimension  respectively. 
The data point $\hat{U}_i^n$ is an  approximation to the true solution of a differential equation
\begin{equation*}
\hat{U}_{i}^{n}  \approx  u(x_i, t^n) \;\;\text{for}\;\; (x_i,t^n) \in \Omega,
\end{equation*}
 at the spatial location $x_i = i\Delta x \in [X_1,X_2]$ and $t^n = n\Delta t \in [0,T]$. Here  $\Delta x = (X_2 - X_1)/{(\mathbb{N}_x-1)}$ and  $ \Delta t = {T}/{(\mathbb{N}_t-1)}$.
In the noisy case, we express the noisy data ${\hat{U}}_{{i}}^n$ in terms of the
clean data ${{U}}_{{i}}^n = u(x_i, t^n)$ as:
\begin{equation}
{\hat{U}}_{i}^n = {U}_{i}^n 
+ {\epsilon}_{i}^n,
\label{e: noisy discretized data}
\end{equation}
where ${\epsilon}_{i,}^n$ represents the noise at $(x_i,t^n)$.   
The objective is to identify a differential equation in the form of \eqref{eqdifferential1} from the given data (\ref{e: given data}).

We assume that the governing equation $f$ in \eqref{eqdifferential1} is a linear combination of linear and nonlinear terms including the derivatives of $u$.  This covers a vast range of ODEs and PDEs in applications, e.g., the  Lorenz equation, the Lotka-Volterra equation, transport equations, Burgers' equation, the heat equation, the KS equation, the KdV equation, and reaction-diffusion  equations.
In this paper, we consider the function $f$ to be a linear combination of different derivatives of powers of $u$:
\begin{align}
      \frac{\displaystyle \partial u}{\displaystyle \partial t}(x, t)= \sum_{l=1}^L c_l F_l
      \quad \text{with} \quad
       F_l = \frac{\partial^{
      \alpha_l}}{\partial x^{\alpha_l}}f_l, \ \text{where} \ f_l = f_l(u) = u^{\beta_l}.
      \label{e: feature formulation}
\end{align}
The $l^{\rm th}$ feature $F_l(u)$  represents the ${\alpha}_l^{\rm th}$ spatial derivative of the monomial $f_l = f_l(u) = u^{\beta_l}$ for some nonnegative integer $\beta_l$.
Let the highest order of derivative be $\bar{\alpha}$ such that $\alpha_l \in \{0,\ldots,\bar{\alpha}\}$, and the highest order of monomial be $\bar{\beta}$ such that $\beta_l \in \{0,\ldots,\bar{\beta}\}$. 
We use $L$ to denote the total number of features in the  dictionary, which depends on $\bar{\alpha}$ and $\bar{\beta}$, since it includes all combinations.  The formulation of (\ref{e: feature formulation}) has the advantage in accurate feature approximation particularly for the weak form, since integration by parts moves the derivatives to the test function.   
When the spatial domain is multi-dimensional, we consider $f_l$ as monomials in the multivariable case, and we allow $F_l$ to be partial derivatives of $f_l$ across different spatial dimensions.

In (\ref{e: feature formulation}), the coefficient can be considered as a sparse vector
\begin{equation}
\boldsymbol{c} = (c_1,...,c_L)^T  \in \mathbb{R}^{L}
\label{eqcoefficientc}
\end{equation} 
which parametrizes the differential equation.
The objective of this paper is to recover  the differential equation from the given noisy data set ${\mathcal D}$
\eqref{e: given data}, by finding a sparse coefficient vector $\boldsymbol{c}$ \eqref{eqcoefficientc} of the linear system \eqref{eqdifferential1}.

\subsection{The Weak Formulation}
\label{subsec: feature matrix with weak form} 

The weak formulation of \eqref{e: feature formulation} is
\begin{equation}
      \int_{\Omega_{h(x_i, t^n)}} \phi_{h(x_i, t^n)}(x,t) \frac{\displaystyle \partial u(x,t)}{\displaystyle \partial t}dx dt = \sum_{l=1}^L c_l  \int_{\Omega_{h(x_i, t^n)}} \phi_{h(x_i, t^n)}(x,t)F_l dxdt,
      \label{eq:weakformulation1}
\end{equation}
where the test function $\phi_h(x,t)$ is locally defined on a region $\Omega_{h(x_i, t^n)}$, which is centered at $(x_i, t^n)$ and indexed by $h$. Specifically, each test function $\phi_h(x,t)$ is a translation of a fixed function $\phi(x,t)$ such that $\phi_{h(x_i, t^n)}(x,t) = \phi(x - x_i, t - t^n)$.
Integration by parts of \eqref{eq:weakformulation1}  gives rise to
\begin{equation}
    -\int_{\Omega_{h(x_i, t^n)}} u(x,t) \frac{\displaystyle \partial \phi_h(x,t)}{\displaystyle \partial t} dxdt =
    \sum_{l=1}^L c_l
    \int_{\Omega_{h(x_i, t^n)}} 
    (-1)^{\alpha_l}  u^{\beta_l}
    \frac{\partial^{\alpha_l} \phi_h}{\partial x^{\alpha_l}} dxdt,
    \label{e: integral form system}
\end{equation}
as long as $\phi_h$ and its derivatives up to order $\bar{\alpha}$ vanish on the boundary of $\Omega_{h(x_i,t^n)}$.
The $l^{\rm th}$ term $$ \int_{\Omega_{h(x_i, t^n)}}  (-1)^{\alpha_l}  u^{\beta_l} \frac{\partial^{\alpha_l} \phi_h(x,t)}{\partial x^{\alpha_l}} dxdt$$ is the $l^{\rm th}$ integral feature with the test function $\phi_{h}$. 
Since the test function is smooth, the numerical integration can be carried out with higher order accuracy. 
With numerical integration, we obtain the following discrete linear system for WeakIdent: 
\begin{equation}
    \boldsymbol{W}\boldsymbol{c} = \boldsymbol{b}
\label{e: Wc=b}
\end{equation}
where 
\begin{equation*}
    \boldsymbol{W} = (w_{h(x_i, t^n),l}) \in \mathbb{R}^{H \times L}, \;\; \boldsymbol{c} = (c_l) \in \mathbb{R}^{L},\;\; \text{ and }\;\; \boldsymbol{b}  = (b_{h(x_i, t^n)})\in \mathbb{R}^{H},
\end{equation*} 
for 
\begin{eqnarray}
 w_{h(x_i, t^n),l} &=&  \sum_{(x_{j},t^{k}) \in \Omega_{h(x_i, t^n)}}
    (-1)^{\alpha_l} \hat{U}_{j}^{k}
    \frac{\partial^{\alpha_l} }{\partial x^{\alpha_l}}\phi_h(x_j,t^k) \Delta x \Delta t
    \quad \text{and} \quad  \nonumber\\
    b_{h_(x_i,t^n)}& =&  -\sum_{(x_j,t^k) \in \Omega_{h(x_i, t^n)}}\hat{U}_{j}^{k} \frac{ \partial \phi_h (x_j,t^k)}{ \partial t}  \Delta x \Delta t.
\label{e: approx of bW}
\end{eqnarray}
Here the numerical integration is computed with the data points $(x_{j},t^{k}) \in \Omega_{h(x_i, t^n)}$, and $w_{h(x_i, t^n),l}$ represents an  approximation of the integral of the feature $F_l$ in the integral region $\Omega_{h(x_i,t^n)}$ centered at $(x_i, t^n)$.  
The numerical integration is computed from $\mathbb{N}_x\mathbb{N}_t$ grid points. 

For the test function, we choose to use $\phi(x,t)$ as in \cite{messenger2021weakPDE}:
\begin{equation}
    \phi(x,t) =  \left(1 - \left(\frac{x}{m_{x}\Delta x}\right)^2\right)^{p_x}\left(1-\left(\frac{t}{m^t\Delta t} \right)^2\right)^{p_t}, \;\; \; 
    (x,t) \in \Omega_{h(x_i, t^n)}
    \label{e: test function}
\end{equation}
for $i = 1,..,\mathbb{N}_x, n = 1,...,\mathbb{N}_t$ where $p_x$ and $p_t$ give the smoothness of $\phi$ in terms of $x$ and  $t$. The test function satisfies  $\int_{\Omega_{h(x_i,t^n)}} \phi(x,t)dx dt  = 1$ and $ \phi(x,t) = 0 \; \text{for} \; (x,t) \in \partial \Omega_{h(x_i, t^n)} $, with $\phi(x,t)$ localized around $(x_i,t^n)$ and is supported on $\Omega_{h(x_i,t^n)}  = [x_{i} - m_x \Delta x, x_{i} + m_x \Delta x ] \times [t^{n} - m_t \Delta t, t^{n} + m_t \Delta t]$ for some positive integers $m_x$ and $m_t$.  
The weak features $ w_{h(x_i,t^n)}$ in \eqref{e: approx of bW} can be written into a convolution form $U * \frac{\partial^{\alpha_l} }{\partial x^{\alpha_l}}\phi$ and calculated through Fast Fourier Transform in terms of $\mathcal{F}^{-1}\left(\mathcal{F}(U) \circ \mathcal{F} \left( \frac{\partial^{\alpha_l} }{\partial x^{\alpha_l}}\phi \right) \right)$, where $\circ$ denotes point-wise multiplication, and  $p_x$, $p_t$, $m_x$ and $m_t$ are carefully chosen to give a denoising effect depending on the frequency of the given data as in \cite{messenger2021weakPDE}.
For the completeness, more details are presented in Appendix \ref{Asec:testfunction}.

The the weak form (\ref{e: Wc=b}) has $\mathbb{N}_x\mathbb{N}_t$ rows. For computational efficiency, we subsample $\boldsymbol{W}$ to 
\begin{equation}\label{eq:H}
H={N}_x{N}_t \leq \mathbb{N}_x \mathbb{N}_t,
\end{equation}
rows by uniformly subsampling  $N_x$ and $N_t$  points  in space and time respectively. Then, we consider highly dynamic regions to further reduce the size of $\boldsymbol{W}$ and $\boldsymbol{b}$ for an improved coefficient  recovery (details in Subsection \ref{ss:region}). In comparison, 
 random subsampling is used in \cite{reinbold2020using} for sparse regression, and  regions with large gradients in time are considered in \cite{messenger2021weak}. 

\subsection{Error Analysis of the Weak formulation} 
\label{subsec: error analysis}

We next analyze the approximation error of the weak formulation in \eqref{e: integral form system}. 
Suppose  the given noisy data $\mathcal{D}$ \eqref{e: given data} has mean-zero i.i.d Gaussian noise,  $\mathbb{E}[\epsilon_i^n] = 0$, and $
{\rm Var}(\epsilon_i^n) = \sigma^2$. Let $c_l$ be the $l^{\rm th}$ true coefficient in the true support $\text{Supp}^*$. 
The associated integral formulation using the test function (\ref{e: test function}) with the true coefficients from the true support becomes
\begin{equation}
    \int_{\Omega_h}u(x,t) \frac{\displaystyle \partial \phi_h(x,t)}{\displaystyle \partial t} dxdt +
    \sum_{l\in \text{Supp}^*} (-1)^{\alpha_l} c_l \int_{\Omega_h}
     f_l(x,t)
    \frac{\partial^{\alpha_l} \phi_h(x,t)}{\partial x^{\alpha_l}}  dxdt = 0.
    \label{e: integral form - cts - approximated system}
\end{equation}

We next analyze the error for the discretized system in \eqref{e: Wc=b} using the noisy data $\{\hat{U}_i^n\}$, approximating the true equation \eqref{e: integral form - cts - approximated system}.  The $h^{\rm th}$ row of the linear system \eqref{e: Wc=b} is obtained from the weak form with the test function $\phi_h$. The error for the discretized system in \eqref{e: Wc=b} is defined as
\begin{align}
    \boldsymbol{e}   =  {\boldsymbol{W}}\boldsymbol{c} -{\boldsymbol b}   
    \label{e: error W_hat c - h_hat}
\end{align}
where the row-wise error is 
$$e_h =\sum_{l \in {\rm Supp}^*} \sum_{(x_{j},t^{k}) \in \Omega_{h(x_i, t^n)}}
    (-1)^{\alpha_l}c_l \hat{U}_{j}^{k}
    \frac{\partial^{\alpha_l} }{\partial x^{\alpha_l}}\phi_h(x_j,t^k) \Delta x \Delta t
    +\sum_{(x_j,t^k) \in \Omega_{h(x_i, t^n)}}\hat{U}_{j}^{k} \frac{\displaystyle \partial \phi_h}{\displaystyle \partial t} (x_j,t^k) \Delta x \Delta t.$$
We decompose the error as 
\begin{align}
    \boldsymbol{e} =    \boldsymbol{e}^{\rm int} + \boldsymbol{e}^{\rm noise}
    \label{eq:errordecomposition}
\end{align}
where 
\begin{align*}
{e}^{\rm noise}_h &= e_h - \left(\sum_{l \in {\rm Supp}^*} \sum_{(x_{j},t^{k}) \in \Omega_{h(x_i, t^n)}}
    (-1)^{\alpha_l}c_l {U}_{j}^{k}
    \frac{\partial^{\alpha_l} }{\partial x^{\alpha_l}}\phi_h(x_j,t^k) \Delta x \Delta t
    +\sum_{(x_j,t^k) \in \Omega_{h(x_i, t^n)}}{U}_{j}^{k} \frac{\displaystyle \partial \phi_h}{\displaystyle \partial t} (x_j,t^k) \Delta x \Delta t\right) 
    \\
 {e}^{\rm int}_h  & = \left(\sum_{l \in {\rm Supp}^*} \sum_{(x_{j},t^{k}) \in \Omega_{h(x_i, t^n)}}
    (-1)^{\alpha_l}c_l {U}_{j}^{k}
    \frac{\partial^{\alpha_l} }{\partial x^{\alpha_l}}\phi_h(x_j,t^k) \Delta x \Delta t
    +\sum_{(x_j,t^k) \in \Omega_{h(x_i, t^n)}}{U}_{j}^{k} \frac{\displaystyle \partial \phi_h}{\displaystyle \partial t} (x_j,t^k) \Delta x \Delta t\right) \\
    &\qquad -  \left(
    \sum_{l \in {\rm Supp}^*} c_l
    \int_{\Omega_{h(x_i, t^n)}} 
    (-1)^{\alpha_l}  u^{\beta_l}
    \frac{\partial^{\alpha_l} \phi_h}{\partial x^{\alpha_l}} dxdt + \int_{\Omega_{h(x_i, t^n)}} u(x,t) \frac{\displaystyle \partial \phi_h(x,t)}{\displaystyle \partial t} dxdt.
    \right).
\end{align*}    
In this decomposition, $\boldsymbol{e}^{\rm int}$
represents the numerical integration error of the noise-free data $U$.
It has been shown in \cite{messenger2021weakPDE} that $\boldsymbol{e}^{\rm int} = \mathcal{O}((\Delta x \Delta t)^{q+1})$, where $q$ is the order of the numerical  integration  as in \cite{messenger2021weakPDE}, if the the decay of test function $\phi$ near the boundary of the test region satisfies $\max\{\phi(1-1/m_x,0) , \phi(0,1-1/m_t) \} \leq   (\frac{2 \max\{m_x, m_t\} -1}{\max\{ m_x, m_t\} ^2})^{q+1}$.

The following Theorem \ref{th:error}  provides an estimate of the error $\boldsymbol{e}^{\rm noise}$ arising from noise.

\begin{theorem} \label{th:error}
Consider a dynamical system $$u_t = \sum_{l\in \text{Supp}^*} c_l \frac{\partial^{\alpha_l}}{\partial x^{\alpha_l}}u^{\beta_l}$$ of one spatial variable where $\text{Supp}^*$ denotes the true support of the underlying differential equation. 
Assume the noise $\epsilon_i^n$ are i.i.d. and  satisfies $\mathbb{E}[\epsilon_i^n]=0$,  ${\rm Var}(\epsilon_i^n) = \sigma^2$, and 
$|\epsilon_i^n| \le \epsilon$ for all $i$ and $n$. 
Each test region $\Omega_{h(x_i,t^n)}$ for $i = 1,...,\mathbb{N}_x, n = 1, ..., \mathbb{N}_t$ has area $|\Omega_h| = m_xm_t\mathbb{N}_x\mathbb{N}_t$.
Then, 
\begin{enumerate}[(a)]
\item{In \eqref{eq:errordecomposition}, the error  from noise $\boldsymbol{e}^{\rm noise}$ for the discretized system  satisfies
\begin{equation}
    \|\boldsymbol{e}^{\rm noise} \|_\infty \le \bar{S}^*   |\Omega_h|   \epsilon + \mathcal{O}\left(\epsilon^2\right)
    \label{enoiseestimate}
\end{equation}
with a constant
\begin{equation} \bar{S}^* =\max_h \sup_{(x_j,t^k) \in \Omega_h}
\bigg|
\sum_{l\in \text{Supp}^*}(-1)^{\alpha_l}{c_l}
{\beta_l}(U_j^k)^{ \beta_l - 1}
\frac{\partial^{\alpha_l} \phi }{\partial x^{\alpha_l}}(x_j, t^k)
-  \frac{\displaystyle \partial \phi}{\displaystyle \partial t}  (x_j, t^k) 
\bigg|.
\label{eq:sbarstar}
\end{equation}
}
\item{The leading error  in $e^{\rm noise}_h$ (that is linear in noise) for the test function $\phi_h$ has mean $0$ and variance $\sigma^2 S_h^*$ where
\begin{equation}
S_h^* = \Delta x \Delta t \sum_{(x_j,t^k) \in \Omega_{h(x_i,t^n)}}
    \left(
     \sum_{l\in \text{Supp}^*}(-1)^{\alpha_l}{c_l}
     {\beta_l}(U_j^k)^{ \beta_l - 1}
    \frac{\partial^{\alpha_l} \phi_h }{\partial x^{\alpha_l}}(x_{j}, t^{k})
    + \frac{\displaystyle \partial \phi_h}{\displaystyle \partial t}  (x_{j}, t^{k}) 
    \right)^2 .
    \label{eqshstart}
    \end{equation}
    }
\end{enumerate}
\end{theorem}

Theorem \ref{th:error} is proved in Appendix \ref{app:proof:th:error}.
In summary, we prove that the error $\boldsymbol{e}$ in \eqref{e: error W_hat c - h_hat} for the discretized linear system under the weak formulation satisfies the following upper bound\begin{equation}\label{eq:errorBound}
\|\boldsymbol{e}\|_{\infty} \leq \mathcal{O} ((\Delta x \Delta t)^{q+1}) + \bar{S}^* |\Omega_h| \epsilon + \mathcal{O}(\epsilon^2) \end{equation}
where $q$ is the order of the numerical  integration  as in \cite{messenger2021weakPDE}.
By comparison, the error for the discretized system under the differential form  \cite{kang2021ident} is on the order of
\begin{equation}
 \mathcal{O}\left(\Delta t + \Delta x^{p+1-r} + \frac{\epsilon}{\Delta t} + \frac{\epsilon}{\Delta x^r}\right),
 \label{eq:identerrorBound}
\end{equation}
where $r$ is the highest order of derivatives for the features in the true support, and the  numerical differentiation is carried by interpolating the data by a $p$th order polynomial.  By comparing \eqref{eq:errorBound} and \eqref{eq:identerrorBound}, we observe that the error for the discretized linear system in the weak form is significantly smaller than the error in the differential form.

\section{WeakIdent Algorithm}
\label{sec: weakIdent Algorithm}

In this section, we present the details of the proposed Weak formulation for Identifying Differential Equation
using Narrow-fit and Trimming (WeakIdent) model.
There are mainly four steps to the algorithm: After the system is set-up as in  (\ref{e: Wc=b}), 
\begin{enumerate}
\item[\textbf{[Step 1]}]{For each sparsity level $k$, we use Subspace Pursuit (SP)\cite{dai2009subspace} to find an initial choice of support  $\mathcal{A}_0^k$ from the dictionary of $L$ features.  SP finds the choice with the minimum residual from a column-wise normalized \eqref{eq:Wdagger} linear system as in \cite{he2020robust}.}
\item[\textbf{[Step 2]}]{Narrow-fit. To recover the coefficient value using the support $\mathcal{A}_j^k$, we (i) identify  highly dynamic regions of certain features of interest; (ii)  normalize the reduced feature matrix according to the leading error term, then (iii) determine a coefficient value vector $\boldsymbol{c}(k,j)$ from this reduced narrow system (We set $j=0$ on the first iteration).}
\item[\textbf{[Step 3]}]{Trimming. With the updated coefficient values $\boldsymbol{c}(k,j)$ in [Step 2], we identify a single feature with the least contribution to $f$.  
If the contribution score is less than a preset trimming parameter $\mathcal{T}$, we trim the  corresponding coefficient.  This trimming yields a new updated support $\mathcal{A}^k_j$.  We iterate [Step 2] and [Step 3], with increment $j$, until no change is made to $\mathcal{A}^k_j$ at $j=J_k$. }
\item[\textbf{[Step 4]}]{Cross Validation. With the final support $\mathcal{A}^k_{J_k}$ and coefficient value vector $\boldsymbol{c}(k,J_k)$ for each different sparsity level $k$, we select the one $\boldsymbol{c}(k^*,J_{k^*})$ with the minimum Cross-Validation error \eqref{eq:CV} as the final result.  } 
\end{enumerate}

\begin{figure}
\centering
\includegraphics[width =0.9\textwidth]{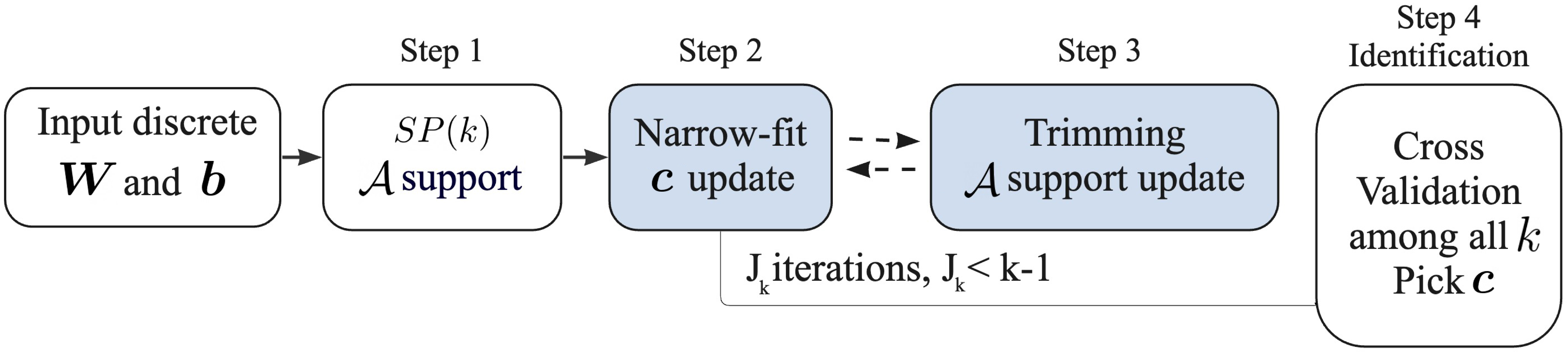}
\caption{WeakIdent flowchart: Input weak formulation $\boldsymbol{W}$ and $\boldsymbol{b}$ in \eqref{e: Wc=b} subsampled as (\ref{eq:H}). [Step 1] SP for a given sparsity $k$ gives the first candidate of coefficient support $\mathcal{A}^k_0$.  [Step 2] Narrow-fit  and [Step 3] Trimming improves the coefficient values $\boldsymbol{c}(k,j)$ and support $\mathcal{A}^k_j$. Steps 2 and 3 are iterated at most $k-1$ times.  Finally, in [Step 4] the result $\boldsymbol{c}(k^*,J_{k^*})$ with the minimum Cross Validation among all different sparsity level $k$ give the identification of the differential equation.   }\label{fig:flowchart}
\end{figure}

A schematic of the algorithm is given in Figure~\ref{fig:flowchart}. From the weak form input $\boldsymbol{W}$ and $\boldsymbol{b}$, for a fixed sparsity level $k$, SP is used to find the initial set of support $\mathcal{A}_0^k$.  Then [Step 2] Narrow-fit and [Step 3] Trimming are iterated  until the support does not change, where the number of iterations is at most  $k-1$.  Here we use $\boldsymbol{c}(k,j)$ to indicate the coefficient vector for the sparsity level $k$ and $j$ iteration.  
The cross validation is used to select the optimal solution $\boldsymbol{c}(K,J_K) $ among all $k\leq L$.
   
We present the details in the following subsections.
In [Step 2], we normalize each column of the feature matrix according to its leading error term, to balance the effect of noise perturbations across the features. The details for this error normalization of the feature matrix are given in Subsection \ref{ss: error normalized matrix}. We detail the implementation of
Narrow-fit using the highly dynamic regions  in Subsection \ref{ss:region}.
In [Step 3], we trim the support removing features with contributions below a threshold, as described in detail in Subsection \ref{ss:supression}. The algorithm is summarized Subsection \ref{ssec:algo}.

\subsection{Column-wise error normalized matrix}
\label{ss: error normalized matrix}

We use least squares for coefficient recovery.  The accuracy of least squares is highly dependent on the conditioning of the feature matrix  \cite{bjorck1990least,bjorck1991error}. 
In this paper, we utilize two types of normalization for the columns of the feature matrix to improve the coefficient recovery. 
For the linear system \eqref{e: Wc=b}, we introduce a diagonal matrix $\boldsymbol{D}= \text{diag}(d_1,...,d_L)$ and solve 
\begin{equation}
    \bW {\bD}^{-1} \bar{\bc} =\bb \quad \text{ and then } \quad \bc = {\bD}^{-1}\bar{\bc}
\end{equation}
instead.

\begin{figure}[t]
\centering
\begin{tabular}{ccc}
(a) & (b) & (c) \\
\includegraphics[width = 0.3\textwidth ]{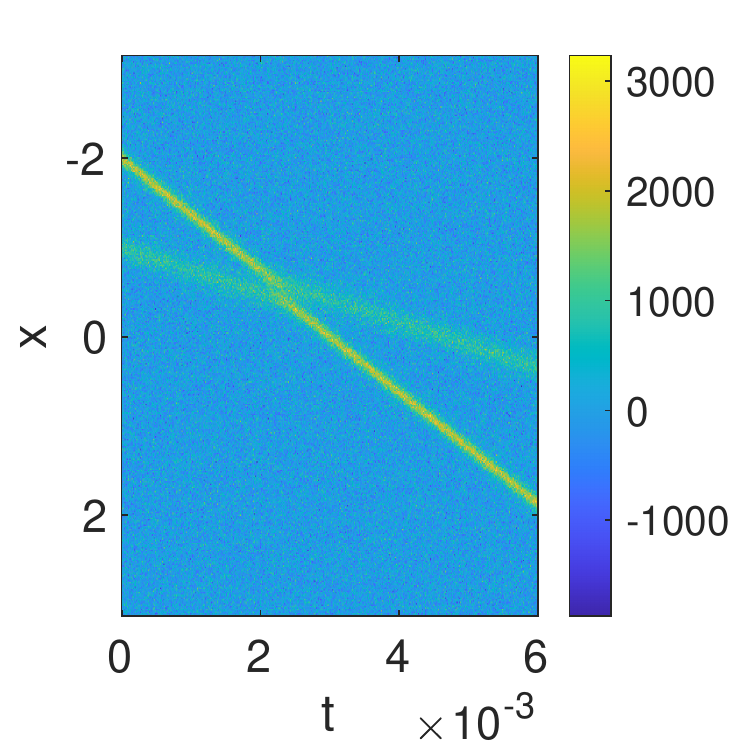} & 
\includegraphics[width = 0.3\textwidth ]{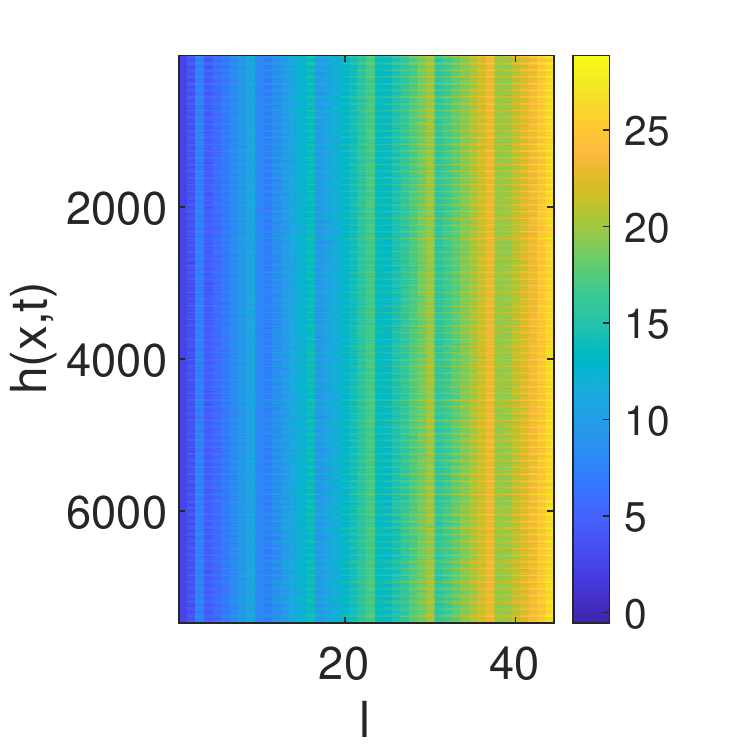}&
\includegraphics[width = 0.3\textwidth ]{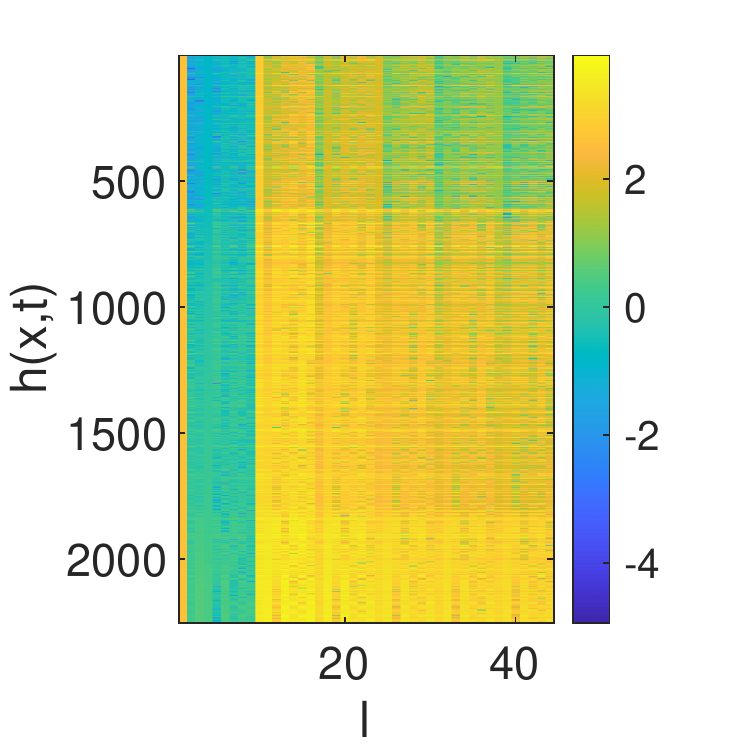}
\end{tabular}
\caption{Error normalization: (a) The given noisy data $\hat{U}$ with $\sigma_{\rm NSR}=0.5$ in $x-t$ plane. (b) The entry-wise magnitude of the matrix $\boldsymbol{W}$. 
(c) The matrix $\tilde{\boldsymbol{W}}_{\rm narrow}$ in (\ref{eq:Wtilde}).  We use $\log10$ scale in (b) and (c). 
The difference in scale  has been reduced approximately from  $10^{29}$ in the unnormalized matrix (b) to $10^6$ after normalization in (c). Our error normalization results in more uniform entry values with less variance across different columns. 
}
\label{fig:error_Normal}
\end{figure}

The first type of normalization we consider  is \textbf{column normalization}, which is applied to the feature matrix as an input to SP in [Step 1]. Denote  $\bW = [\bw_1 \ \bw_2 \ \ldots \ \bw_L]$. We let  $\boldsymbol{D}= \text{diag}(\|\boldsymbol{w}_1\|,...,\|\boldsymbol{w}_L\|)$ and  each column of $\boldsymbol{W}$ is normalized by its own norm:
\begin{equation} \label{eq:Wdagger}
  \boldsymbol{W}^\dagger = \displaystyle{ \left[\frac{\boldsymbol{w}_1}{\|\boldsymbol{w}_1\|}, \frac{\boldsymbol{w}_2}{\|\boldsymbol{w}_2\|}, \dots, \frac{\boldsymbol{w}_L}{\|\boldsymbol{w}_L\|}\right]}.
\end{equation}
We observe that the scale of the columns in the feature matrix usually varies  substantially  from column to column, which negatively affects the SP step.
This column normalization helps to prevent a large difference in the scale among the columns.  For example, in Figure \ref{fig:error_Normal} (b) shows that the magnitude of the entries in $\bW$ vary from $0$ to $10^{29}$. 

In [Step 2], we introduce our second normalization --  \textbf{error normalization}, which is particularly effective for coefficient recovery. 
The columns in $\boldsymbol{W}$ are given by certain  derivatives of a monomial of $u$. When we compute the feature matrix with noisy data, the noise has different effects on different features. For the feature $ \frac{\partial^{\alpha}}{\partial x ^{\alpha}} \left( u^{\beta}\right)$,
the noisy data with noise $\epsilon$ in   \eqref{e: noisy discretized data} give rise to the following integral feature:
\begin{equation*}
\int_{\Omega_h} (-1)^{\alpha}(u + \epsilon)^{\beta} \frac{\partial^{\alpha}}{\partial x ^{\alpha}} \left( \phi_h(x,t)\right) dxdt = \sum_{k= 0}^{\beta}(-1)^{\alpha} {\beta \choose k}   \epsilon^{\beta - k} \int_{\Omega_h} u^k \frac{\partial^{\alpha}}{\partial x ^{\alpha}}   \phi_h(x,t) dxdt.
\end{equation*}
The leading coefficient in the error (that is linear in $\epsilon$) in this integral feature is  obtained for $k=\beta -1$:

\begin{equation}
s(h,l) = \beta \left| \int_{\Omega_h} u^{\beta -1} \frac{\partial^{\alpha}}{\partial x ^{\alpha}}   \left(\phi_h(x,t)\right) dxdt \right|, \;\; h = 1,2,..., H, \; \beta \ge 1.
\label{e: leading coefficients}
\end{equation}
When $\alpha = \beta = 0$, we set $s(h,l)=1$. This leading coefficient $s(h,l)$ depends on the row index $h$ and the column index $l$. For the $l^{\rm th}$ column, we define
$$\langle s(h,l)\rangle_{h} = \frac{1}{H} \sum_{h=1}^H s(h,l)$$
as an average of these leading coefficients over the rows.

By error normalization, we normalize $\bW$ with the diagonal matrix $\bD = {\rm diag}(\langle s(h,1)\rangle_{h}, \ldots,\langle s(h,L)\rangle_{h})$ such that $\bW$ is normalized to
\begin{equation} 
  \tilde{\boldsymbol{W}} = \displaystyle{ \left[\frac{\boldsymbol{w}_1}{\langle s(h,1)\rangle_{h}}, \frac{\boldsymbol{w}_2}{\langle s(h,2)\rangle_{h}}, \dots, \frac{\boldsymbol{w}_L}{\langle s(h,L)\rangle_{h} }\right]}
  \label{eq:Wtilde}
\end{equation}
Figure \ref{fig:error_Normal} shows an example, with the given noisy data in (a) and the unnormalized feature matrix $\boldsymbol{W}$ in (b).    
Figure \ref{fig:error_Normal} (c) shows the normalized matrix $\tilde{\boldsymbol{W}}$ after the error normalization. We use $\log10$ scale in Figure \ref{fig:error_Normal}.  The difference in scale  has been reduced approximately from  $10^{29}$ in the unnormalized matrix (b) to $10^6$ after normalization in (c). Our error normalization results in more uniform entry values with less variation across different columns.

In the following Subsection, we further discuss how error normalization is used to select the highly dynamic regions .

\subsection{Highly dynamic regions: choice of the domain $\Omega_{h(x_i, t^n)}$} \label{ss:region}

One of the benefits of using the weak form is to consider the influence of different regions on the integral computation. We take advantage of this and choose a subset of test functions indexed by $\{ h| h = 2, \dots, H\}$ to improve the coefficient recovery.  We propose the following Narrow-fit procedure: (i) define the features of interest, (ii) determine the highly dynamic regions of the chosen features, and then (iii) use the subsampled matrix based on the highly dynamic regions for the  coefficient  recovery. This Narrow-fit procedure focuses on  the regions with higher dynamical behaviors for the features of interest, so that these regions  play a larger role in the coefficient recovery.  

\textbf{Features of interest:} 
We focus on a small group of features which give the variation information for the differential equation, thus highlighting which rows to choose for the coefficient recovery. 
In this paper, we choose the features of interest as follows: In 1D, we choose the features with $(\alpha, \beta) = (1,2)$ for the case of one variable in 1D and $(\alpha,\beta_u,\beta_v) = (1,2,0),(1,0,2)$ for the case of 2 variables ($u$ and $v$) in 1D. For a scalar equation in 1D,  the features of interest correspond to $\displaystyle{\frac{\partial}{\partial x} u^2}$. This term is $u u_x$,  giving combined information about $u$ and $u_x$. For a system with two variables $u,v$ in  1D, there are two features of interest: $\displaystyle{\frac{\partial}{\partial x} u^2}$ and $\displaystyle{\frac{\partial}{\partial x} v^2}$.
In 2D, we choose the features with  $(\alpha_x, \alpha_y, \beta) = \{(1,0,2), (0,1,2),(1,1,3) \} $ for a scalar equation in 2D, i.e., the features of interest are $\displaystyle{\frac{\partial}{\partial x} u^2}$, $\displaystyle{\frac{\partial}{\partial y} u^2}$ and $\displaystyle{\frac{\partial^2}{\partial x \partial y} u^3}$. 
For the case of 2 variables ($u$ and $v$) in 2D, $(\alpha_x, \alpha_y, \beta_u, \beta_v) = \{(1,0,2,0), (0,1,2,0), (1,0,0,2), (0,1,0,2),(1,0,2,1), (0,1,1,2) \}$, that is there are six features of interest:   $\displaystyle{\frac{\partial}{\partial x} u^2}$, $\displaystyle{\frac{\partial}{\partial x} v^2}$,
$\displaystyle{\frac{\partial}{\partial y} u^2}$, $\displaystyle{\frac{\partial}{\partial y} v^2}$,
$\displaystyle{\frac{\partial}{\partial x} u^2v}$, and 
$\displaystyle{\frac{\partial}{\partial y} uv^2}$. 
We explored including other terms as features of interest, but did not provide consistent improvements.

For each feature of interest, we utilize  the leading coefficient error \eqref{e: leading coefficients} to select highly dynamic regions.  
For multiple features of interest with indices $l=l_1,l_2,...,l_{\cal L}$,  we take the average  over $l$, and let
\[\bar{s}(h) =  \frac{1}{{\cal L}}  \sum_{i = 1}^{\cal L} |s(h,l_i)|, 
\]
with $s=\bar{s}$ for ${\cal L}=1$. 

\begin{figure}[t]
    \centering
    \begin{tabular}{ccc}
(a) & (b) & (c)\\    \includegraphics[width = 0.3\textwidth ]{figures_/KdVNoisyData0.5.pdf} & 
\includegraphics[width = 0.29\textwidth ]{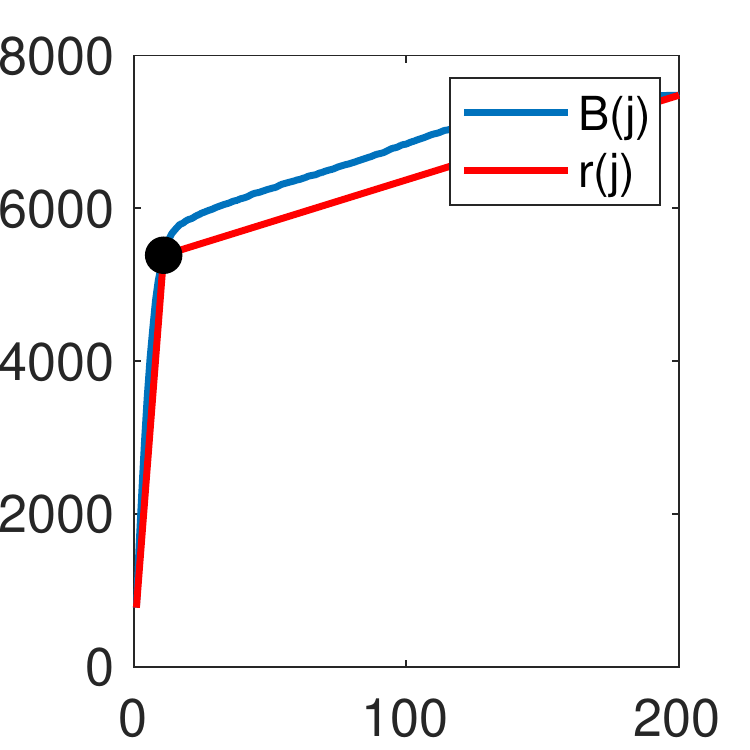} 
&
\includegraphics[width = 0.29\textwidth ]{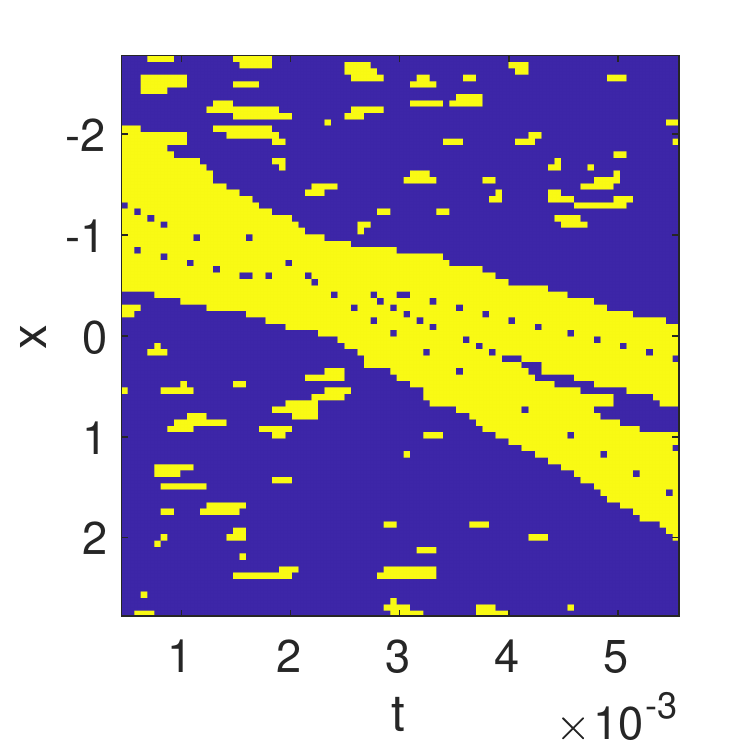}
\end{tabular}
\caption{Highly dynamic regions  for an experiment using the KdV equation \eqref{e: pde kdv} with $\sigma_{\rm NSR}=0.5$. (a) The given noisy data $\hat{U}$ with $\sigma_{\rm NSR}=0.5$ in $x-t$ plane.  (b) The separation point $\Gamma$(black) for $\mathbb{H}$ (\ref{eq:Hset}) is found, from the accumulated function $B(j)$ (blue) and the fitted piecewise linear function $r(j)$ with one junction at $\Gamma$ (red).  (c) The location of highly dynamic regions in the $x-t$ plane.}
\label{fig:highDynamic}
\end{figure}
\textbf{Highly dynamic regions:} We consider the set $\mathcal{S}=\{\bar{s}(h) | h=1, \dots, H\}$, which is the collection of averaged leading coefficient errors over the features of interest.
We  divide the set $\mathcal{S}$ into  mildly and highly dynamic regions, automatically identifying the transition point $\Gamma$ between these two types of dynamics as follows. 

After partitioning the histogram of $\mathcal{S}$ into $N_{\cal S}$ bins $(b_1,b_2,...,b_{N_{\cal S}})$, we consider the cumulative sum of the bins $B(j) = \sum_{i=1}^j b_i$. We used $N_{\cal S} = 200$ for PDEs and $N_{\cal S} = 100$ for ODEs in this paper.  We fit the  function
$B(j)$ with a piecewise linear function $r(j)$ with one junction point,  using the cost function $ \sum_j (B(j)-r(j))^2/B(j)^2$.  The junction point $\Gamma$ separates the highly dynamic   and  mildly dynamic regions.  Any $h$ with  $\bar{s}(h) \geq \Gamma$ gives the highly dynamic region  $\Omega_h$ which  we  include for the coefficient recovery.  Let the collection of the row indices of highly dynamic regions  be an ordered set:
\begin{equation} \label{eq:Hset}
    \mathbb{H} =\{ h_i \;| \; \bar{s}(h_i) \geq \Gamma, \; h_i< h_j \text{ for } i <j \}.
\end{equation}

Figure \ref{fig:highDynamic}  illustrates
 how the transition point $\Gamma$ is computed in (b) from the given data in (a). Figure \ref{fig:highDynamic} (c) shows  the locations in $x-t$ plane of the  highly dynamics regions with the index set $\mathbb{H}$.

\textbf{Narrow-fit:} We consider a submatrix  using only the ordered rows from  the highly dynamic region  $\mathbb{H}$, indicated by a subscript $\mathbb{H}$, for both $\boldsymbol{W}$ and $\boldsymbol{b}$ :
\[
\boldsymbol{W}_{\rm narrow} := \boldsymbol{W}_{\mathbb{H}}  \quad \text{ and } \quad 
\boldsymbol{b}_{\rm narrow}: = \boldsymbol{b}_{\mathbb{H}}.
\]
We also error normalize this matrix, using the rows in $\mathbb{H}$:
\begin{equation} \label{eq:Wtildenarrow}
  \tilde{\boldsymbol{W}}_{\rm narrow} = \displaystyle{ \left[\frac{\boldsymbol{w}_{1 \mathbb{H}}}{\langle s(h,1)\rangle_{\mathbb{H}}}, \frac{\boldsymbol{w}_{2 \mathbb{H}}}{\langle s(h,2)\rangle_{\mathbb{H}}}, \dots, \frac{\boldsymbol{w}_{L  \mathbb{H}}}{\langle s(h,L)\rangle_{\mathbb{H}}}\right], }
\end{equation}
where $\boldsymbol{w}_{i \mathbb{H}}$ represents the  $i^{\rm th}$ column with  the rows indexed by $\mathbb{H}$, and $\langle s(h,l)\rangle_{\mathbb{H}}$ takes  the average of $s(h,l)$ for $h \in \mathbb{H}$.  This matrix is represented in Figure \ref{fig:error_Normal} (c). 
Let $\bar{b}=\langle\boldsymbol{b}_{narrow}\rangle$ be the average of the entries of $\boldsymbol{b}_{\rm narrow}$. After narrow-fitting, We solve: 
\begin{equation}\label{eq:narrow}
\tilde{\boldsymbol{W}}_{\rm narrow}\tilde{\boldsymbol{c}} = \tilde{\boldsymbol{b}}_{\rm narrow} \quad \text{ where } \quad 
\tilde{\boldsymbol{b}}_{\rm narrow} = \boldsymbol{b}_{\rm narrow}/\bar{b}. 
\end{equation}
We then compute the coefficient $\boldsymbol{c}$ by rescaling $\tilde{\boldsymbol{c}}$ back: 
\begin{equation}\label{eq:cUpdate}
\boldsymbol{c} = \bar{b} \; \tilde{\boldsymbol{c}} \; \text{diag} \left\{ \frac{1}{\langle s(h,1)\rangle_{\mathbb{H}}}, \frac{1}{\langle s(h,2)\rangle_{\mathbb{H}}}, \dots, \frac{1}{\langle s(h,L)\rangle_{\mathbb{H}}} \right\}.
\end{equation}

\subsection{Trimming the support 
} \label{ss:supression}

After the coefficient values in $\boldsymbol{c}$ are recovered, some  features give very small contributions to $u_t$.   We further trim the support by eliminating these features corresponding to small contributions. 

From the solution $\tilde{\boldsymbol{c}}$ of the linear equation (\ref{eq:narrow}), we define a  \textbf{contribution score} $a_i$ of each feature as 
\begin{equation} \label{eq:CScore}
a_i =  \frac{ n_i }{\max_{i \leq L} n_i} 
\quad \text{where} \quad 
n_i =  || \tilde{\boldsymbol{w}}_i ||_2 |\tilde{c}_i|, \quad  i=1,2,\dots, L. \end{equation} 
Here $\tilde{\boldsymbol{w}}_i$ denotes the  $i^{\rm th}$column of $\tilde{\boldsymbol{W}}_{\rm narrow}$. We consider the $L_2$ norm of this column multiplied by the coefficient value of the  $i^{\rm th}$ component of $\tilde{c}$.
Since $a_i$ is normalized by the maximum value of $n_i$, $a_i$ gives the score of the contribution of the $i^{\rm th}$ feature relative to the contribution of the feature with the largest contribution. 

We trim the coefficient, thus the feature,  when the contribution score of that feature is below $\mathcal{T}$, i.e. $a_i<  \mathcal{T}$.
Typically, we set $\mathcal{T}=0.05$ to trim the features with contributions less than $5\%$ of $u_t$.  Each time [Step 3] is called to trim the support set $\mathcal{A}^k_j$ to the new support set  $\mathcal{A}^k_{j+1}$, and [Step 2] narrow-fit is called to find the updated coefficient value $c(k,j+1)$. 

\begin{figure}[t]
\centering
\begin{tabular}{cc}
    (a) $\sigma_{\rm NSR}=0.1$ &  
    (b) $\sigma_{\rm NSR}=0.5$ \\
    \includegraphics[width = 0.48 \textwidth]{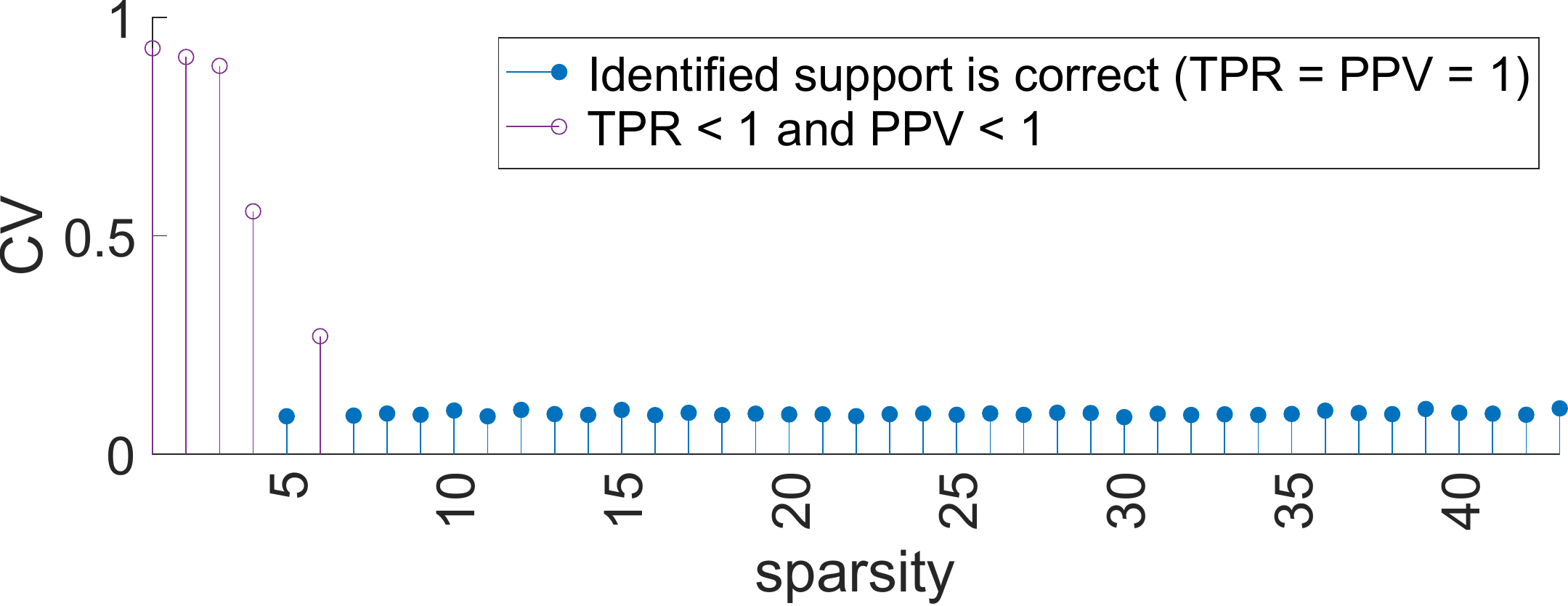} &
    \includegraphics[width = 0.48 \textwidth]{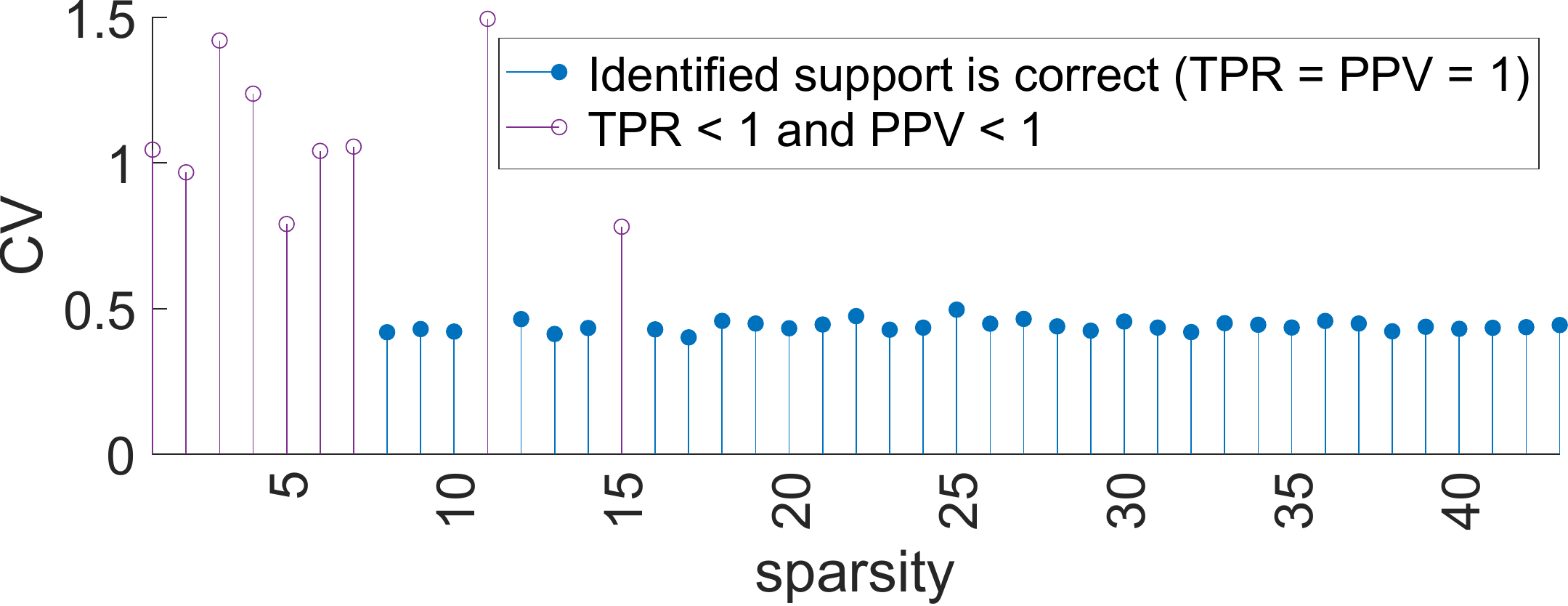} 
\end{tabular}
\caption{ Trimming is demonstrated in an experiment using the KS equation \eqref{e: pde KS}.  For each sparsity level $k$ in $x$-axis, the bar shows the cross validation (\ref{eq:CV}) of the recovered coefficient $c(k^*,J_{k^*})$.  Notice for most sparsity levels 5 and above the correct support is found.   After SP finds $k$ supports, the trimming step reduces the support until only the correct ones are left. 
Here $\sigma_{NSR}$ is the noise-to-signal ratio \eqref{e: signal noise ratio}, TPR is true positive rate \eqref{e: err tpr} and PPV is positive prediction value \eqref{e: err ppv}.}
\label{F:cScore}
\end{figure}

Figure \ref{F:cScore} shows the effect of trimming.  For each sparsity level $k$ in $x$-axis, the bar shows the cross validation value (\ref{eq:CV}) of the recovered coefficient $c(k^*,J_{k^*})$.  For a large sparsity level, thanks to the trimming step, the correct support and coefficient values are found.  

\subsection{Algorithms} \label{ssec:algo}

Our WeakIdent algorithm is summarized in Algorithm \ref{alg: WeakIdent}. 
From the  linear system in (\ref{e: Wc=b}) 
$$  \boldsymbol{W}\boldsymbol{c} = \boldsymbol{b},$$
we input $\boldsymbol b$ and $\boldsymbol{W}$ computed through  (\ref{e: approx of bW}), with subsampling in \eqref{eq:H}.  For each sparsity level $k =1, 2, \dots, K$, 

\begin{enumerate}
\item[\textbf{[Step 1]}]{First, Subspace Pursuit (SP)\cite{dai2009subspace} is applied to find $\mathcal{A}^k_o = \text{supp}\{SP(\boldsymbol{W} \dagger, \tilde{\boldsymbol{b}}, s) \}$ using the column normalized matrix $\boldsymbol{W} \dagger$ in (\ref{eq:Wdagger}) and $\tilde{\boldsymbol{b}} = \boldsymbol{b}/||\boldsymbol{b}||$. 
}   
\item[\textbf{[Step 2]}]{Narrow-fit. To recover the coefficient values using the support $\mathcal{A}_j^k$, we find the row index set $\mathbb{H}$ of highly dynamic regions in (\ref{eq:Hset}), and solve 
\[\tilde{\boldsymbol{W}}_{\rm narrow}\tilde{\boldsymbol{c}} = \tilde{\boldsymbol{b}}_{\rm narrow}\]
in (\ref{eq:narrow}) and get  $\boldsymbol{c}(k,j)$ in (\ref{eq:cUpdate}). }
\item[\textbf{[Step 3]}]{Trimming. Update to  $\mathcal{A}_{j+1}^k$, if there is any column with the contribution score in (\ref{eq:CScore}) below $\mathcal{T}$, i.e.  $a_i <\mathcal{T} $ .  If trimmed, move to [Step 2] to get a new updated $\boldsymbol{c}(k,j+1)$. If no column is trimmed, move to [Step 4] and set $J_k=j$.  }
\item[\textbf{[Step 4]}]{Cross Validation. With the  support $\boldsymbol{c}(k,J_k)$ computed for each sparsity level $k =1, \dots, K$, we select the final support by finding the $k^*$ which gives the minimum cross-validation error. For a sparsity level $k$, we randomly sample regions  from the $N_xN_t$ regions and equally partition these regions into two sets indexed by $\mathbb{A}$ and $\mathbb{B}$ respectively. %
We consider the linear system in \eqref{eq:narrow}: 
\begin{equation}
\tilde{\boldsymbol{W}} = \displaystyle{ \left[\frac{\boldsymbol{w}_{1 }}{\langle s(h,1)\rangle_{\mathbb{H}}}, \frac{\boldsymbol{w}_{2 }}{\langle s(h,2)\rangle_{\mathbb{H}}}, \dots, \frac{\boldsymbol{w}_{L}}{\langle s(h,L)\rangle_{\mathbb{H}}}\right]},\; \quad \text{and}\quad \Tilde{\boldsymbol{b}} = \boldsymbol{b}/\bar{b}
    \label{e: Wtilda big}
\end{equation}
utilizing the highly dynamic region error normalization for the large full matrix.
Here $ \mathbb{H}$ indicates   ordered row index from the set $\mathbb{H}$, and $\langle s(h,l)\rangle_{\mathbb{H}}$ taking the  average of $s(h,l)$ for $h \in \mathbb{H}$. 
We solve least square problems  $\boldsymbol{\Tilde{W}}_{\mathbb{A}}\boldsymbol{\Tilde{c}}_{\mathbb{A}} = \boldsymbol{\Tilde{b}}_{\mathbb{A}}$ and $\boldsymbol{\Tilde{W}}_{\mathbb{B}} \boldsymbol{\Tilde{c}}_{\mathbb{B}} = \boldsymbol{\Tilde{b}}_{\mathbb{B}}$, where $\tilde{\boldsymbol{W}}_{\mathbb{A}}$ and $\mathbb{B}$ contain the  rows of $\tilde{\boldsymbol{W}}$  indexed by  $\mathbb{A}$ and  $\mathbb{B}$ respectively. 
Then, we compute the \textbf{cross validation} (CV) error 
\begin{equation} \label{eq:CV}
{\rm CV}(k)=
\lambda ||\boldsymbol{\Tilde{W}}_{\mathbb{A}} \boldsymbol{\Tilde{c}}_{\mathbb{B}} - \boldsymbol{\Tilde{b}}_{\mathbb{A}} ||_2 + (1-\lambda) ||\boldsymbol{\Tilde{W}}_{\mathbb{B}} {\boldsymbol{\Tilde{c}}}_{\mathbb{A}} - \boldsymbol{\Tilde{b}}_{\mathbb{B}} ||_2, 
\end{equation}
where we set $\lambda = 1/100$.
In practice,  for each $k$, we generate 30 different random partitions of $\mathbb{H}$ to $\mathbb{A}$ and $\mathbb{B}$,  then select the minimum:
\begin{equation}\label{eq:cFinal}
 \boldsymbol{c}(k^*,J_{k^*}) = \argmin_k \{ {\rm CV}(k) | k = 1,2,...,L \}.   
\end{equation}
Here  $K \leq L$, since $L$ is the total number of features in the dictionary.  In practice, a small $K$ is needed.   Figure \ref{F:cScore} illustrates that for  (small) values of  $K$ around $K=10$ and below, the correct coefficients are found, thanks to the trimming step.  
}
\end{enumerate}

\begin{algorithm}[t]
\caption{WeakIdent Algorithm} \label{alg: WeakIdent}
\begin{algorithmic}
\State \textbf{Input:} $\boldsymbol{W} \in \mathbb{R}^{H \times L}, \boldsymbol{b} \in \mathbb{R}^{H}$, from \eqref{e: Wc=b} uniformly subsampled as (\ref{eq:H}); Parameter $\mathcal{T}=0.05$\\
\For{k = 1,2,...,K}
{
\State{[Step 1] $\mathcal{A}_0^k = \text{supp}\{SP(\boldsymbol{W} \dagger, \boldsymbol{\Tilde{b}}, s)\}$\, use SP \cite{dai2009subspace} and set $j=0$;} 
\State{[Step 2] Find $c(k,j)$ by narrow-fit (\ref{eq:narrow})\;}

\While{there exists $a_i < \mathcal{T}$ as in (\ref{eq:CScore})}
{
\State{[Step 3] Trim as in Subsection \ref{ss:supression} and set $j=j+1$ \;} 
\State{[Step 2] Find $c(k,j+1)$ by Narrow-fit (\ref{eq:narrow})\;}
}
}
\State{Among $k=1,\dots, K$, find $\boldsymbol{c}(k^*, J_{k^*})$ by Cross Validation in (\ref{eq:cFinal}).}
\State{\textbf{Output:} ${\boldsymbol{c}}= \boldsymbol{c}(k^*, J_{k^*}) \in \mathbb{R}^L$ such that $\boldsymbol{W} {\boldsymbol{c}} \approx \boldsymbol{b}$.}
\end{algorithmic}
\end{algorithm}

\section{WeakIDENT results and comparions}
\label{sec: numerical experiments}

\begin{table}[t!]
\setlength\abovedisplayskip{0pt}
\setlength\belowdisplayskip{0pt}
\centering
\begin{tabularx}{\textwidth}{CC}
\toprule
Equation  & Parameters \\
\midrule
\multirow{3}{=}{
Transport equation
\begin{equation}
\displaystyle{\frac{\partial u}{\partial t}} = - \frac{\partial u}{\partial x} + 0.05 \frac{\partial^2 u}{\partial x^2} \quad \quad  \quad \quad \quad
\label{e: pde transport}
\end{equation}
}   & $L=43$,  $\bar{\alpha}=6$, $\bar{\beta}=6$ , $[X_1, X_2]=[0,1]$, $\Delta x=0.039$, $T=0.3$, $\Delta t=0.001$ \\
& $u(x,0) = \sin(4\pi / (1-T)x)^3\cos(\pi/(1-T)x)$\\
& \hspace{1cm} for $x < 1-T$, and 0 otherwise\\
\midrule
 
\multirow{3}{=}{
Korteweg-de Vires (KdV)
\begin{equation}\displaystyle{\frac{ \partial u }{\partial t}} = -0.5u\frac{\partial u }{ \partial x} - \frac{\partial^3 u }{\partial x^3} \quad \quad \quad \quad \quad
 \label{e: pde kdv}
 \end{equation}} 
 &
$L= 43$,  $\bar{\alpha}= 6$, $\bar{\beta}= 6$ ,$[X_1, X_2]= [-\pi, \pi]$, $\Delta x= 0.0157$,  $T= 0.006$, $\Delta t= 10^{-5}$ \\
& $u(x,0) = 3.0 \times 25 ^2 * \text{sech} ( 0.5 \times ( 25 \times ( x + 2.0 ) ) ) ^ 2$ \\
& \hspace{1cm} $+ 3.0 \times 16^2 * \text{sech} ( 0.5 \times ( 16 * ( x + 1.0 ) ) ) ^ 2$\\
\midrule
\multirow{3}{=}{
Kuramoto-Sivashinsky (KS) 
\begin{equation}
\displaystyle{\frac{\partial u}{\partial t}} = - u - \displaystyle{\frac{\partial^2 u }{\partial x^2}} - \displaystyle{\frac{\partial^4 u }{\partial x^4}} \quad \quad  \quad \quad \;
\label{e: pde KS}
\end{equation}} 
 &
$L= 43$,  $\bar{\alpha}= 6$, $\bar{\beta}= 6$,
$[X_1, X_2]= [0,100.53]$, $\Delta x= 0.3927$, $T= 150$, $\Delta t= 0.5$ \\
& $u(x,0) = \cos(x/16) (1+\sin(x/16)). $\\
\midrule
 \multirow{4}{=}{
Nonlinear Schrodinger (NLS) (1D) 
\begin{equation}
 \begin{cases}
 \displaystyle{\frac{  \partial u }{ \partial t}} & = 0.5\displaystyle{ \frac{ \partial^2 v }{  \partial x^2}+u^2v +v^3}\\ 
  \displaystyle{\frac{ \partial v }{ \partial t}} &= -0.5\displaystyle{\frac{   \partial^2 u }{  \partial x^2}} -uv^2 -u^3 \quad \;
 \end{cases}
 \label{e: pde NLS}
\end{equation}} \\
 & $L= 190$,  $\bar{\alpha}= 6$, $\bar{\beta}= 6$ \\
& $[X_1, X_2]= [-5,5]$, $\Delta x= 0.0391$\\
& $T= 3.1416$, $\Delta t= 0.0126$ \\
\\
\midrule
\multirow{3}{=}{
Anisotropic Porous Medium (PM) (2D) 
\begin{equation}
\displaystyle{\frac{\partial u }{\partial t}}  =  +0.3\displaystyle{\frac{\partial^2 u^2 }{\partial^2 y}} -0.8 \displaystyle{\frac{\partial^2 u^2 }{\partial x \partial y}} +\displaystyle{\frac{\partial^2 u^2 }{\partial^2 x}}  \quad \quad
\label{e: pde PM}
\end{equation}} 
& $L= 65$,  $\bar{\alpha}= 4$, $\bar{\beta}= 4$\\
& $[X_1, X_2]= [-5,5]$, $\Delta x= 0.0503$\\
& $T= 5$, $\Delta t= 0.0503$ \\
\midrule

\multirow{2}{=}{
Reaction-Diffusion (2d) 
\begin{equation}
\begin{cases}\displaystyle{\frac{ \partial u }{\partial t}} & =     0.1  \displaystyle{\frac{\partial^2 u }{\partial^2 x}} +0.1  \displaystyle{\frac{\partial^2 u }{\partial^2 t}} +u + v^3  -uv^2   +u^2v - u^3 \\ 
\displaystyle{\frac{\partial v }{\partial t}} & = 0.1 \displaystyle{\frac{\partial^2 v }{\partial^2 y}}   +0.1 \displaystyle{\frac{\partial v^2 }{\partial^2 x}} + v -v^3 - uv^2 - u^2v - u^3
\end{cases}
\label{e: pde RD}
\end{equation}} 
 &  
$L= 155$,  $\bar{\alpha}= 4$, $\bar{\beta}= 5$, $[X_1, X_2]= [-10,10]$ \\
& $\Delta x= 0.0781$, $T= 9.9219$, $\Delta t= 0.0781$ \\ 
& $u(x,y,0) = \text{tanh}(\sqrt{x^2 + y^2} \cos ( \theta (x+iy) - \pi \sqrt{x^2+y^2} )$, \\
& $v(x,y,0) = \text{tanh}(\sqrt{x^2 + y^2} \sin(\theta (x+iy) - \pi \sqrt{x^2 + y^2})$\\
& \\
\bottomrule
\end{tabularx}
\caption{A list of PDEs considered in this paper.  Here $L$ is the total  number of features, $\bar{\alpha}$ is the  highest order of partial derivative,  $\bar{\beta}$ is the highest degree used in $f_l$ in \eqref{e: feature formulation}, $[X_1,X_2]$ is the range of the spatial domain,  $T$ is the final time for simulation. $\Delta x$ and $\Delta t$ are the spatial and temporal increment of the given data.  The set up of \eqref{e: pde kdv}, \eqref{e: pde KS}, \eqref{e: pde NLS}, \eqref{e: pde PM}, and \eqref{e: pde RD} are identical to \cite{messenger2021weakPDE}. } \label{T:PDE}
\end{table}

In this section, we provide detailed experimental results.  
We summarize a list of PDEs and ODE systems in Table \ref{T:PDE} and  Table \ref{T: odes}. 
For the systems of ODEs, we consider features with polynomial order between 3 and 5 ,  with $L\leq 21$ for all the cases. 
For the systems of PDEs, we consider features with both polynomial order and derivative order between 4 and 6, which gives a dictionary of size $L\leq 65$ for the  1 spatial dimension   and $L \leq 190$ for 2 spatial dimensions. 
Simulation and feature details are presented in Table \ref{T:PDE} and \ref{T: odes} for each experiment.

For PDEs, $N_x$ and $N_t$ are chosen such that $N_xN_t\in (1,000,3,000)$ to reduce the computational cost. In particular, we set 
\begin{equation} \label{eq:NxNt}
N_x =\lceil \frac{\mathbb{N}_x - 2m_x-1}{\lfloor \mathbb{N}_x/\mathbf{N} \rfloor} +1\rceil \;\; \text{ and } \;\; 
N_t =\lceil \frac{\mathbb{N}_t - 2m_t-1}{\lfloor \mathbb{N}_t/\mathbf{N} \rfloor} +1 \rceil,
\end{equation}
with $\mathbf{N}=50$ as a default choice. Here $\lceil \cdot \rceil$ and $\lfloor \cdot \rfloor$ denotes the ceiling and floor operator. 
In Table \ref{T:PDE}, \eqref{eq:NxNt} is used for the transport question \eqref{e: pde transport}, the KS equation \eqref{e: pde KS}   and the nonlinear Schrodinger equation  \eqref{e: pde NLS}. 
For certain cases such as  the KdV equation \eqref{e: pde kdv} where $|\mathbb{H}|$ is very small, we increase $N_x$ and $N_t$, e.g., using $\mathbf{N}=70$, such that $|\mathbb{H}|>800$.
For the spatially 2 dimensional  cases, we use $\mathbf{N}=(25,25)$ for the anisotropic porous medium equation (PM) \eqref{e: pde PM}, and $\mathbf{N}=(19,16)$ for the 2D reaction-diffusion equation \eqref{e: pde RD} to reduce the time of computation.  
For the ODEs listed in Table \ref{T: odes}, we choose $N_t\approx 1000$ by default with $\mathbf{N}=1000$. Since we use different subsampling, we present additional comparisons in Section \ref{sss:subsample} to demonstrate that the effect of subsampling on the result is minimal.

\begin{table}[t]
\centering
\setlength\abovedisplayskip{0pt}
\setlength\belowdisplayskip{0pt}
\centering
\begin{tabularx}{\textwidth}{sbs}
\toprule
 Name & \multicolumn{1}{c}{Equation} & parameters \\ 
\midrule
2D Linear System & \multirow{3}{=}{\begin{equation}
    \displaystyle{\frac{d}{dt}} \begin{bmatrix}
    x \\ y
    \end{bmatrix} = \begin{bmatrix}
    -0.15 & 2.5 \\ -2.5 & -0.15
    \end{bmatrix}
    \begin{bmatrix}
    x \\ y
    \end{bmatrix}
    \label{e: ode, linear 2d}
    \end{equation}}
     & 
      $(x_0, y_0)=(2,50)$ ,\\
      && $\Delta t =  0.01 , T = 10$ \\
      && $L = 21, \bar{\beta}= 5$  \\
      \midrule
 2D Nonlinear& \multirow{2}{=}{\begin{equation}\displaystyle{\frac{d}{dt}} \begin{bmatrix}
    x \\ y
    \end{bmatrix} = \begin{bmatrix}
    0 & 1 & 0 \\
    4 & -1 & -4 \\
    \end{bmatrix}
    \begin{bmatrix}
    x \\ y \\ x^2y
    \end{bmatrix}
    \label{e: ode, van-der-pol}
    \end{equation} }  &
    $(x_0, y_0) = (0,1)$\\
     (Van der Pol)& & $\Delta t =  0.001, T =  15  $\\
    & & $L =  21, \bar{\beta}= 5$\\
    \\
    \midrule
 2D Nonlinear & \multirow{3}{=}{\begin{equation}\displaystyle{\frac{d}{dt}} \begin{bmatrix}
    x \\ y
    \end{bmatrix} = \begin{bmatrix}
    0 & 1 & 0 \\
    -0.2 & -0.05 & -1  \\
    \end{bmatrix}
    \begin{bmatrix}
    x \\ y \\ x^3
    \end{bmatrix}
    \label{e: ode, Duffing}
    \end{equation}}& 
    $(x_0,y_0) = (0,2)$\\
    (Duffing) & & $\Delta t= 0.01, T = 10 $\\
    & & $L = 21, \bar{\beta }=5$\\
    \\
    \midrule
 2D Nonlinear& \multirow{4}{=}{\begin{equation}
\displaystyle{\frac{d}{dt}} \begin{bmatrix}
    x \\ y
    \end{bmatrix} = \begin{bmatrix}
    0.67 & 0 & -1.33 \\
    0 & -1 & 1  \\
    \end{bmatrix}
    \begin{bmatrix}
    x \\ y \\ xy
    \end{bmatrix}\label{e: ode, Lotka}\end{equation}}& 
    $(x_0, y_0) = (10,10)$\\
    & \\
    (Lotka-Volterra) & & $ \Delta t = 0.05, T = 50$ \\
    & & $L = 21,\bar{\beta}= 5$\\
    \midrule
 & \multirow{6}{=}{\begin{equation}\displaystyle{\frac{d}{dt}} \begin{bmatrix}
    x \\ y \\ z
    \end{bmatrix} = \begin{bmatrix}
    -10.2 & 10.2 & 0 & 0 & 0 \\
    29 & -1 &0 & 0 & -1 \\
    0 & 0 & -2 & 1 & 0\\
    \end{bmatrix}
    \begin{bmatrix}
    x \\ y \\ z \\ xy \\ xz \\
    \end{bmatrix}\label{e: ode, lorenz}\end{equation}} & \\
     3D Nonlinear & & 
  $(x_0, y_0, z_0) = (-8 , 7 , 10)
  $\\
  (Lorenz) & & $\Delta t = 0.001 , T = 15$ \\
  & & $L = 20 , \bar{\beta} = 3 $    \\
  & \\
  & \\
\bottomrule
\end{tabularx}
\caption{A list of ODEs considered in this paper. This table includes the initial condition, the temporal increment $\Delta t$, the total simulation time $T$, the total number of  features $L$ and the highest degree of polynomials $\bar{\beta}$ in \eqref{e: feature formulation} for each equation. The Solution is simulated with RK45 with tolerance $10^{-10}$.}
\label{T: odes}
\end{table}

The experiments are performed on both clean data and noisy data with various Noise-to-Signal Ratio,  $\sigma_{\rm NSR}$ defined as follows:
\begin{equation} \sigma_{\rm NSR} = \frac{\displaystyle \epsilon_i^n}{\frac{1}{\mathbb{N}_t\mathbb{N}_x}\displaystyle \sum_{{i},n}|{U}_i^n - (\max_{i, n} {U}_i^n +\min_{i, n} {U}_i )/2 |^2} 
\label{e: signal noise ratio}
\end{equation}
for $i = 1,...,\mathbb{N}_x$, $n = 1,...,\mathbb{N}_t$.  Note that our definition of NSR reflects the local variation of the given data.  This is different from the absolute variation (absolute root mean squared of ${U}_i^n$ ) $\sigma_{\rm NR}$ used in \cite{messenger2021weakPDE}, and this $\sigma_{\rm NSR}$ value tends to be smaller than  the $\sigma_{\rm NR}$ value.  We also mention the $\sigma_{\rm NR}$ value  in the following experiments when it is relevant.  
We use Gaussian noise, such that $\epsilon_i^n \sim \mathcal{N}(0, \sigma_{\rm NSR})$ for $\epsilon_i^n$, and $\hat{U}_i^n$ in \eqref{e: noisy discretized data}.
For the case of multiple variables, we compute \eqref{e: signal noise ratio} for each variable.

\textbf{Error measures}:  To quantify the quality of the recovery, we utilize different error measurements  listed in Table \ref{T: different errors}. 
\begin{table}[]
\setlength\abovedisplayskip{0pt}
\setlength\belowdisplayskip{0pt}
\centering
\begin{tabularx}{\textwidth}{p{5cm} p{10cm}}
\hline 
Relative coefficient Error $l_2$ &\begin{equation}
E_{2}  = {|| \boldsymbol{c}^* - \boldsymbol{c} ||_2}/{||\boldsymbol{c}^*||_2} \hfill
\label{e: err l2}\end{equation}\\
Relative coefficient Error $l_\infty$ & \begin{equation}
E_{\infty} = \max_l\{|\boldsymbol{c}^*(l) - \boldsymbol{c}(l)|/ |\boldsymbol{c}^*(l)|:  \boldsymbol{c}^*(l) \neq 0\} \hfill 
\label{e: err linfty}\end{equation}\\
True Positive Rate &
\begin{equation}
\text{TPR} = \displaystyle{| \{l: \boldsymbol{c}^*(l) \neq 0, \; \boldsymbol{c}(l) \neq 0 \}| / | \{ l: \boldsymbol{c}^*(l) \neq 0\}|} 
\hfill \label{e: err tpr} \end{equation}\\
Positive Predictive Value & \begin{equation}
\text{PPV} = \displaystyle{| \{l: \boldsymbol{c}^*(l) \neq 0, \; \boldsymbol{c}(l) \neq 0 \}| / | \{ l: \boldsymbol{c}(l) \neq 0\}|} 
\hfill \label{e: err ppv}
\end{equation}\\
Residual Error & \begin{equation}
E_{\rm res}= {||\boldsymbol{W}\boldsymbol{c}  - \boldsymbol{b}||_2}/{||\boldsymbol{b}||_2} \hfill \label{e: err e_res} \end{equation} \\
Dynamic Error & \begin{equation} E_{\rm dyn} =  \sum_{1 \leq i \leq \mathbb{N}_x,1 \leq n \leq \mathbb{N}_t }(|U^n_{\text{i,forward}} - U^n_{\textbf{i},clean}|^2) / (\mathbb{N}_x\mathbb{N}_t) \hfill \label{e: err e_dyn} \end{equation} \\
\bottomrule
\end{tabularx}
\caption{Error measurements used for comparisons.}
\label{T: different errors}
\end{table}
The relative coefficient errors $E_2$ in \eqref{e: err l2}  and $E_\infty$ in \eqref{e: err linfty} measure the accuracy of the recovered coefficients $\boldsymbol{c}$ against the true coefficients $\boldsymbol{c}^*$ in terms of the $l_2$ and the infinity norm, respectively. 
We introduce two new measures to quantify the accuracy of the support recovery.  The  True Positive Rate (TPR)  \footnote{The definition of TPR in \eqref{e: err tpr}  is different from that used in \cite{messenger2021weakPDE}} \eqref{e: err tpr} measures the fraction of features that are found out of all features in the true equation, and is defined as the ratio of the cardinality of the correctly identified support over the cardinality of the true support. The TPR is $1$ if all the true features are found. The Positive Predictive Value  (PPV) \eqref{e: err ppv} indicates the presence of false positives: it is
the ratio of the cardinality  of  the correctly identified support over the total cardinality of the identified support.  The PPV is $1$ if the recovered support is also in the true support.  
The residual error $E_{\rm res}$ in \eqref{e: err e_res}, which is also used in \cite{he2020robust},  measures the relative difference between the learned differential equation and the given data. 
To show the effectiveness of WeakIdent in the recovery of the dynamics, we define the dynamical error $E_{\rm dyn}$ in \eqref{e: err e_dyn} to measure the difference between the true dynamics and the  expected dynamics simulated from the recovered equation. In \eqref{e: err e_dyn}, we use $U^n_{\text{\textbf{i},\rm forward}}$ and $U^n_{\textbf{i},\rm clean}$ to denote the simulated data and the true data without noise. We simulate ODEs using RK45 with the relative error tolerance to be $10^{-10}$.  This is measured for ODEs only, due to restricted stability conditions for PDEs.  
If  the identified equation blows up before the final time $T$ is reached, we compare $U^n_{\text{\textbf{i},\rm forward}}$ and $U^n_{\textbf{i},\rm clean}$ just before the blow-up.

\subsection{WeakIdent results and comparisons for PDEs }

We present the WeakIdent results, and compare with existing methods, such as the IDENT in  \cite{kang2021ident}, the Robust Ident, with Subspace pursuit Cross validation (SC) and Subspace pursuit Time evolution (ST) in \cite{he2020robust}, 
SINDy \cite{brunton2016discovering}, and methods using the weak form such as  RGG \cite{reinbold2020using}, Weak SINDy for first order dynamical systems (WODE) \cite{messenger2021weak}, and Weak SINDy for PDEs (WPDE) \cite{messenger2021weakPDE}.    
We note that RGG in \cite{reinbold2020using} uses a subset of features (e.g. 8-14 features), which is different from other methods which use the full feature matrix ($L =$ 21 to 190 features). For each experiment in the comparison, we specify which features are used for RGG. 
In many of the PDE experiments in this section, we show comparisons only between our proposed WeakIdent and WPDE \cite{messenger2021weakPDE}, since these two methods give the best results compared to others,  based on the error measures in Table \ref{T: different errors}.  

\subsubsection{Transport equation} 
\begin{figure}[t!]
 \begin{center}
\begin{tabular}{cc}
  (a) \\
\includegraphics[width = 0.17\textwidth]{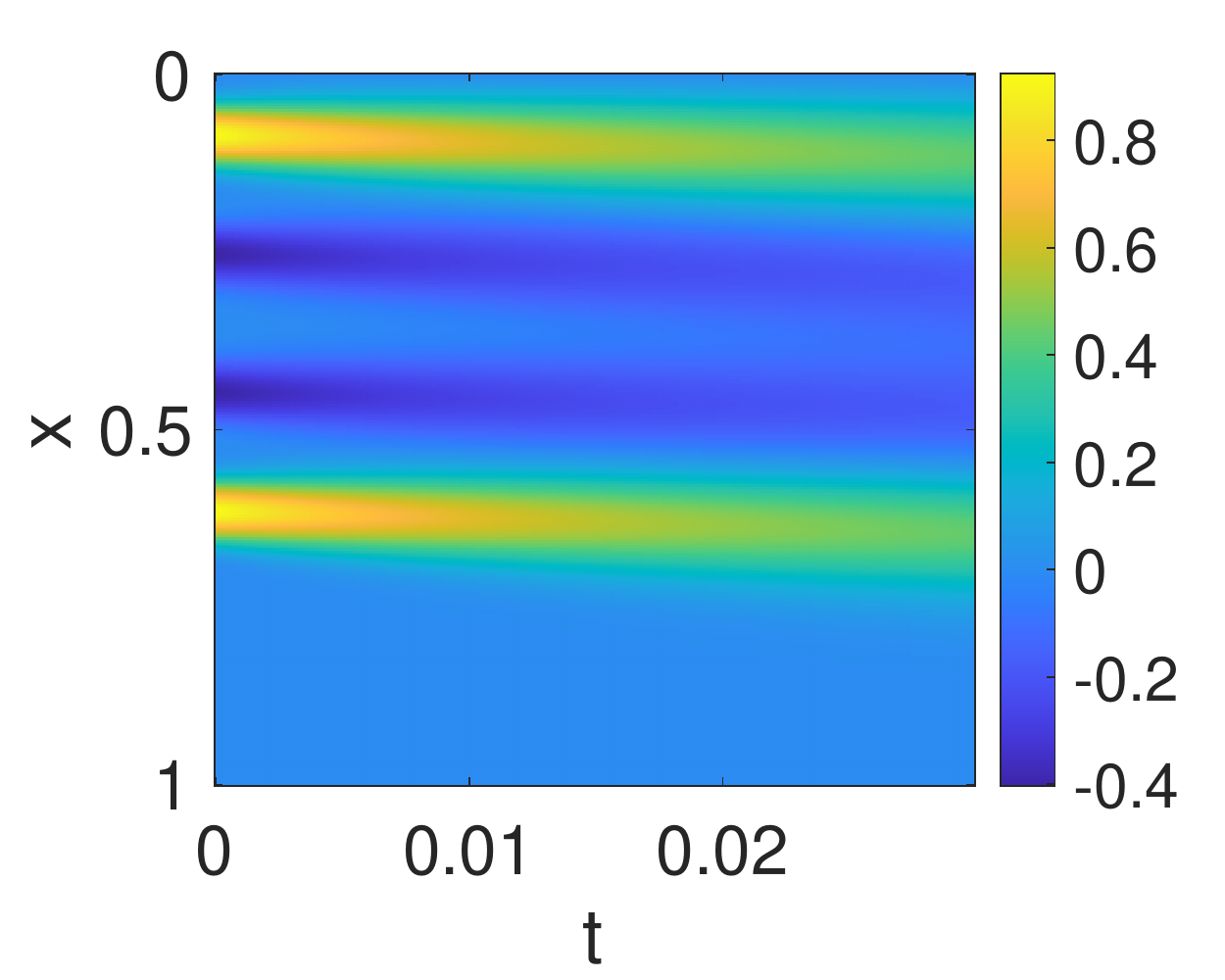}
\end{tabular}
\begin{tabular}{|l|r|r|r|r|r|r|}
\multicolumn{7}{c}{(b) $\sigma_{\rm NSR}=0$}  \\
\toprule 
                        & \textbf{WeakIdent} & WPDE \cite{messenger2021weakPDE} & RGG\cite{reinbold2020using} & IDENT\cite{kang2021ident} & SC\cite{he2020robust} & ST\cite{he2020robust} \\ \midrule
$E_2$                   & \textbf{0.001}     & \textbf{0.001}                   & \textbf{0.001}                        & -                       & 2.26           & 2.24           \\
$E_{\infty}$            & \textbf{0.001}     & \textbf{0.001}                   & \textbf{0.001}                        & -                       & -            & 3.08           \\
$E_{\rm res}$               & \textbf{0.001}     & \textbf{0.001}                   & \textbf{0.001}                       & 0.98                      & 0.03           & 0.03           \\
TPR                     & \textbf{1.0}       & \textbf{1.0}                     & \textbf{1.00}               & -                       & 0.00           & 0.50           \\
PPV                     & \textbf{1.0}       & \textbf{1.0}                     & \textbf{1.00}               & 0.00                      & 0.00           & 0.20           \\ \bottomrule
\end{tabular}
\begin{tabular}{|l|l|}
\multicolumn{2}{c}{(c) $\sigma_{\rm NSR} =0$ } \\
\toprule
True equation & $u_t = -1.00000 u_x + 0.05000 u_{xx}$ \\
\hline
\textbf{WeakIdent} & $\bf u_t    =  -1.00145 u_{x} +0.04999   u_{xx}$  \\
WPDE \cite{messenger2021weakPDE}  & $\bf u_t  =    -1.00144 u_{x} +0.05000 u_{xx}$ \\
RGG \cite{reinbold2020using} &  $\bf u_t =  -1.00119 u_x +0.04999 u_{xx}$\\
IDENT\cite{kang2021ident} & $u_t =    -0.0006
    +0.0036u
    +0.0244u^2
  -0.9992u_x
    +0.0004(u_x)^2+...$
  \\
SC\cite{he2020robust} & $ u_t  =    +1.74039 u^2 -1.03236 u_x +0.05168 u_{xx} +0.00298 uu_{xx}  $ \\
ST\cite{he2020robust} & $ u_t  =    +1.73061 u^2 -1.01121 u_x -0.10390 uu_x +0.05167 u_{xx} +0.00298   uu_{xx}  $ \\

\bottomrule
\end{tabular}

\begin{tabular}{cc}
 (d) \\
\includegraphics[width = 0.17\textwidth]{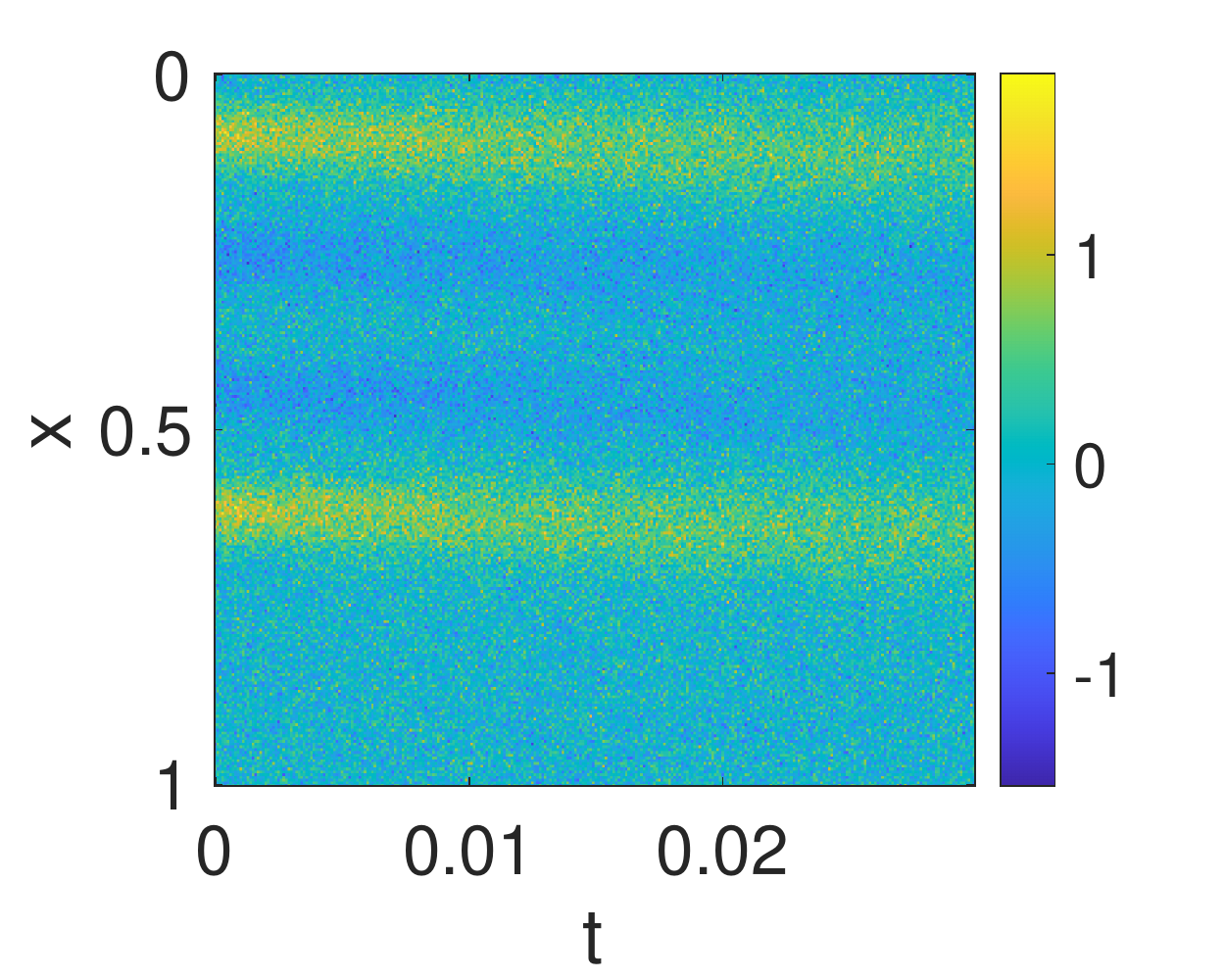}
\end{tabular}
\begin{tabular}{|l|r|r|r|r|r|r|}
\multicolumn{7}{c}{(e) $\sigma_{\rm NSR}=1$}  \\
\toprule
 & \textbf{WeakIdent} & WPDE \cite{messenger2021weakPDE} & RGG \cite{reinbold2020using} & IDENT \cite{kang2021ident} & SC
 \cite{he2020robust} & ST\cite{he2020robust} \\ \midrule
$E_2$        & \textbf{0.008} & 0.184        & 135.33 & -  & 17.43 & 20.32 \\
$E_{\infty}$ & \textbf{0.008} & 1.129        & 0.13   & -  & -   & 18.23 \\
$E_{\rm res}$    & \textbf{0.811} & 0.830        & 0.95   & 0.82 & 0.91  & 0.89  \\
TPR          & \textbf{1.0}   & \textbf{1.0} & 0.50   & 0.00  & 0.00  & 0.50  \\
PPV          & \textbf{1.0}   & 0.5          & 0.25   & 0.00 & 0.00  & 0.20  \\ 
\bottomrule
\end{tabular}

\begin{tabular}{|l|l|}
\multicolumn{2}{c}{(f) $\sigma_{\rm NSR} =1$ } \\
\toprule

True equation & $u_t = -1.00000 u_x + 0.05000 u_{xx}$ \\
\hline
\textbf{WeakIdent} & $\bf u_t    =  -1.00792 u_{x} +0.05029   u_{xx}$ \\
WPDE \cite{messenger2021weakPDE} & $u_t  =    -1.02983 u_{x} +0.10647 u_{xx} -0.15741 (u^3)_{xx} +0.07197   (u^6)_{xx}$ \\
RGG \cite{reinbold2020using} & $u_t  = -0.00003 u_{xxxx} -0.87531 u_x +44.70146 u^2 -127.91622 u^3  $ \\
IDENT\cite{kang2021ident} & $u_t =        -0.2710
    +6.1120u
   -3.4346u^2
   -0.0000u_x
    +0.0000(u_x)^2+...$
   \\
SC\cite{he2020robust} & $ u_t  =    -3.35704 u -17.09628 u^2 -0.26659 u_x    $ \\
ST\cite{he2020robust} & $ u_t  =    +0.60609  -4.44906 u -19.79311   u^2 -0.17997 u_x -0.86156 uu_x  $ \\
\bottomrule
\end{tabular}
\end{center}
\footnotesize{For RGG \cite{reinbold2020using}, we use 8 default features $\{ uu_x, u_{xx}, u_{xxxx}, u, u_{x}, u_{xxx}, u^2, u^3 \}$ and the parameters $p_x = 4, p_t = 3, N_d = 100, D = (40,20)$ are used.  For IDENT \cite{kang2021ident}, we use $\lambda = 200$ for the sparse regression algorithm, and the dictionary is set to be $\{ 1, u, u^2, u_x, u_x^2, uu_x, u_{xx}, u_{xx}^2, uu_{xx}, u_xu_{xx}\}$. SC and ST \cite{he2020robust} use  the same dictionary as IDENT. For SC, we use $\alpha = 100$ and for ST, we use $s = 20$ and $n = 5$. These parameters are from the original papers.  }
\caption{Transport equation with diffusion \eqref{e: pde transport}: clean data case in (a), (b) and (c), and noisy data with $\sigma_{\rm NSR}=100\%$ in (d), (e) and (f).  WeakIdent is compared with WPDE \cite{messenger2021weakPDE}, RGG  \cite{reinbold2020using}, IDENT\cite{kang2021ident}, SC\cite{he2020robust}, and  ST\cite{he2020robust}. The error measures are in Table (b) and (e) and the recovered equations are in (c) and (f).}
\label{f: transport equation}
\end{figure}

The first set of results in Figure \ref{f: transport equation} shows results for the transport equation \eqref{e: pde transport} with clean and noisy data.  (a), (b) and (c) compare the recovery results  with clean data, and (d), (e) and (f) compare the results with highly corrupted data where $\sigma_{\rm NSR} = 100\%$. 
For the case of clean data, RGG \cite{reinbold2020using},  WPDE \cite{ messenger2021weakPDE} and the proposed WeakIdent find the correct support $u_x, u_{xx}$, while the latter two methods have higher accuracy.  In the noisy case of $\sigma_{\rm NSR} = 100\%$, only WeakIdent is able to identify the correct support with the $E_2$ value as low as 0.008.

\begin{figure}[t!]
    \begin{center}
    \begin{tabular}{cccc}
    \multicolumn{4}{c}{\textbf{WeakIdent}} \\
    (a1) $E_2$ & (a2) $E_\infty$  & (a3) TPR & (a4) PPV \\
    \includegraphics[width = 0.20\textwidth]{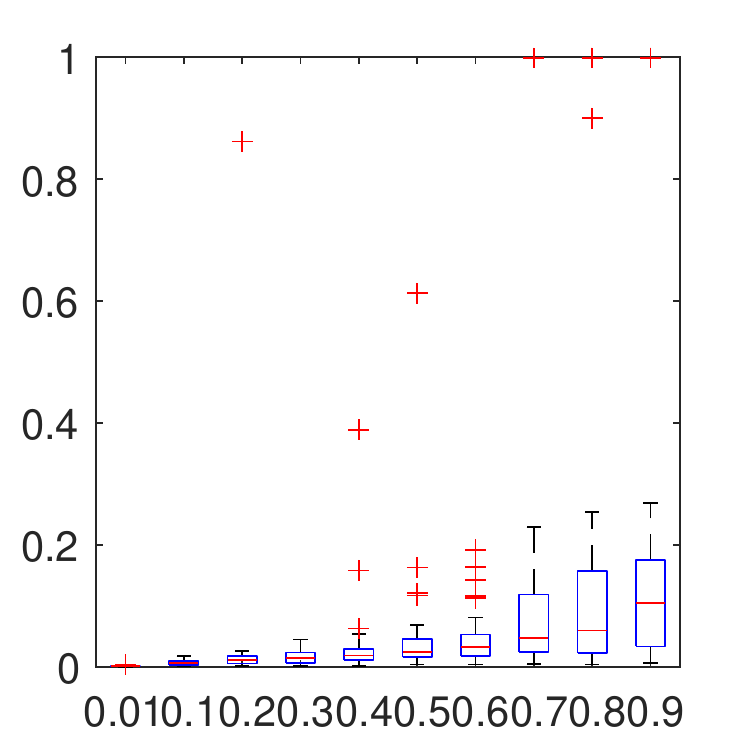}
    &\includegraphics[width = 0.20\textwidth]{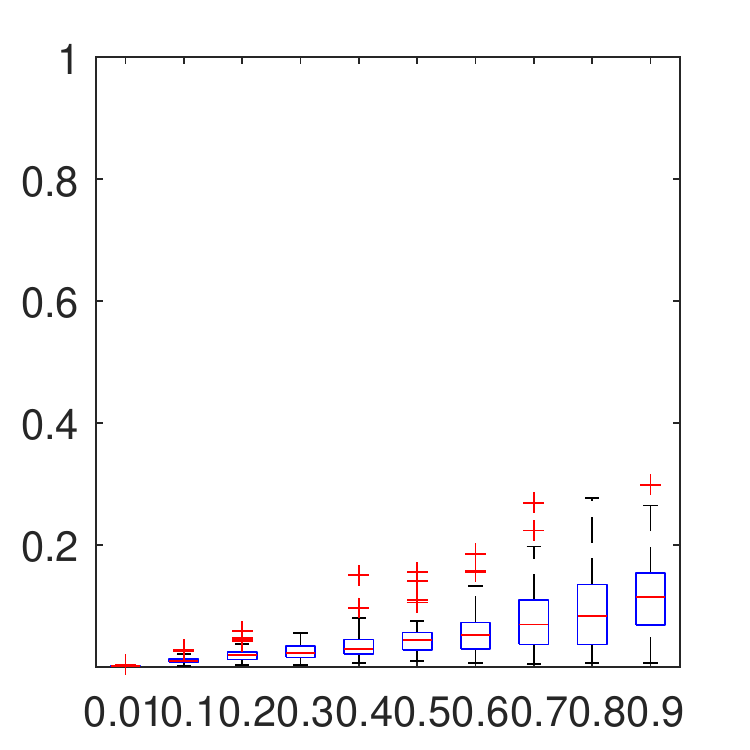}& 
    \includegraphics[width = 0.20\textwidth]{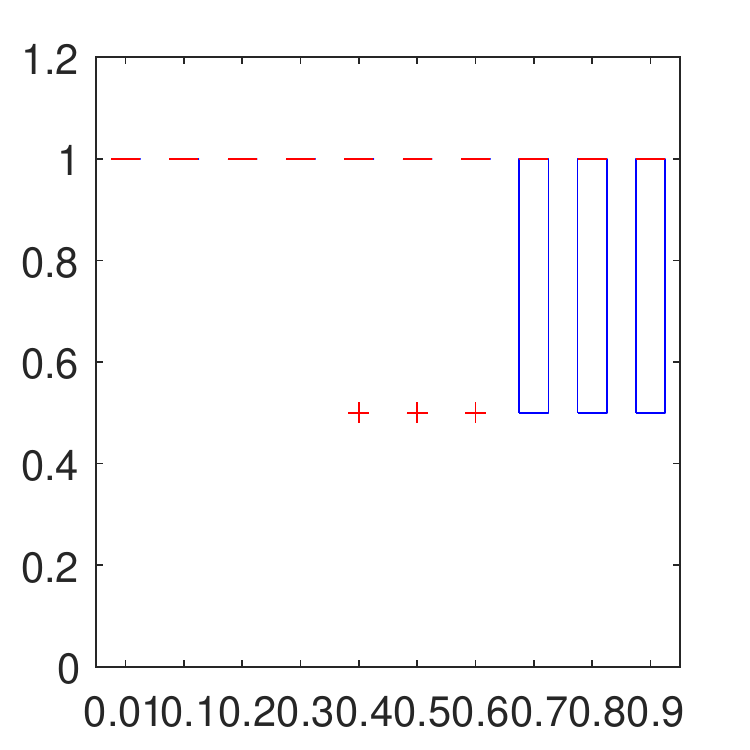} & 
    \includegraphics[width = 0.20\textwidth]{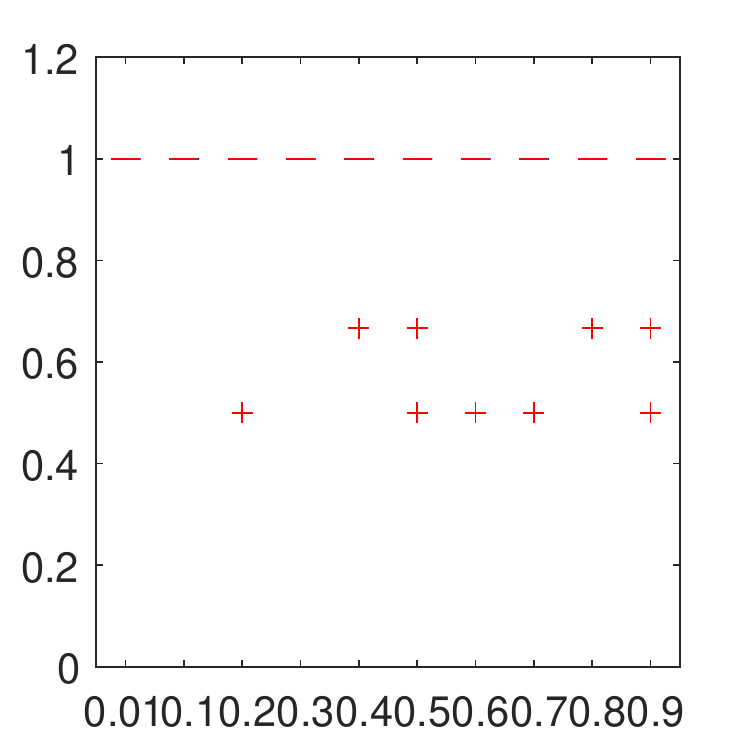}
    \\
    \multicolumn{4}{c}{WPDE \cite{messenger2021weakPDE}} \\
    (b1) $E_2$ & (b2) $E_\infty$ & (b3) TPR & (b4) PPV \\
    \includegraphics[width = 0.20\textwidth]{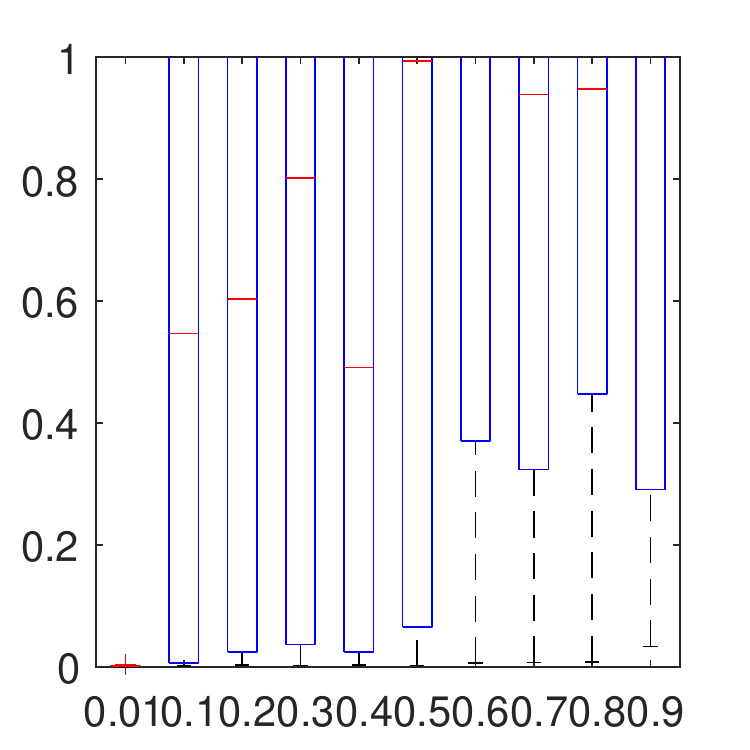} 
    & 
    \includegraphics[width = 0.20\textwidth]{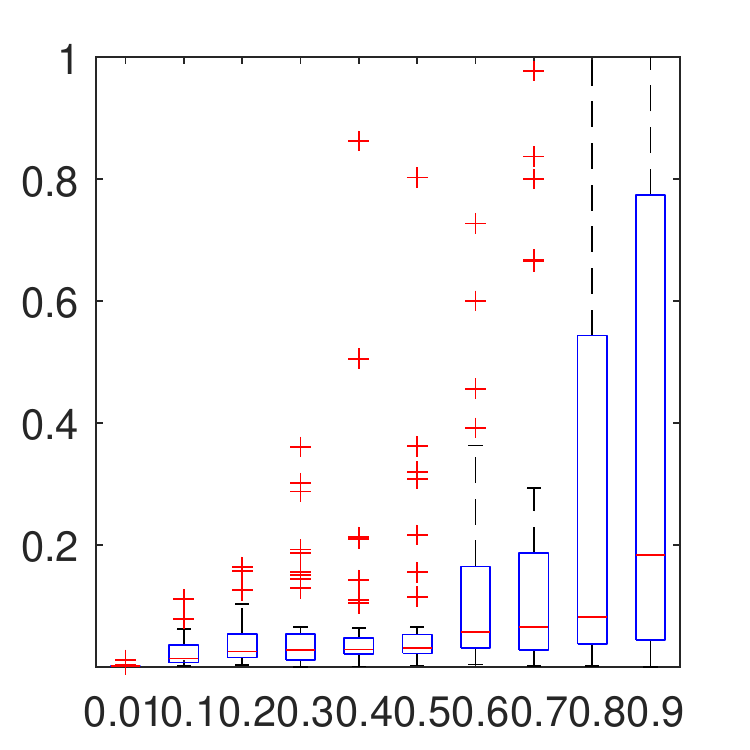}& 
    \includegraphics[width = 0.20\textwidth]{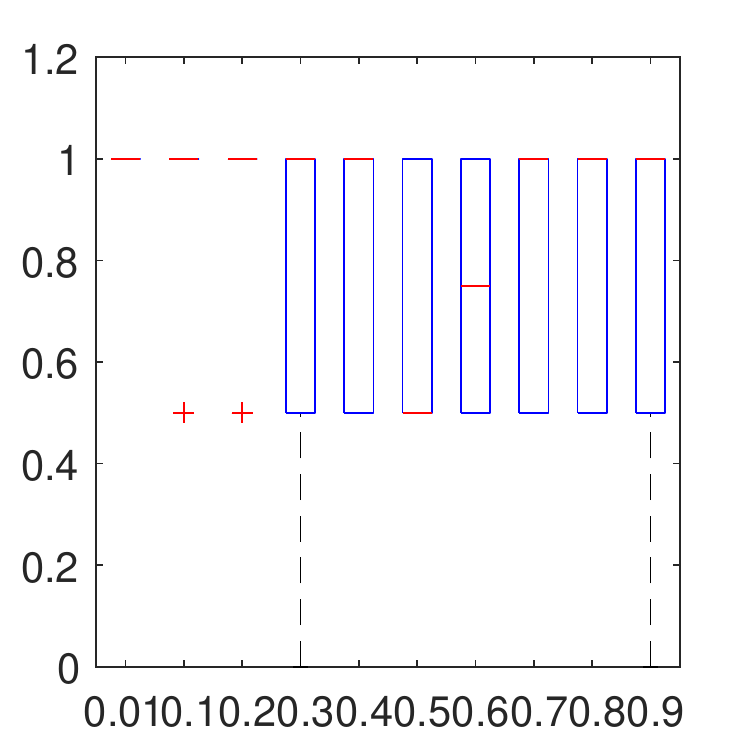} & 
    \includegraphics[width = 0.20\textwidth]{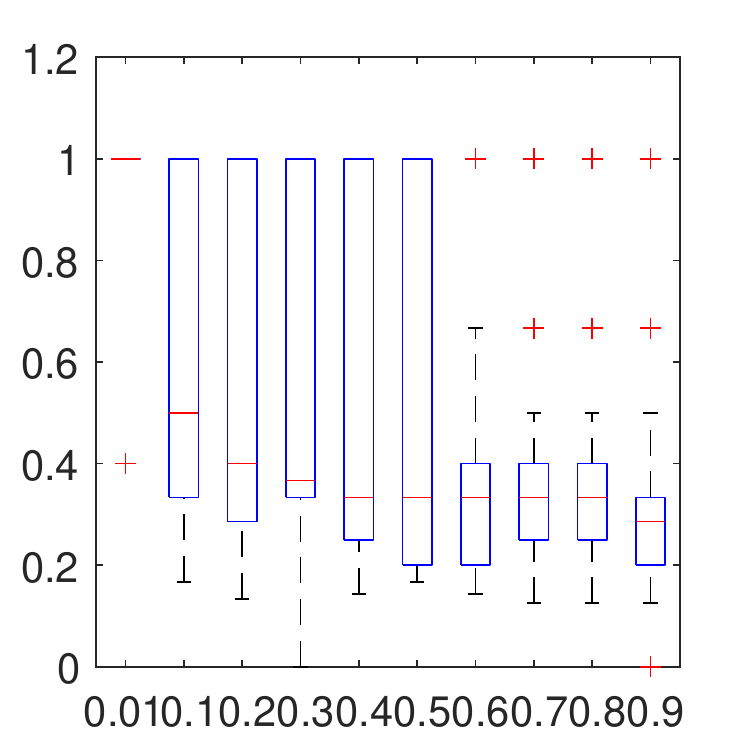}
    \end{tabular}
    \end{center}
    \footnotesize{In each box-plot, the red line is the median, the lower bound is the $25\%$ quantile, the upper bound is the $75\%$ quantile, and $+$ signs represent outliers of each identification error. We use the same criteria for the box-plots in the rest of the figures.}
    \caption{Transport equation (\ref{e: pde transport}), statistical comparison between WeakIdent (the top row) and  WPDE \cite{messenger2021weakPDE} (the second row).
    The errors $E_2, E_{\infty}$, TPR and PPV are shown from 50 experiments for each  $\sigma_{\rm NSR}\in \{ 0.01,0.1, 0.2, ,...,0.9\}$  using box-plots.  The  $E_2$ and $E_\infty$ errors by WeakIent are lower than the errors of WPDE, with less variations. The TPR and PPV by WeakIdent are closer to 1 with less variations as well.}
    \label{fig: box plot transport equation}
\end{figure}

In Figure \ref{fig: box plot transport equation}, we provide statistical comparisons between our proposed  WeakIdent and  WPDE \cite{messenger2021weakPDE} applied to the transport equation  \eqref{e: pde transport} for different levels of $\sigma_{\rm NSR}$.  
We show box-plots for the distribution of the identification errors $E_2, E_{\infty}$,
TPR and PPV over 50 experiments for each level of $\sigma_{\rm NSR}\in \{ 0.01,0.1, 0.2, ,...,0.9\}$. The WeakIdent results are robust even for large noise levels: Panels (a3) and (a4) show that in the  majority($>75\%$) of the cases, a correct support is found by WeakIdent with low $E_2$ error  in Panel (a1). 

\subsubsection{Anisotropic Porous Medium (PM) equation }
\begin{figure}[]
\begin{center}
\begin{tabular}{cc}
 (a) $\hat{U}(\boldsymbol{x},0)$ & (b) $\hat{U}(\boldsymbol{x},T)$ \\
\includegraphics[width = 0.2 \textwidth]{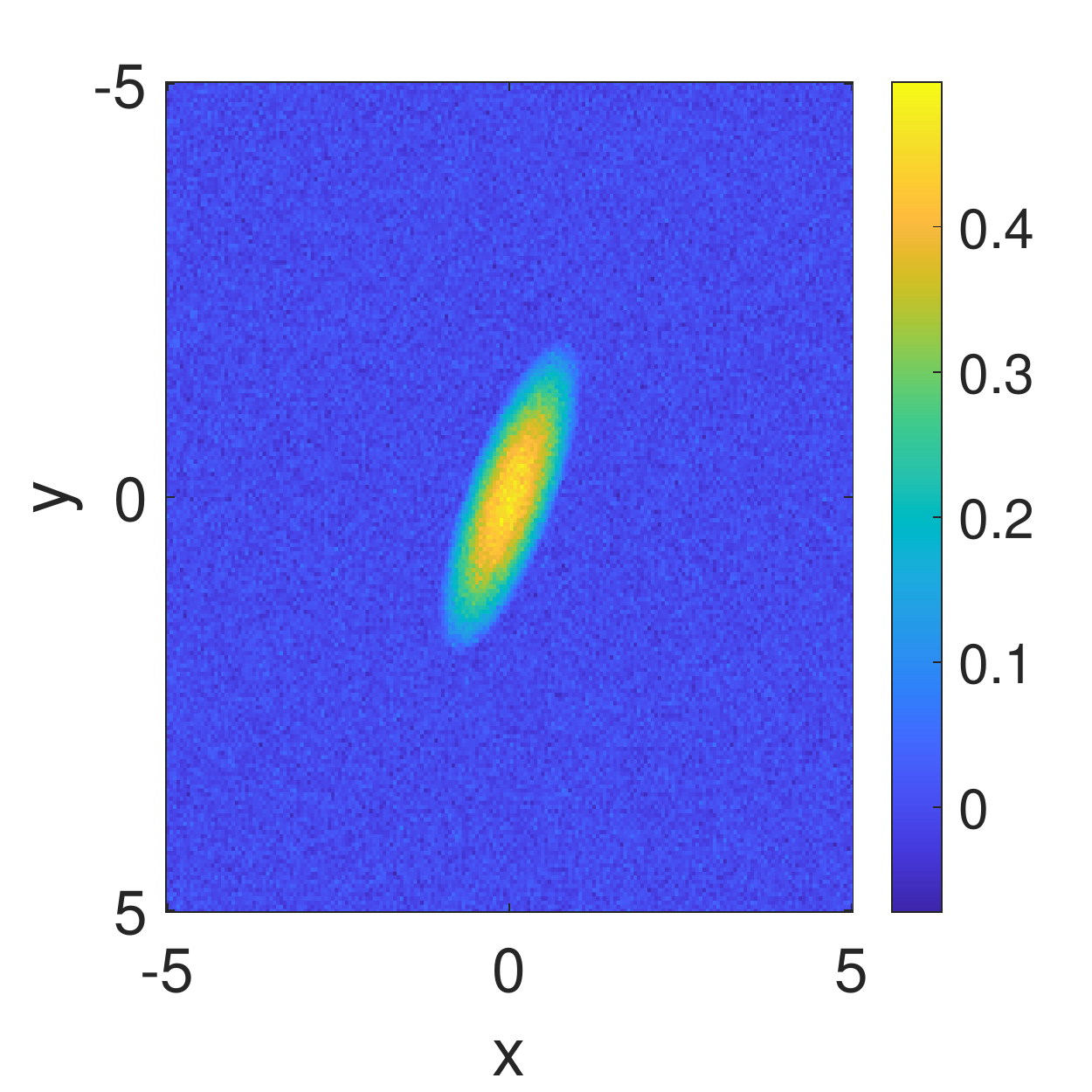} & 
\includegraphics[width = 0.2 \textwidth]{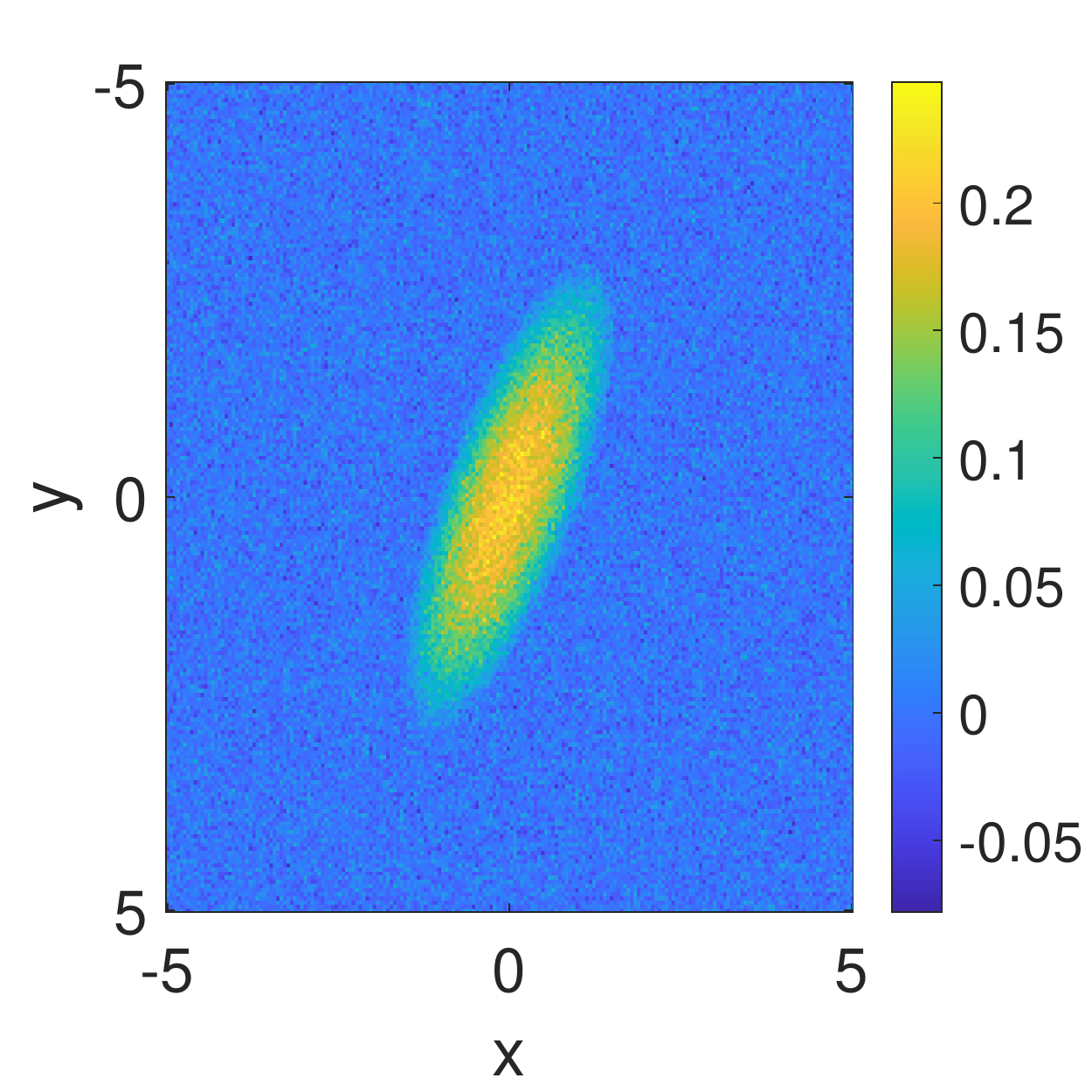}
\end{tabular}
\begin{tabular}{|l|p{10cm}|r|}
\multicolumn{3}{c}{(c) $\sigma_{\rm NSR} =0.08$ } \\
\midrule
True equation & $u_t  =  +0.30000 (u^2)_{yy}
- 0.80000 (u^2)_{xy} + 1.00000 (u^2)_{xx}$ & \\
\hline
\textbf{WeakIdent} & $\bf u_t  =  +0.29912 (u^2)_{yy} -0.79416 (u^2)_{xy} +0.99568 (u^2)_{xx}$
     & \textbf{$E_2= $ 0.0056} \\
     \hline
WPDE  & $\bf u_t  =  +0.29928 (u^2)_{yy} -0.79362 (u^2)_{xy} +0.99512 (u^2)_{xx}$
& $E_2= $ \textbf{0.0061} \\
\hline
RGG\cite{reinbold2020using} & $  u_t  =  +0.00028   +0.23457 u -6.28949 u^2 +21.31151 u^3 -0.36341 (u^2)_x +0.48573 (u^2)_y -0.12914 (u^2)_{xx} +1.34146 (u^2)_{yy} -1.18102 (u^2)_{xy} +0.03796 u_x -0.03006 u_y +0.06174 u_{xx} -0.04775 u_{yy} +0.06238 u_{xy}   $ & $E_2= $ 8.9376   \\
     \bottomrule
\end{tabular}\\
\end{center}
\footnotesize{For RGG \cite{reinbold2020using}, we use a dictionary of 14 features $\{1, u, u^2, u^3, (u^2)_x, (u^2)_y, (u^2)_{xx}, 
    (u^2)_{yy}, (u^2)_{xy}, u_x, u_y,u_{xx}, u_{yy}, u_{xy} \}$ adding the true features, and the parameters $p_x = 2, p_t = 1, N_d = 100$, and $D = (20,10)$.}
\caption{Anisotropic Porous Medium (PM) equation \eqref{e: pde PM} on a 2-D spatial domain with cross derivative feature. We set  $\sigma_{\rm NSR}=0.08$, which is equivalent to  $\sigma_{\rm NR}=0.4139$ in WPDE \cite{messenger2021weakPDE}. 
(a) Given noisy data $\hat{U}(\boldsymbol{x},0)$ and (b) $\hat{U}(\boldsymbol{x},T)$. (c) Identified equations with the $E_2$ error. }
\label{f: recovered equations - PM equation}
\end{figure}

In  Figure \ref{f: recovered equations - PM equation}, we compare the recovery results for the 2D anisotropic porous medium equation (PM) \eqref{e: pde PM}, which includes a feature with the cross-dimensional derivative $u_{xy}$.  
Figure \ref{f: recovered equations - PM equation} (a) shows $\hat{U}(\boldsymbol{x},0)$ and (b) shows $\hat{U}(\boldsymbol{x},T)$, where the  given noisy data has noise-to-signal ratio $\sigma_{\rm NSR}=0.08$.  This noise level is equivalent to $\sigma_{\rm NR}=0.4139$ as defined in WPDE \cite{messenger2021weakPDE}. We show different  recovered equations with the identification error $E_2$ in (c). 
WeakIdent is able to identify the correct support with the  coefficient error $E_2= 0.0056$, demonstrating WeakIdent's capability to identify features across multiple dimensions on 2D spatial domain.

\begin{figure}[]
\begin{center}
\begin{tabular}{cc}
 (a) $\hat{U}(\boldsymbol{x},0)$ & (b) $\hat{U}(\boldsymbol{x},T)$ \\
\includegraphics[width = 0.2  \textwidth]{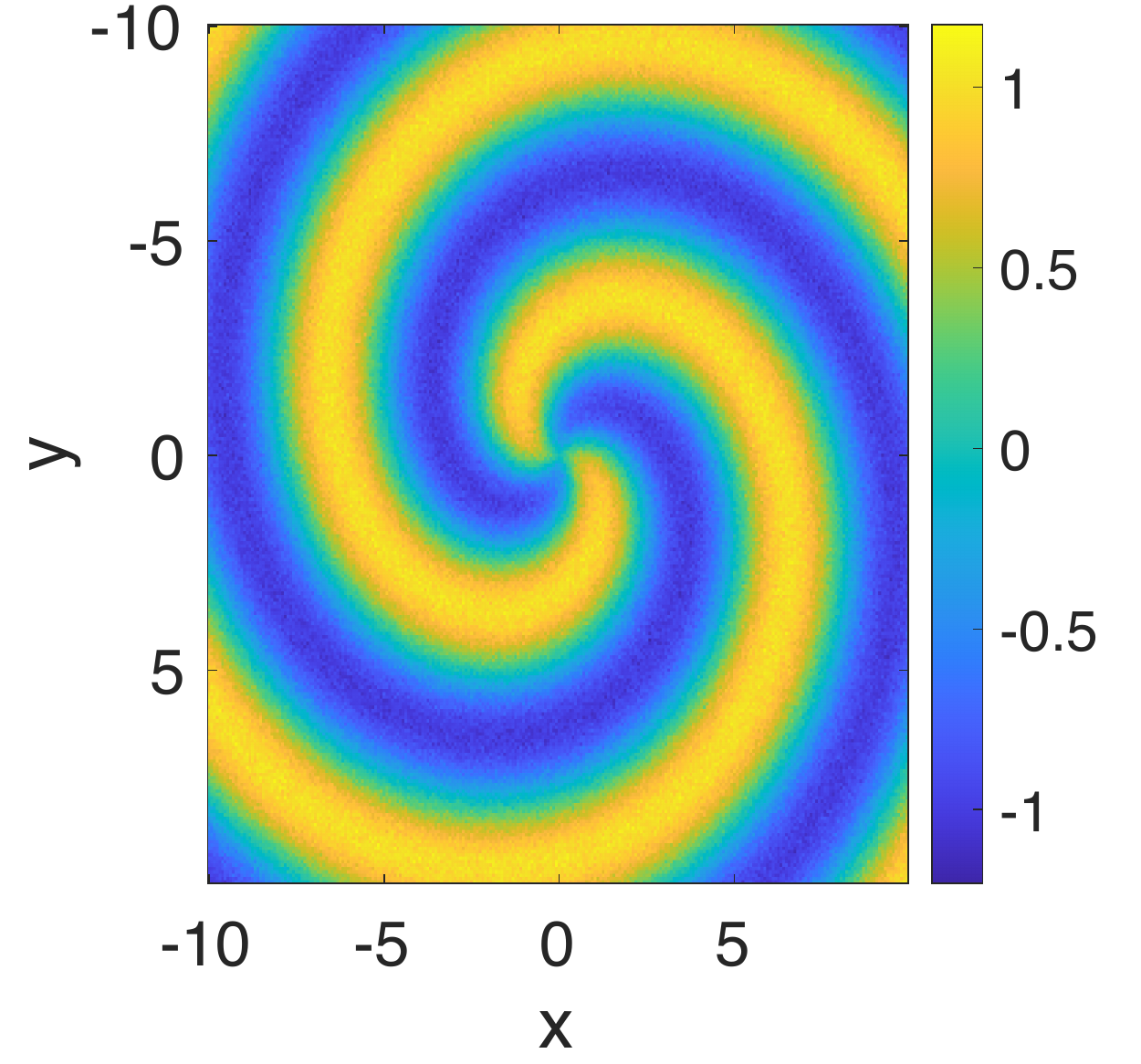} &
\includegraphics[width = 0.2  \textwidth]{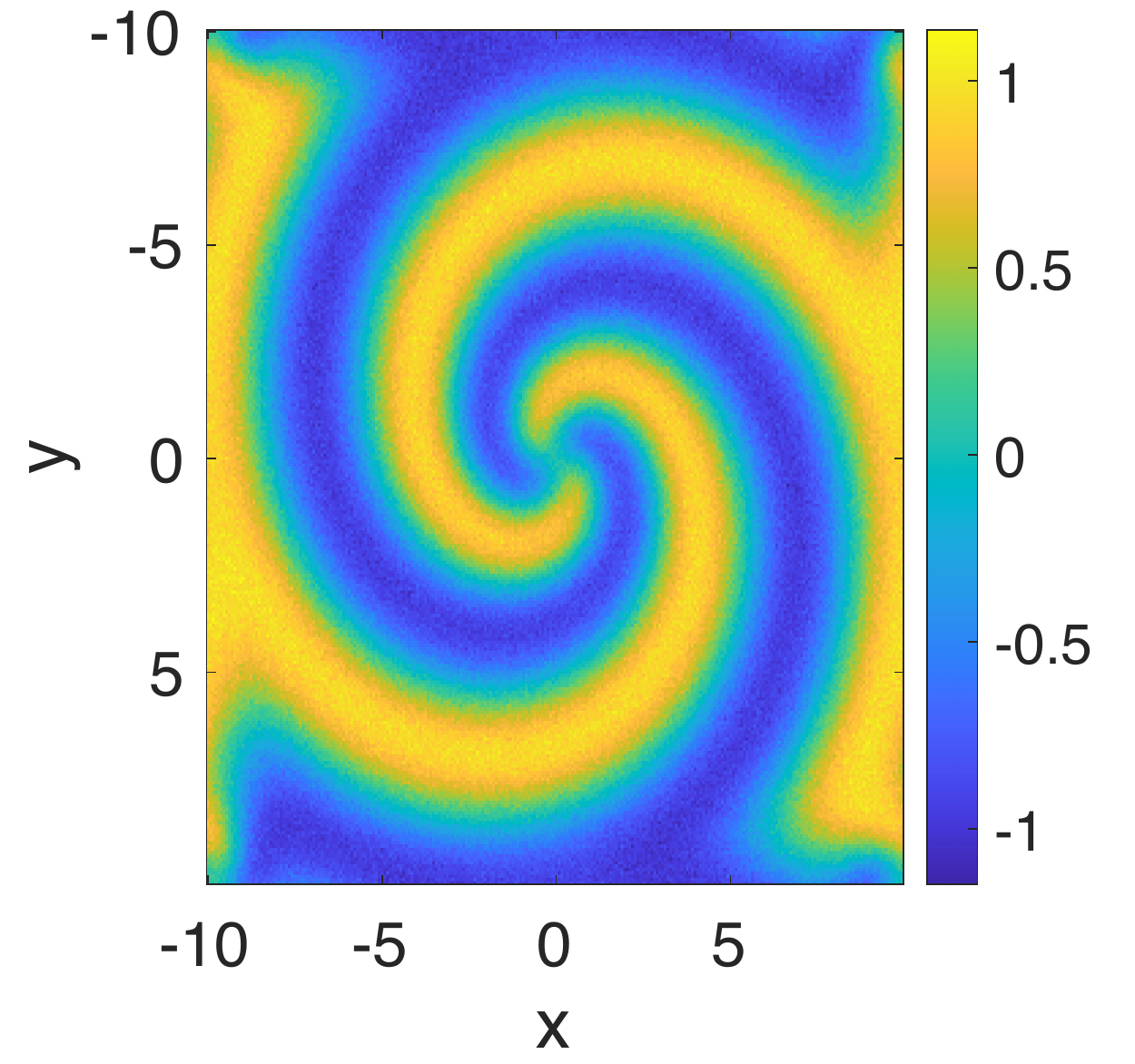}
\end{tabular}
\begin{tabular}{|l|p{10cm}|r|}
\multicolumn{3}{c}{(c) $\sigma_{\rm NSR}=0.08$} \\
\toprule
True equation & $u_t  =  + v^3 + u +0.1 u_{yy} +0.1u_{xx} - uv^2 +u^2v - u^3$ & \\
& $v_t  =   v +0.1v_{yy} +0.1v_{xx} - v^3 -uv^2 - u^2v -u^3$ & \\
\hline
\textbf{WeakIdent} & $u_t  =  +0.99213 v^3 +0.98572 u +0.09660 u_{yy} +0.09695 u_{xx} -0.93229 uv^2 +0.97678 u^2v -0.99018 u^3$ & \textbf{$E_2 = $0.0316}\\
& 
 $v_t  =  +0.97792 v +0.09662 v_{yy} +0.09636 v_{xx} -0.97161 v^3 -0.96468 uv^2 -0.95572 u^2v -0.99605 u^3$
    &  \\
     \hline
WPDE  & $u_t  =  +1.34525 v^3$
& $E_2 = $0.9081\\
&  $v_t  =  -1.34499 u^3$ & \\
\hline
RGG\cite{reinbold2020using} & $ u_t  =  +0.10204 \nabla u +1.02296 u -1.01966 u^3 +1.01341 v^3 +1.03003 u^2v -1.01767 uv^2  $ & $E_2 = $0.0793\\
& $ v_t  =  +0.09244 \nabla v -0.07400 u -0.93640 u^3 +0.95099 v -0.95370 v^3 -0.95450 u^2v -0.93750 uv^2  $ & \\
     \bottomrule
\end{tabular}\\
\end{center}
\footnotesize{For RGG \cite{reinbold2020using}, the provided default features for reaction-diffusion type equation in \cite{reinbold2020using} is used: for $u$, the dictionary is $\{\nabla u, u, u^2, u^3, v,v^2, v^3, uv, u^2v, uv^2\}$ and for $v$, the dictionary is $\{\nabla v, u, u^2, u^3, v,v^2, v^3, uv, u^2v, uv^2\}$,  and parameters $p_x = 2, p_t = 1, N_d = 100, D = (20,10)$.}
\caption{Reaction-diffusion equation  \eqref{e: pde RD} on a 2D spatial domain with $\sigma_{\rm NSR}=0.08$ (equivalent to $\sigma_{\rm NR} = 0.08$  defined in \cite{messenger2021weakPDE}). 
(a) Given noisy data $\hat{U}(\boldsymbol{x},0)$ and (b) $\hat{U}(\boldsymbol{x},T)$. (c) The identified equations and the $E_2$ errors. WeakIdent finds the  correct terms with a small coefficient error.}
\label{f: recovered equations - RD equation}
\end{figure}

\begin{figure}[]
    \centering
    \begin{tabular}{cccc}
    (a1) $E_2$ & (a2) $E_\infty$  & (a3) TPR & (a4) PPV \\
    \includegraphics[width = 0.20\textwidth]{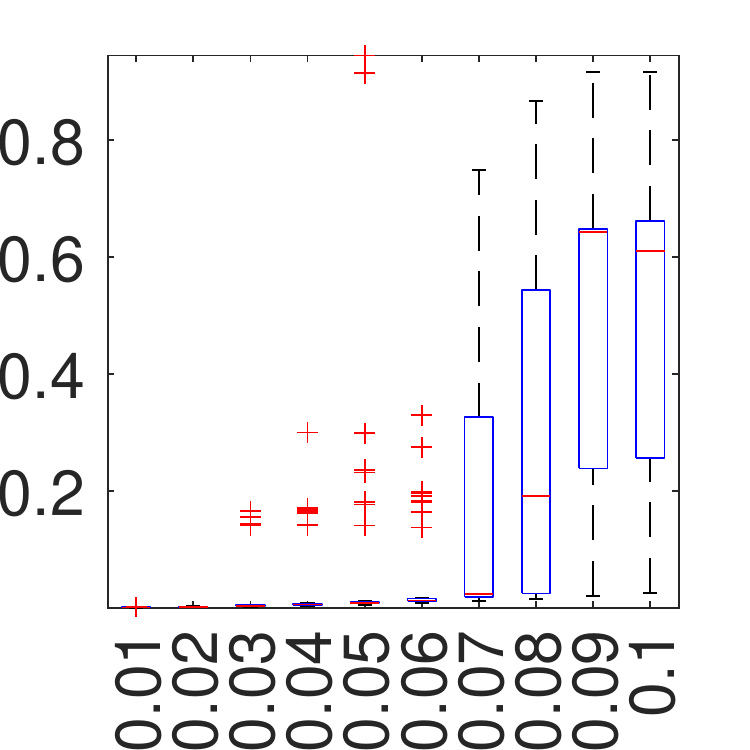}
    &\includegraphics[width = 0.20\textwidth]{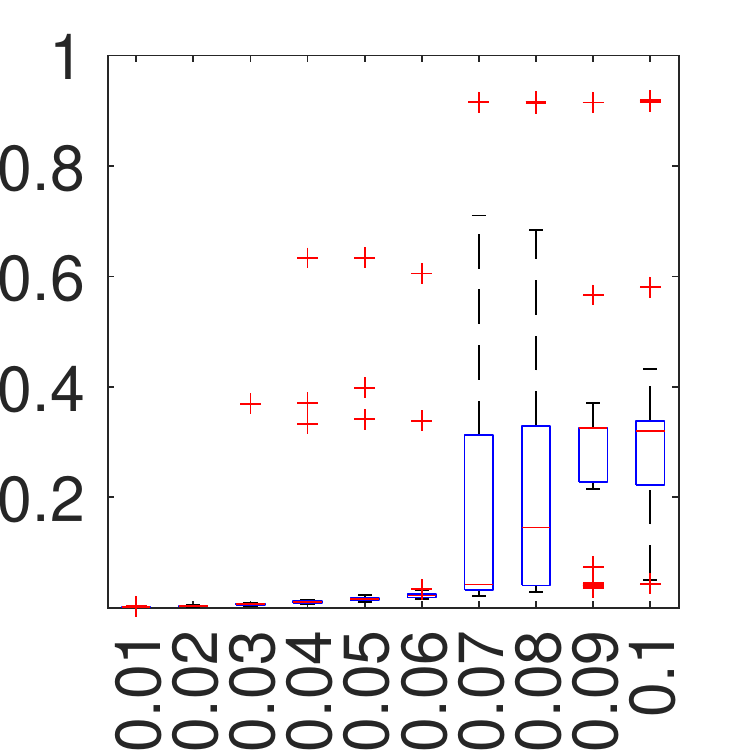}& 
    \includegraphics[width = 0.20\textwidth]{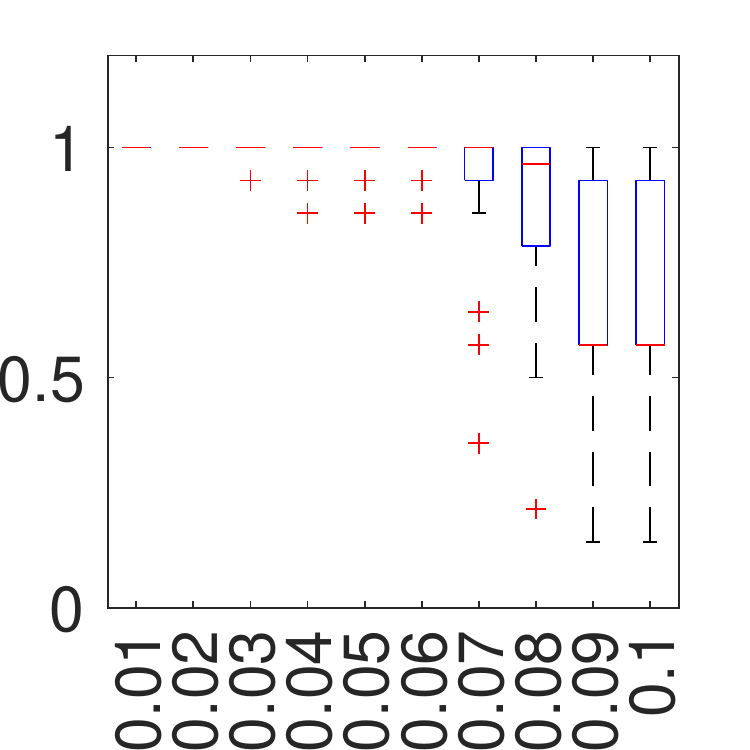} & 
    \includegraphics[width = 0.20\textwidth]{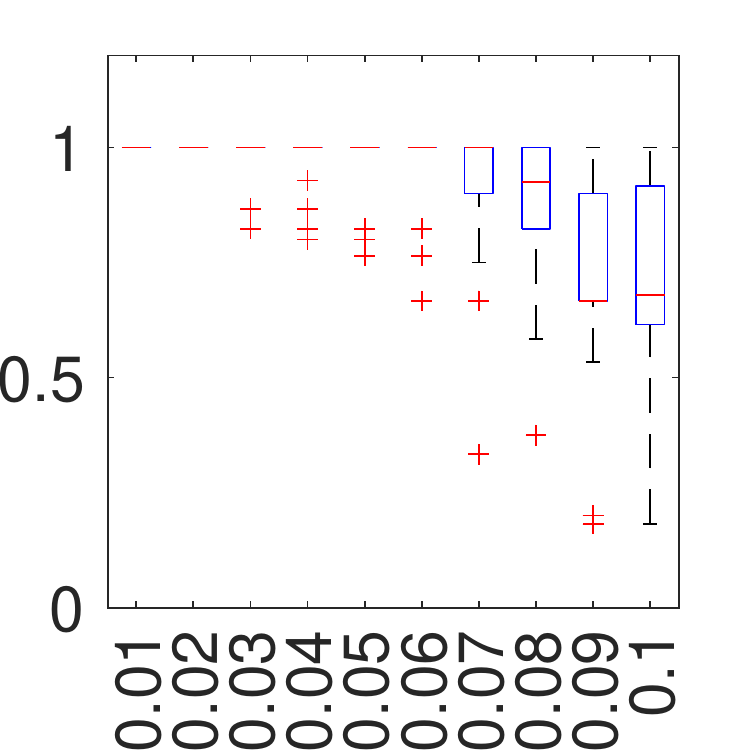}
    \end{tabular}
    \caption{
    The Identification results from WeakIdent for the reaction diffusion  equation (\ref{e: pde RD}):
    The  $E_2, E_{\infty}$ errors , TPR and PPV are shown from 50 experiments for each  $\sigma_{\rm NSR}\in \{ 0.01,0.02,  ,...,0.1\}$  using box-plots.  
    }
    \label{fig: box plot reaction diffusion equation}
\end{figure}

\subsubsection{Reaction-diffusion equation }
In  Figure \ref{f: recovered equations - RD equation}, we compare the recovery results for the 2D reaction-diffusion equation \eqref{e: pde RD}. These systems can generate a variety of patterns such as dots, strips, waves and hexagons. 
The Laplacian (diffusion) features $\Delta u, \Delta v$ in this equation may be difficult to identify in general, particularly in the case where  the diffusion coefficients are small compared to those of other features, and accumulated noise can be emphasized.
We use the spiral pattern data set from \cite{messenger2021weakPDE}.  
Figure \ref{f: recovered equations - RD equation} (a) shows $\hat{U}(\boldsymbol{x},0)$ and (b) shows $\hat{U}(\boldsymbol{x},T)$, where the  given noisy data has  $\sigma_{\rm NSR}=0.08$ (equivalent to $\sigma_{\rm NR} = 0.08$  defined in \cite{messenger2021weakPDE}).  We show different  recovered equations with the $E_2$ identification error  in Figure \ref{f: recovered equations - RD equation} (c). 
WeakIdent finds the  correct terms with a small coefficient error.

In Figure \ref{fig: box plot reaction diffusion equation}, we present the statistical results of WeakIdent over 50 experiments for the 2D reaction-diffusion equation \eqref{e: pde RD}. 

\subsubsection{PDEs and sytems of PDEs with higher order features }

\begin{figure}
    \centering
    \begin{tabular}{|p{1.6cm}|c|c|c|}
    \toprule 
    & $E_2$ & TPR & PPV \\
    \midrule
     Transport equation \eqref{e: pde transport} \; $\tilde{\sigma} =1.19$ & \includegraphics[width = 0.27\textwidth, align = c]{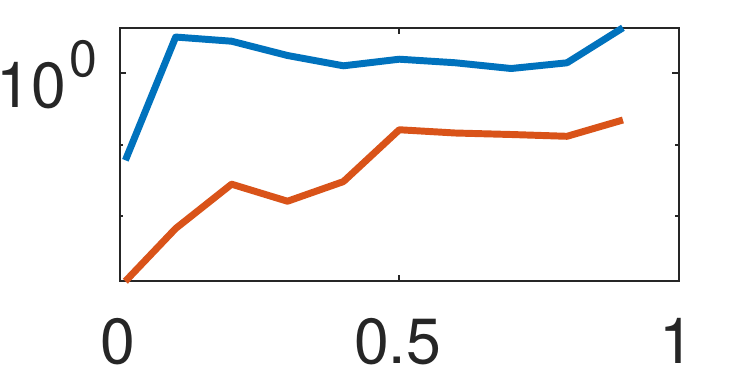} & \includegraphics[width = 0.27\textwidth,  align =c]{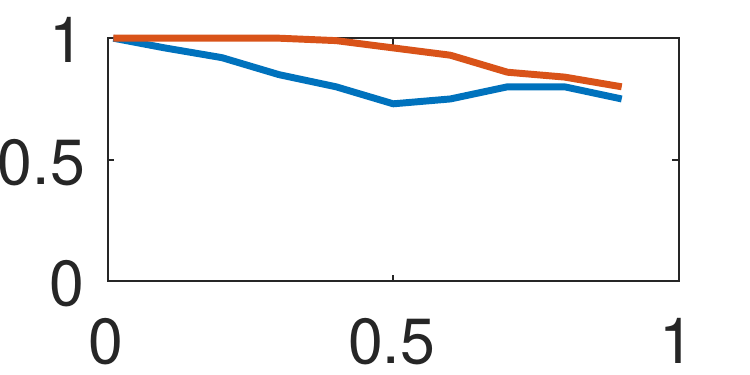} &
     \includegraphics[width = 0.27\textwidth,  align =c]{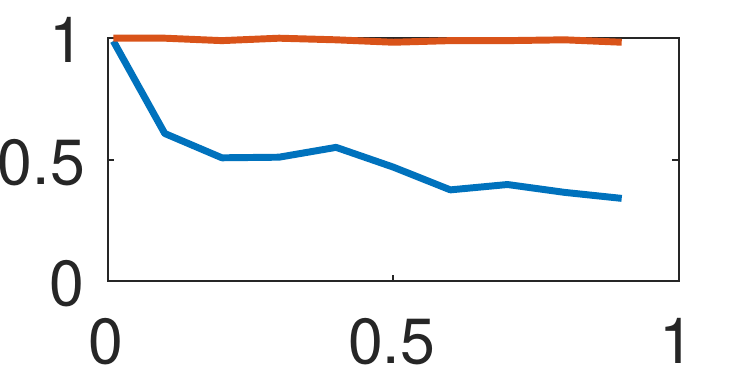}   \\
     \midrule
     KdV \eqref{e: pde kdv} \; $\tilde{\sigma} =3.27$ & \includegraphics[width = 0.27\textwidth,  align =c]{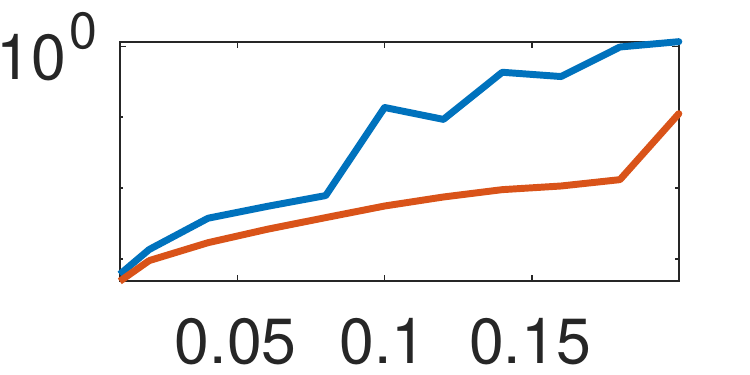} & \includegraphics[width = 0.27\textwidth,  align =c]{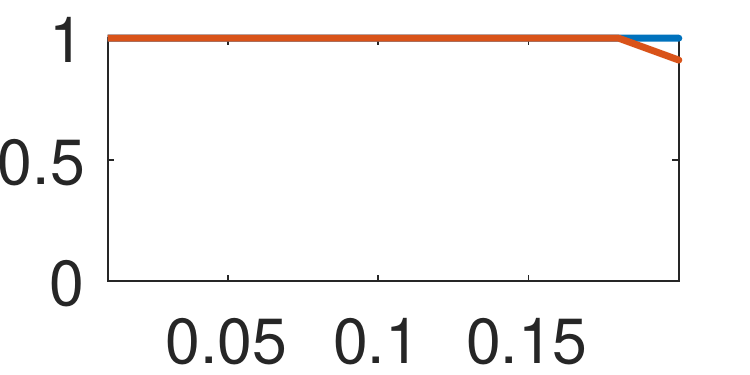} & \includegraphics[width = 0.27\textwidth,  align =c]{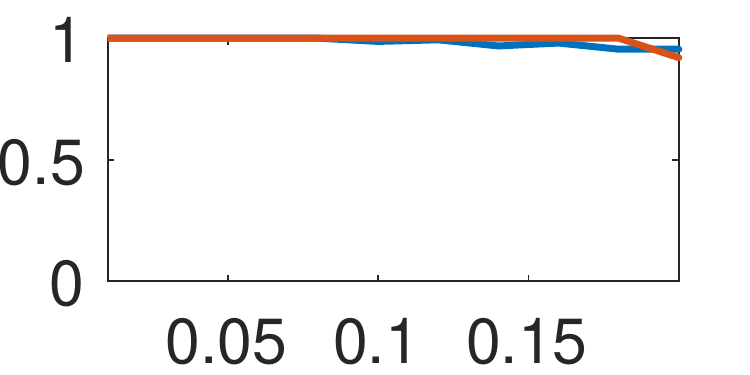}  \\
     \midrule
     KS \eqref{e: pde KS} \; $\tilde{\sigma}=1$  & \includegraphics[width = 0.27\textwidth,  align =c]{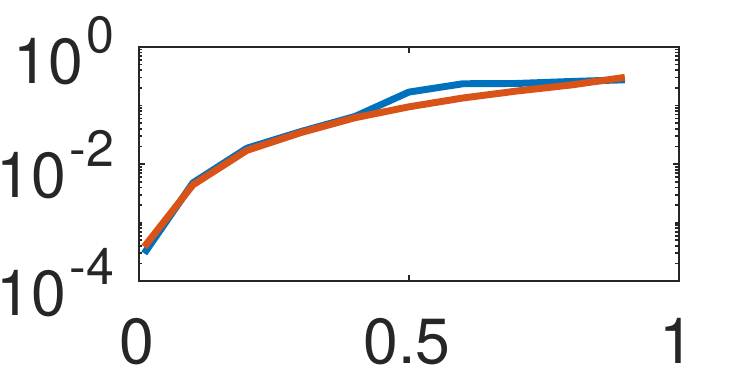} 
     & \includegraphics[width = 0.27\textwidth,  align =c]{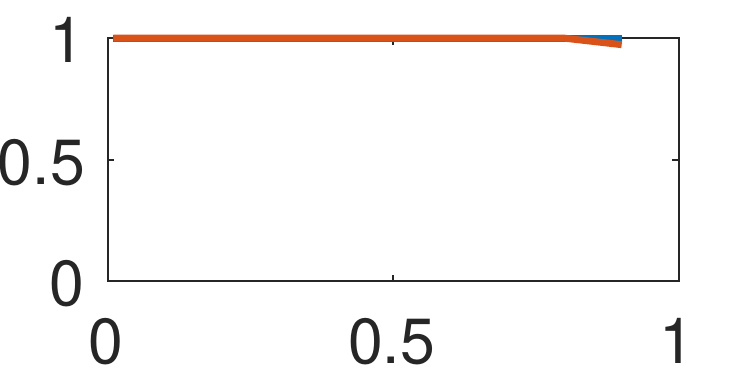} 
     & \includegraphics[width = 0.27\textwidth,  align =c]{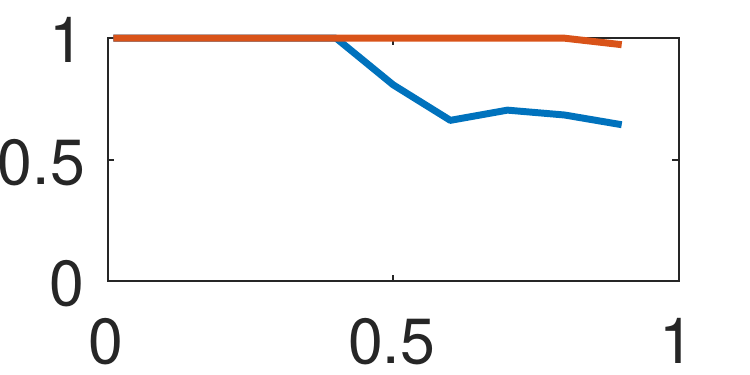} \\
     \midrule
      NLS \eqref{e: pde NLS}\;\;\; \; $\tilde{\sigma}=(1.62,1.62)$ & \includegraphics[width =0.27\textwidth,  align =c]{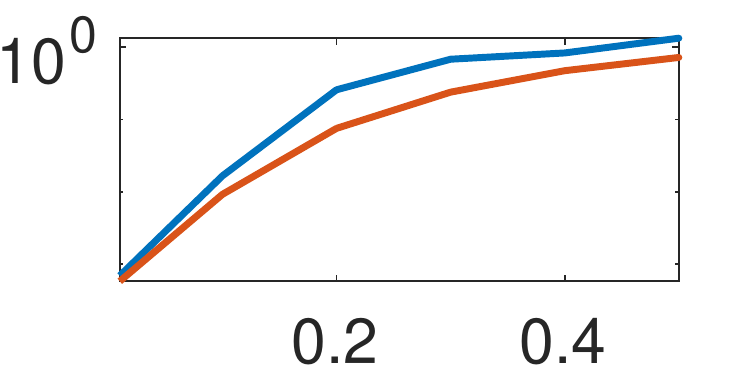} 
     & \includegraphics[width =0.27\textwidth,  align =c]{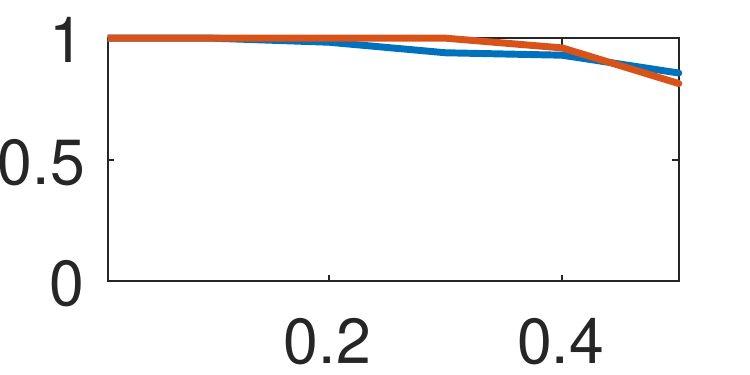} 
     & \includegraphics[width =0.27\textwidth,  align =c]{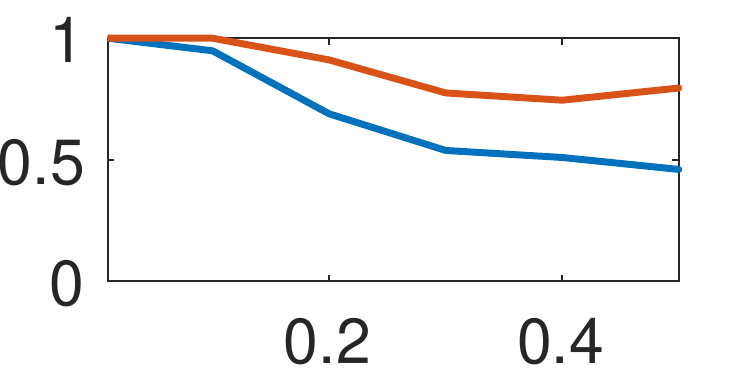} \\
     \midrule
     PM \eqref{e: pde PM}\;\; $\tilde{\sigma}=5.17$ & \includegraphics[width =0.27\textwidth,  align =c]{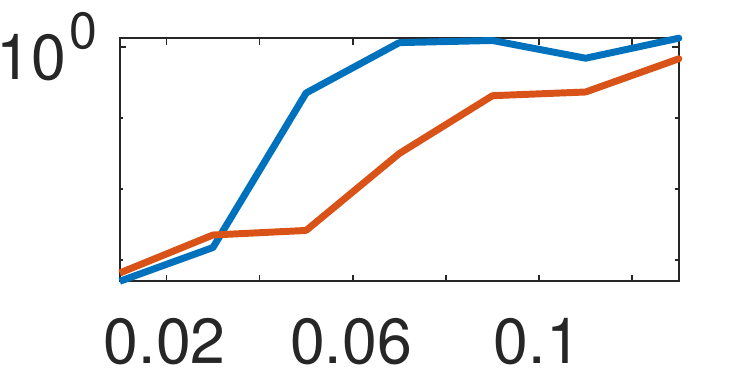} & \includegraphics[width =0.27\textwidth,  align =c]{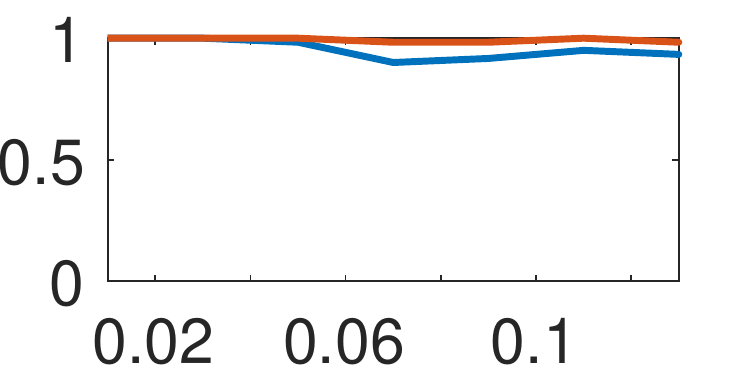} 
     & \includegraphics[width =0.27\textwidth,  align =c]{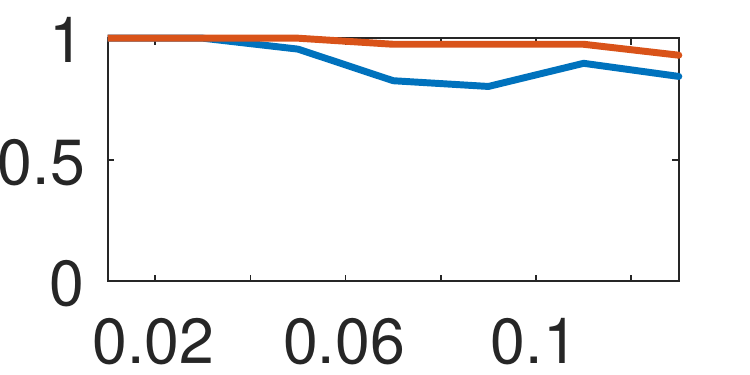} \\
     \midrule
     Reaction-Diffusion \eqref{e: pde RD}\;\; $\tilde{\sigma}=(1,1)$  & \includegraphics[width =0.27\textwidth,  align =c]{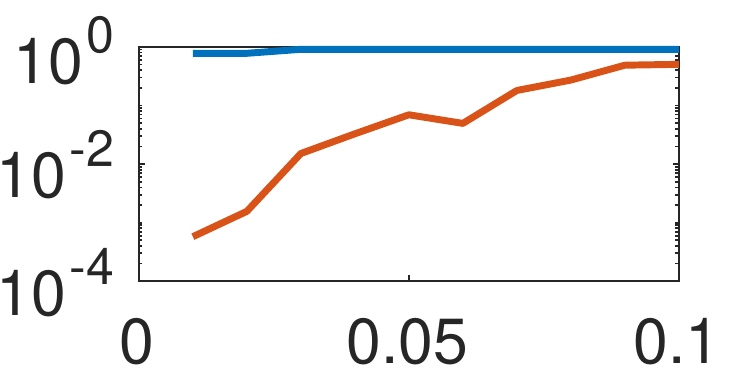} 
     & \includegraphics[width =0.27\textwidth,  align =c]{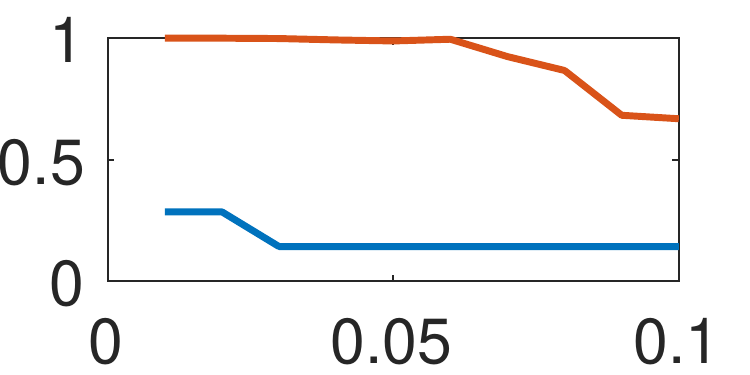} 
     & \includegraphics[width =0.27\textwidth,  align =c]{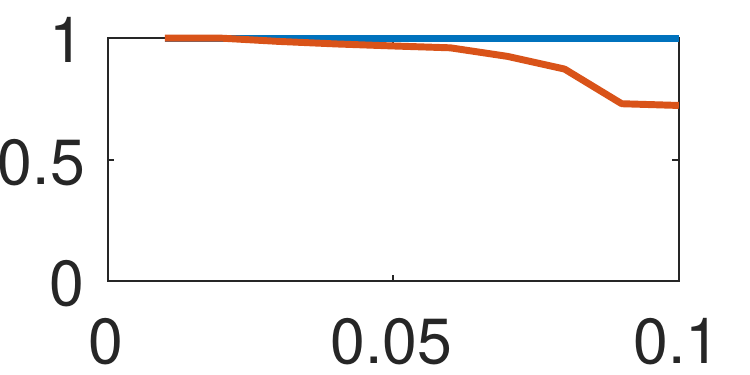} \\
     \bottomrule
    \end{tabular}
    \caption{ The identified PDEs in Table \ref{T:PDE} for different noise levels. We compare WeakIdent (Red) and WPDE (Blue).  The $x$-axis is $\sigma_{\rm NSR}$, while  the $y$-axis is the average  $E_2$ error, TPR and PPV over 50 experiments.  The relative noise ratio $\tilde{\sigma}=\sigma_{\rm NSR}/\sigma_{\rm NR}$ compares our noise level $\sigma_{\rm NSR}$ vs. $\sigma_{\rm NR}$ in \cite{messenger2021weakPDE}. We present results for  the transport equation \eqref{e: pde transport},  the KdV equation  \eqref{e: pde kdv},   the KS equation   \eqref{e: pde KS}, the NLS equation \eqref{e: pde NLS}, the PM equation \eqref{e: pde PM}, and  the reaction-diffusion (2D) equation \eqref{e: pde RD}. The noise-to-signal ratio $\sigma_{\rm NSR}$ ranges  in $\{0,0.1,0.2,...,0.9\}$, $\{0.01,0.02,0.04,...,0.24 \}$,
$\{0,0.1,0.2,...,0.9\}$, $\{0.01,0.1,0.2,...,0.5  \}$, 
$\{0.01,0.03,0.05,...,0.15\}$, and $\{0.01,0.02,...,0.1\}$ for each equation respectively.
     }
    \label{f: a summary of recovering pde examples v.s wsindy}
\end{figure}

In Figure \ref{f: a summary of recovering pde examples v.s wsindy}, we show the average errors of WeakIdent and WPDE over 50 experiments on the PDEs and systems of PDEs in Table \ref{T:PDE} with different noise levels.
Each column gives the $ E_2$ error , TPR and PPV  respectively. In each row, we present the results from the transport equation \eqref{e: pde transport}, KdV   equation  \eqref{e: pde kdv}, the KS \eqref{e: pde KS}, the nonlinear Schrodinger \eqref{e: pde NLS}, the  anisotropic PM equation \eqref{e: pde PM}, and the 2D reaction-diffusion  equation \eqref{e: pde RD}. 
In the first column, we present the ratio $\tilde{\sigma}=\sigma_{\rm NSR}/\sigma_{\rm NR}$ where $\sigma_{\rm NR}$ denotes the noise ratio in WPDE\cite{messenger2021weakPDE}. (The upper bounds of the noise ratio $\sigma_{\rm NR}$ \cite{messenger2021weakPDE} are 1.07, 0.78,  0.9,  0.81, 0.78, 0.1 for each equation.) 
Here the KdV  \eqref{e: pde kdv} and KS equations \eqref{e: pde KS} include higher order derivative features $u_{xxx}$ and $u_{xxxx}$.  These features are in general difficult to recover, especially from highly corrupted noisy data. 
Each plot gives comparisons between WeakIdent (Red) and WPDE (blue), with $\sigma_{\rm NSR}$ on  the $x$-axis. The $y$-axis is the $E_2$ error, TPR, or PPV averaged  over 50 experiments for a given $\sigma_{\rm NSR}$. 
According to the  $E_2$ error shown in the first column, WeakIdent has smaller $E_2$ errors than other methods, showing that WeakIdent is more accurate in the coefficient recovery. According to the TPR and PPV  in the  second and third column,  WeakIdent is more accurate in support recovery since the TPR and PPV values of WeakIdent are closer  to 1.  

\subsection{WeakIdent results and comparisons for ODEs}

\begin{figure}[h]
    \centering
    \begin{tabular}{ccccc}
    (a) Equ.  \eqref{e: ode, linear 2d} & (b) Equ. \eqref{e: ode, Duffing}   & (c) Equ. \eqref{e: ode, van-der-pol}  & (d) Equ. \eqref{e: ode, Lotka}    & (e) Equ. \eqref{e: ode, lorenz}   \\
     \includegraphics[width = 0.18\textwidth ]{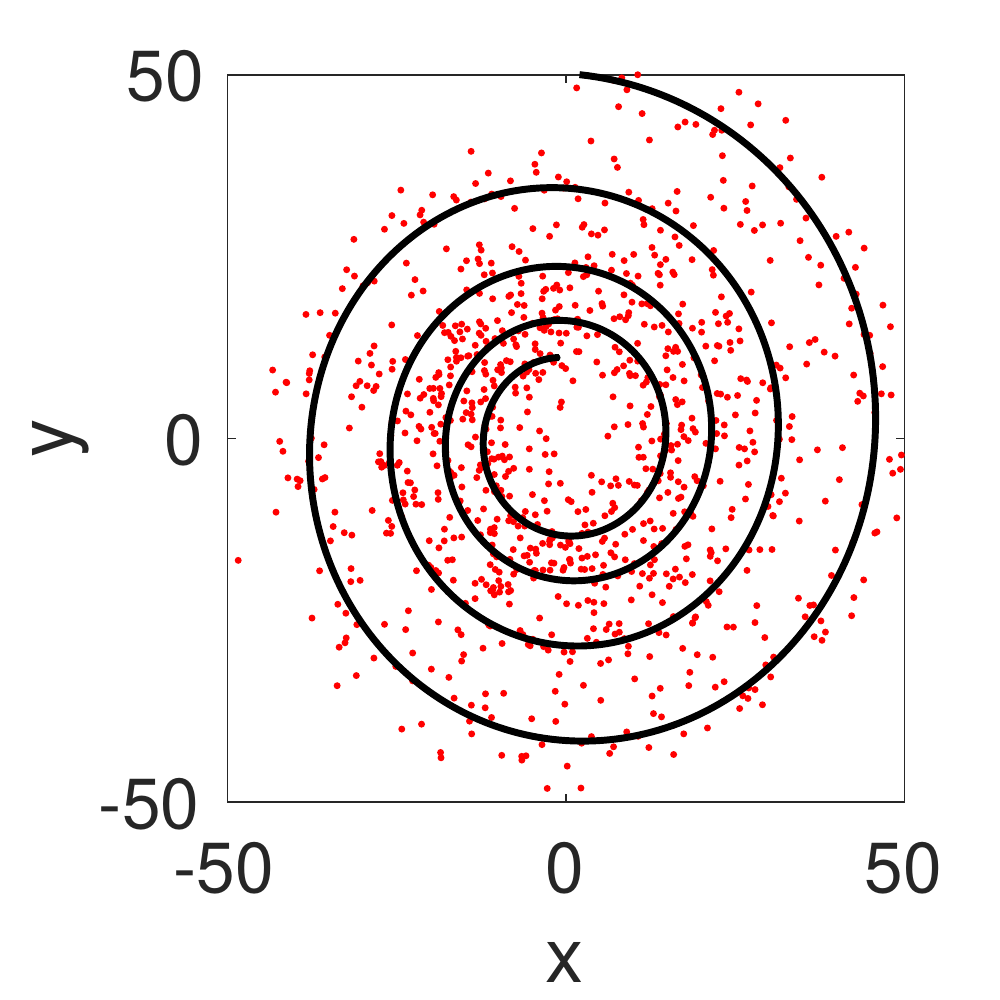} &
     \includegraphics[width = 0.18\textwidth ]{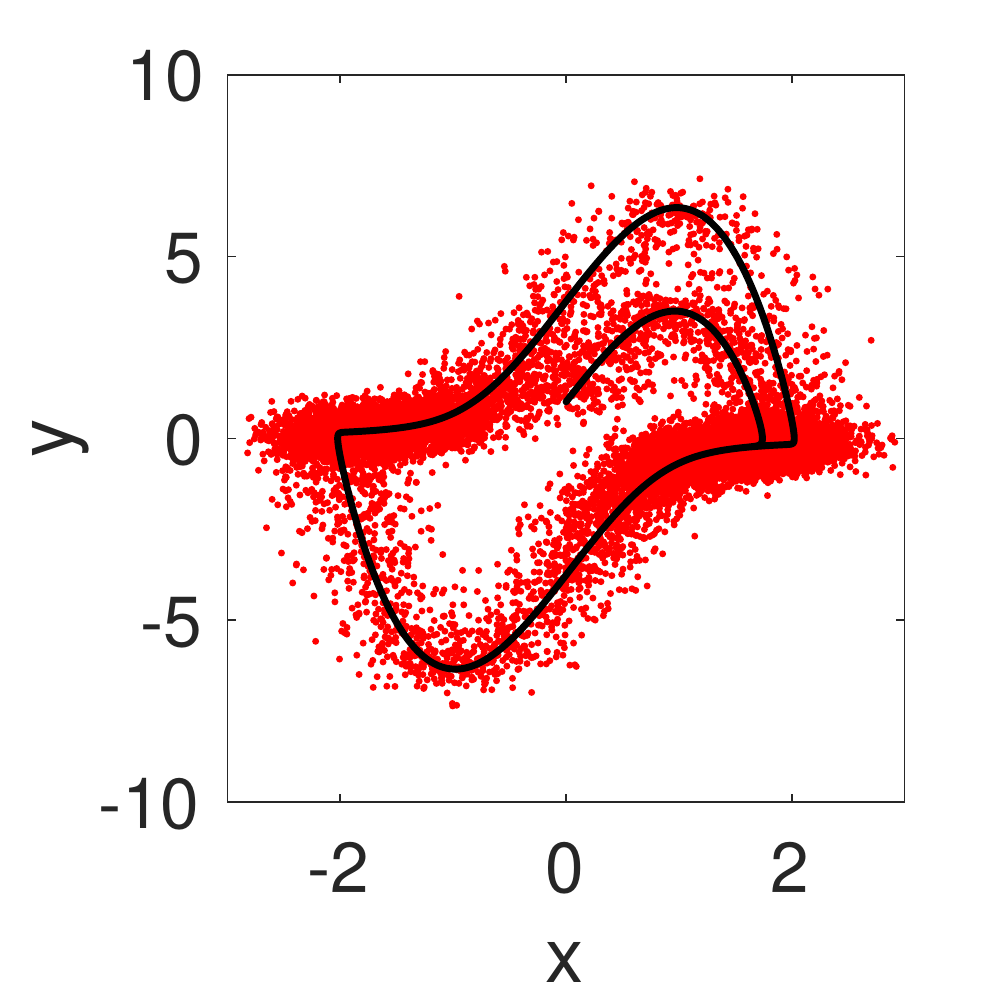} & 
     \includegraphics[width = 0.18\textwidth ]{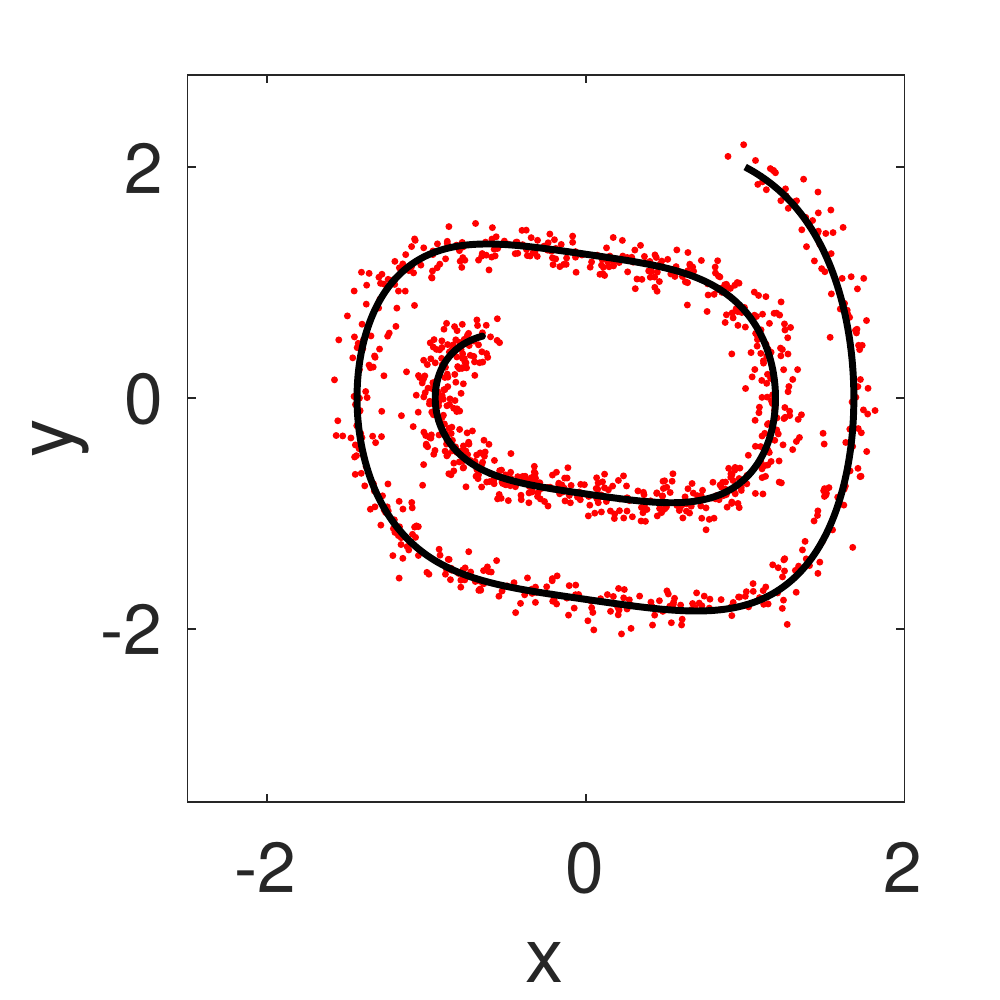}
     &
     \includegraphics[width = 0.18\textwidth ]{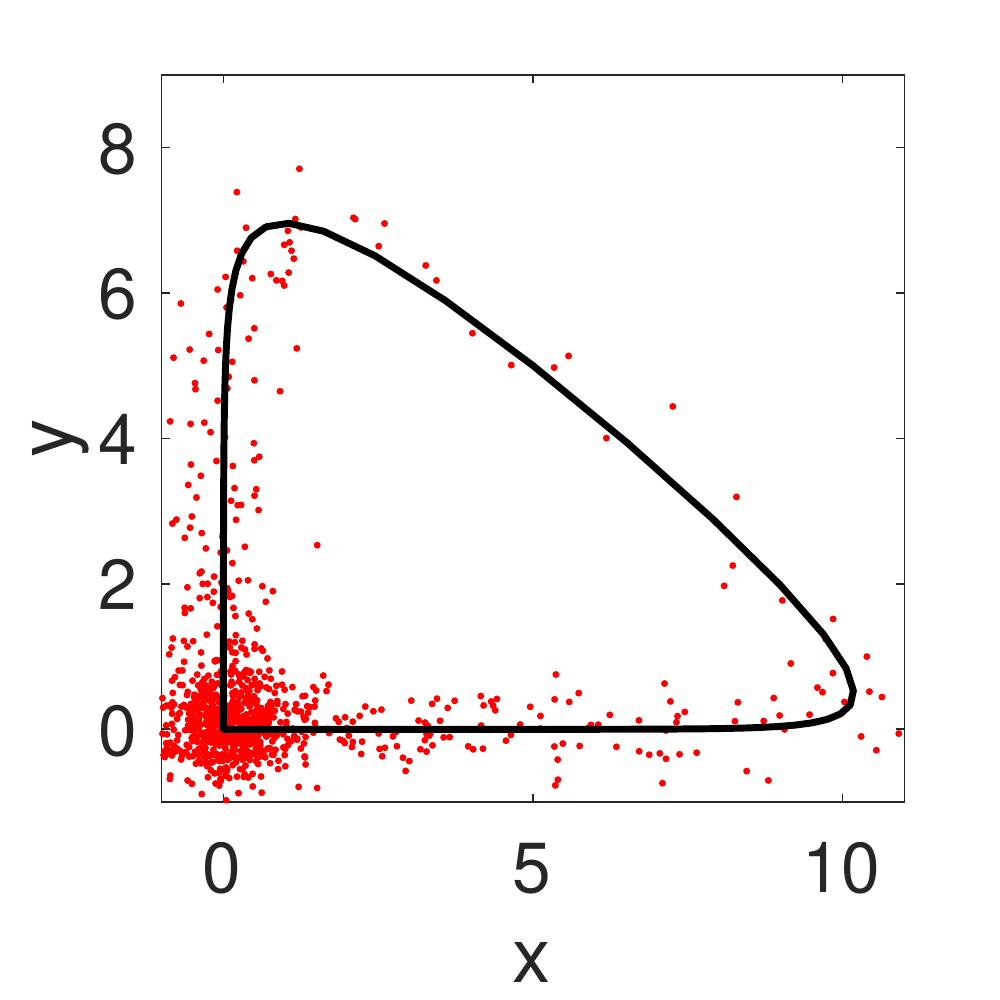} &
     \includegraphics[width = 0.19 \textwidth ]{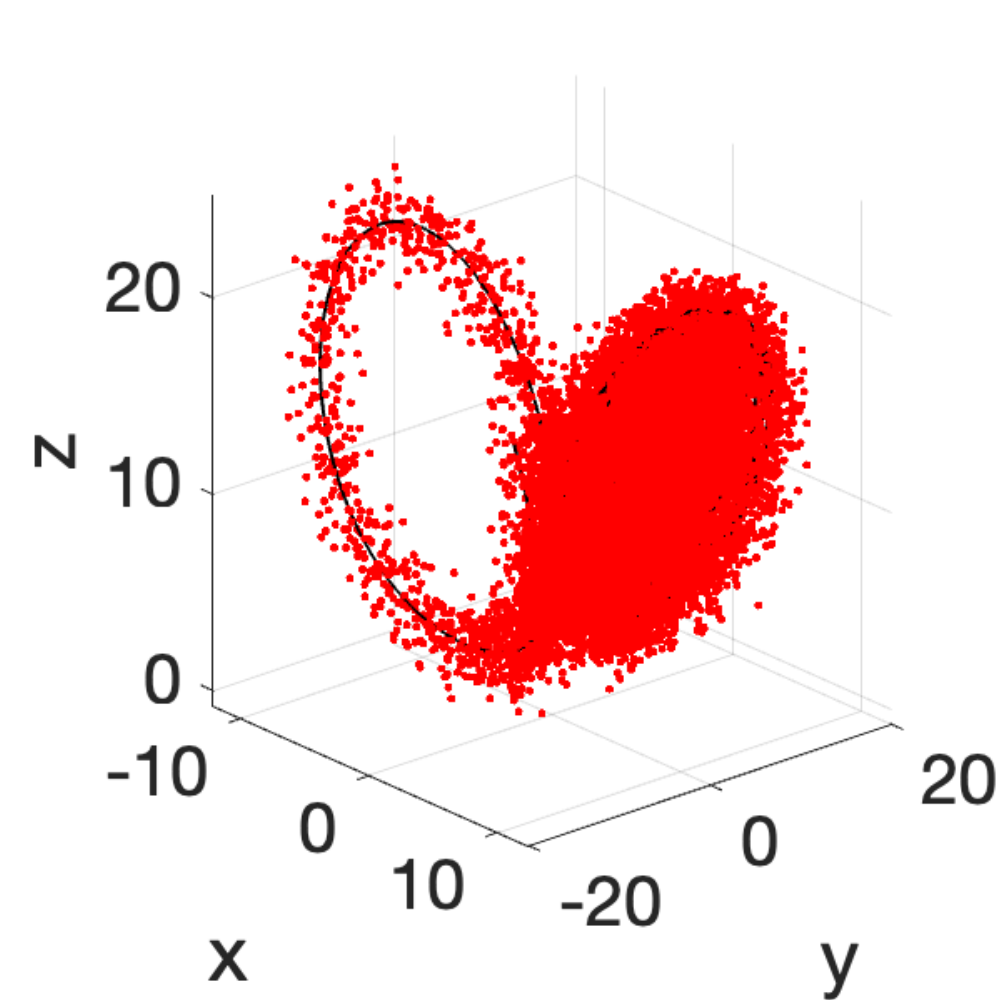}\\
    (f) Equ.  \eqref{e: ode, linear 2d} & (g) Equ. \eqref{e: ode, Duffing}   & (h) Equ. \eqref{e: ode, van-der-pol}  & (i) Equ. \eqref{e: ode, Lotka}    & (j) Equ. \eqref{e: ode, lorenz}  \\
    \includegraphics[width = 0.18\textwidth ]{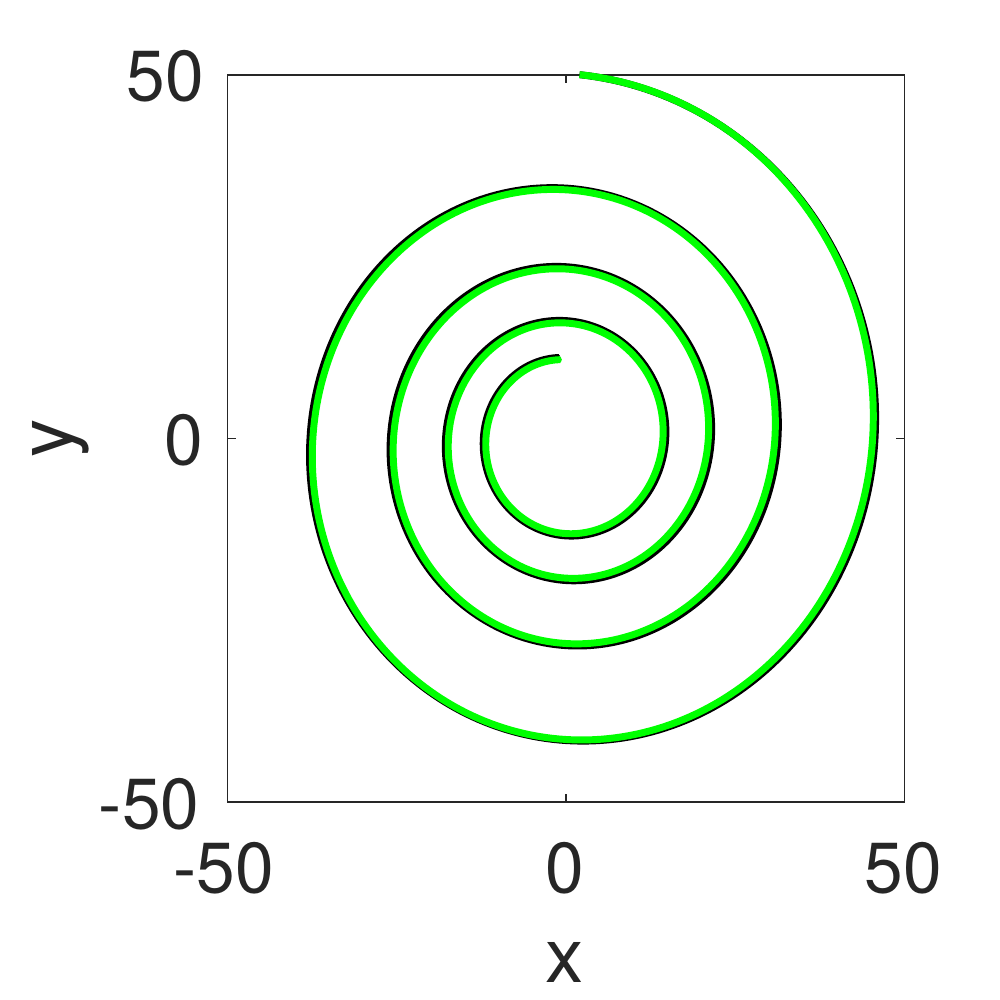}
         &  
     \includegraphics[width = 0.18\textwidth ]{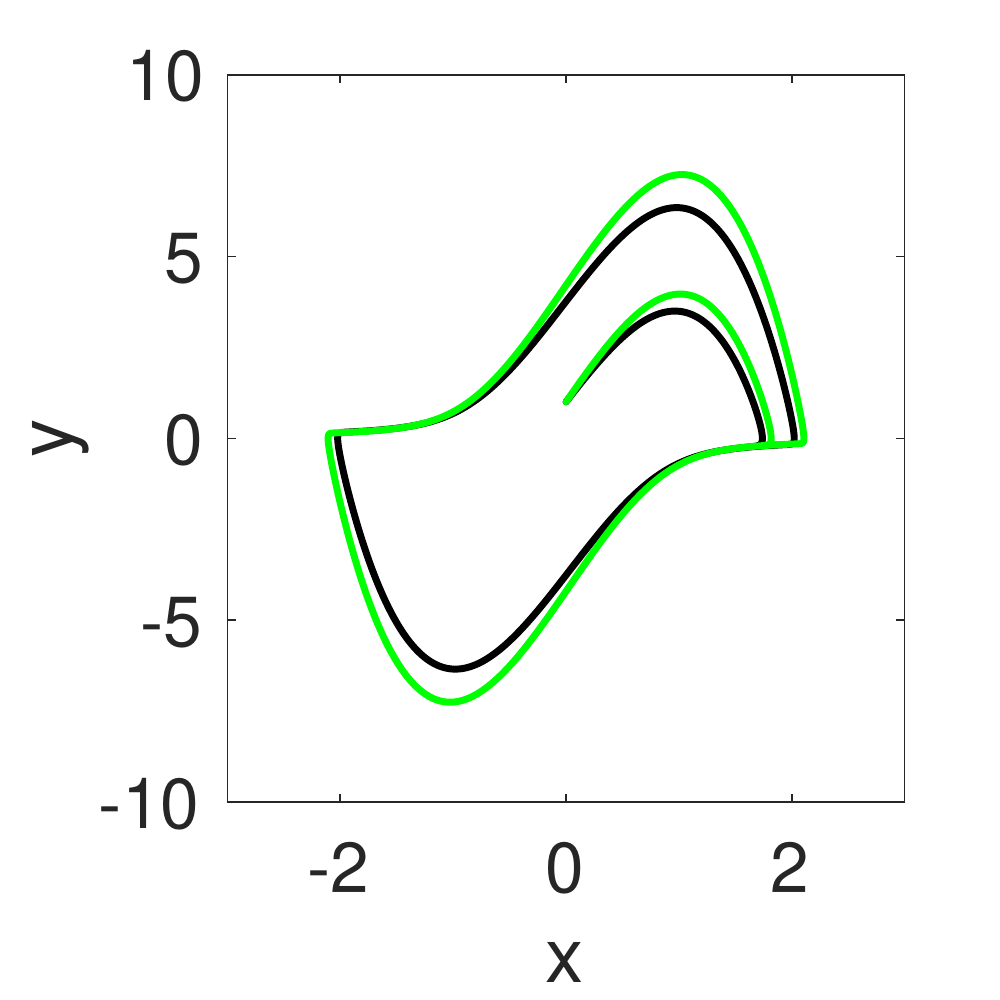}
         & 
    \includegraphics[width = 0.18\textwidth ]{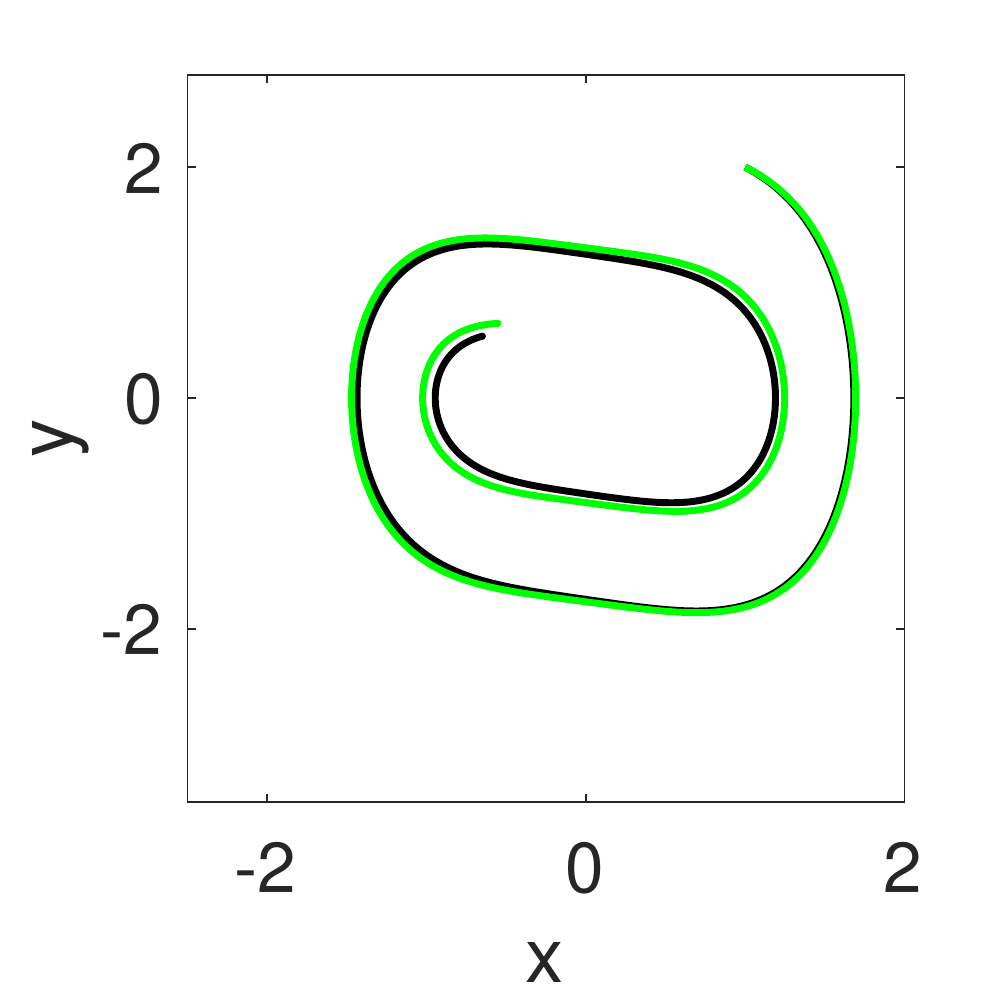} 
    &
    \includegraphics[width = 0.18\textwidth ]{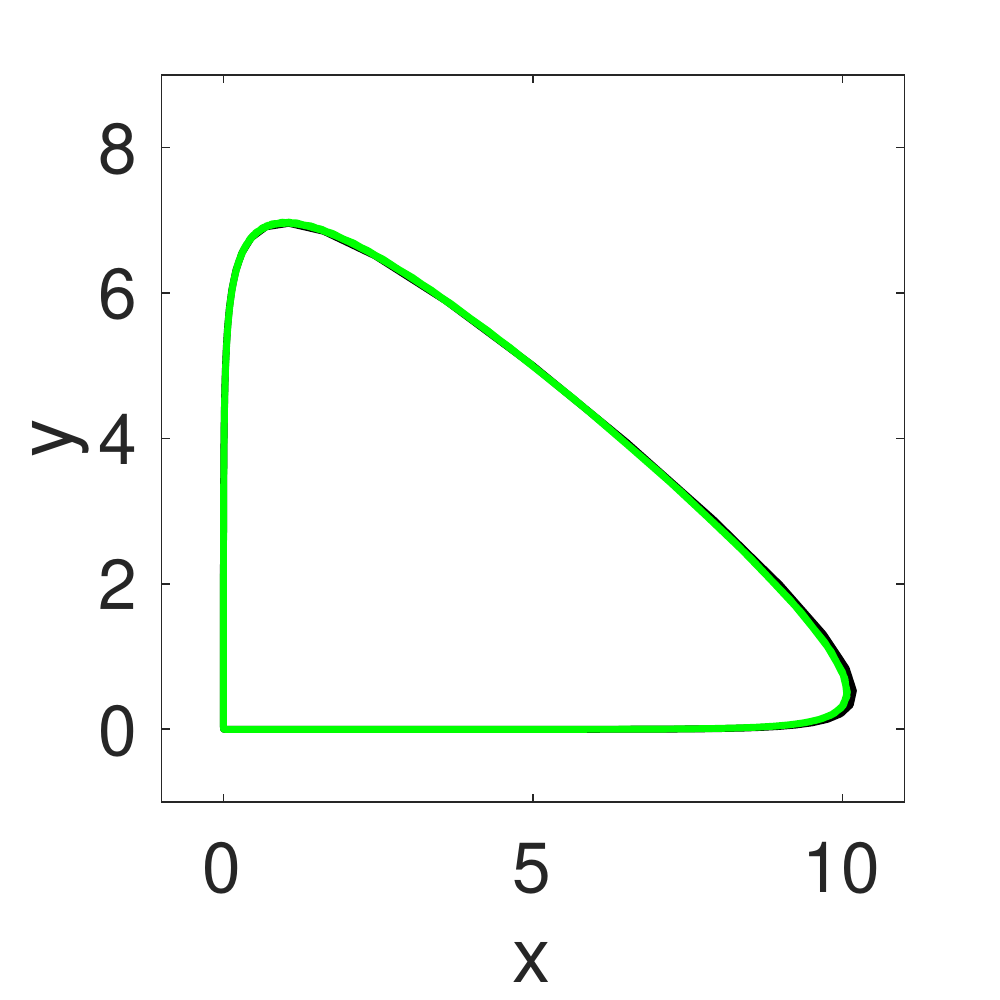} &
    \includegraphics[width = 0.19 \textwidth ]{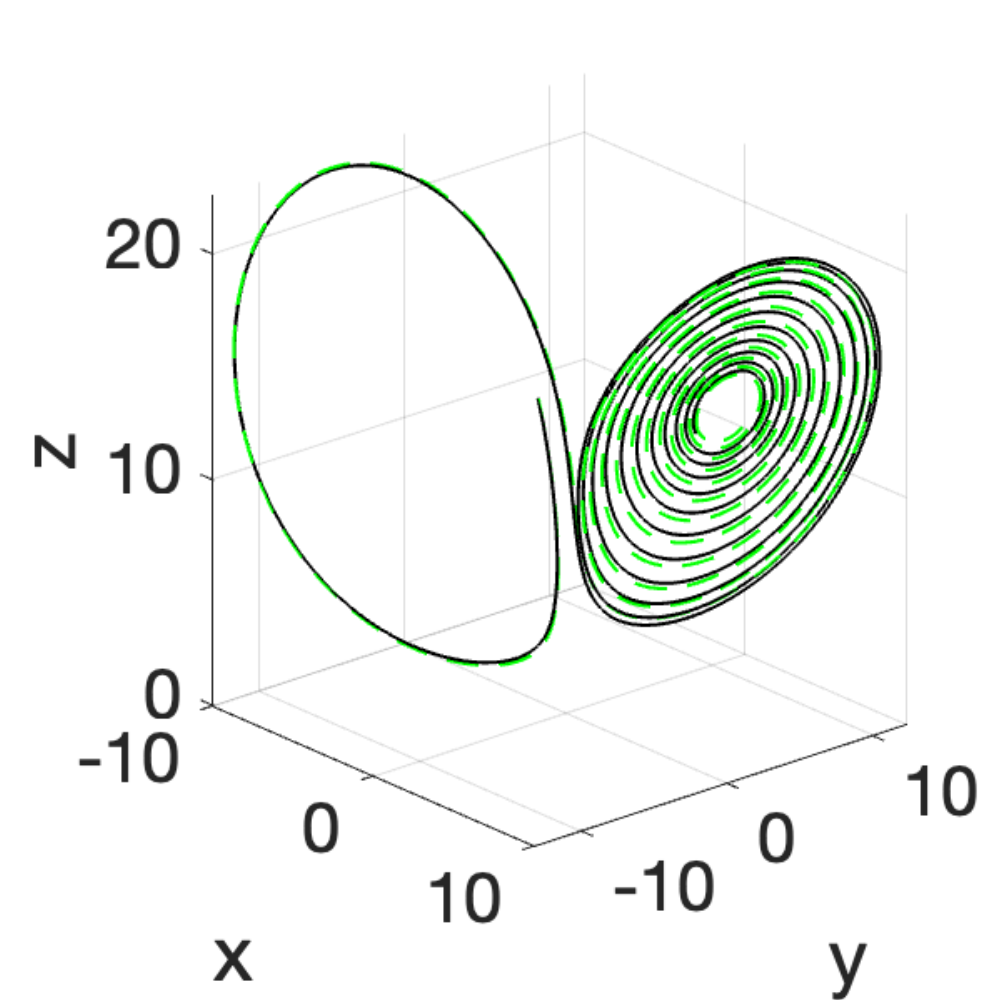} \\
    \end{tabular}
    \begin{tabular}{lrccccrr}
    \toprule
     & $\sigma_{\rm NSR}$ 
     & $E_2$ & $E_{\infty}$ & $E_{\rm res}$ & $E_{\rm dyn}$ & TPR & PPV  \\ \hline
     (f) 2D Linear System  \eqref{e: ode, linear 2d} & 0.2  
     & 0.006 & 0.063 & 0.061 & 0.484 & 1 & 1 \\
     (g) Van der Pol \eqref{e: ode, Duffing}  &0.2 
     & 0.055& 0.080 & 0.042 & 0.120 &1 & 1\\
     (h) Duffing \eqref{e: ode, van-der-pol} &0.2 
     & 0.088 & 0.204 & 0.048 & 0.156 & 0.75& 1 \\
    (i) Lotka-Volterra \eqref{e: ode, Lotka} &0.1 
    & 0.013 & 0.023 & 0.120 & 0.365 & 1 & 1\\
    (j) Lorenz \eqref{e: ode, lorenz} &0.2 
     & 0.011 & 0.041 & 0.032 & 2.380 & 1 & 1 \\ 
     \bottomrule
    \end{tabular}
    \caption{WeakIdent results for  ODE systems in Table \ref{T: odes}. (a)-(e): Given noisy data compared to the true dynamics. (f)-(j): Recovered systems via WeakIdent using true initial conditions.
    WeakIdent  recovers the dynamics close to the true dynamics with a small identification error. }
    \label{fig: vis dynamic of odes}
\end{figure}

Since ODE systems do not include spatial derivatives, they have lower computational cost in feature computation.
We consider polynomial terms with the highest order being 5. Table \ref{T: odes} presents details of the parameters used for simulation.
In Figure \ref{fig: vis dynamic of odes}, we show the identified dynamics and various identification errors obtained from WeakIdent on the 5 ODE systems listed in Table \ref{T: odes}. The noise-to-signal ratio is $\sigma_{\rm NSR} =0.2$ for  the linear system \eqref{e: ode, linear 2d}, the Van der Pol nonlinear system \eqref{e: ode, van-der-pol} and $\sigma_{\rm NSR} = 0.1$ for the rest of the systems.
Figure \ref{fig: vis dynamic of odes} (a)-(e) show the phase portraits of the given noisy data for the different ODEs (red) superimposed on the simulated true data (black). 
Figure \ref{fig: vis dynamic of odes} (f)-(j) show the WeakIdent results (green) compared to the true solution (black).  WeakIdent is able to find the correct support in the majority of the cases with $E_2\leq 0.088$. 

\begin{figure}
    \centering
\begin{tabular}{lcccc}
    & $E_2$ & $E_{\rm res}$ & TPR & PPV \\
    \toprule
    & (a1) & (a2) & (a3) & (a4)\\
     \textbf{WeakIdent} & 
     \includegraphics[width = 0.17\textwidth, align = c]{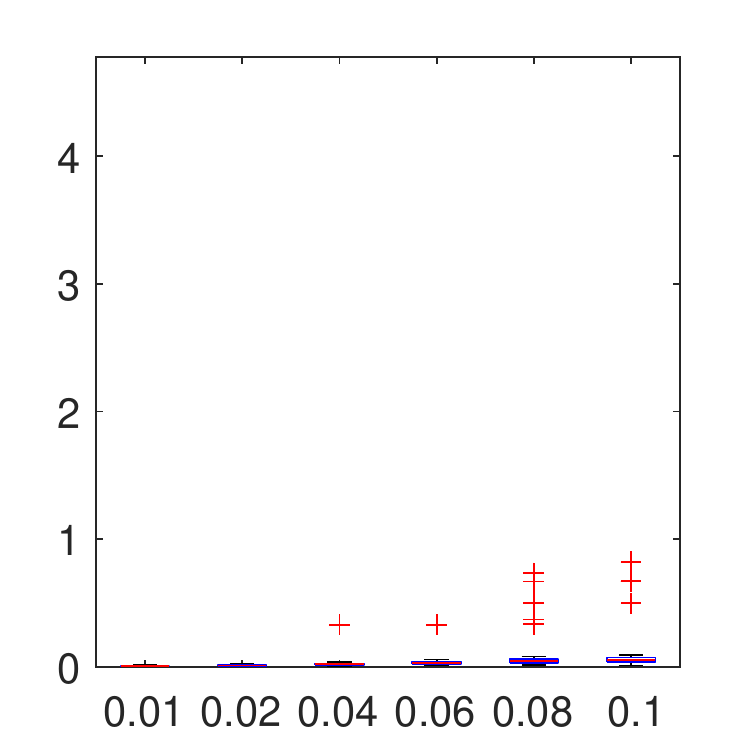} &
     \includegraphics[width = 0.17\textwidth, align = c]{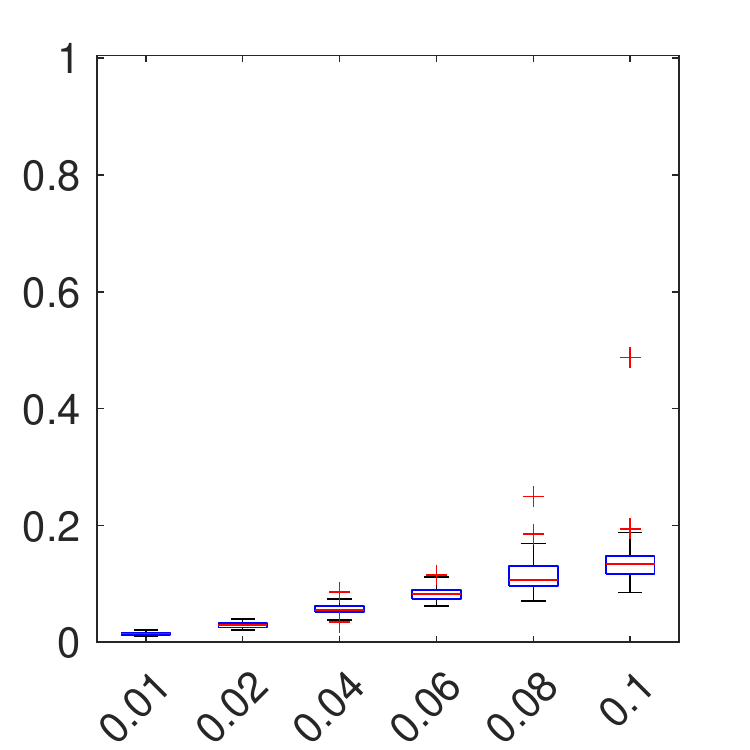} & 
     \includegraphics[width = 0.17\textwidth, align = c]{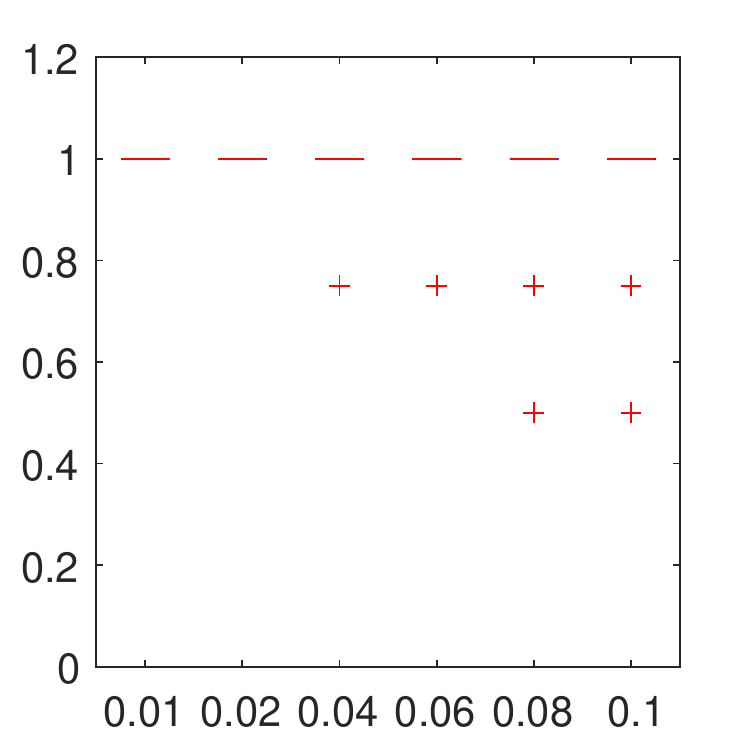} & 
     \includegraphics[width = 0.17\textwidth, align = c]{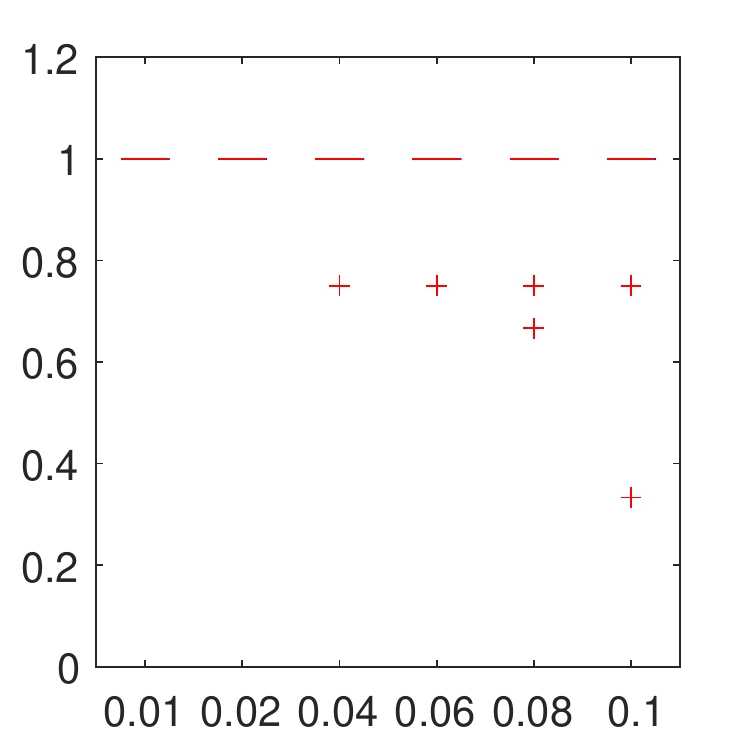} \\  
     & (b1) & (b2) & (b3) & (b4) \\
     WODE\cite{messenger2021weak} & 
     \includegraphics[width = 0.17\textwidth, align = c]{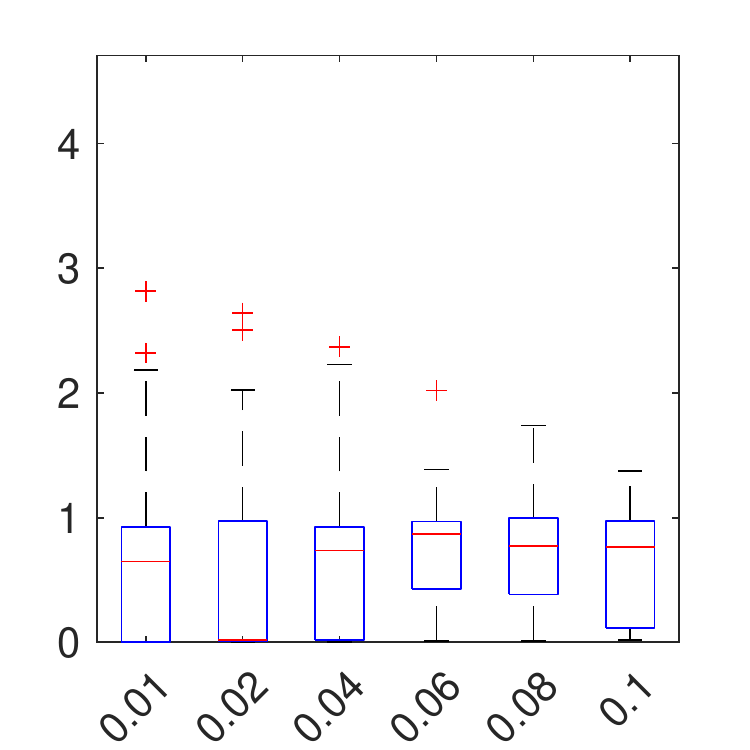} &
     \includegraphics[width = 0.17\textwidth, align = c]{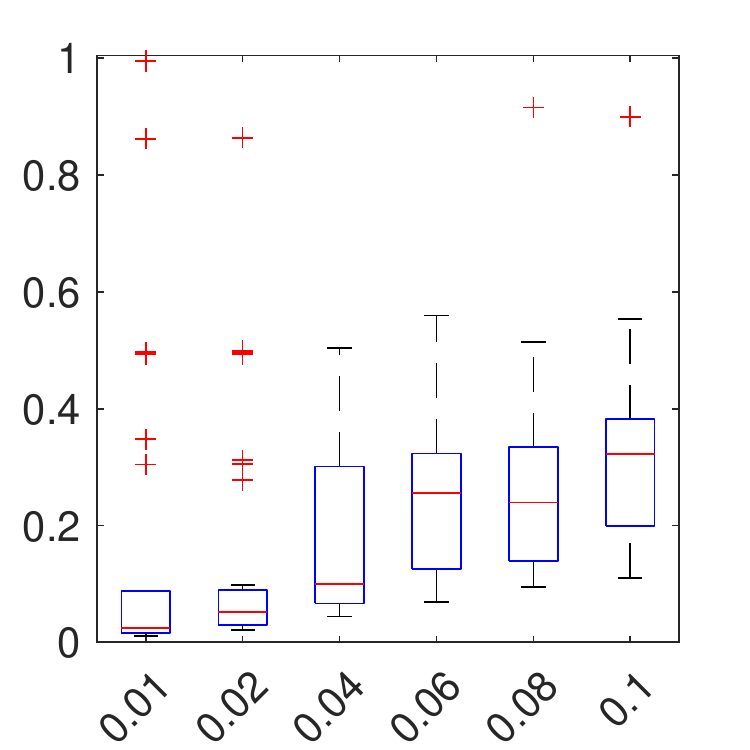} & 
     \includegraphics[width = 0.17\textwidth, align = c]{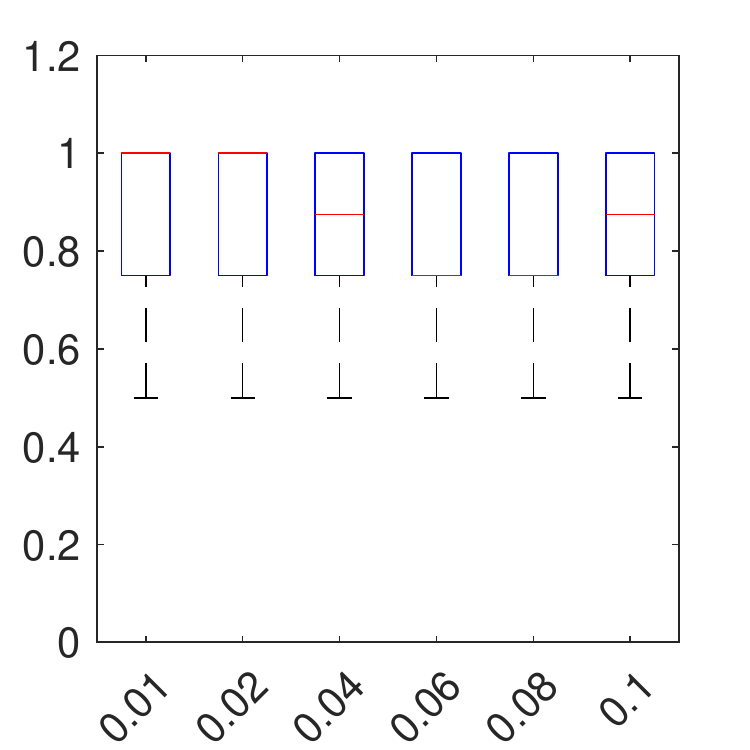} & 
     \includegraphics[width = 0.17\textwidth, align = c]{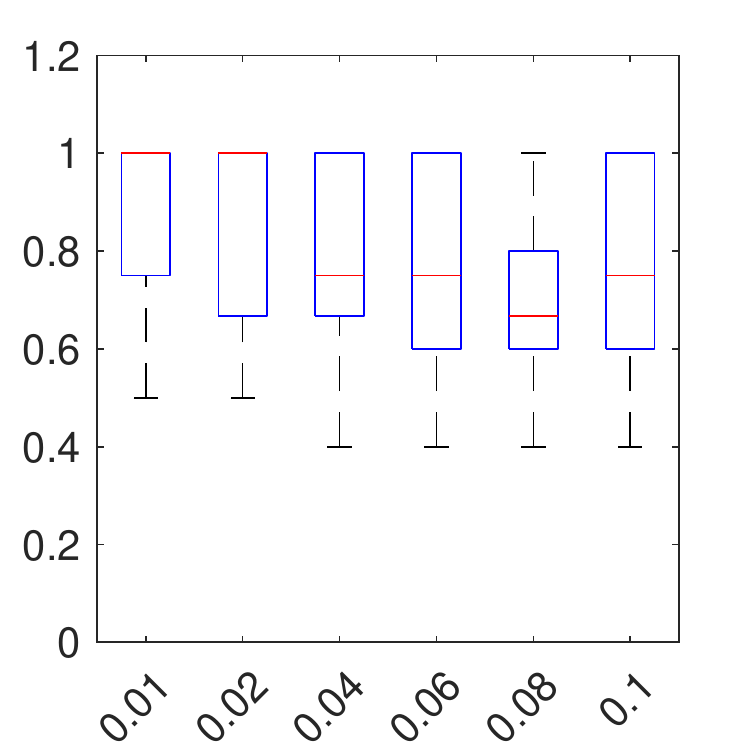} \\ 
    & (c1) & (c2) & (c3) & (c4) \\
    SINDy \cite{brunton2016sparse} &
     \includegraphics[width = 0.17\textwidth, align = c]{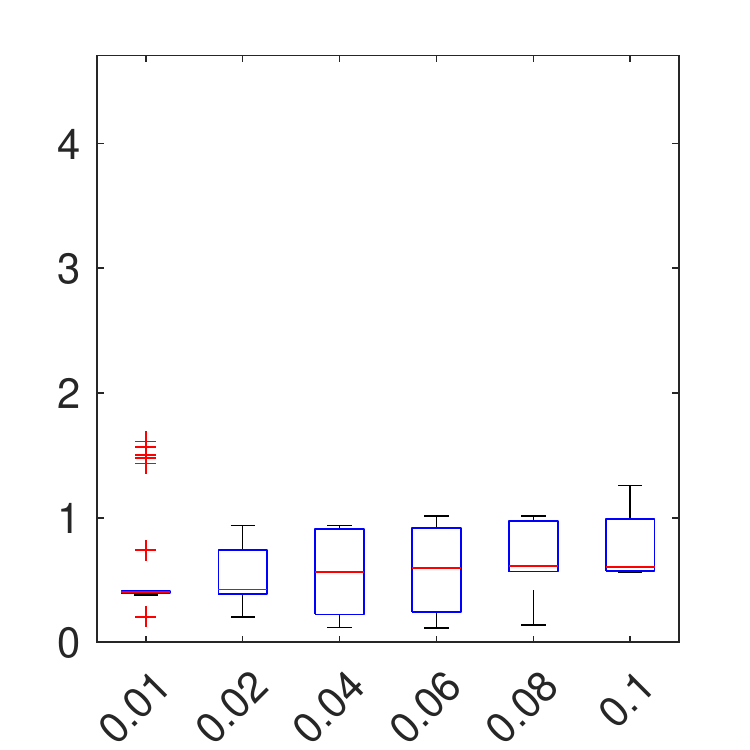} &
     \includegraphics[width = 0.17\textwidth, align = c]{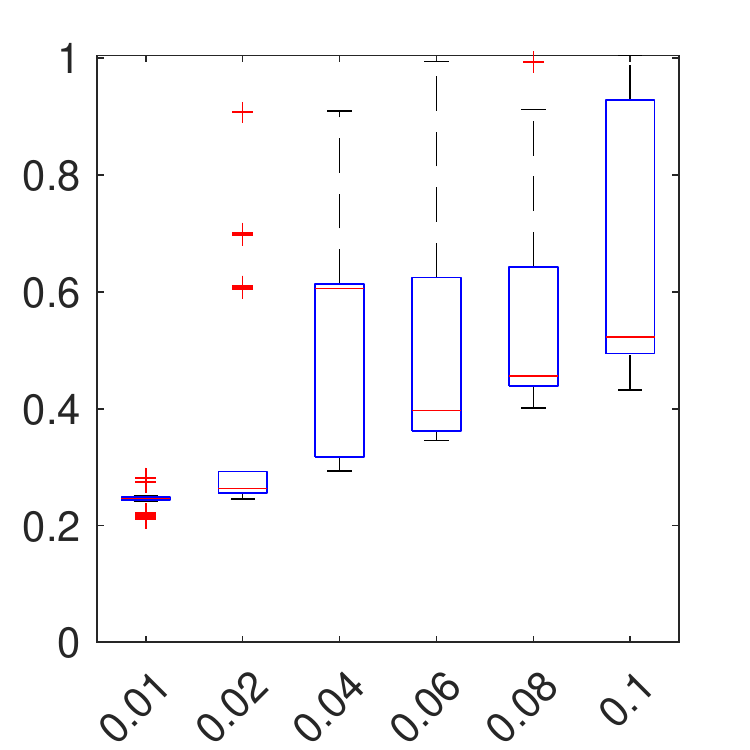} & 
     \includegraphics[width = 0.17\textwidth, align = c]{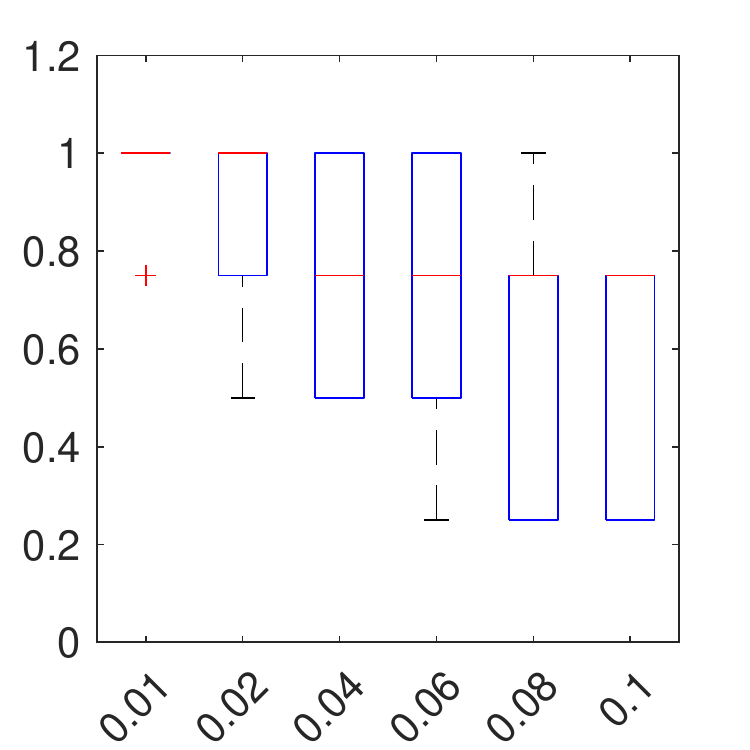} & 
     \includegraphics[width = 0.17\textwidth, align = c]{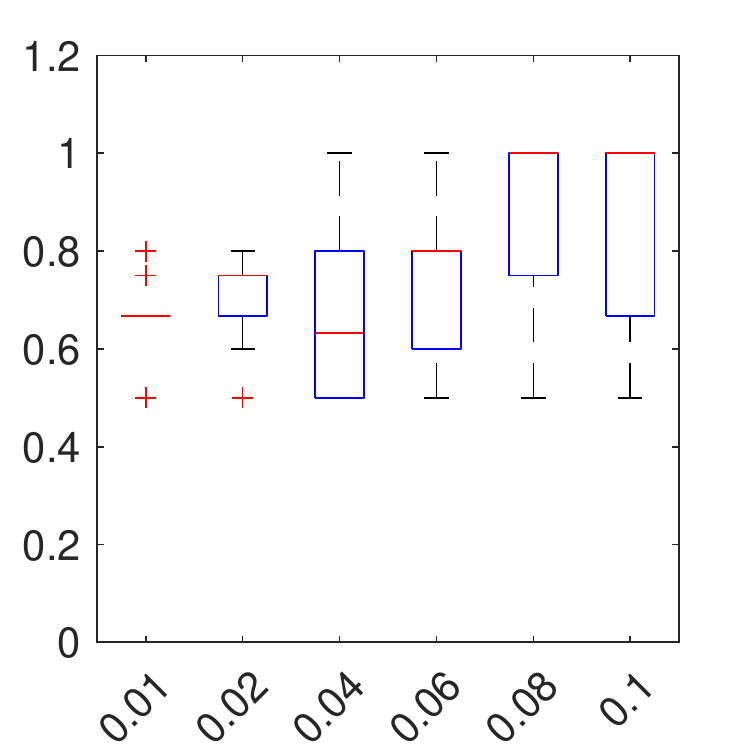} \\
     & (d1) & (d2) & (d3) & (d4) \\
     SC \cite{he2020robust} & \includegraphics[width = 0.17\textwidth, align = c]{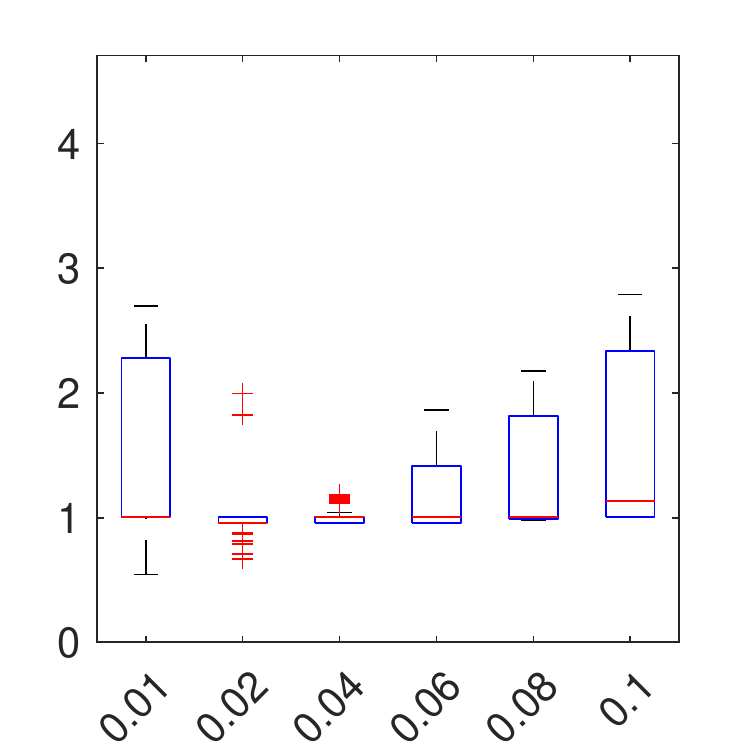} &
     \includegraphics[width = 0.17\textwidth, align = c]{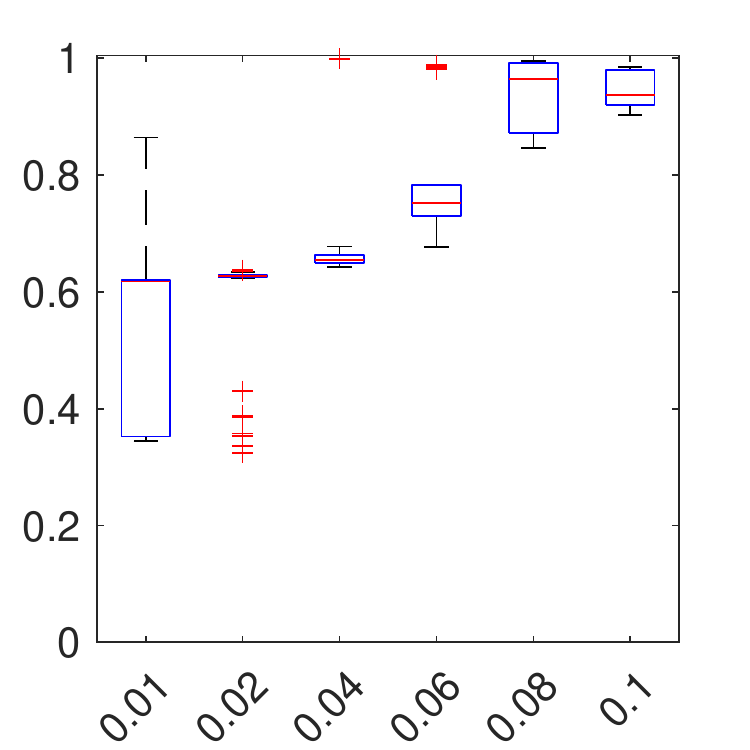} & 
     \includegraphics[width = 0.17\textwidth, align = c]{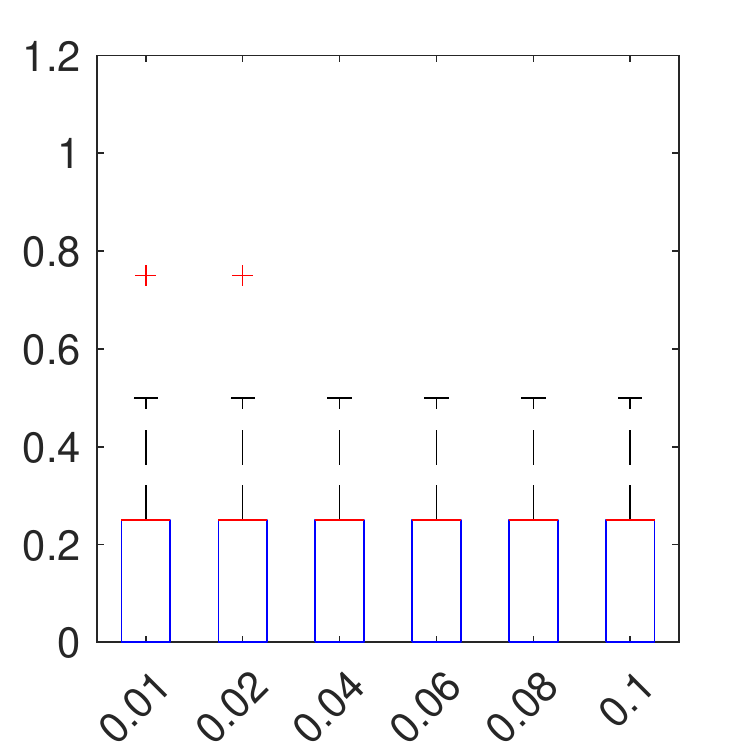} & 
     \includegraphics[width = 0.17\textwidth, align = c]{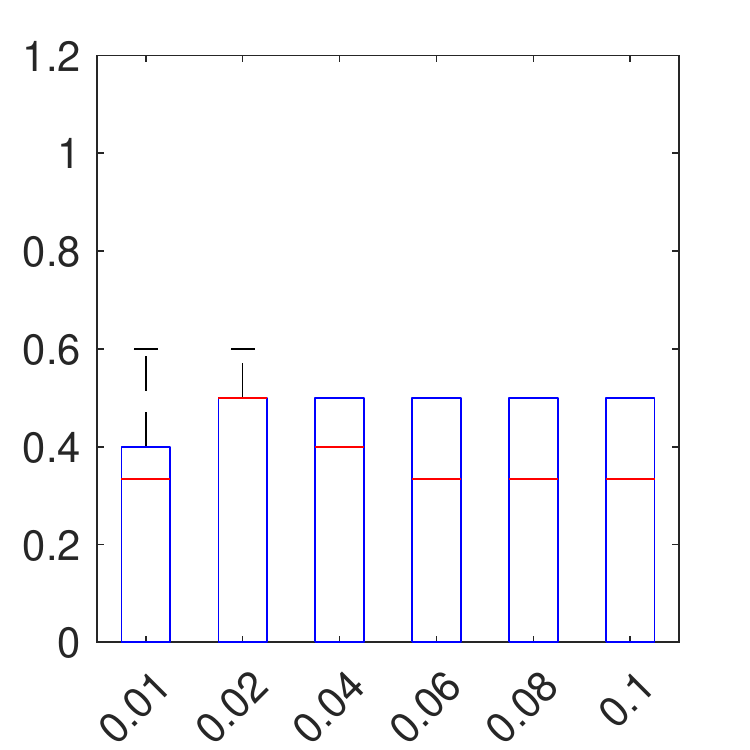} \\
     & (e1) & (e2) & (e3) & (e4) \\
     ST \cite{he2020robust} &
     \includegraphics[width = 0.17\textwidth, align = c]{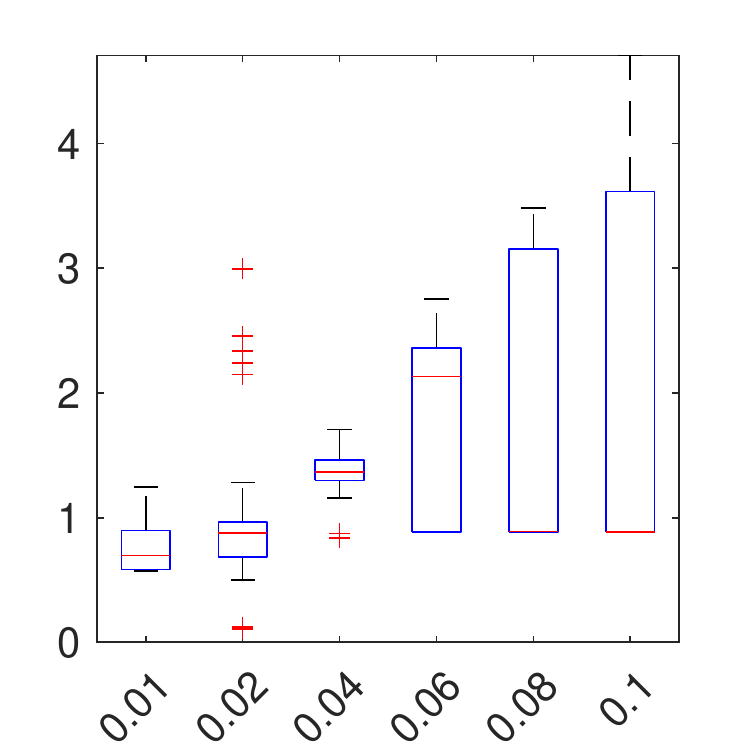} &
     \includegraphics[width = 0.17\textwidth, align = c]{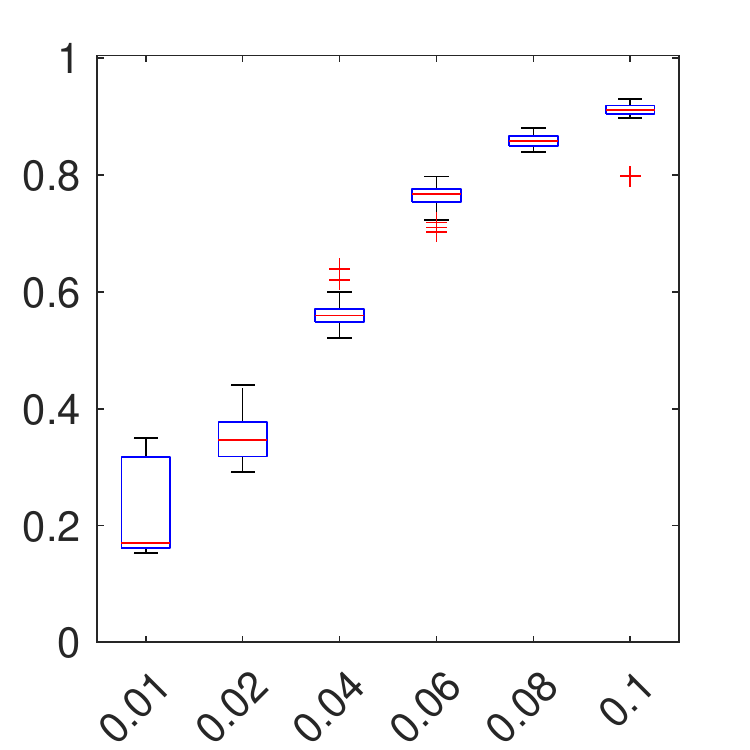} & 
     \includegraphics[width = 0.17\textwidth, align = c]{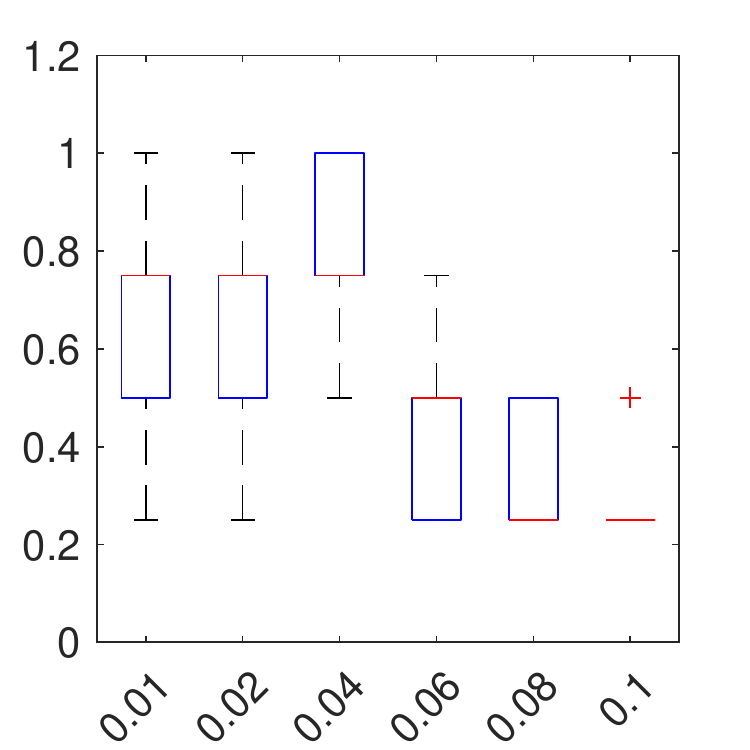} & 
     \includegraphics[width = 0.17\textwidth, align = c]{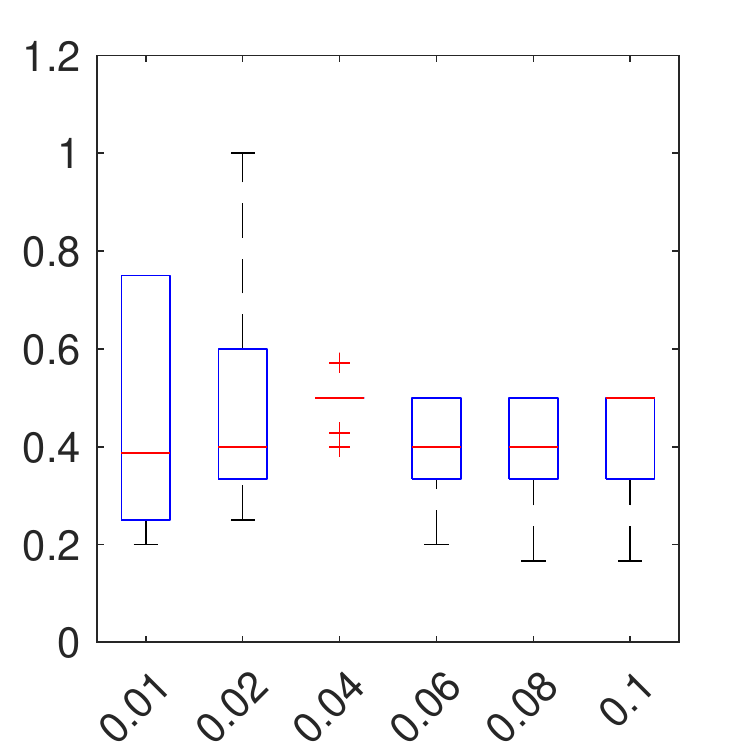} \\     
     \bottomrule
    \end{tabular}
    \caption{
     The Lotka-Volterra equation (\ref{e: ode, Lotka}). Statistical comparisons between (a1)-(a4) WeakIdent, (b1)-(b4) WODE \cite{messenger2021weak}, (c1)-(c4) SINDy\cite{brunton2016sparse}, (d1)-(d4) SC\cite{he2020robust} and (e1)-(e4) ST\cite{he2020robust}. 
    The  $E_2, E_{\rm res}$ errors, TPR and PPV are shown from 50 experiments for each  $\sigma_{\rm NSR}\in \{ 0.01,0.02,  ,...,0.1\}$  using box-plots.  
    Notice that for WeakIdent, the $E_2$ error  is  lower with less variations, and the TPR and PPV are closer to 1 as compared with that obtained from other methods.
    }
    \label{fig: ode - box plot - Lotka-Volterra system}
\end{figure}

Figure \ref{fig: ode - box plot - Lotka-Volterra system} compares the recovery results for the Lotka-Volterra (LV) system \eqref{e: ode, Lotka} across different methods, showing results for the given data sets with various noise levels. 
The methods we compare include
WODE\cite{messenger2021weak},
SINDy\cite{brunton2016discovering},
Robust IDENT SC\cite{he2020robust} and ST\cite{he2020robust}. 
Each column is associated with an error type and each row gives results from one method. WeakIdent is able to capture the correct support with a low coefficient error in the last rows. 
WODE, SINDy, SC and ST has larger coefficient errors with incorrect support in many cases.  
A similar statistical comparison between these methods on the Lorenz system \eqref{e: ode, lorenz} is shown in Figure \ref{fig: ode - box plot - Lorenz system} in the Appendix \ref{AS:ode}. 
We refer to Table \ref{F: an ode example - recovering noisey data - Lotka-Volterra} and Table \ref{F: an ode example - recovering noisey data - Nonlinear Lorenz}  for the recovery results of the Lotka-Volterra system \eqref{e: ode, Lotka} and  the Lorenz system \eqref{e: ode, lorenz} from two noisy data sets with $\sigma_{S\rm NR}= 0.1$. 
We also provide a comprehensive comparison on all ODE systems listed in Table \ref{T: odes} in Appendix \ref{AS:ode} (See Figure \ref{fig: summary recovering odes} for the details).

\subsection{Influence of the initial condition in WeakIdent}
\begin{figure}[t]
 \begin{center}
 \begin{tabular}{ccccccc}
 IC (1) & IC (2) & IC (3) & IC (4) & IC (5) \\
 \includegraphics[width = 0.18 \textwidth ]{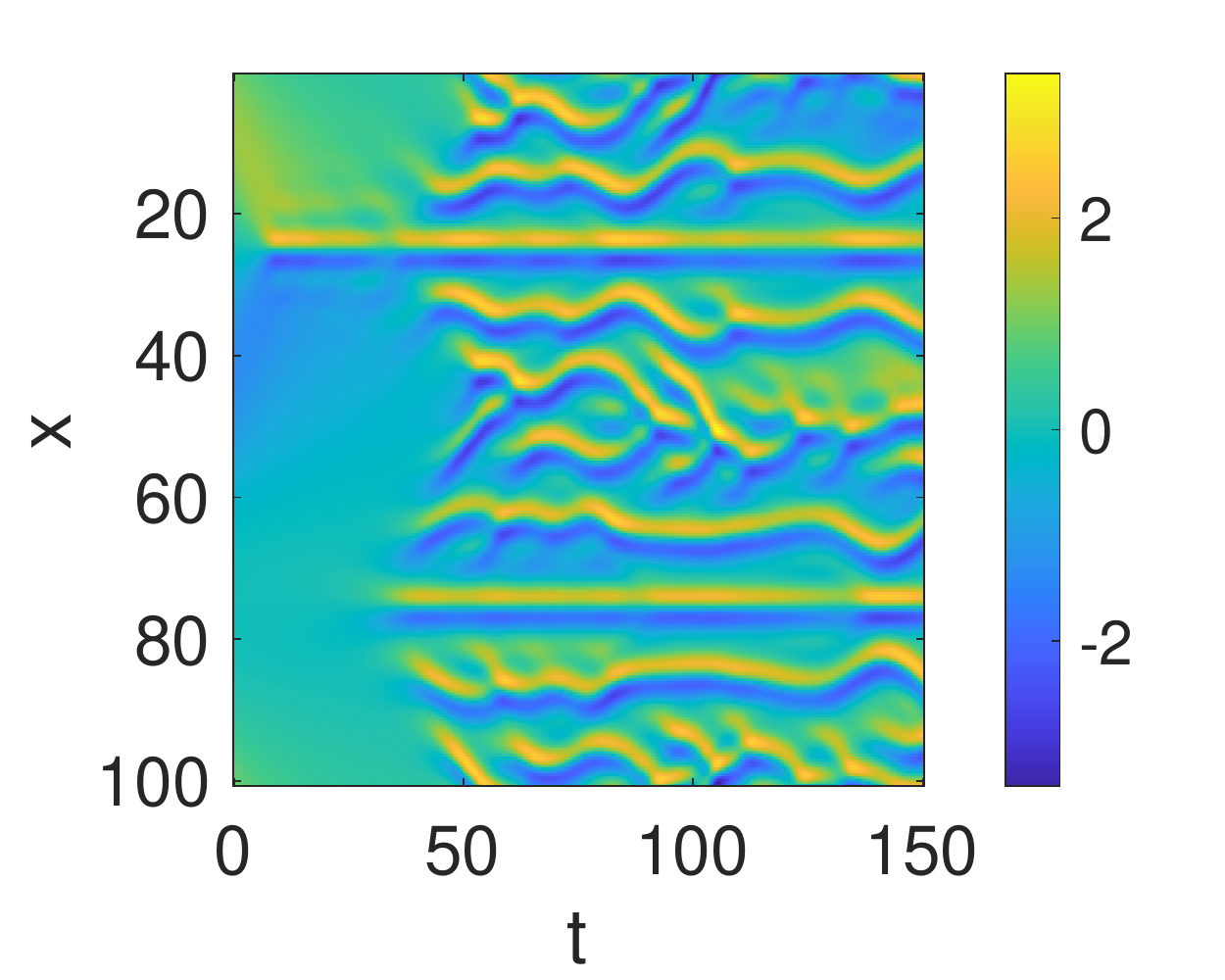} & 
 \includegraphics[width = 0.18 \textwidth ]{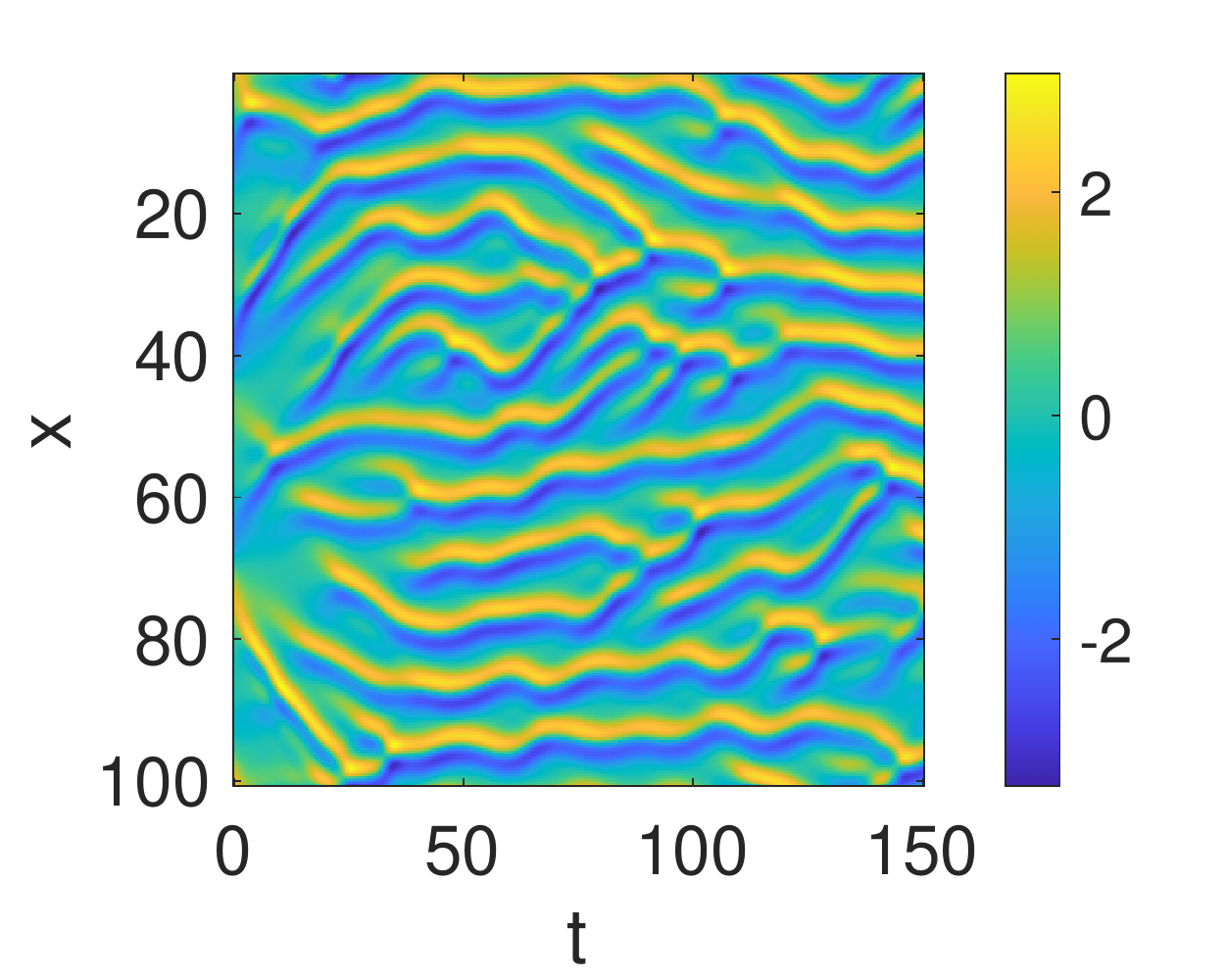} & 
 \includegraphics[width = 0.18 \textwidth ]{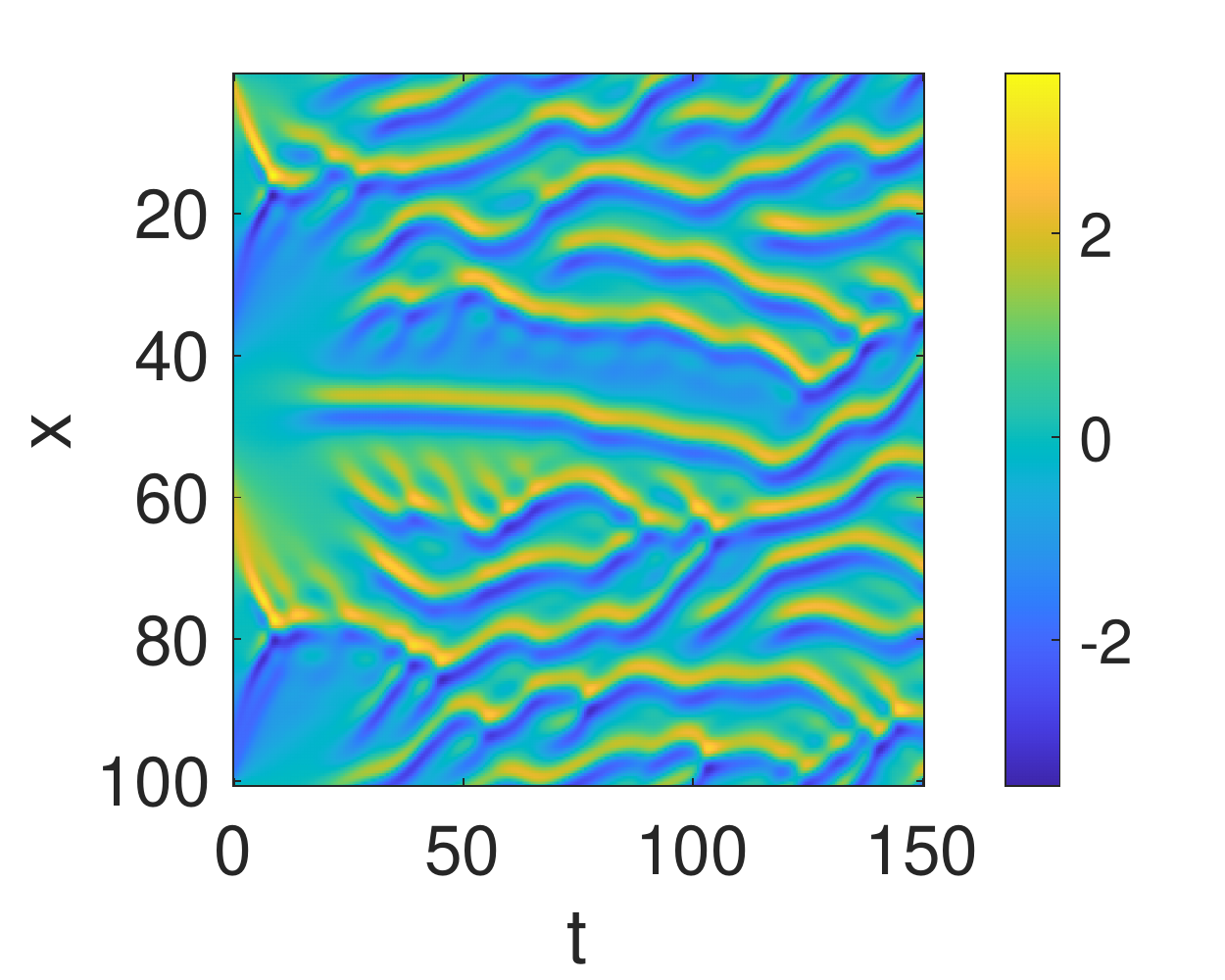} & 
 \includegraphics[width = 0.18 \textwidth ]{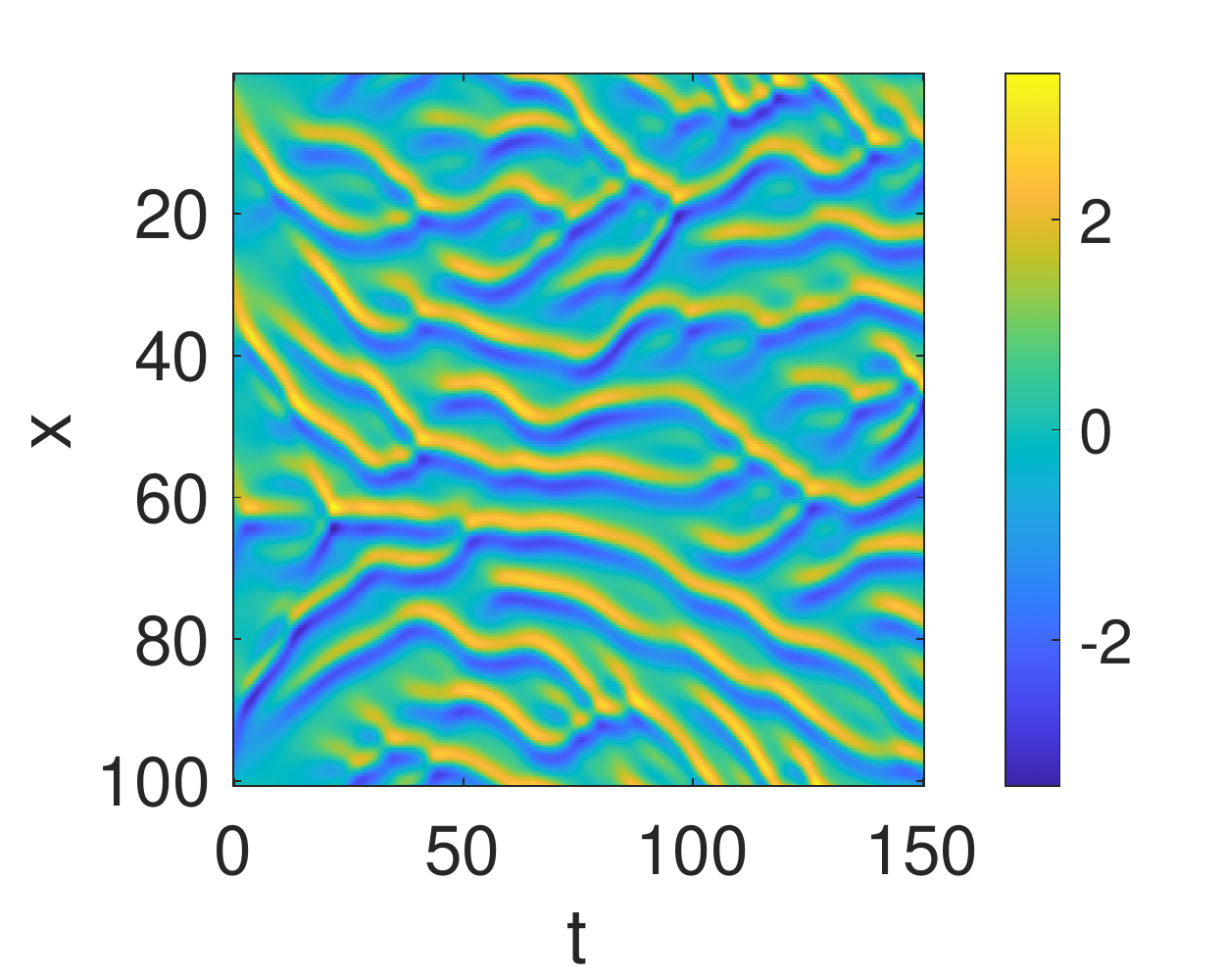} & 
 \includegraphics[width = 0.18 \textwidth ]{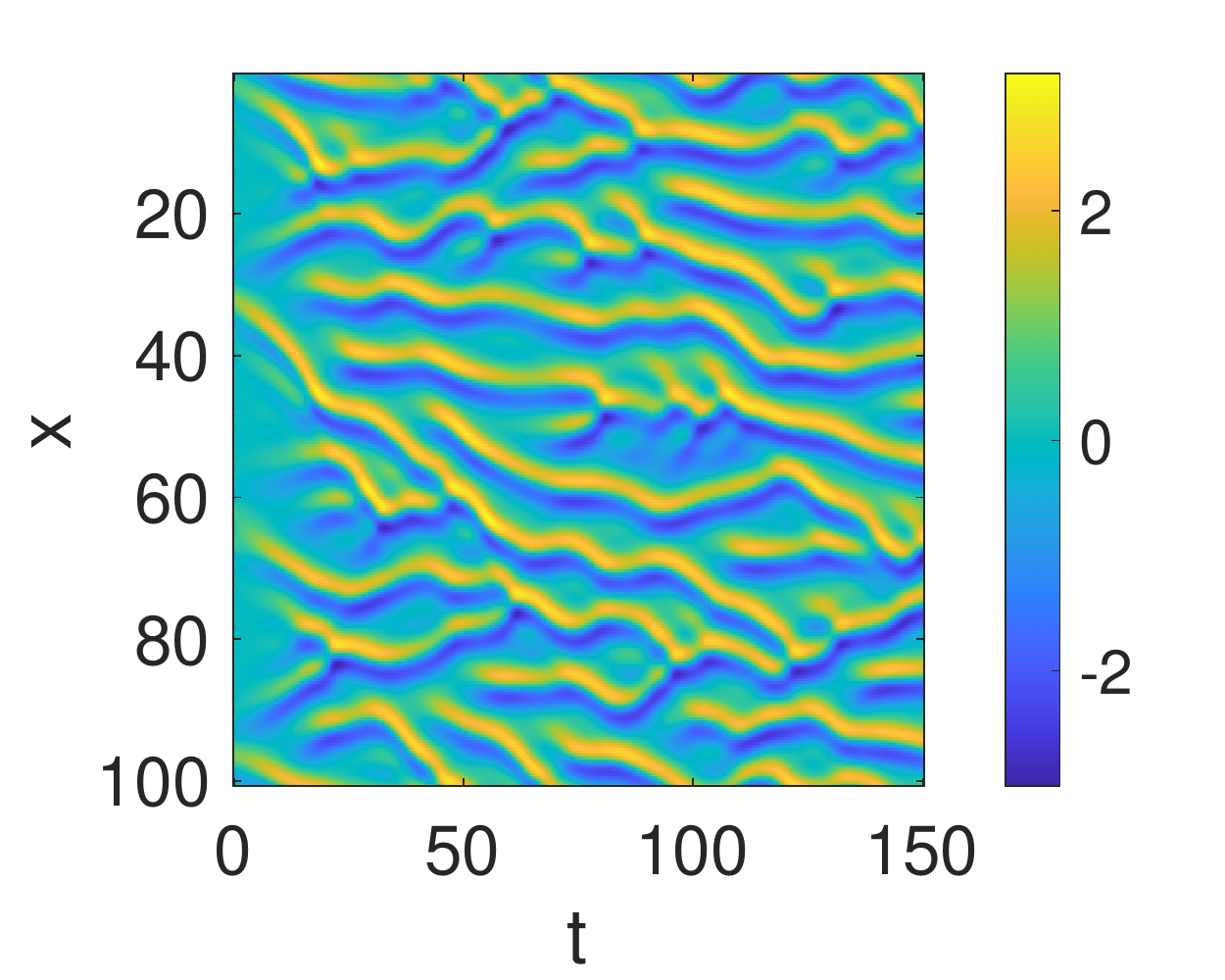} \\
 \end{tabular}
\begin{tabular}{cc|cc|cc}
\toprule
\multicolumn{2}{c|}{$E_2$} & 
\multicolumn{2}{c|}{TPR} &
\multicolumn{2}{c}{PPV} \\ 
\midrule
 (a) \textbf{WeakIdent} & (b) WPDE &  (c) \textbf{WeakIdent} & (d) WPDE & (e) \textbf{WeakIdent} & (f) WPDE \\
 \includegraphics[width = 0.14 \textwidth ]{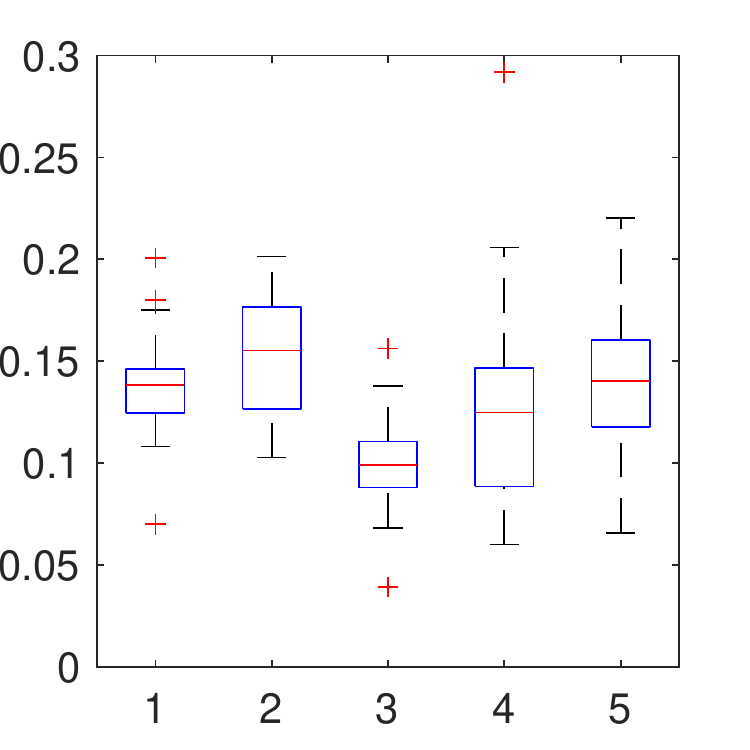} &
  \includegraphics[width = 0.14 \textwidth ]{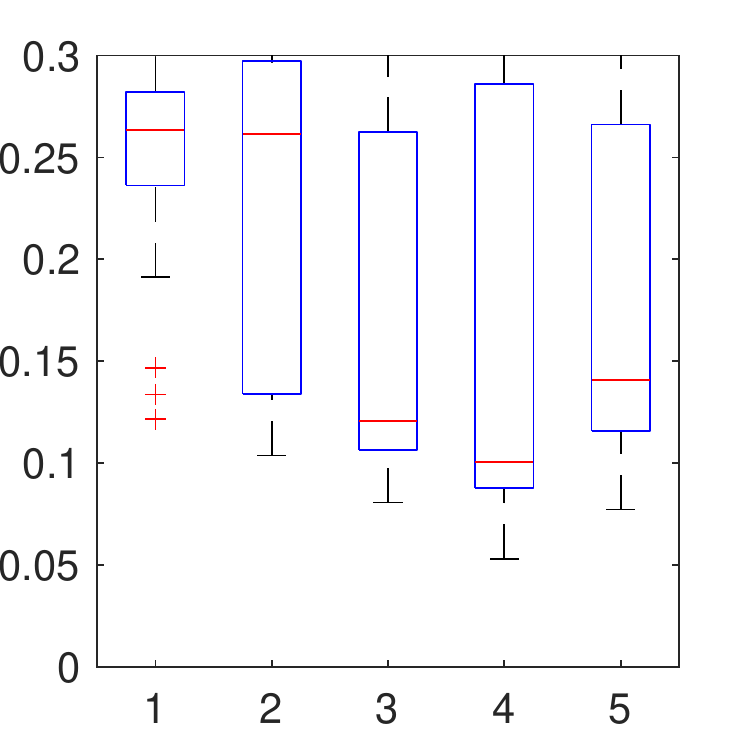} & 
 \includegraphics[width = 0.14 \textwidth ]{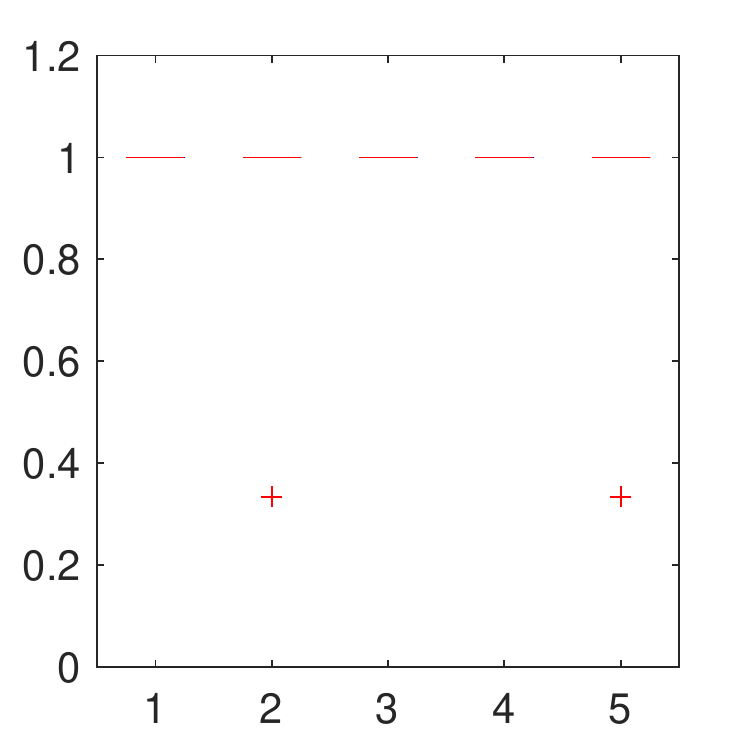} &
  \includegraphics[width = 0.14 \textwidth ]{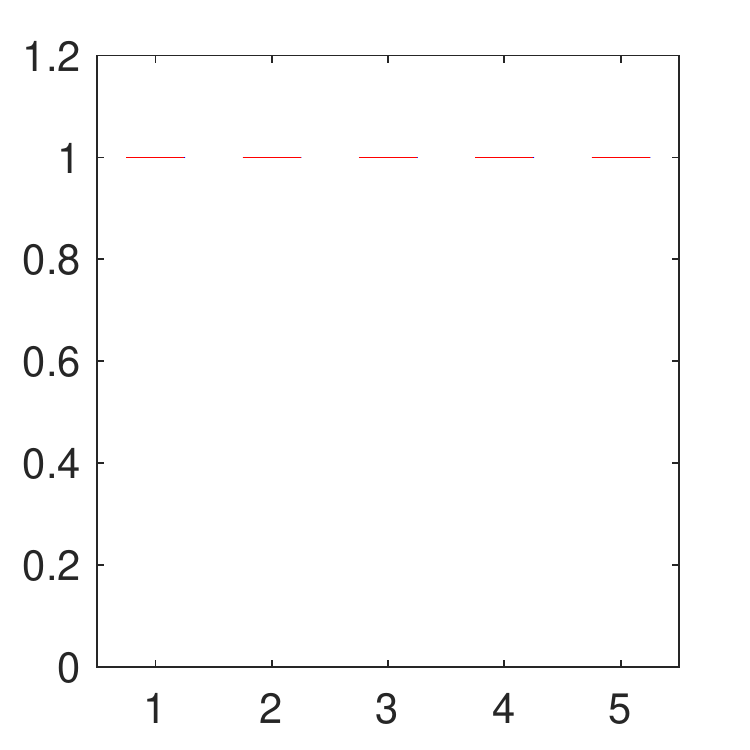} & 
 \includegraphics[width = 0.14 \textwidth ]{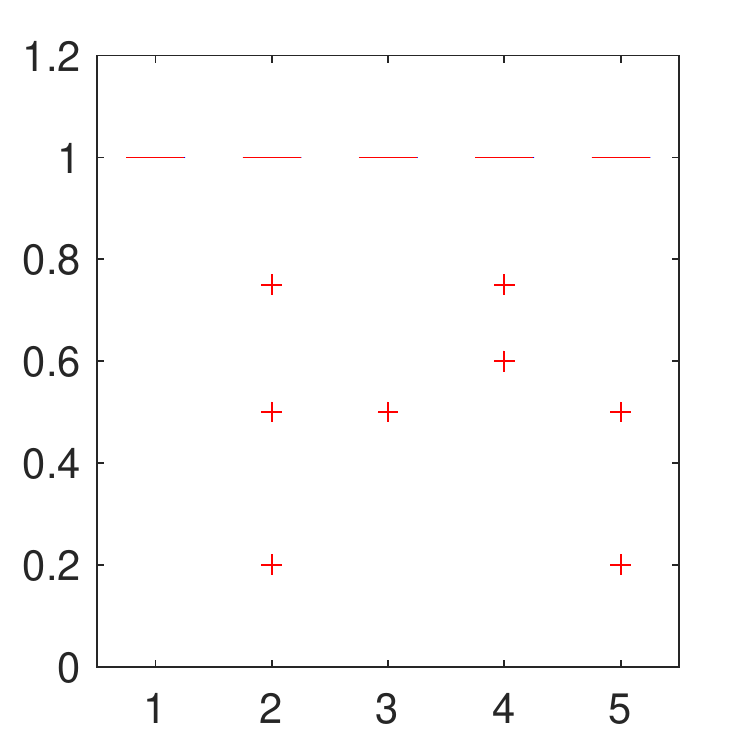} & \includegraphics[width = 0.14 \textwidth ]{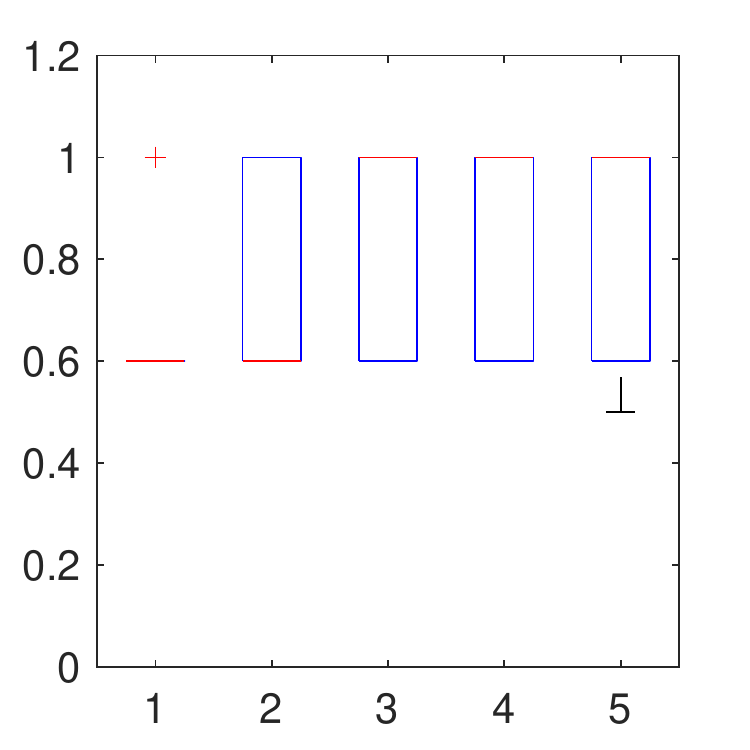} \\
\bottomrule
 \end{tabular}
 \end{center}
\caption{The  KS equation \eqref{e: pde KS} using five different  initial conditions (1)-(5)  with the noisy level of $\sigma_{\rm NSR}=0.6$. In (a)-(f) the $x$-axis is the index of initial conditions (1)-(5). For each initial condition, the box plot represents the statistical results over 50  experiments.  WeakIdent gives a smaller $E_2$ error, and PPV is  closer to 1 with less variations. } 
\label{fig: KS different initial conditions}
 \end{figure}

Figure \ref{fig: KS different initial conditions} shows comparisons of WeakIdent and WPDE for the KS equation \eqref{e: pde KS} on noisy data with $\sigma_{\rm NSR}=0.6$, using 5 different initial conditions:  (1)  $u(x,0)= \cos(x/16).*(1+\sin(x/16))$, (2) 
$u = \cos(x/4).*(1+\sin(x/5))$,
(3) $u = \cos(x/10).*(1+\cos(x/5))$,
(4) $u = \sin(x/4).*(1+\cos(x/5))$,
(5) $u = \sin(x.^2/4)$.
The top row illustrates the given clean data from the different initial conditions yielding  different pattern evolution. In each box plot, the $x$-axis gives  the indices of the initial condition (1)-(5). WeakIdent recovery is robust across these different patterns  in recovering this system with higher order features.

\subsection{The choice of the trimming parameter $\mathcal{T}$} 
\begin{figure}
    \centering
    \begin{tabular}{cccc}
    (a) The KdV equation \eqref{e: pde kdv} & (b) The KS equation \eqref{e: pde KS} \\
    \includegraphics[width = 0.4\textwidth]{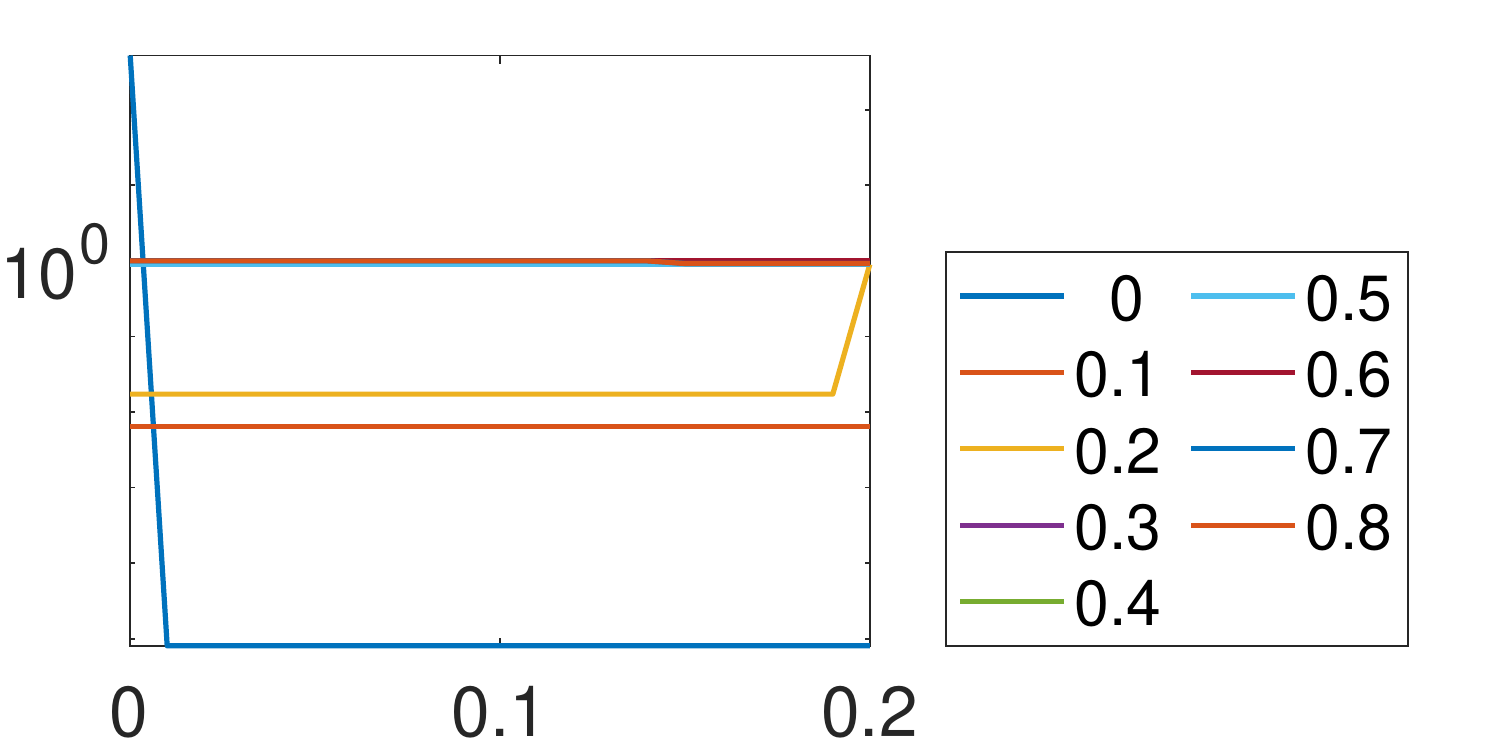}   
    & 
      \includegraphics[width = 0.4\textwidth]{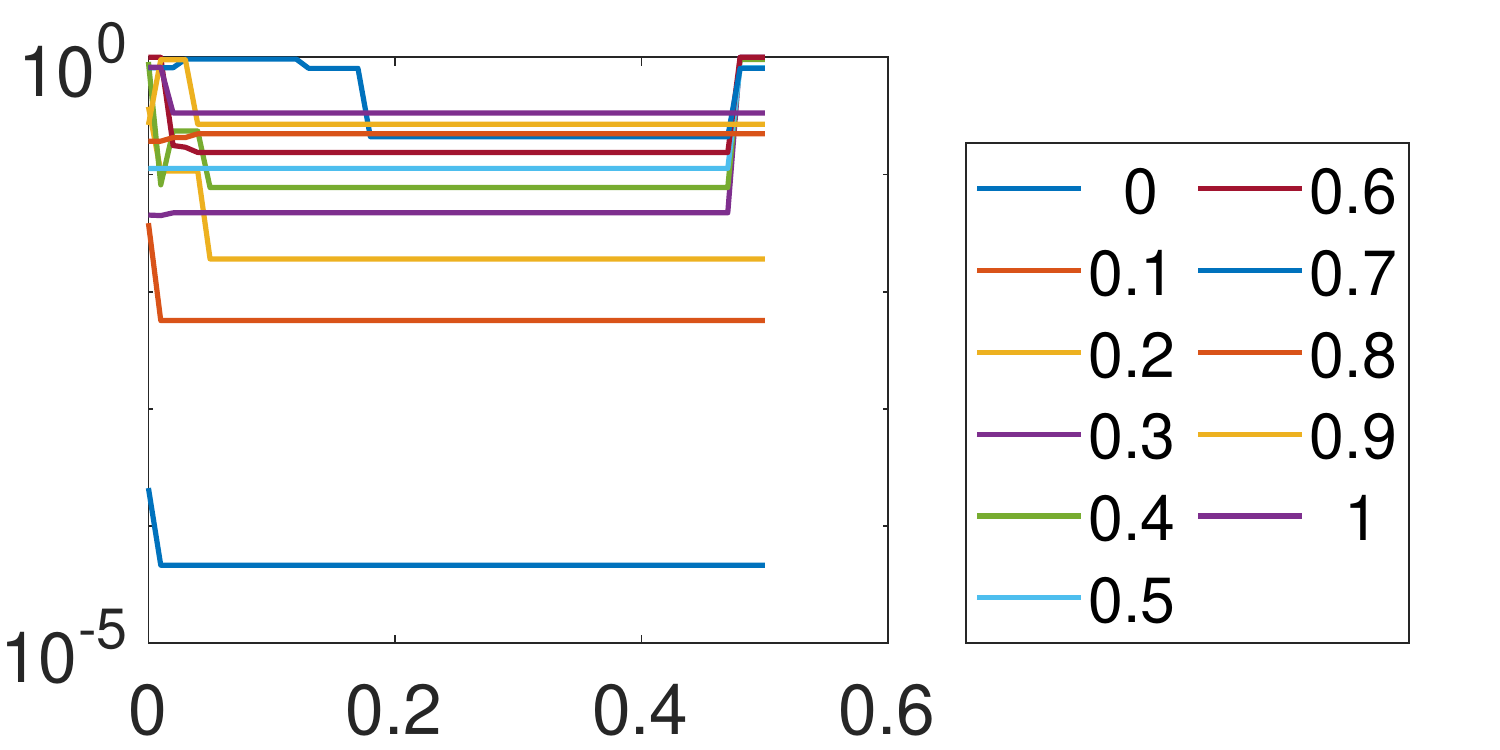}   
    \end{tabular}
\caption{The coefficient $E_2$ error ($y$-axis) versus the trimming parameter $\mathcal{T}$ ($x$-axis) for the identification of (a) the KdV equation \eqref{e: pde kdv} and (b) the KS equation \eqref{e: pde KS}.  Different color curves represent results for various noise-to-signal ratios $\sigma_{\rm NSR}\in \{ 0,0.1,...,1\}$. Notice a  wide range  of $\mathcal{T}$ gives the same recovery.  }
\label{f: effect of importance score}
\end{figure}

In Figure \ref{f: effect of importance score}, we present the coefficient $E_2$ error ($y$-sxis) against different values of the trimming parameter $\mathcal{T}$ ($x$-axis) for different noise-to-signal ratios (different color curves)  for (a) the KdV equation \eqref{e: pde kdv} and (b) the KS equation \eqref{e: pde KS}.  In general, we use $\mathcal{T} = 0.05$ as a default for all equations in Table \ref{T:PDE} and Table \ref{T: odes}, except for the KS equation \eqref{e: pde KS} and the PM equation \eqref{e: pde PM} for which we use $\mathcal{T} = 0.2$. 
Our experiments use the same distribution of seeds {for the noise} with different variances.  Different color curves represent the different values of noise-to-signal ratio $\sigma_{\rm NSR}\in \{ 0,0.1,...,1\}$.  For example, when there is no noise,  $\sigma_{\rm NSR} = 0$ (the lowest blue curve), it gives the lowest recovery error (compared to other colored curves) over the widest range of allowable $\mathcal{T}$.
There is a wide range of $\mathcal{T}$ that yields the same recovery.  We use  $\mathcal{T}=0.2$ for the KS and PM equations, by choosing a value of $\mathcal{T}$ from a large plateau.  This makes the algorithm more robust.  In general, since the colored curves are decreasing functions in terms of $\mathcal{T}$, if the given data is highly corrupted by noise, using a larger $\mathcal{T}$ can help with the identification. 

\subsection{Effects of subsampling in data acquisition and the feature matrix $W$} \label{sss:subsample}

In Figure \ref{fig: effect of T, dt, dx}, we show  the effects of changing the final time  $T$ (the top row), and of changing $\Delta x$ and $\Delta t$ for $\mathbb{N}_x\mathbb{N}_t$ (the second row), and of changing the uniform subsampling in (\ref{eq:H}), i.e., $\Delta t^*$ and $\Delta x^*$  for the generation of the feature matrix (the third row).
We compare for the KS equation \eqref{e: pde KS}, the 2D linear ODE system \eqref{e: ode, linear 2d},  the Van der Pol equation \eqref{e: ode, van-der-pol}, and the Duffing equation \eqref{e: ode, Duffing} to illustrate the effects. The noise level is $\sigma_{\rm NSR}=0.1$ for each example. 
We present the average of the $E_2$ error , the TPR and PPV values from 20 independent experiments for one varying variable among the variables $\{ T, \Delta t, \Delta x, \Delta t^*, \Delta x^* \}$ while fixing the rest. 
The first row shows that the recovery by WeakIdent is robust as long as $T$ is above a sufficiently large value (e.g. 100 or 10), which indicates that  there is a time $T$ such that the solution of the differential equations contains enough dynamics up to time $T$.  The second row shows that WeakIdent gives  a smaller error with smaller $\Delta x$ and $\Delta t$.
The bottom row shows that the  size of uniform subsampling in space and time of the feature matrix does not affect the recovery. 

\begin{figure}[h]
    \centering
    \begin{tabular}{cc|ccc}
    \toprule
         \multicolumn{2}{c|}{ KS\eqref{e: pde KS}}
        & 2D Linear system.\eqref{e: ode, linear 2d} & Duffing \eqref{e: ode, Duffing} &  Lotka-Volterra \eqref{e: ode, Lotka} 
        \\
    \midrule
    \multicolumn{2}{c|}{Effect of $T$ } & \multicolumn{3}{c}{Effect of $T$ } \\
     \multicolumn{2}{c|}{ \includegraphics[width = 0.17\textwidth, align = c]{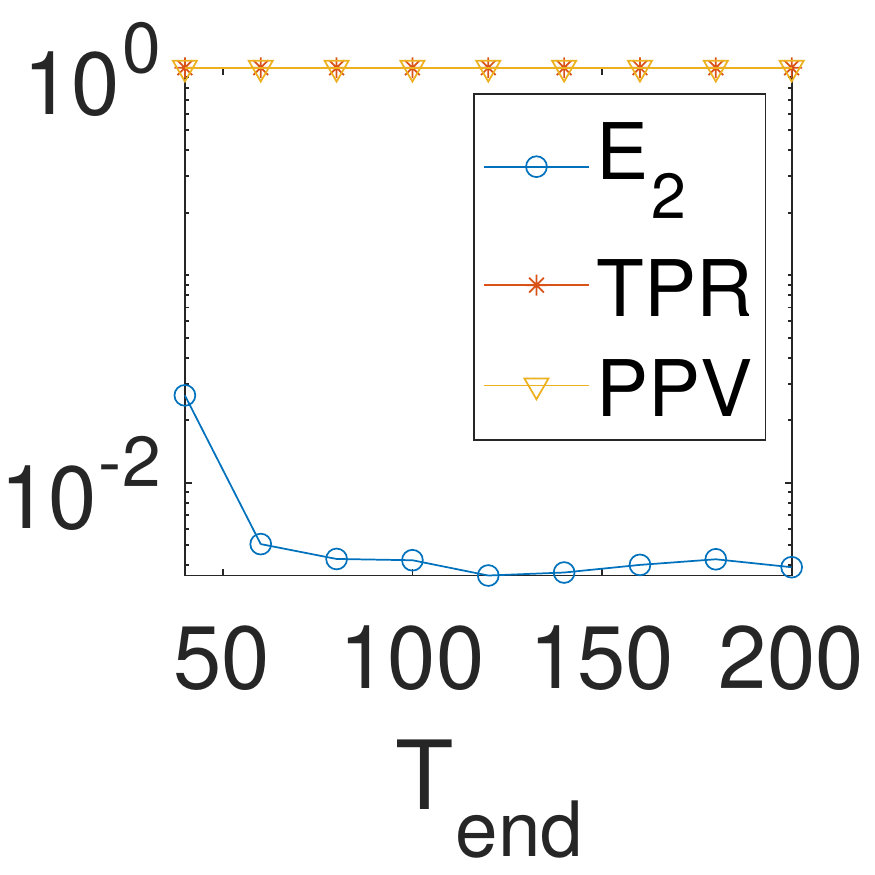}}
    & \includegraphics[width = 0.17\textwidth, align = c]{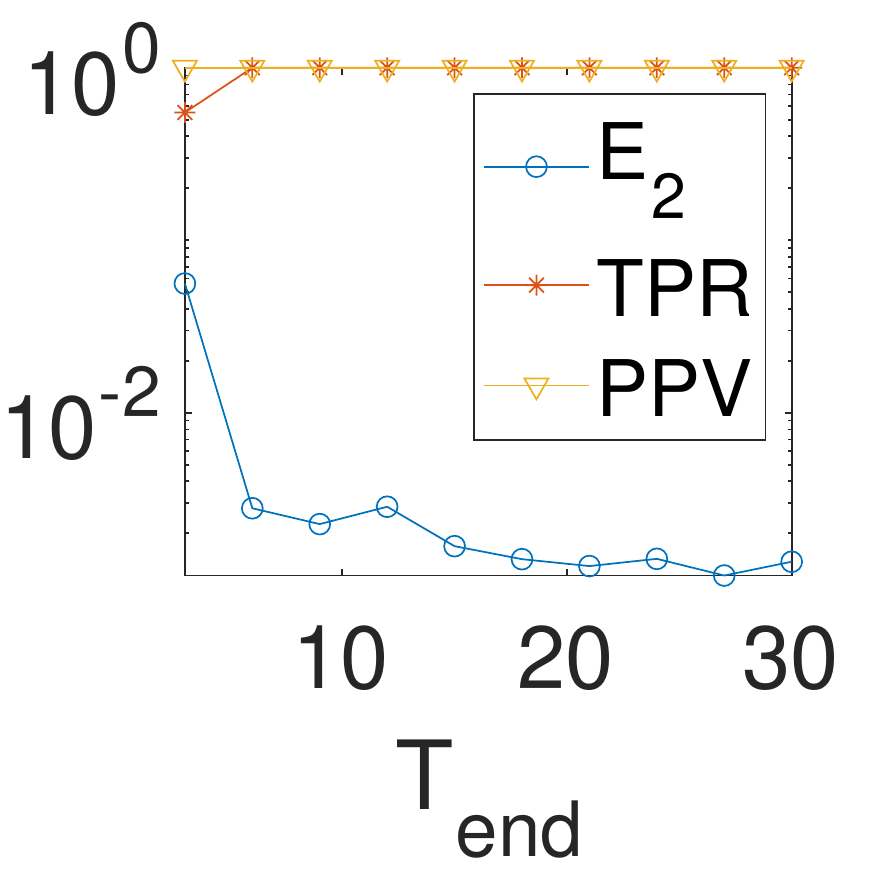}
    & \includegraphics[width = 0.17\textwidth, align = c]{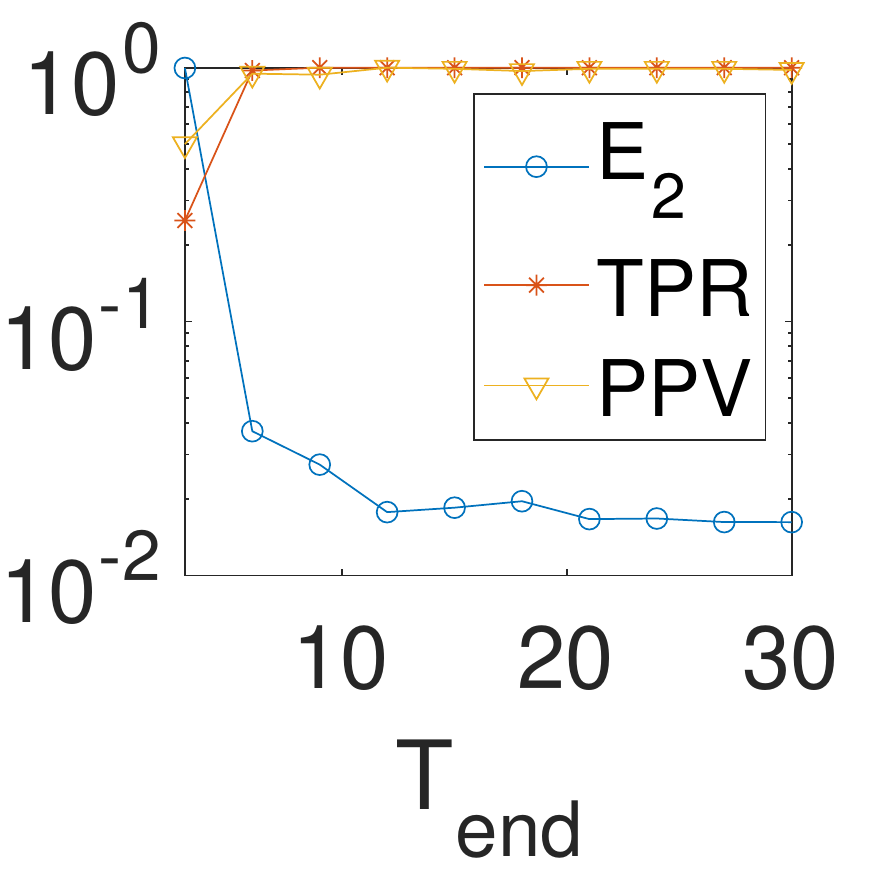}    &  \includegraphics[width = 0.17\textwidth, align = c]{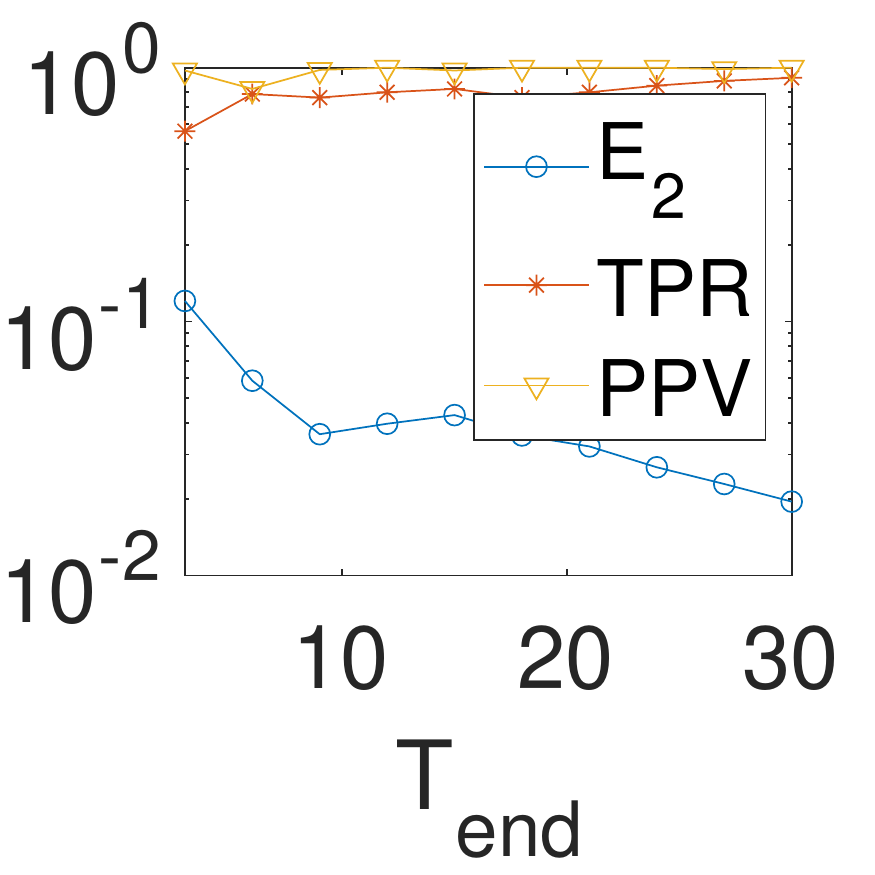} 
    \\
    \midrule
    \multicolumn{1}{c}{Effect of $\Delta t$} & Effect of $\Delta x$ 
    & \multicolumn{3}{c}{Effect of $\Delta t$}\\
    \midrule
    \includegraphics[width = 0.17\textwidth, align = c]{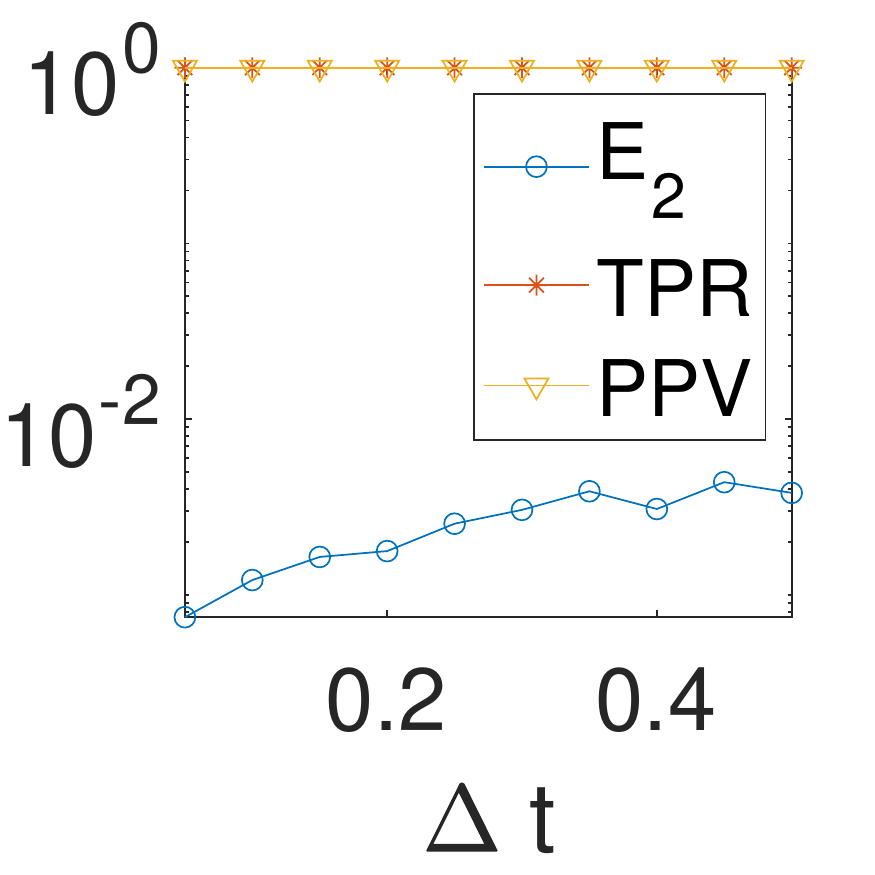} 
    & \includegraphics[width = 0.17\textwidth, align = c]{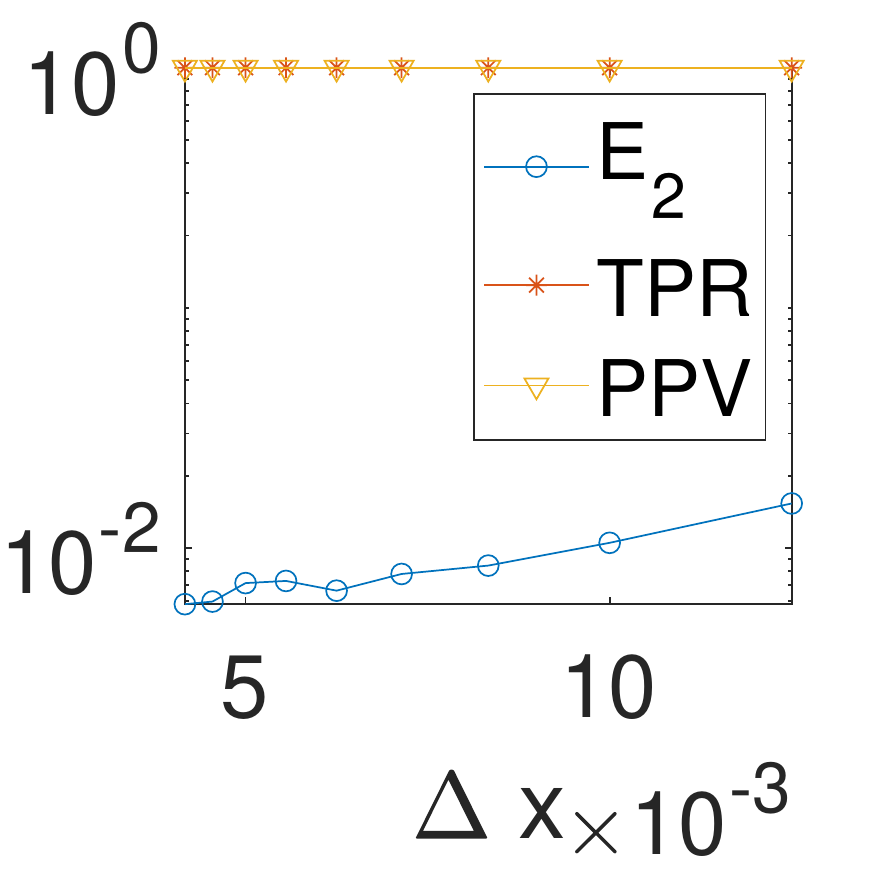}
    & 
    \includegraphics[width = 0.17\textwidth, align = c]{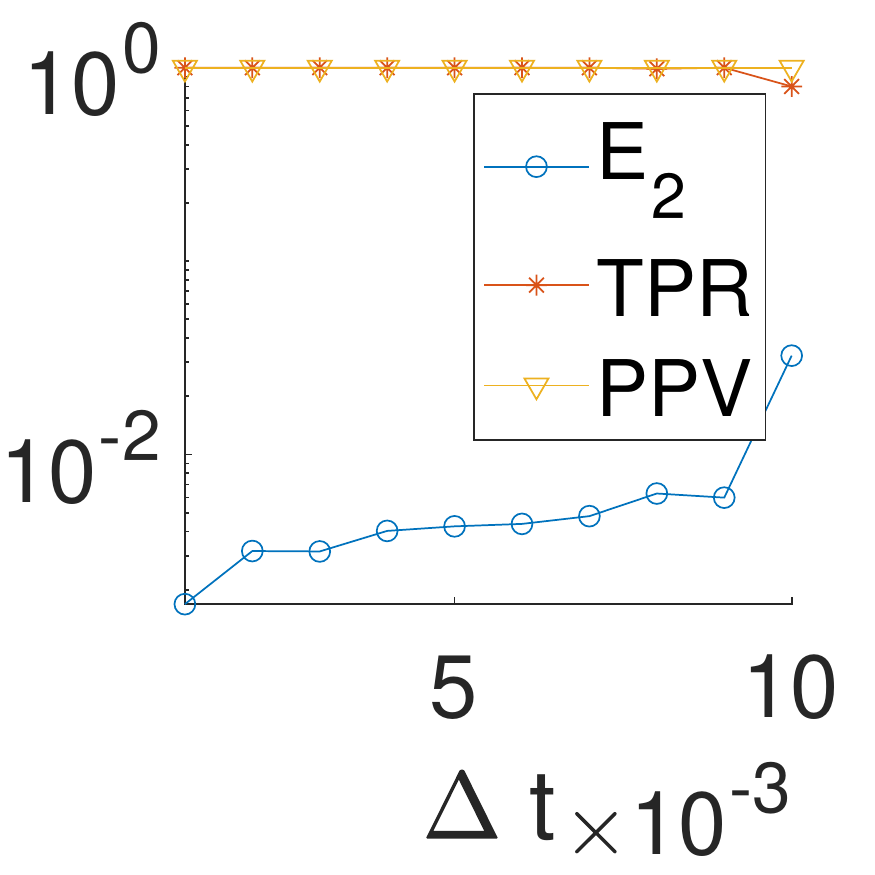} &
    \includegraphics[width = 0.17\textwidth, align = c]{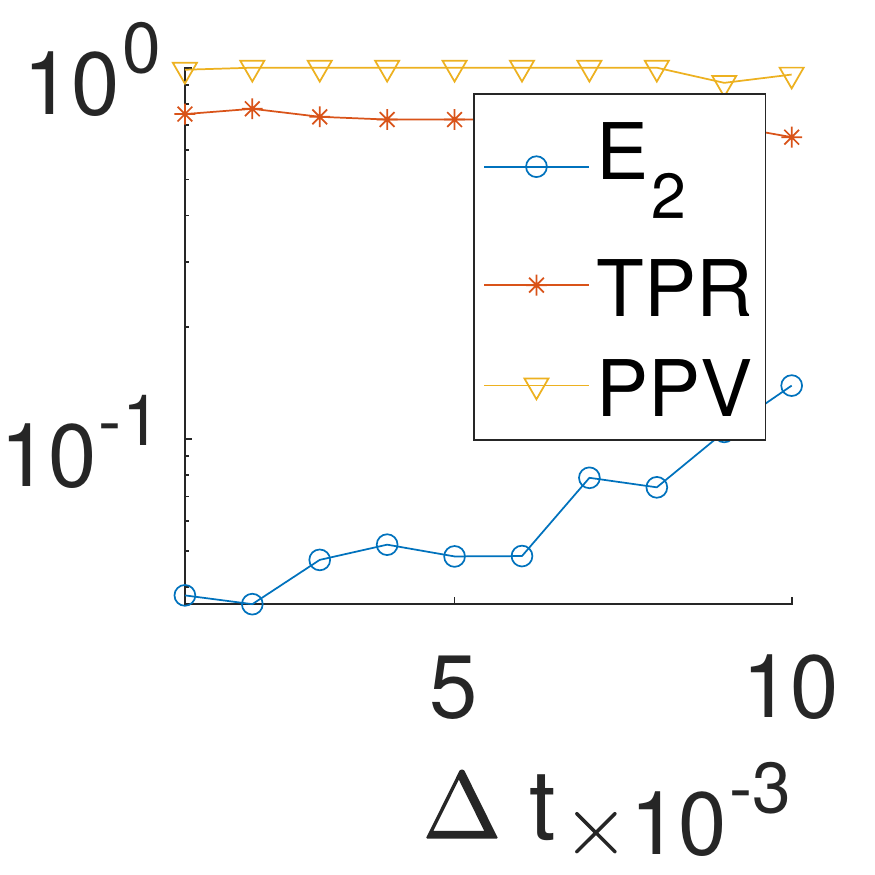} & 
    \includegraphics[width = 0.17\textwidth, align = c]{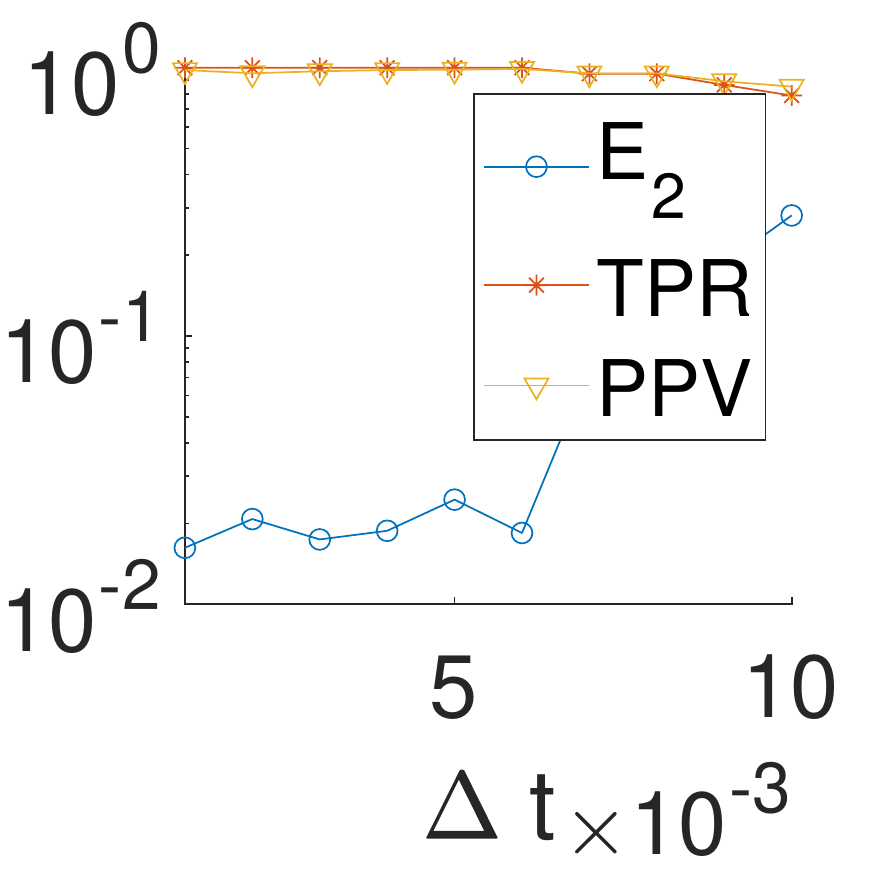}    \\
     \midrule
    \multicolumn{1}{c}{Effect of $\Delta t^*$} &
    \multicolumn{1}{c|}{Effect of $\Delta x^*$}
    &  \multicolumn{3}{c}{Effect of $\Delta t^*$}\\
\midrule
    \includegraphics[width = 0.17\textwidth, align = c]{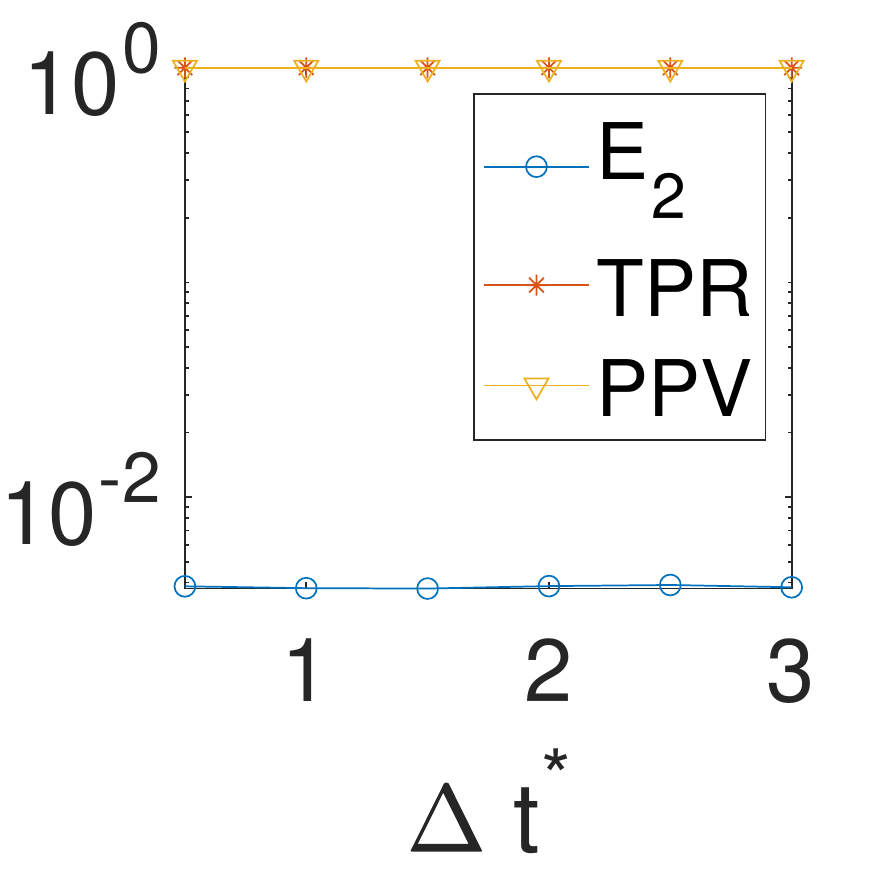}
    & 
    \includegraphics[width = 0.17\textwidth, align = c]{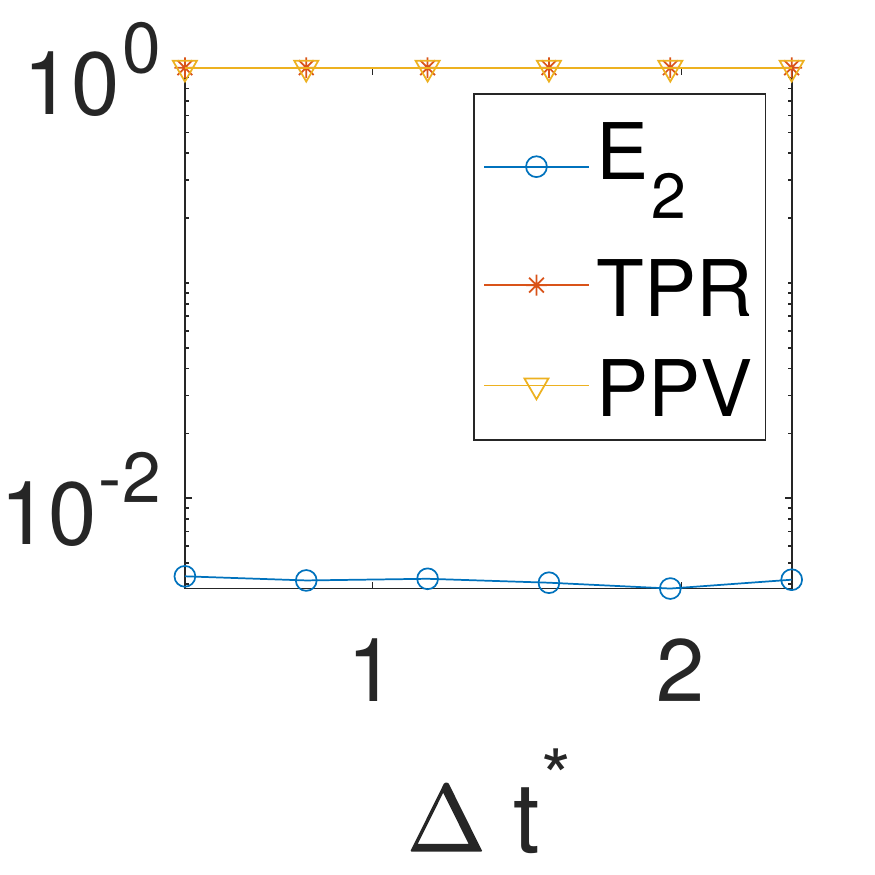} &
    \includegraphics[width = 0.17\textwidth, align = c]{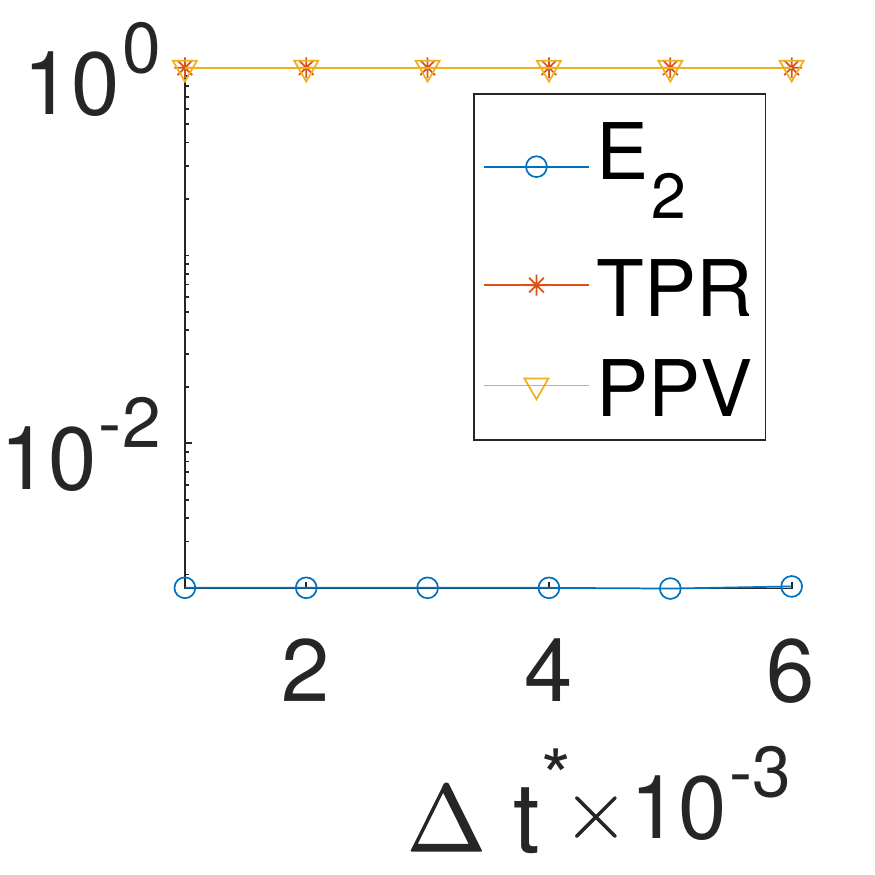} &
    \includegraphics[width = 0.17\textwidth, align = c]{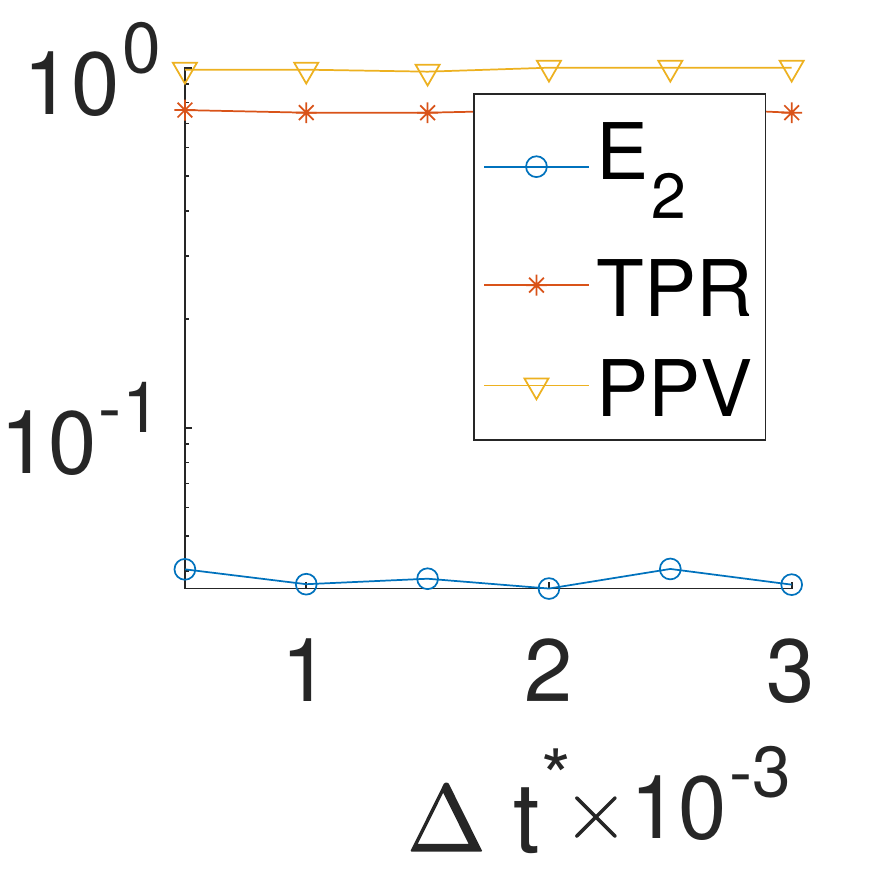} & 
    \includegraphics[width = 0.17\textwidth, align = c]{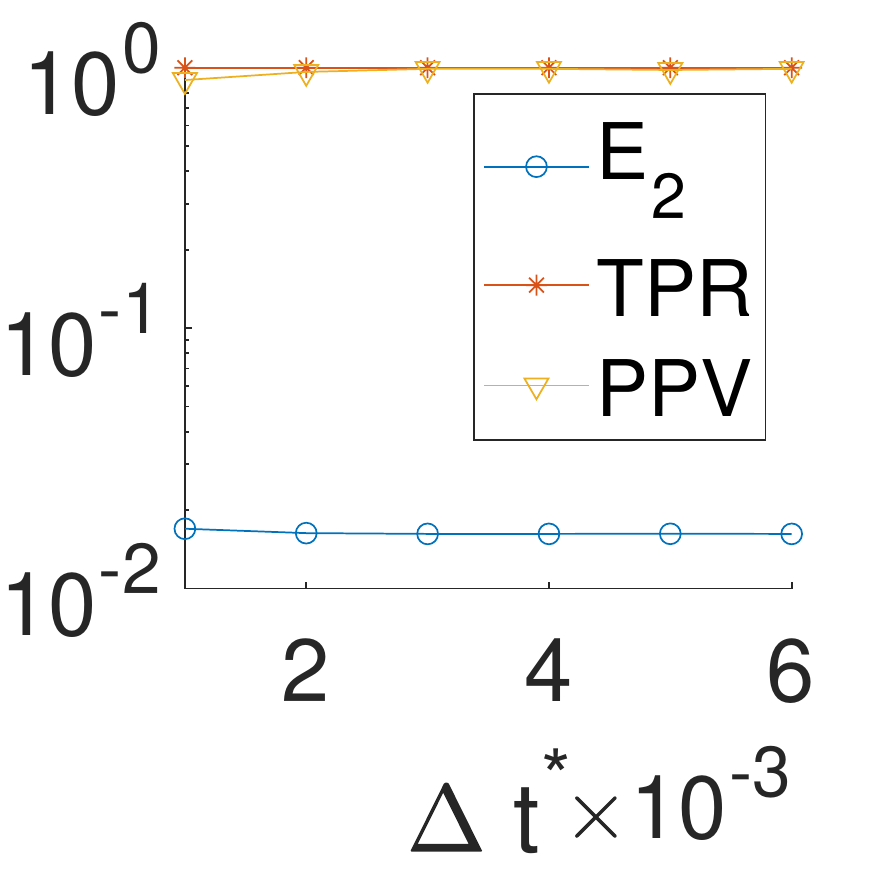}   \\   
    \bottomrule
    \end{tabular}
    \caption{Effects of the final time $T$ (the top row), $\Delta t, \Delta x$ for $\mathbb{N}_x\mathbb{N}_t$ of the given data (the second row), and subsampling $\Delta t^*, \Delta x^*$ in (\ref{eq:H}) in the third row.  Each graph shows the average of the $E_2$ error, the TPR and PPV values over 20 experiments for one varying variable among the variables in $\{ T, \Delta x, \Delta t, \Delta x^*, \Delta t^* \}$ while the rest is fixed. The noise level is $\sigma_{\rm NSR} = 0.1$. The left column gives  the PDE results for  the KS equation while  both $\Delta t, \Delta x$ are shown. The right columns  show  $\Delta t$ only for  ODEs, including the 2D Linear system \eqref{e: ode, linear 2d},  the Duffing equation \eqref{e: ode, Duffing} and the Lotka-Volterra equation\eqref{e: ode, Lotka}. There is a transition point in $T$ such that the given data up to $T$ contain enough dynamics. The recovery is in general better with smaller $\Delta t$ and $\Delta x$, and the rate of uniform subsampling has a minimal effect on the results.  }
    \label{fig: effect of T, dt, dx}
\end{figure}

In Table \ref{t: size of given data and feature matrix}, we show an example of the size reduction from $\mathbf{W}$ to $\mathbf{W}_{\rm narrow}$for the PDEs and ODEs considered in this paper. We use $\sigma_{\rm NSR}=0.1$ for the RD equation \eqref{e: pde RD} and $\sigma_{\rm NSR}=0.2$ for the rest of the equations.  The given data is of size $\mathbb{N}_x\mathbb{N}_t$ and it is subsampled to $H=N_xN_t$ number of rows for $\boldsymbol{W}$.  The narrow-fit further reduces the feature matrix to $\tilde{\mathbf{W}}_{\rm narrow}$ for computational accuracy. 
 \begin{table}[!]
     \centering
     \begin{tabular}{lrrrrrrr}
       \toprule
       Equation   & $\mathbb{N}_x$ & $\mathbb{N}_t$ & ${N}_x$ & ${N}_t$
       & $\boldsymbol{W} \text{ size} (H \times L)$ & $\tilde{\boldsymbol{W}}_{\rm narrow} \text{ size}$
       \\
       \midrule
       Linear Equ.\eqref{e: pde transport} & 257 & 300 & 36 & 39 & 1404 $\times$43  & 824 $\times$43\\
       KdV Equ.\eqref{e: pde kdv}& 400 & 601 & 71 & 65 &   4615 $\times$43 & 1367 $\times$43\\  
       KS Equ.\eqref{e: pde KS}& 256 & 301 & 46 & 43 & 1935$\times$43 & 916$\times$43\\
       NLS Equ.\eqref{e: pde NLS}& 256 & 251 & 39 & 42 & 1225$\times$190 & 159$\times$190\\
       PM Equ.\eqref{e: pde PM}& 200$\times$200 & 128 & 14$\times$14 & 16 & 3136$\times$65 & 1349$\times$65  \\
       RD Equ.\eqref{e: pde RD} & 256$\times$256 & 201 & 13$\times$13 & 14 & 2366$\times$155 &  2271$\times$155\\
       \midrule
       Linear Equ.\eqref{e: ode, linear 2d} & - & 1001 & - & 851  & 877$\times$21 & 127$\times$21\\
       VdP Equ.\eqref{e: ode, van-der-pol} & -  & 15001 & - &  958 & 958$\times$21 & 295$\times$21 \\
       Duffing Equ.\eqref{e: ode, Duffing} & - & 1001 & - &  915 & 915$\times$21 &  57$\times$21\\ 
       LV Equ.\eqref{e: ode, Lotka} & - & 1001 & - & 947 & 947$\times$10 & 338$\times$10\\
       Lorenz Equ.\eqref{e: ode, lorenz} & - & 15001 & - & 983 & 983$\times$20 & 930$\times$20\\
       \bottomrule
     \end{tabular}
     \caption{Typical examples of the feature matrix size and the reduction in narrow-fit. The given data is of size $\mathbb{N}_x\mathbb{N}_t$ and it is subsampled to $H=N_xN_t$ rows for $\boldsymbol{W}$.  We use $\sigma_{\rm NSR}=0.1$ for the RD equation \eqref{e: pde RD} and $\sigma_{\rm NSR}=0.2$ for the rest of the equations. For systems of equations, the size of the feature matrix for each dependent variable is identical. }\label{t: size of given data and feature matrix}
 \end{table}

\subsection{Speed of WeakIdent} 
We perform experiments  using Matlab  on  the Apple M1 processor with 8-core CPU and 16GB of RAM.
The computational cost of WeakIdent is typically about 1-5 seconds for an ODE system or  a PDE with one dependent variable in  a 1D spatial domain. For example, the cpu times to recover the Lotka-Volterra system \eqref{e: ode, Lotka} and the KdV equation \eqref{e: pde kdv} are 1.11  and 0.63 seconds, respectively.  For the cases in  2D spatial domains,  such as the anisotropic PM equation \eqref{e: pde PM} with one variable, and the 2D reaction-diffusion equation  \eqref{e: pde RD} with two variables, the recovery can take about 3 and 35 seconds, respectively.  The speed is comparable with WPDE\cite{messenger2021weakPDE}, which takes 16 and 75 seconds for these two examples. 

\bigskip
In Appendix \ref{sec: appendix}, we present additional results and more comparisons. The additional results for PDEs are in Subsection \ref{AS:pde} and additional results for ODEs are in Subsection \ref{AS:ode}. Details about how to construct test functions are given in  Appendix  \ref{Asec:testfunction}.

\section{Conclusion and discussion}
\label{sec: conclusion}

We propose a new method WeakIdent for identifying both PDEs and ODE systems from noisy data using  a weak formulation.
The proposed WeakIdent does not require  prior knowledge of the governing features, but uses all features up to certain polynomial  order, and up to certain order of derivatives. 
We first use Subspace Pursuit to find a candidate support, then propose two novel techniques called narrow-fit and trimming to improve both the support identification and the coefficient recovery.  A careful design of the test functions helps with the recovery, and a proper normalization of the columns in the feature matrix improves the results in the implementation of least-squares. 
The proposed WeakIdent requires at most $L$ sparsity iterations (or including the sub-iteration of narrow-fit and trimming, at most $\frac{L^2}{2}$  iterations),  where $L$ is the number of features.  At the same time the trimming step  improves the recovery and gives good results after a fraction of $L$ is used  to identify the correct support, as shown in Figure \ref{F:cScore}.  Narrow-fit  based on  highly dynamic regions  also makes the computation more efficient, and with error normalization  of the feature matrix, the coefficient recovery is improved. 
Comprehensive numerical experiments on various equations/systems are provided, showing the robust performance  of WeakIdent compared to other  state-of-the-art  methods. The  Weak form in general is effective when the noise level is high,  At the same time, to take advantage of the weak form, the possible features in the differential equation must be in a specific form for the integration of parts.

\bibliographystyle{plain} 
\bibliography{cite_Wident} 

\appendix

\section{Proof of Theorem \ref{th:error}}
\label{app:proof:th:error}

\begin{proof}
(a)   Using the noisy data in the form ${\hat{U}}_{i}^n = {U}_{i}^n  + {\epsilon}_{i}^n$ in 
(\ref{e: noisy discretized data}), the $h$th entry of $\boldsymbol{e}^{\rm noise}$ can be expressed as \begin{align}
{e}^{\rm noise}_h
= & \;\; \Delta x \Delta t \sum_{(x_j,t^k) \in \Omega_{h(x_i,t^n)}} \left(
\sum_{l\in \text{Supp}^*}(-1)^{\alpha_l}{c_l}
\left( ({U}_j^k+\epsilon_j^k)^{\beta_l}  - (U_{j}^{k})^{\beta_l}  \right)
\frac{\partial^{\alpha_l} \phi_h }{\partial x^{\alpha_l}}(x_{j}, t^{k})
+ (\hat{U}_{j}^{k} - U_{j}^{k} )\frac{\displaystyle \partial \phi_h}{\displaystyle \partial t}  (x_{j}, t^{k})
\right)\nonumber \\
= & \;\; \Delta x \Delta t \sum_{(x_j,t^k) \in \Omega_{h(x_i,t^n)}}
    \left(
     \sum_{l\in \text{Supp}^*}(-1)^{\alpha_l}{c_l}
     \epsilon_{j}^{k}
     \left( \sum_{r=1}^{\beta_l} {\beta_l \choose r} ( \epsilon_{j}^{k})^{r-1} (U_{j}^{k})^{\beta_l -r} \right)
    \frac{\partial^{\alpha_l} \phi_h }{\partial x^{\alpha_l}}(x_{j}, t^{k})
    + \epsilon_{j}^{k}\frac{\displaystyle \partial \phi_h}{\displaystyle \partial t}  (x_{j}, t^{k})
    \right) \nonumber\\
= & \;\; \Delta x \Delta t \sum_{(x_j,t^k) \in \Omega_{h(x_i,t^n)}}
    \left(
     \sum_{l\in \text{Supp}^*}(-1)^{\alpha_l}{c_l}
     {\beta_l}(U_j^k)^{ \beta_l - 1}
    \frac{\partial^{\alpha_l} \phi_h }{\partial x^{\alpha_l}}(x_{j}, t^{k})
    + \frac{\displaystyle \partial \phi_h}{\displaystyle \partial t}  (x_{j}, t^{k}) 
    \right) \epsilon_j^k
    + \mathcal{O}\left( (\epsilon_j^k)^2\right).\label{eqenoiseh3}
\end{align}
Hence, \begin{align*}
\|\boldsymbol{e}^{\rm noise}\|_\infty & = \max_h |{e}^{\rm noise}_h|
 \le\max_h\left[\epsilon \bar{S}^*_h |\Omega_h|  
\right] + \mathcal{O}\left( \epsilon^2\right)  \end{align*}
where
\[\bar{S}^*_h =\sup_{(x_j,t^k) \in \Omega_h}
\bigg|
\sum_{l\in \text{Supp}^*}(-1)^{\alpha_l}{c_l}
{\beta_l}(U_j^k)^{ \beta_l - 1}
\frac{\partial^{\alpha_l} \phi_h }{\partial x^{\alpha_l}}(x_j, t^k)
-  \frac{\displaystyle \partial \phi_h}{\displaystyle \partial t}  (x_j, t^k) 
\bigg|  .
\]
Setting $\bar{S}^* =\max_h \bar{S}^*_h$ as in \eqref{eq:sbarstar} gives rise to our estimate in \eqref{enoiseestimate}.

(b) From \eqref{eqenoiseh3}, the leading term in $e_h^{\rm noise}$ is 
\begin{align*}
\Delta x \Delta t \sum_{(x_j,t^k) \in \Omega_{h(x_i,t^n)}}
    \left(
     \sum_{l\in \text{Supp}^*}(-1)^{\alpha_l}{c_l}
     {\beta_l}(U_j^k)^{ \beta_l - 1}
    \frac{\partial^{\alpha_l} \phi_h }{\partial x^{\alpha_l}}(x_{j}, t^{k})
    + \frac{\displaystyle \partial \phi_h}{\displaystyle \partial t}  (x_{j}, t^{k}) 
    \right) \epsilon_j^k.
\end{align*}
Based on our noise assumption, this leading term $e_h^{\rm noise}$ has mean $0$, and variance $\sigma^2 S_h^*$
with $S_h^*$ given in \eqref{eqshstart}.

\end{proof}

\section{Additional results and comparisons}
\label{sec: appendix}

\subsection{Additional results and comparisons for PDEs}
\label{AS:pde}

\begin{figure}[h]
\begin{center}
\begin{tabular}{cc}
 (a) $\hat{U}(\boldsymbol{x},t)$ \\
\includegraphics[width = 0.3 \textwidth]{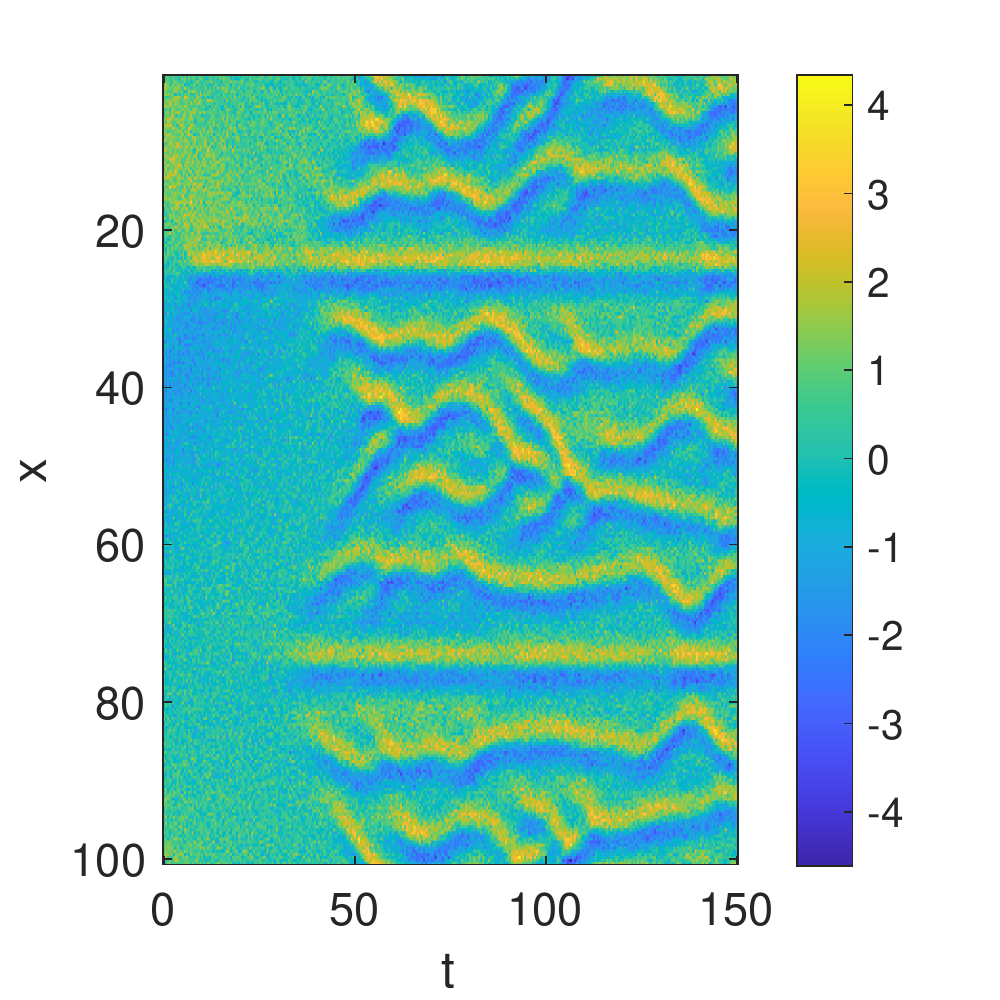} 
\end{tabular}
\begin{tabular}{|l|p{10cm}|c|}
\multicolumn{3}{c}{(b) $\sigma_{\rm NSR} =0.5$ } \\
\toprule
True equation & 
$u_t  =  -1.00000 u_{xx} -1.00000 u_{xxxx} -0.50000 (u^2)_{x}$ & \\  
\hline
\textbf{WeakIdent}&$u_t  =  -0.91387 u_{xx} -0.90906 u_{xxxx} -0.45686 (u^2)_{x}$ & $E_2 =$ \textbf{0.08831} \\
\hline
WPDE\cite{messenger2021weakPDE} & $u_t  =  -1.03383 u_{xx} -1.27316 u_{xxxx} -0.51246 (u^2)_{x} $ &  $E_2 = $ 0.25635  \\  
& \hspace{1cm} $-0.24197 u_{xxxxxx}  -0.11572 (u^2)_{xxx}$
&\\
\hline
RGG \cite{reinbold2020using}& $ u_t  =  0.90654 u_{xx} +1.18313 u_{xxxx} +0.58744 (u^2)_x  $ & $E_2 = $ 0.33155\\
& \hspace{0.5cm} $+0.00906 u -0.02481 u_x +0.34133 u_{xxx} +0.00175 u^2 -0.00952 u^3  $&
 \\
\bottomrule
\end{tabular}\\
\end{center}
\footnotesize{For RGG \cite{reinbold2020using},  8 default features $\{ uu_x, u_{xx}, u_{xxxx}, u, u_{x}, u_{xxx}, u^2, u^3 \}$ and parameters $p_x = 4, p_t = 3, N_d = 100, D = (40,20)$ are used, as in the case for transport equation \eqref{e: pde transport} in  Figure \ref{f: transport equation}. }
\caption{KS equation \eqref{e: pde KS} with $\sigma_{\rm NSR}=0.5$.  (a) Given noisy data $\hat{U}(\boldsymbol{x},t)$. (b) The identified equations using WeakIdent, WPDE\cite{messenger2021weakPDE} and RGG \cite{reinbold2020using} where the $E_2$ error is given in the right column.}
\label{f: recovered equations - KS equation}
\end{figure}

In Figure \ref{f: recovered equations - KS equation}, we experiment on the KS equation (\ref{e: pde KS}) with $\sigma_{\rm NSR}=0.5$.  We compare WeakIdent with WPDE and RGG.  Figure \ref{f: recovered equations - KS equation} (a) shows noisy data $\hat{U}(x,t)$ and (b) gives the recovered equation with the  $E_2$ error.  WeakIdent finds correct support with a small error $E_2=0.08831$.

\begin{figure}
\begin{center}
\begin{tabular}{cc}
 (a) $\hat{U}(\boldsymbol{x},t)$ & (b) $\hat{V}(\boldsymbol{x},t)$ \\
\includegraphics[width = 0.25 \textwidth]{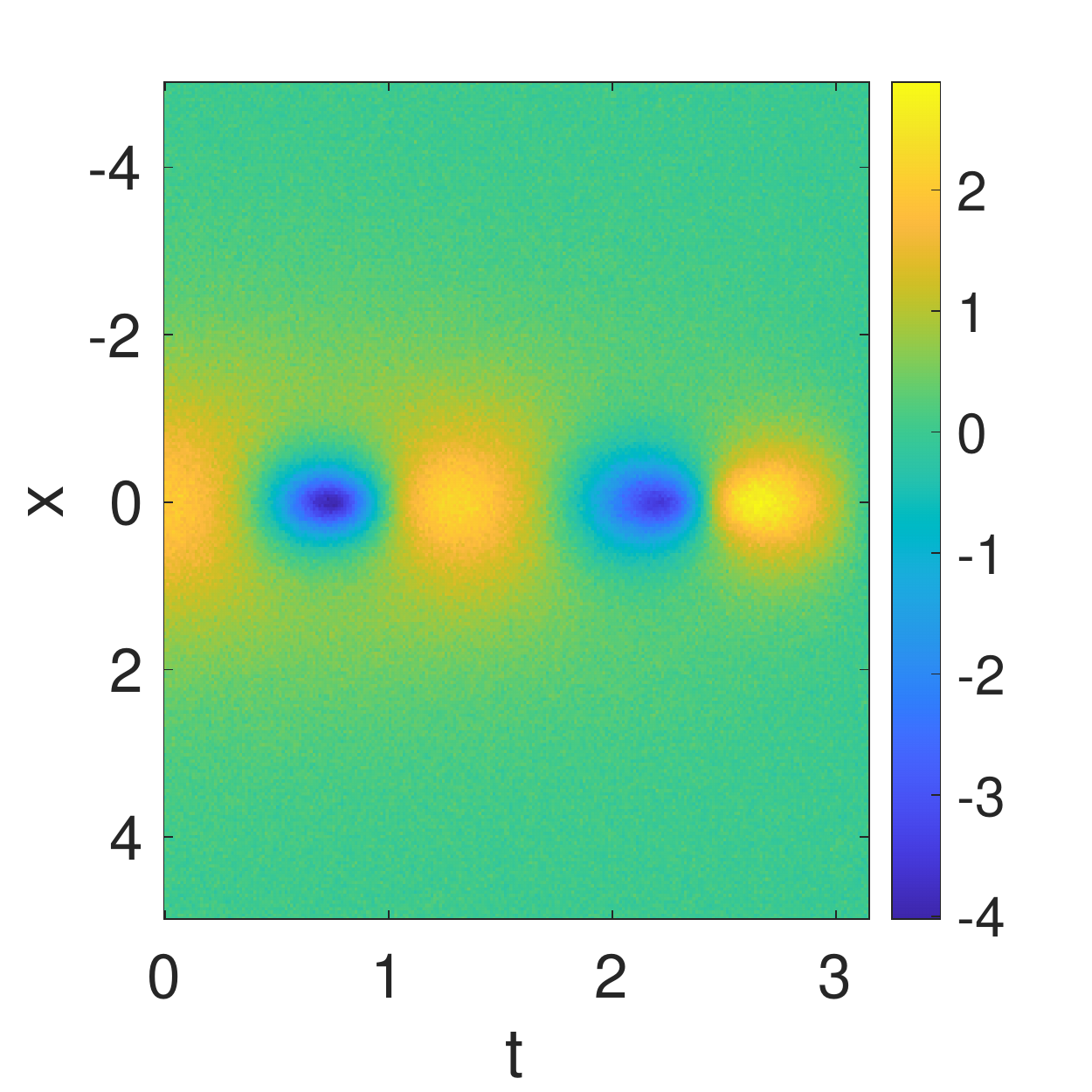} & 
\includegraphics[width = 0.25 \textwidth]{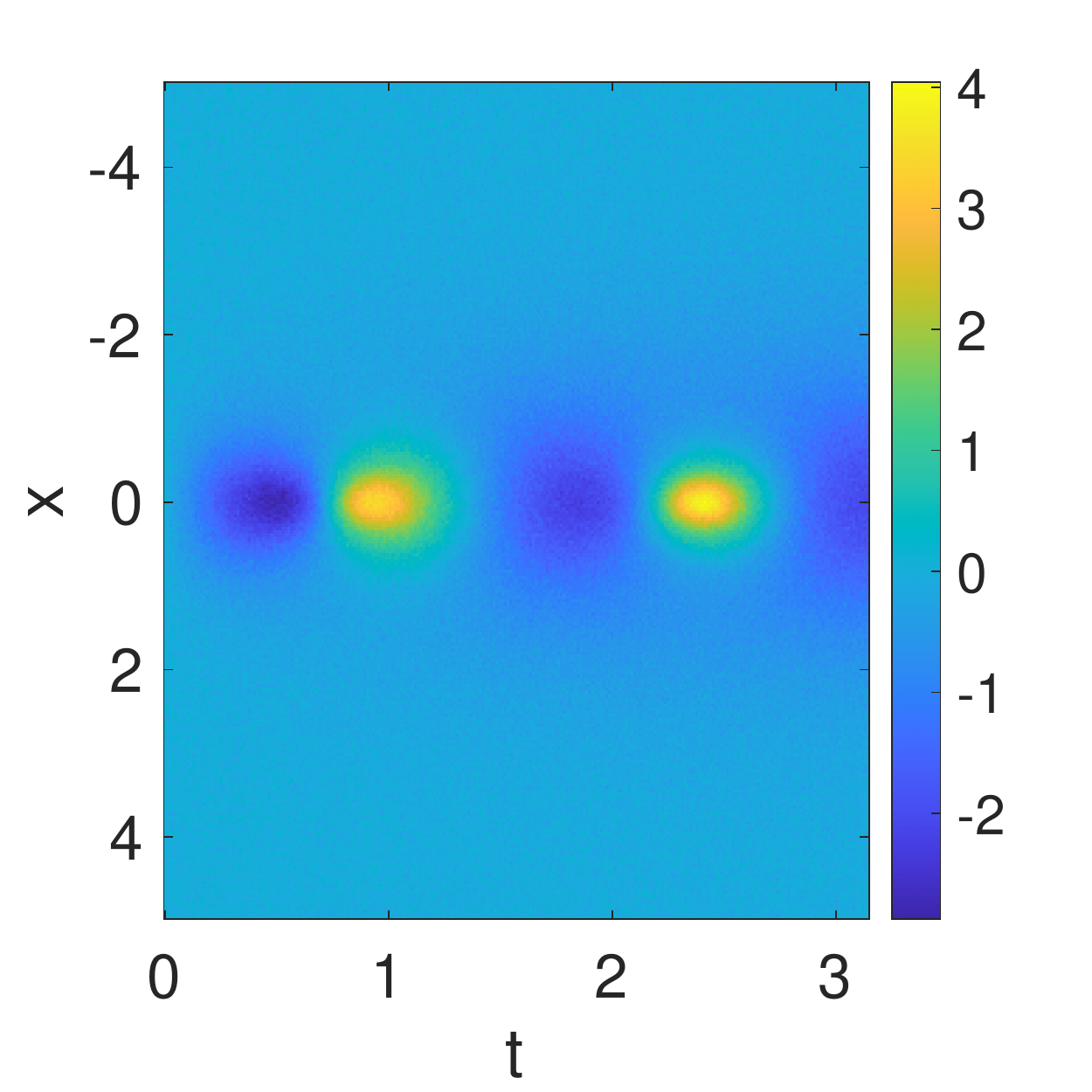} 
\end{tabular}
\begin{tabular}{|l|p{10cm}|c|}
\multicolumn{3}{c}{(c) $\sigma_{\rm NSR} =0$ } \\
\toprule
True equation & 
$u_t  =  +0.50000 v_{xx} +1.00000 v^3 +1.00000 u^2v$ & \\  
& $v_t  =  -0.50000 u_{xx} -1.00000 uv^2 -1.00000 u^3$ & \\
\hline
\textbf{WeakIdent}&$u_t  =  +0.50000 v_{xx} +1.00000 v^3 +1.00000 u^2v$ & $E_2 = \textbf{9.4887e-08}$\\  
& $v_t  =  -0.50000 u_{xx} -1.00000 uv^2 -1.00000 u^3$ & \\
\hline
WPDE\cite{messenger2021weakPDE} & $u_t  =  +0.50000 v_{xx} +1.00000 v^3 +1.00000 u^2v$ & $E_2 = 1.3254e-07$\\  
& $v_t  =  -0.50000 u_{xx} -1.00000 uv^2 -1.00000 u^3$ & \\
\hline
RGG \cite{reinbold2020using}&$ u_t  =  +0.49996 v_{xx} +0.99986 v^3 +0.99988 u^2v  $ & $E_2 = 1.2391e-4$\\ & $ u_t  =  -0.50004 u_{xx} -0.99990 u^3 -0.99988 uv^2  $ &
 \\
\midrule
\multicolumn{3}{c}{(d) $\sigma_{\rm NSR} =0.1$ } \\
\toprule
\textbf{WeakIdent}&$u_t  =  +0.49872 v_{xx} +0.98737 v^3 +1.00001 u^2v$&  $E_2 = $\textbf{0.011} \\
& $v_t  =  -0.49977 u_{xx} -1.01550 uv^2 -0.98785 u^3$&\\
\hline
WPDE\cite{messenger2021weakPDE}&$u_t  =  +0.49868 v_{xx} +0.98722 v^3 +1.00068 u^2v$ &  $E_2 = $0.388    \\
& $v_t  =  -0.06604  +0.10805 v^2 +0.14935 v^2_{xx} -0.04445 u_{xxxx} +0.04995 uv_{xx} -0.98964 uv^2 -0.14989 uv^2_{xx} -0.01231 uv^3_{xx} +0.08142 u^2 -0.03497 u^2v^2 +0.06269 u^2v^2_{xx} -0.01519 u^2v^4_{xx} -0.97427 u^3 -0.09657 u^3_{xx} -0.01753 u^3v^2 +0.01367 u^3v^2_{xx}$ & \\
\hline
RGG \cite{reinbold2020using}& $ u_t  =  +0.09242 (u^2)_x +0.13394 u_{xx} +0.59688 u +0.20176 u_x -3.29871 u^2 +2.49062 u^3 -0.42728 (v^2)_x +0.17625 v_{xx} -2.52277 v +0.27970 v_x -11.82022 v^2 +5.53290 v^3 -5.43222 uv +10.31760 u^2v +16.92442 uv^2$ & $E_2 = $17.668 \\
& $ u_t  =  +0.09242 (u^2)_x +0.13394 u_{xx} +0.59688 u +0.20176 u_x -3.29871 u^2 +2.49062 u^3 -0.42728 (v^2)_x +0.17625 v_{xx} -2.52277 v +0.27970 v_x -11.82022 v^2 +5.53290 v^3 -5.43222 uv +10.31760 u^2v +16.92442 uv^2  $ &
 \\
\bottomrule
\end{tabular}\\
\end{center}
\footnotesize{For RGG \cite{reinbold2020using}, we add an additional dictionary to include the correct features. We use a dictionary of 19 features: $\{(u^2)_x, u_{xx}, u_{xxxx}, u, u_x, u_{xxx},  u^2, u^3, (v^2)_x, v_{xx}, v_{xxxx}, v, v_x, v_{xxx},$  $v^2, v^3,uv, u^2v, uv^2\}$.} 
\caption{Nonlinear Schrodinger equation \eqref{e: pde NLS} with two variables. 
The given noisy data $\hat{U}(\boldsymbol{x},t)$ and $\hat{V}(\boldsymbol{x},t)$ are shown in (a) and (b) respectively. Table (c) and (d) show the identified equations using WeakIdent, WPDE and RGG with  $\sigma_{\rm NSR}=0$ and $\sigma_{\rm NSR}=0.1$. }
\label{f: recovered equations - NLS equation}
\end{figure}

In Figure \ref{f: recovered equations - NLS equation}, we show the identification results  for the  nonlinear Schrodinger equation \eqref{e: pde NLS}.   Table (c) shows the results from noise-free data with $\sigma_{\rm NSR}=0$, and Table (d) shows the noisy case with  $\sigma_{\rm NSR}=0.1$. 
Figure \ref{f: recovered equations - NLS equation} (a) and (b) show the noisy data  $\hat{U}(x,t)$ and $\hat{V}(x,t)$, and the tables  (c) and (d) show the identified equations for WeakIdent, WPDE and RGG.  WeakIdent finds the correct support with small errors in both the noise-free and noisy cases. 

\subsection{Additional results and comparisons for ODEs} \label{AS:ode}

\begin{figure}
    \centering
    \begin{tabular}{|c|c|c|c|c|}
    \toprule
   Eqn & $E_2$ & TPR & PPV\\
    \midrule
     \eqref{e: ode, linear 2d} &
    \includegraphics[width = 0.2\textwidth,align = c]{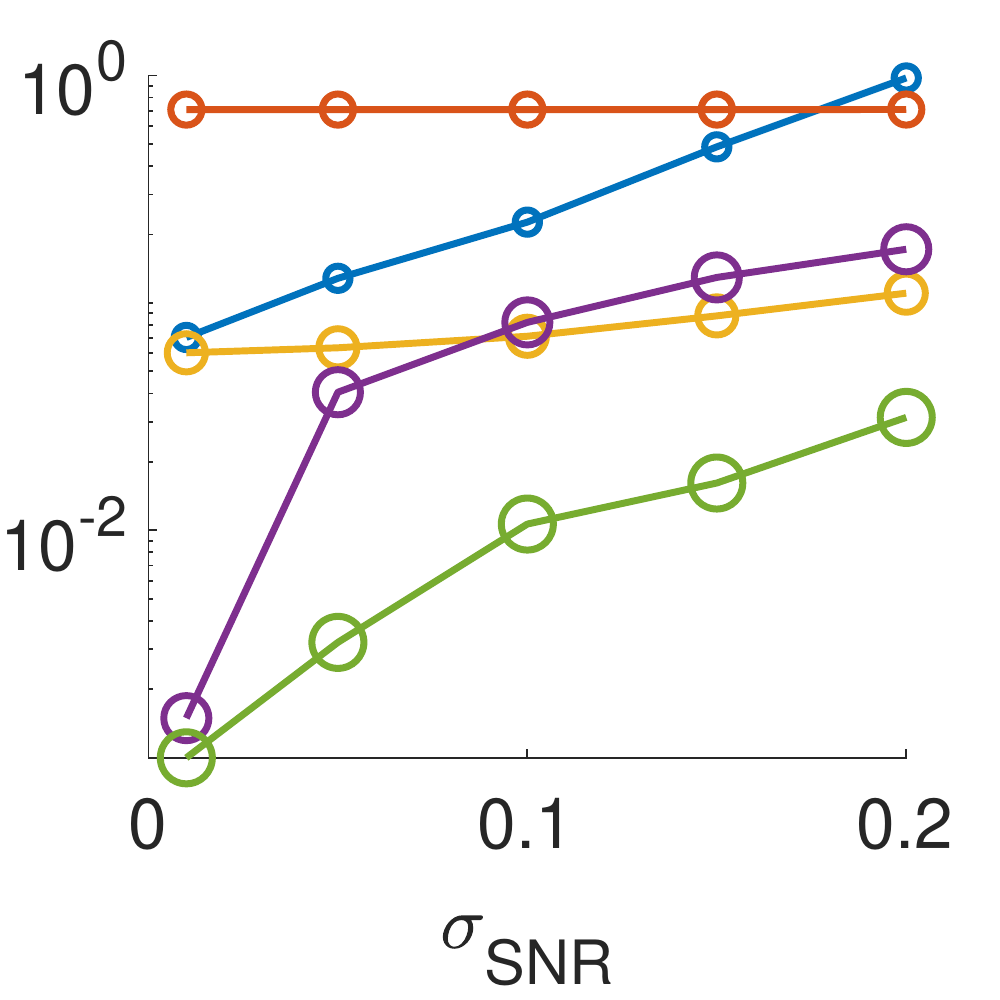} &  \includegraphics[width = 0.2\textwidth,align = c]{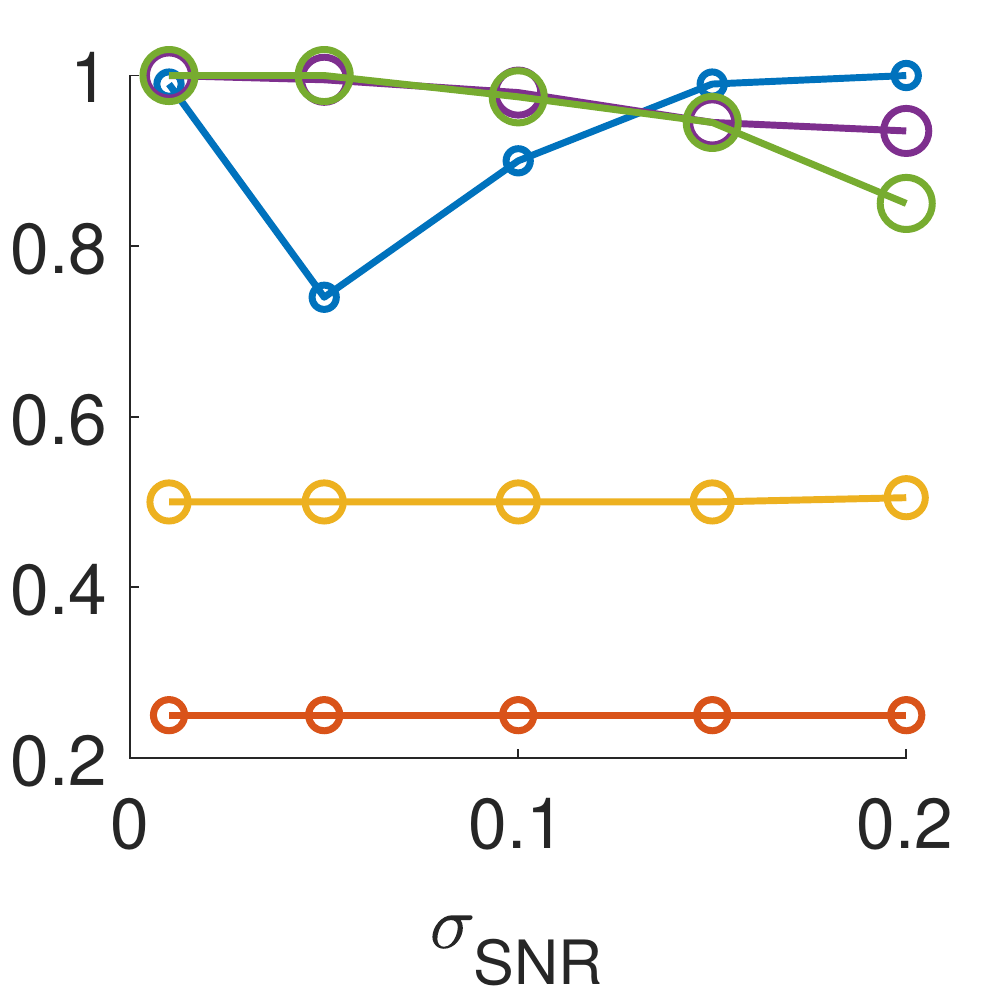}
    &  \includegraphics[width = 0.2\textwidth,align = c]{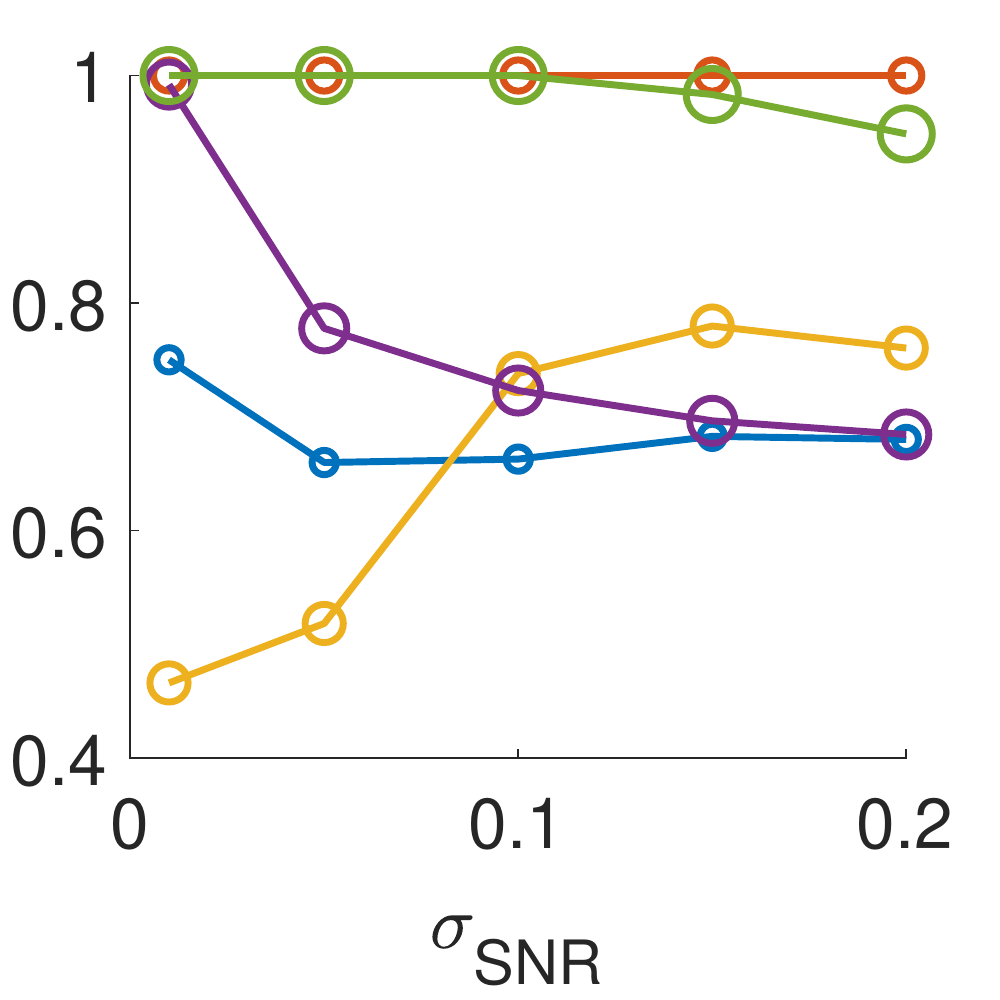} & \includegraphics[width = 0.1\textwidth]{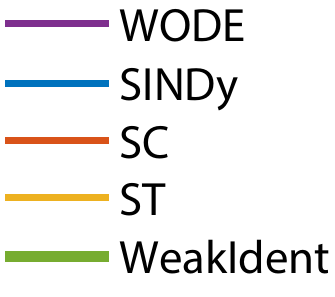}  \\
    \midrule
    \eqref{e: ode, van-der-pol} &
    \includegraphics[width = 0.2\textwidth,align = c]{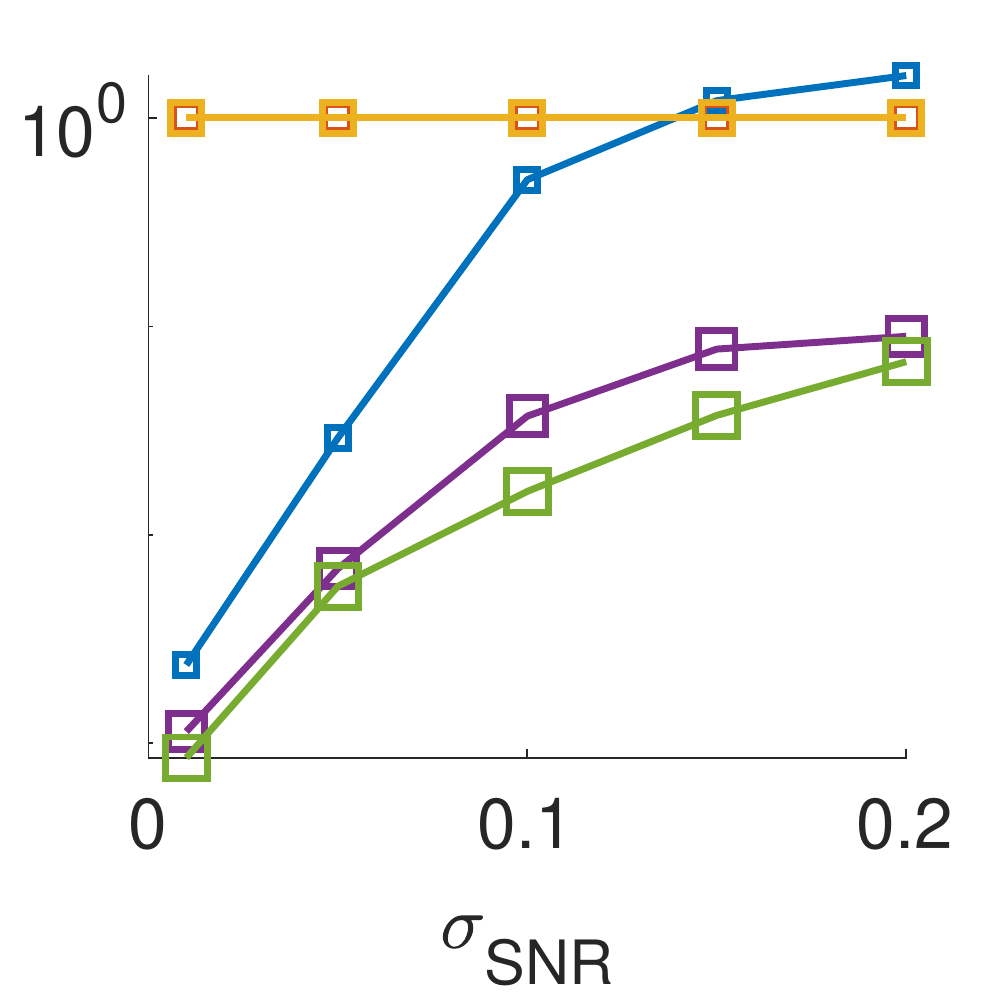} &  \includegraphics[width = 0.2\textwidth,align = c]{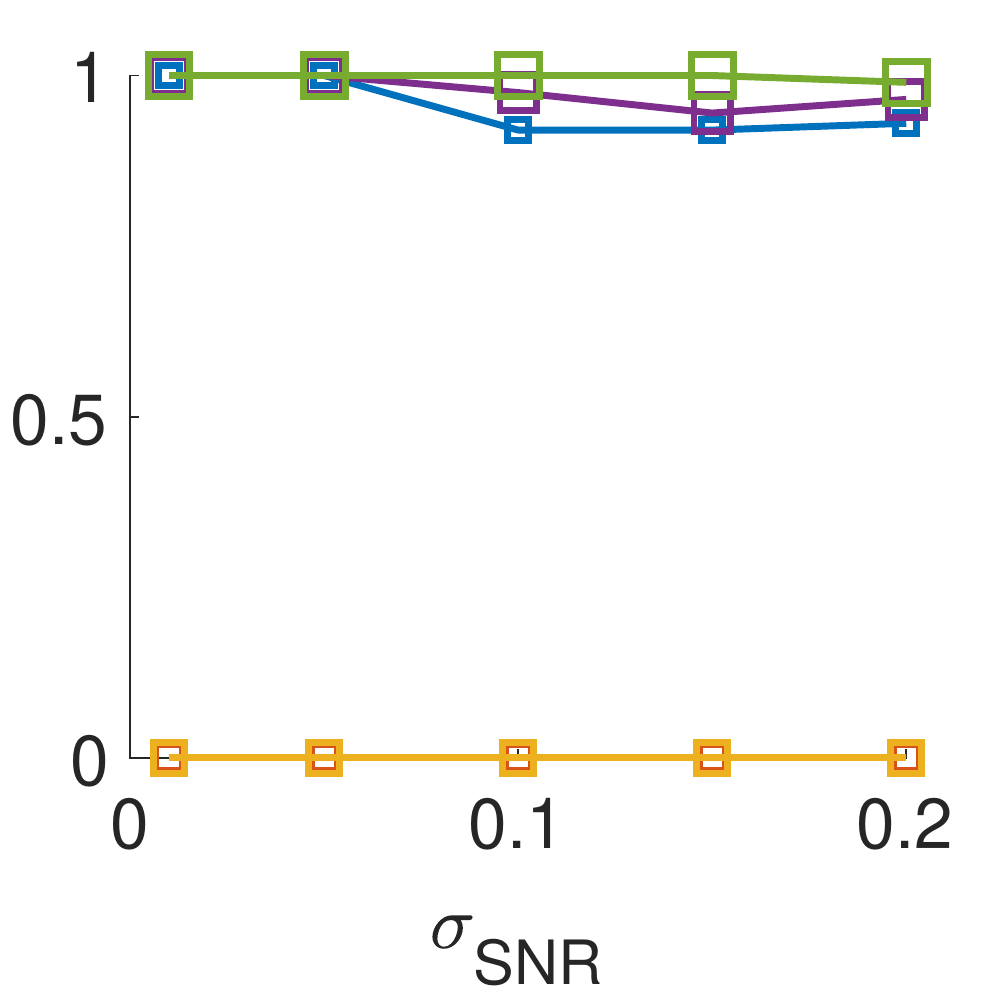}
    &  \includegraphics[width = 0.2\textwidth,align = c]{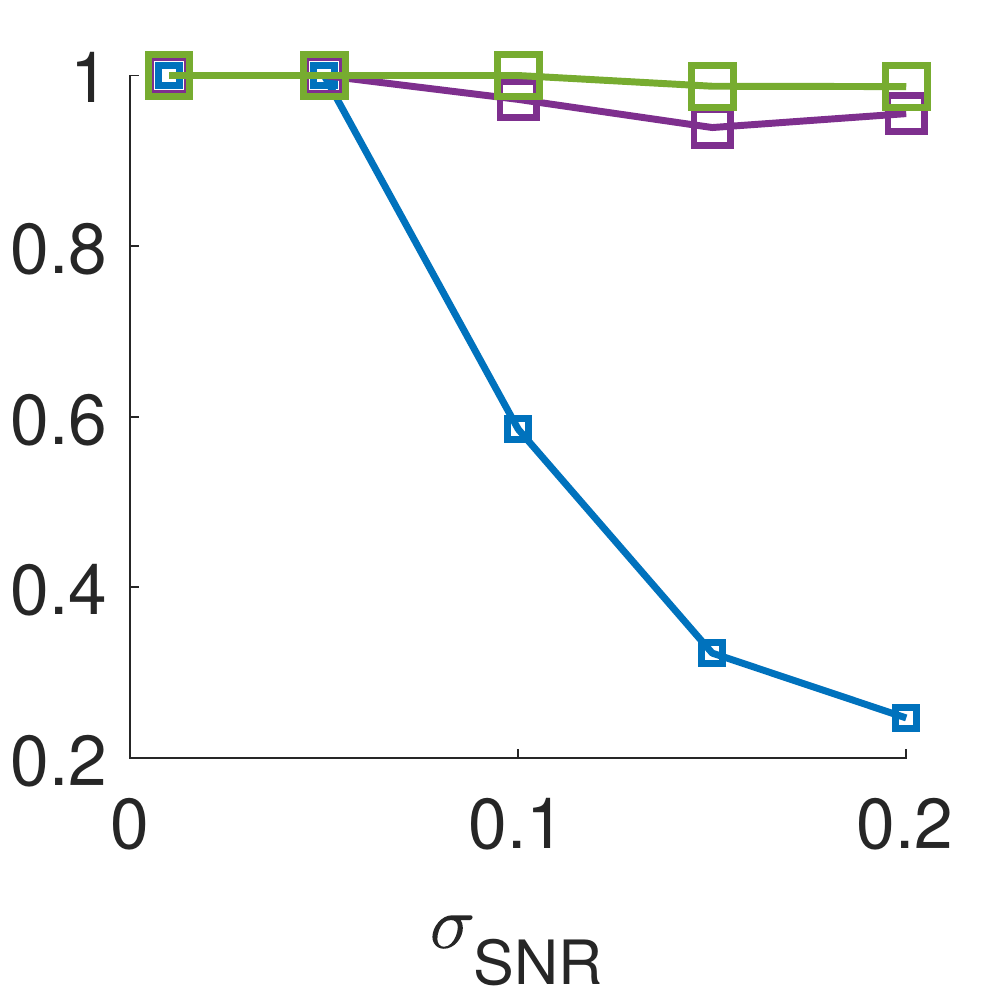} & \includegraphics[width = 0.1\textwidth]{figures_/legendv3.pdf} \\
    \midrule
    \eqref{e: ode, Duffing} &
    \includegraphics[width = 0.2\textwidth,align = c]{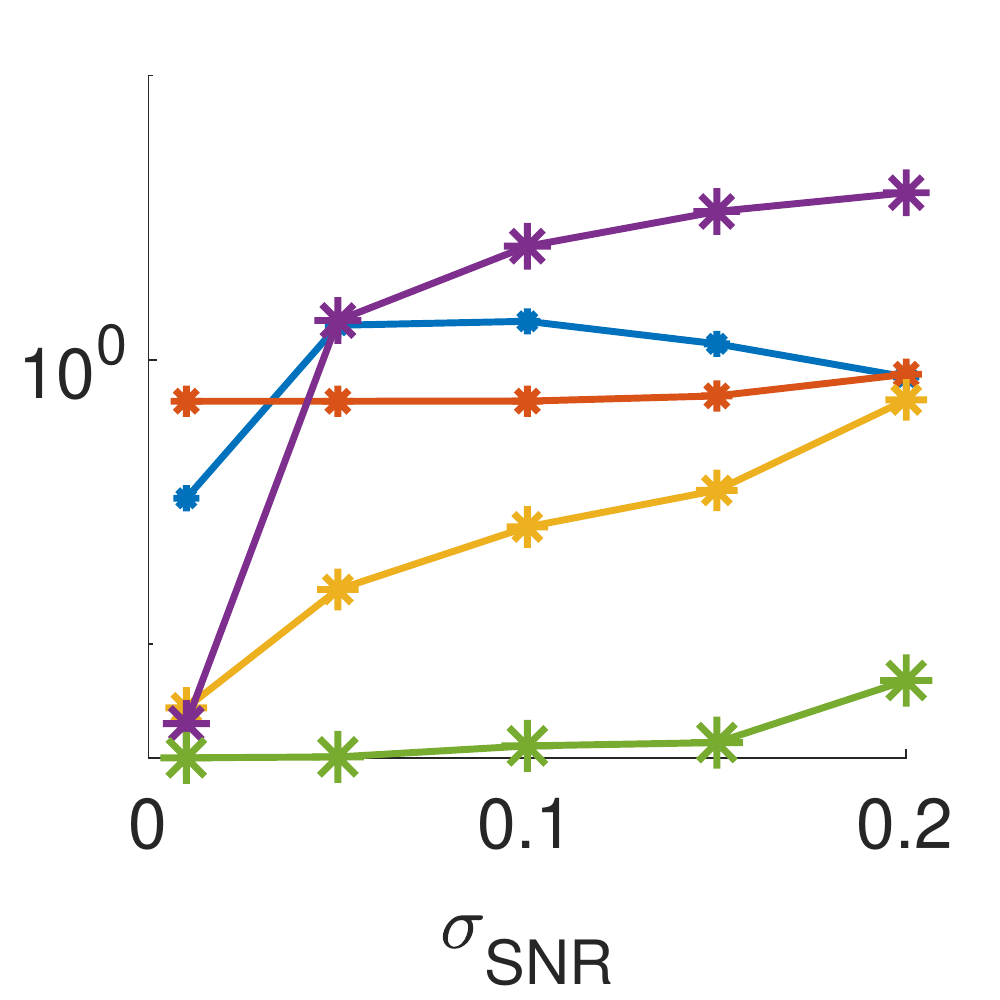} &  \includegraphics[width = 0.2\textwidth,align = c]{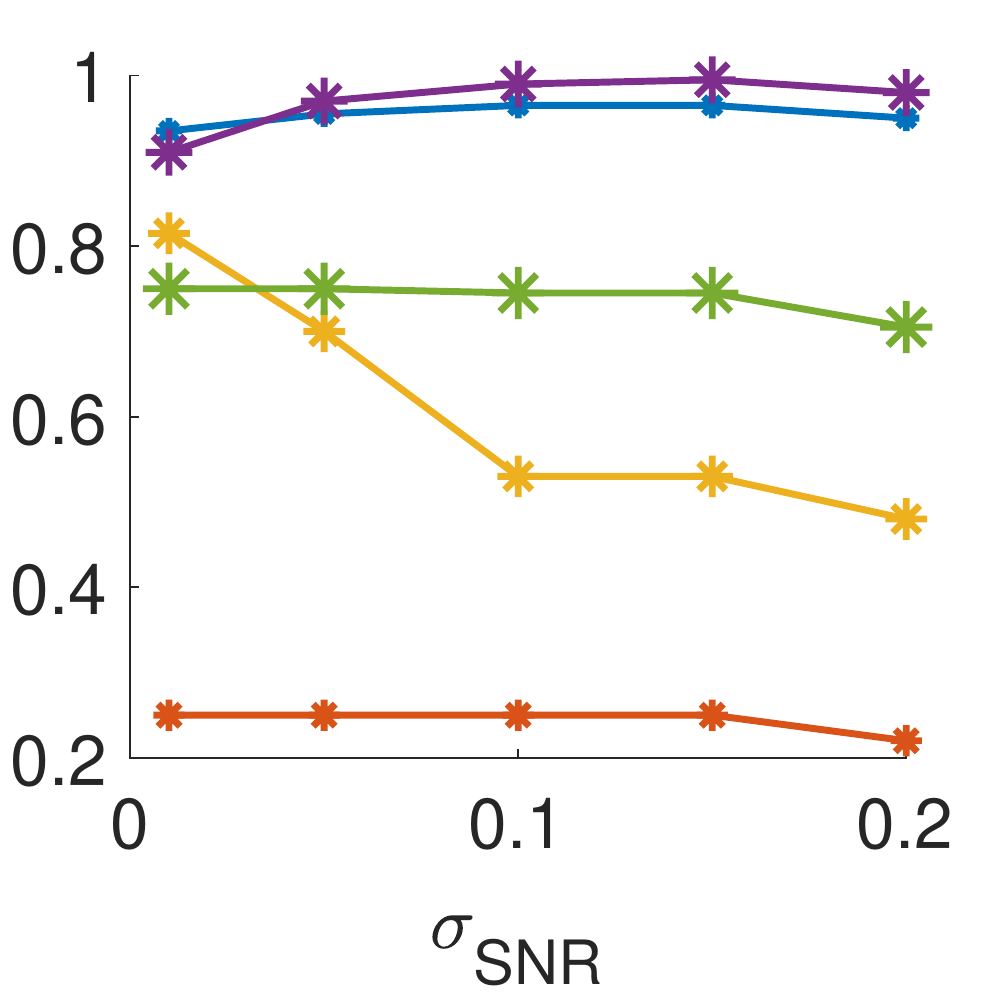}
    &  \includegraphics[width = 0.2\textwidth,align = c]{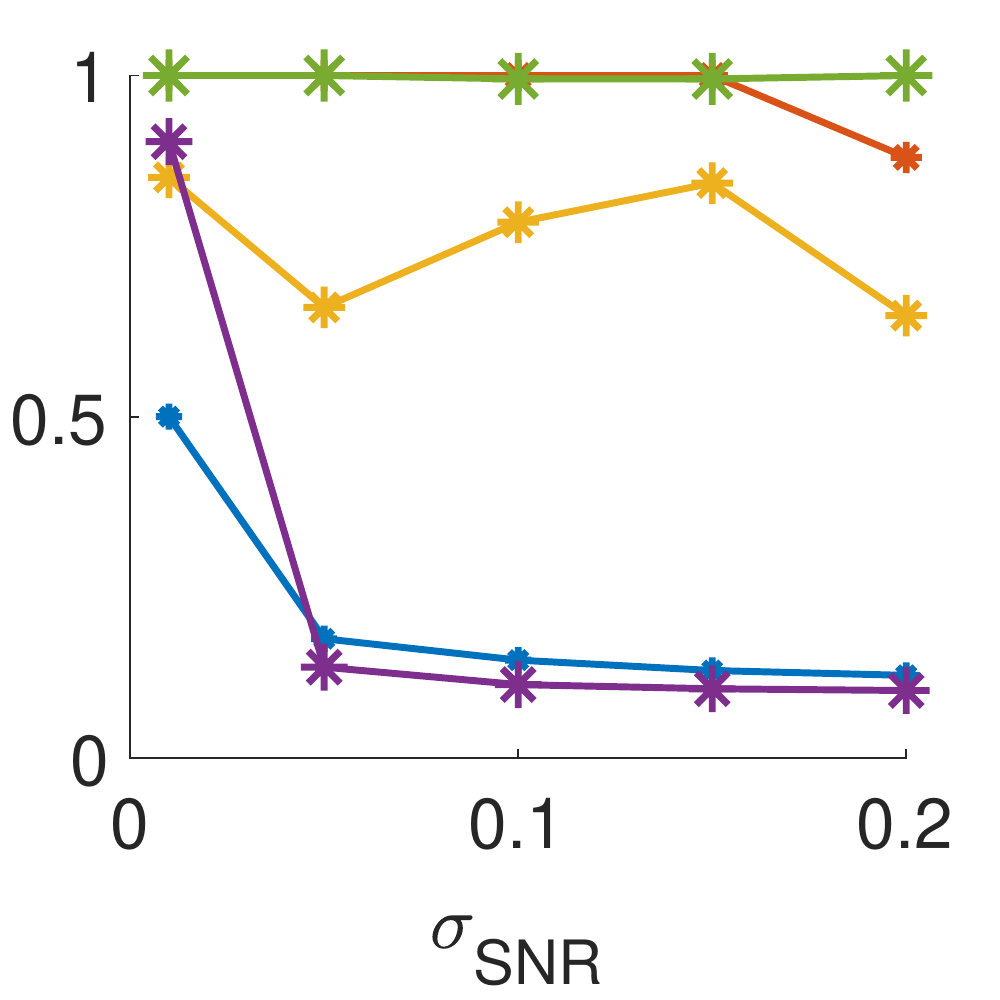} & \includegraphics[width = 0.1\textwidth]{figures_/legendv3.pdf}\\
    \midrule
    \eqref{e: ode, Lotka} &
    \includegraphics[width = 0.2\textwidth,align = c]{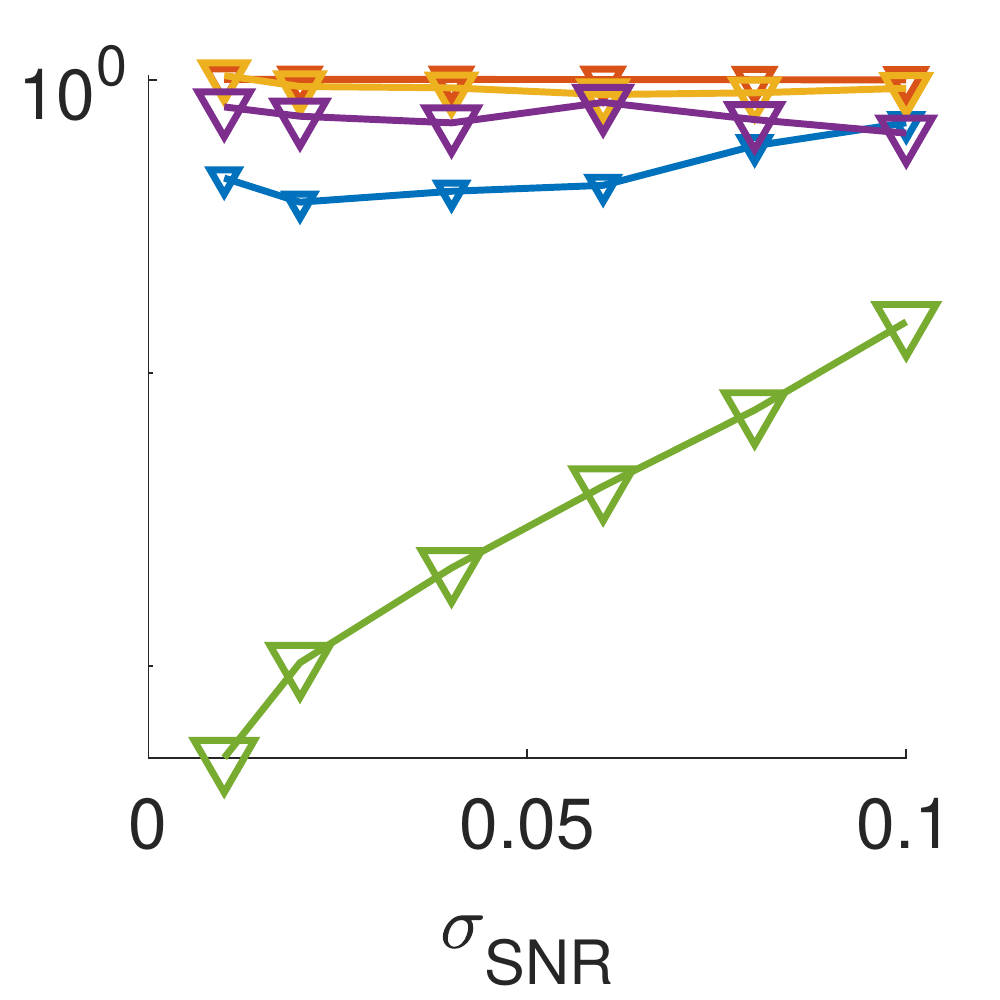} &  \includegraphics[width = 0.2\textwidth,align = c]{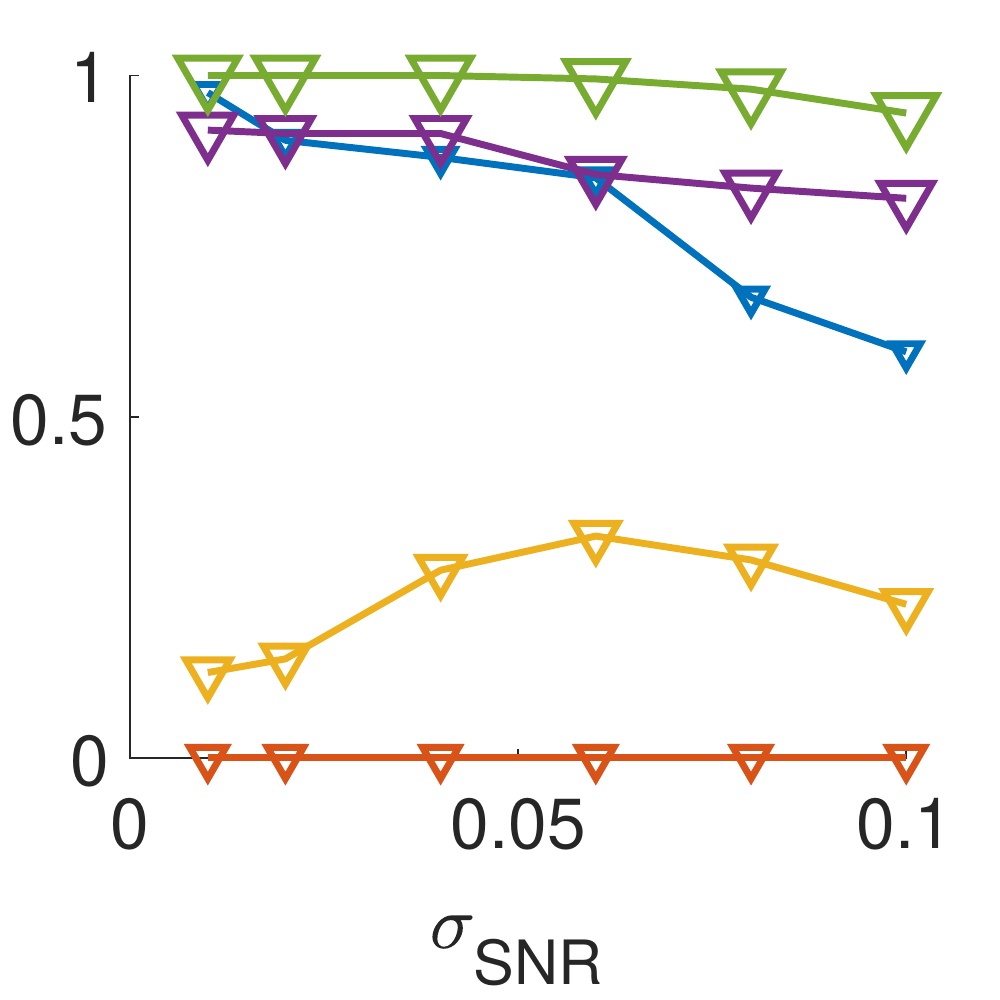}
    &  \includegraphics[width = 0.2\textwidth,align = c]{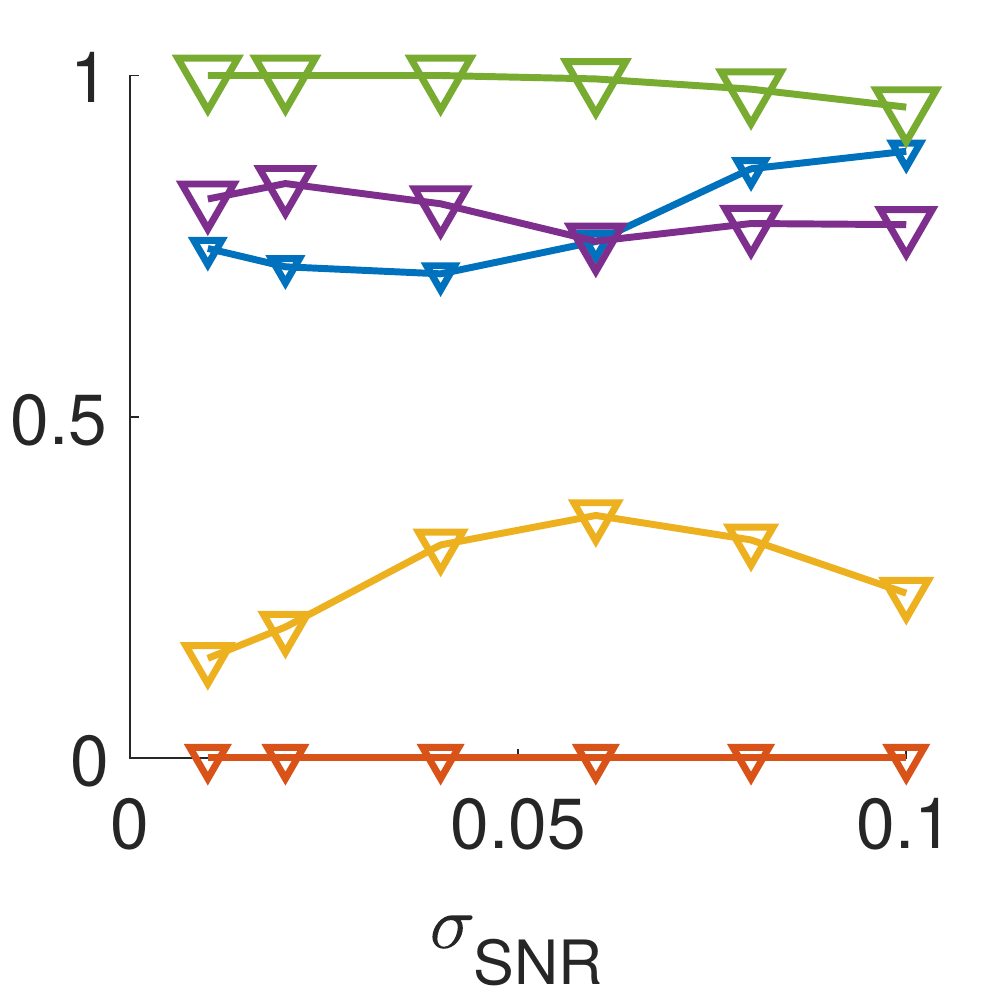} & \includegraphics[width = 0.1\textwidth]{figures_/legendv3.pdf}\\
    \midrule
    \eqref{e: ode, lorenz} &
    \includegraphics[width = 0.2\textwidth,align = c]{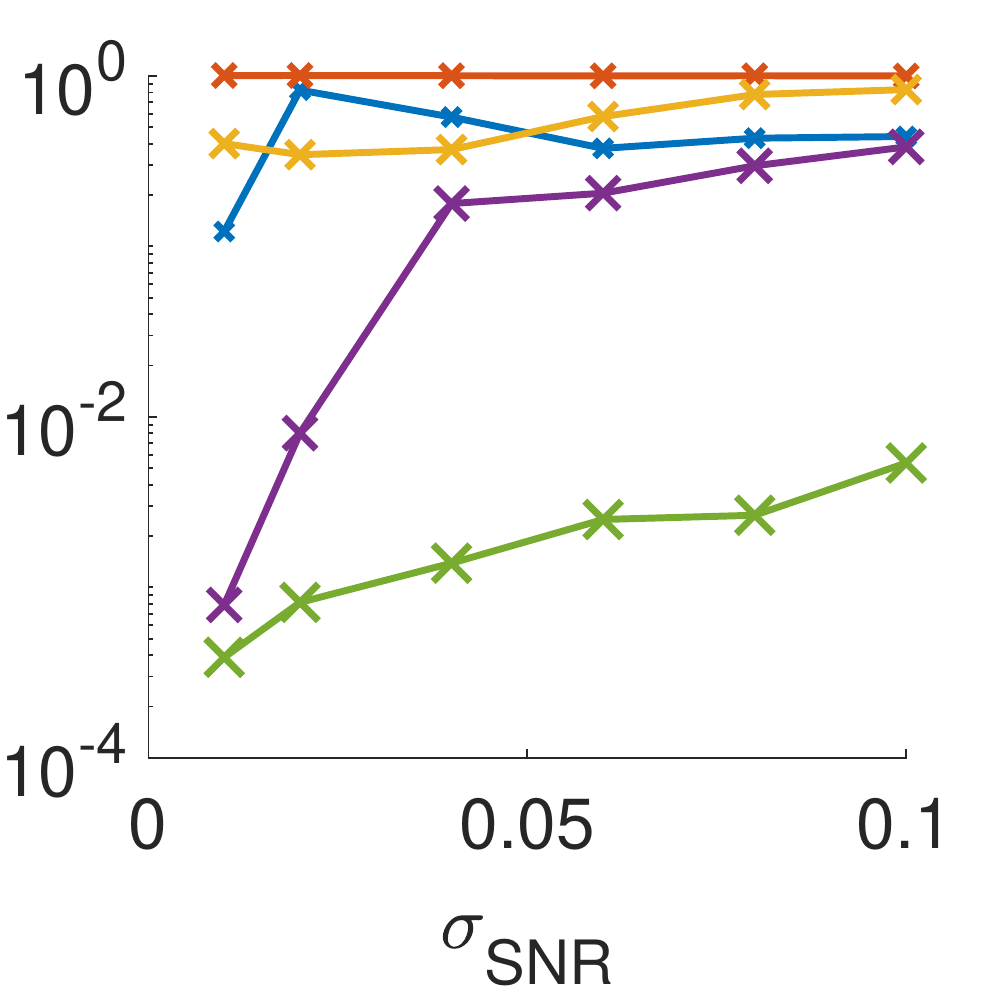} &  \includegraphics[width = 0.2\textwidth,align = c]{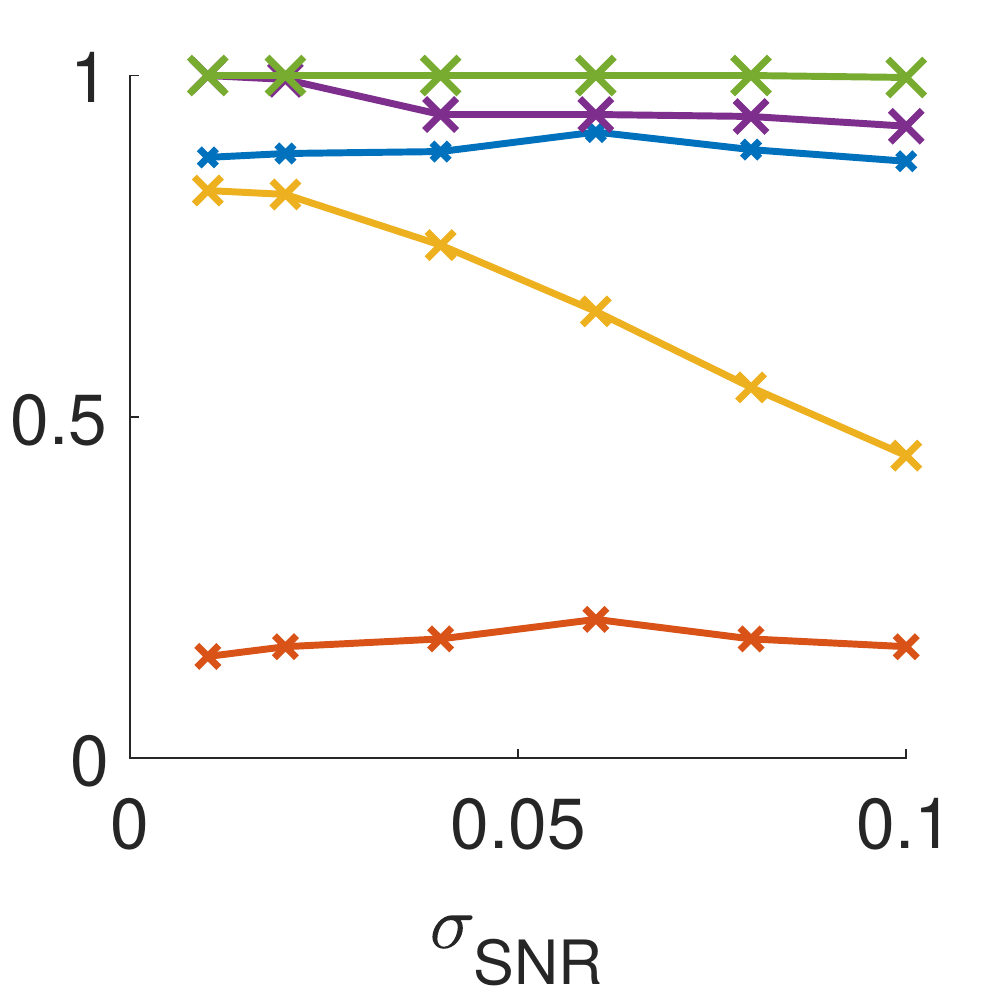}
    &  \includegraphics[width = 0.2\textwidth,align = c]{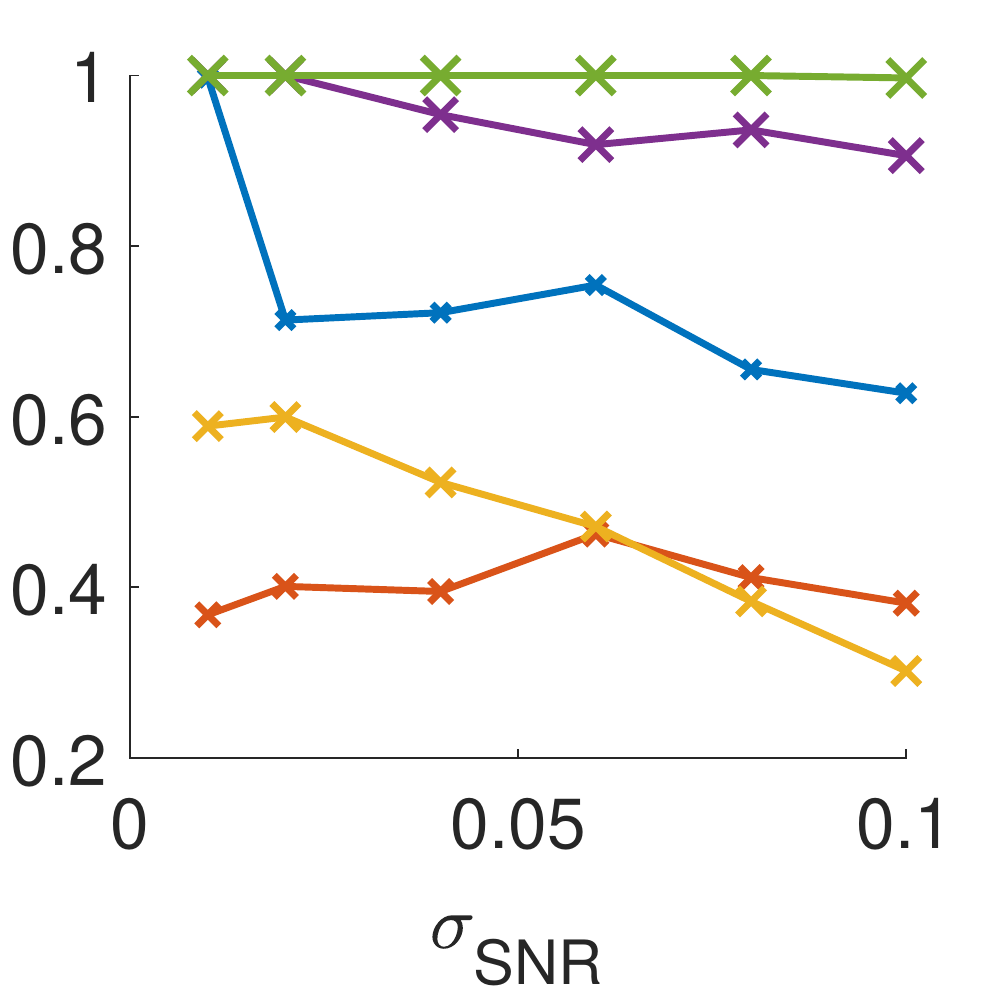} & \includegraphics[width = 0.1\textwidth]{figures_/legendv3.pdf}\\
    \bottomrule
    \end{tabular}
    \caption{WeakIdent results for the identification of ODEs listed in Table \ref{T: odes}, with the noise level  $\sigma_{\rm NSR}$ from 0 to 0.1.  Each graph shows the median over 50 experiments on each equation using WODE (purple), SINDy (blue), SC (red), ST (yellow),  and WeakIdent (green).  Each column shows the $E_2$ error, the TPR and PPV values. The green curve of WeakIdent gives the lowest $E_2$ error and the TPR and PPV values are close to 1.}
    \label{fig: summary recovering odes}
\end{figure}

In Figure \ref{fig: summary recovering odes}, we present  the identification results for the ODEs in Table \ref{T: odes}: the linear system \eqref{e: ode, linear 2d},  Van der Pol \eqref{e: ode, van-der-pol}, Duffing \eqref{e: ode, Duffing},  Lotka-Volterra \eqref{e: ode, Lotka}, and Lorenz \eqref{e: ode, lorenz}.  We experiment with different noise levels with different methods, including WeakIdent, WODE\cite{messenger2021weak},  SINDy\cite{brunton2016discovering}, SC \cite{he2020robust}, and ST \cite{he2020robust}.
Figure \ref{fig: summary recovering odes} shows the median of the $E_2$ error, TPR and PPV over 50 experiments for each equation. 
Overall WeakIdent (light green curves) yields the lowest $E_2$ error in the first column, and the TPR and PPV values near 1, which demonstrates a good support recovery.  

In Table \ref{F: an ode example - recovering noisey data - Lotka-Volterra}, we present the detailed results for the data in Figure \ref{fig: vis dynamic of odes}(d).   The noise level is $\sigma_{\rm NSR}=0.1$.
Table \ref{F: an ode example - recovering noisey data - Lotka-Volterra} (a) shows the dynamics. The noisy data for each of the dependent variables $\hat{X}$ and $\hat{Y}$ are shown in (b) and (c) respectively. Table \ref{F: an ode example - recovering noisey data - Lotka-Volterra} shows the identified systems by WeakIdent,  WODE, SINDy, SC and ST with  the $E_2$ and $E_{dyn}$ errors, TPR and PPV.  WeakIdent gives the  most accurate  recovery. 

\begin{table}[h]
 \begin{center}
\begin{tabular}{ccc}
(a) Equ.\eqref{e: ode, Lotka} & (b)  $\hat{X}$ & (c) $\hat{Y}$ \\
 \includegraphics[width = 0.25\textwidth ]{figures_/phase_portrait__Lotka_Volterra_Noisy.pdf} &
 \includegraphics[width = 0.25\textwidth]{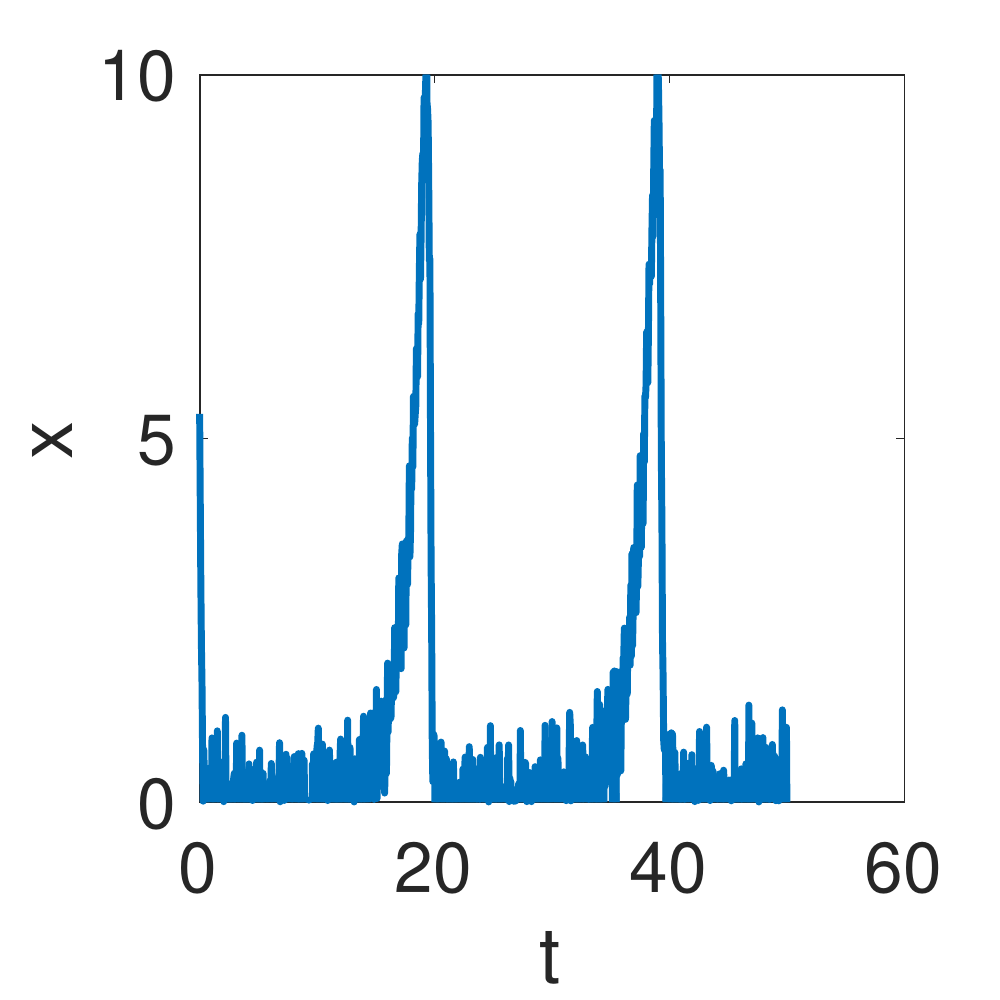}    &
  \includegraphics[width = 0.25\textwidth]{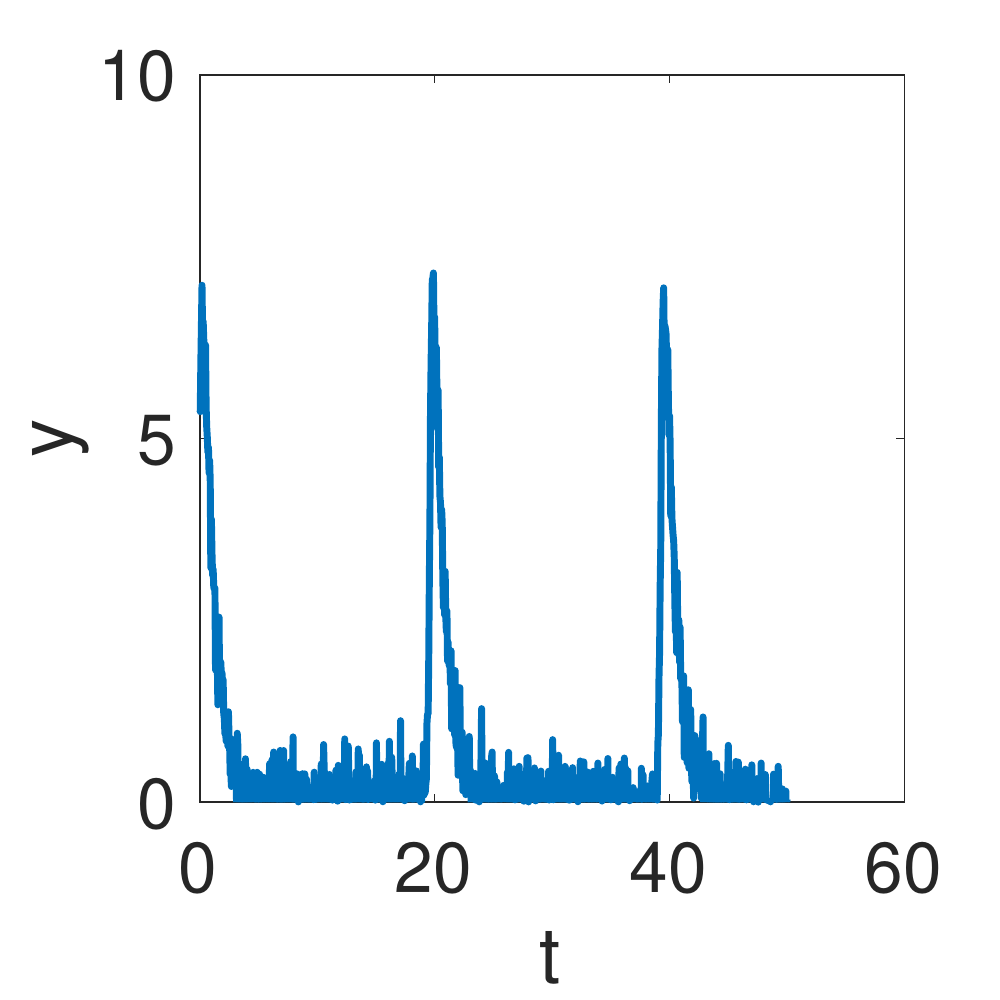} 
\end{tabular}
\begin{tabular}{|l|p{7cm}| r|r|r|r|}
\toprule
Method                         & Equation(s)   &                                             $E_2$  & $E_{\rm dyn}$
& TPR & PPV        \\ \hline
True equation &	$\dot{x}  =  +0.66667 x -1.33333 xy$		& & & & \\					
& $\dot{y}  =  -1.00000 y +1.00000 xy$ & & &&	\\		
\hline
\textbf{WeakIdent}&	$\dot{x}  =  +0.61271 x -1.21729 xy$			&	\textbf{0.07}&	\textbf{0.65}&	\textbf{1.00}&	\textbf{1.00}\\
 &	$\dot{y}  =  -0.97002 y +0.92339 xy$& & & &\\	
 \hline
WODE\cite{messenger2021weak}	&$\dot{x}  =  +0.62699 x -0.20740 x^2y$				&0.73&	2.34&	0.75&	0.75\\
 &	$\dot{y}  =  -0.91948 y +0.91515 xy$& & & &	\\	
 \hline
SINDy\cite{brunton2016sparse}&	$\dot{x}  =  +0.65909 x -1.05803 xy$			&	0.58&	5.57&	0.75	&1.00\\
 	&$\dot{y}  =  +0.61390 xy$	& & & &	\\			
 	\hline
SC\cite{he2020robust}&	$\dot{x}  =  -0.05204 y^3$	&			1.00 & 	112.17&	0.00&	0.00\\
&  	$\dot{y}=0$	& & & &		\\		\hline		
ST\cite{he2020robust}&	$\dot{x}  =  -0.05204 y^3$		&		0.88	&3.86	&0.25	&0.25\\
 &	$\dot{y}  =  +0.91700 xy +0.01157 xy^2 -0.05996 x^2y$	& & & &	\\	
 \bottomrule
\end{tabular}
\end{center}
\caption {The Lotka-Volterra equation  \eqref{e: ode, Lotka} with $\sigma_{\rm NSR}=0.1$.  We use the same data as in Figure \ref{fig: vis dynamic of odes}(d).   We present the comparisons between WeakIdent and WODE \cite{messenger2021weak},  SINDy \cite{brunton2016sparse},  SC, ST \cite{he2020robust}.}
\label{F: an ode example - recovering noisey data - Lotka-Volterra} 
\end{table}

\begin{figure}
    \centering
\begin{tabular}{lcccc}
    & $E_2$ & $E_{\rm res}$ & TPR & PPV \\
    \toprule
     & (a1) & (a2) & (a3) & (a4) \\
     \textbf{WeakIdent} & 
     \includegraphics[width = 0.17\textwidth, align = c]{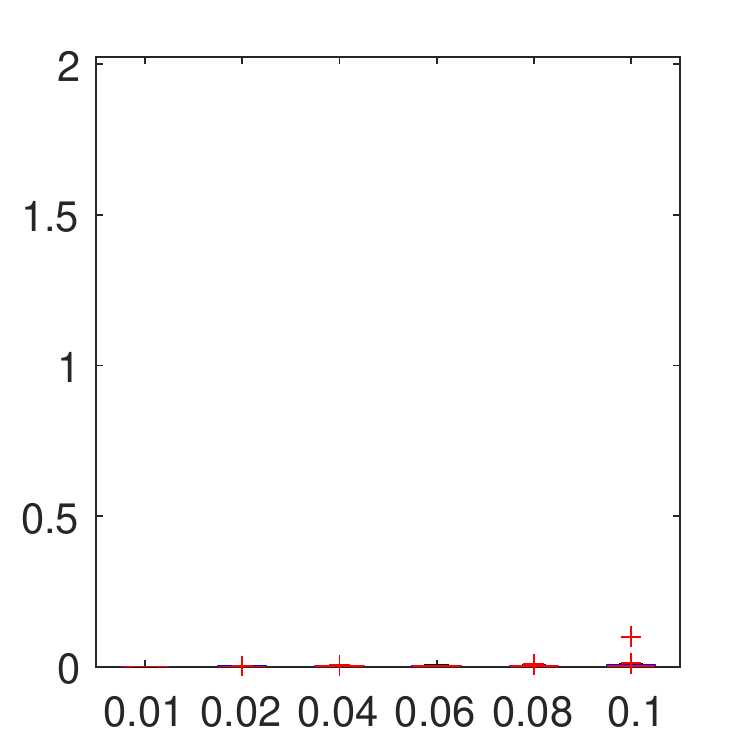} &
     \includegraphics[width = 0.17\textwidth, align = c]{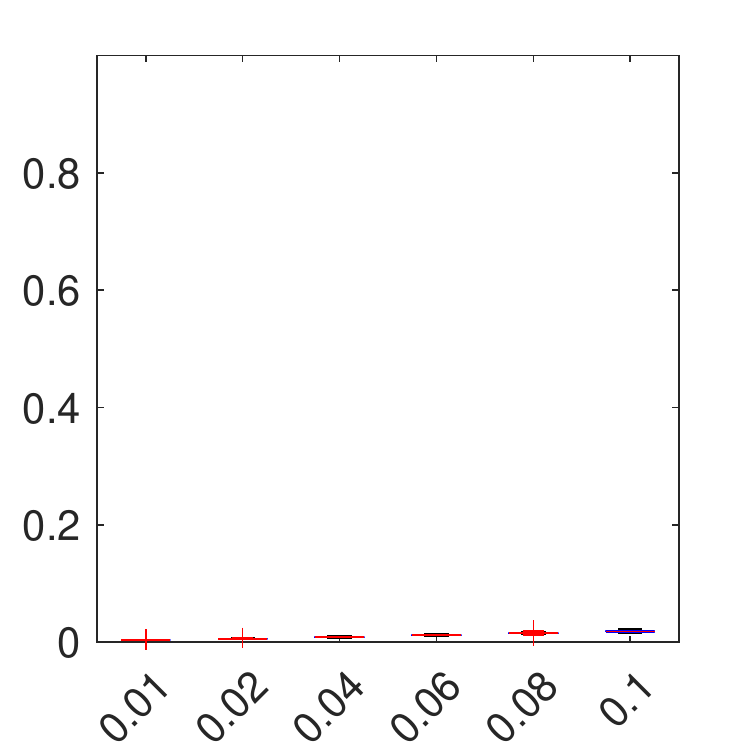} & 
     \includegraphics[width = 0.17\textwidth, align = c, align = c]{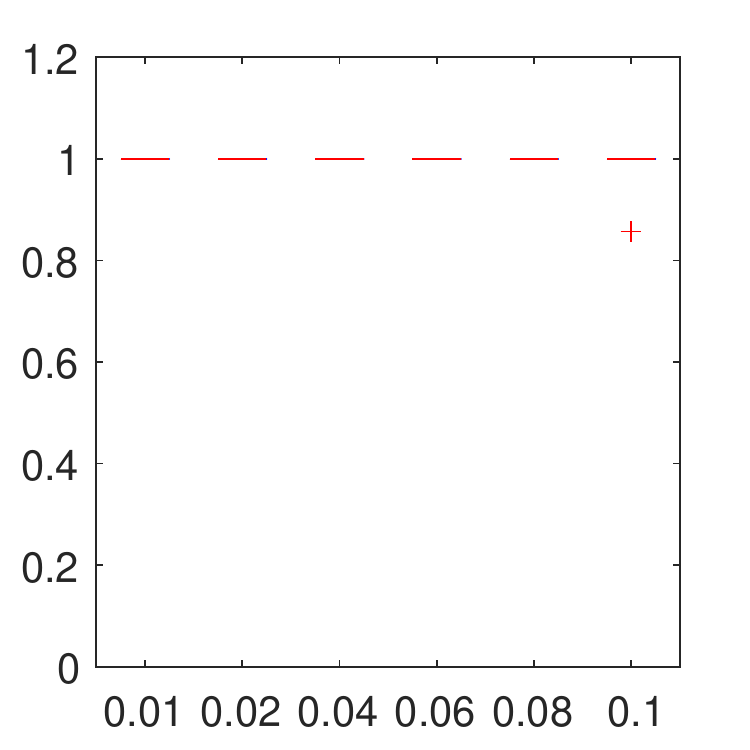} & 
     \includegraphics[width = 0.17\textwidth, align = c]{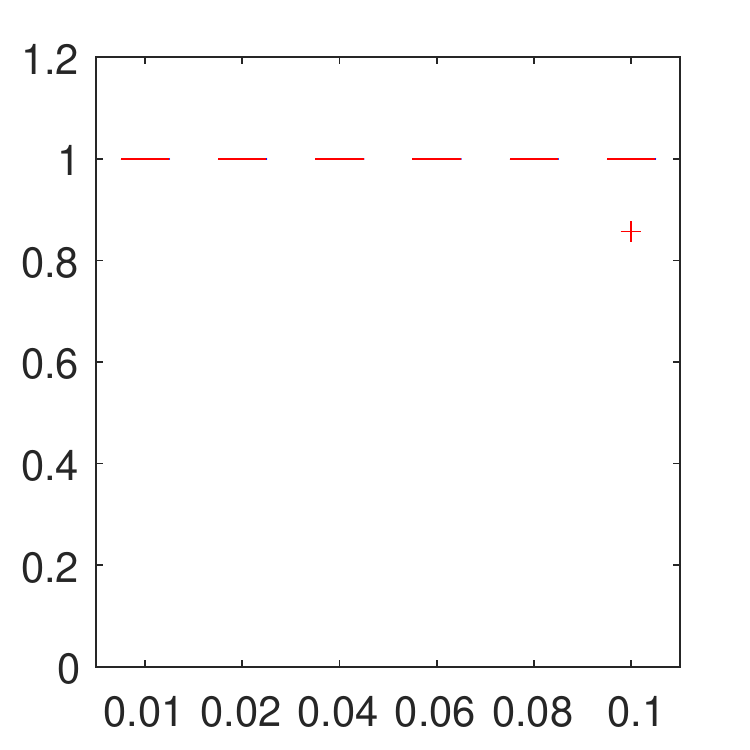} \\  
    & (b1) & (b2) & (b3) & (b4) \\
     WODE\cite{messenger2021weak} & 
     \includegraphics[width = 0.17\textwidth, align = c]{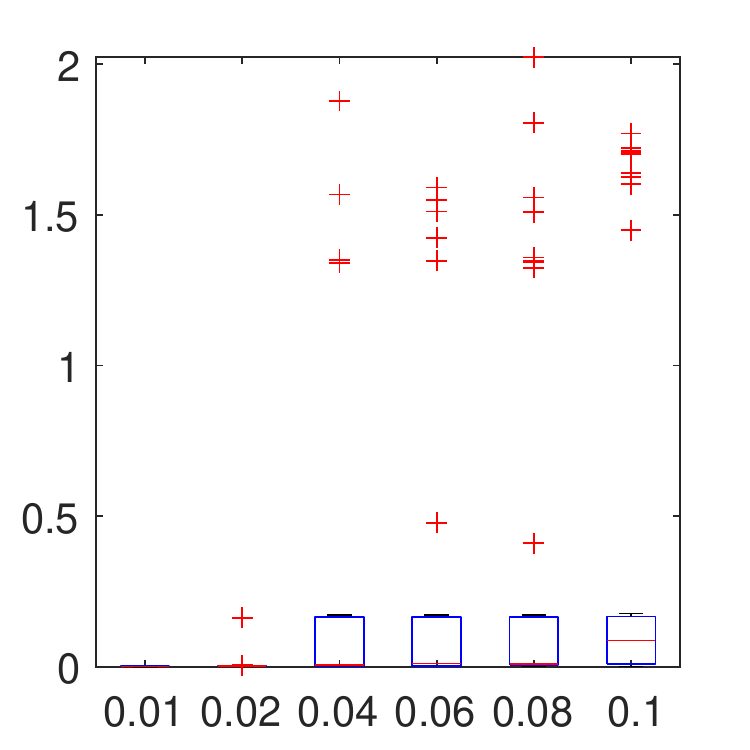} &
     \includegraphics[width = 0.17\textwidth, align = c]{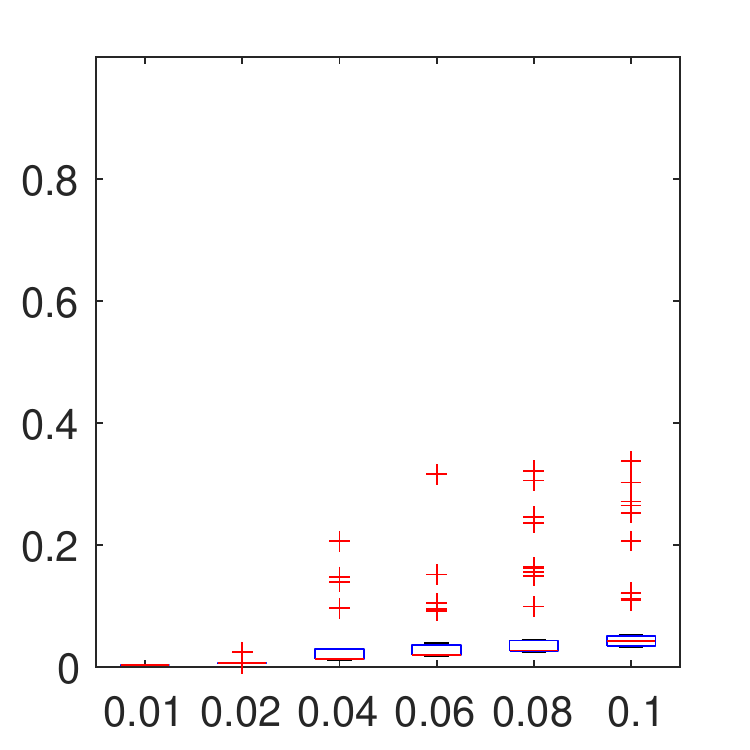} & 
     \includegraphics[width = 0.17\textwidth, align = c]{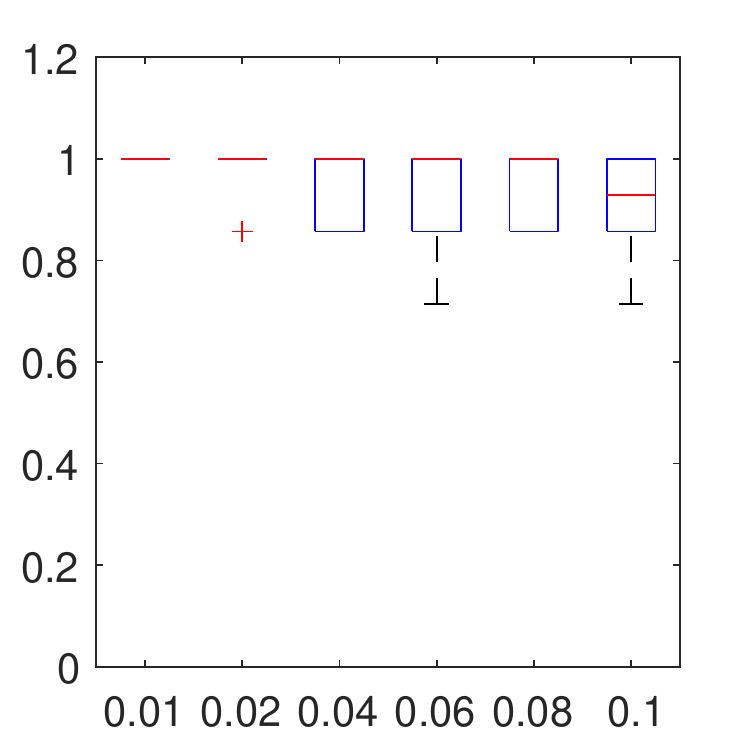} & 
     \includegraphics[width = 0.17\textwidth, align = c]{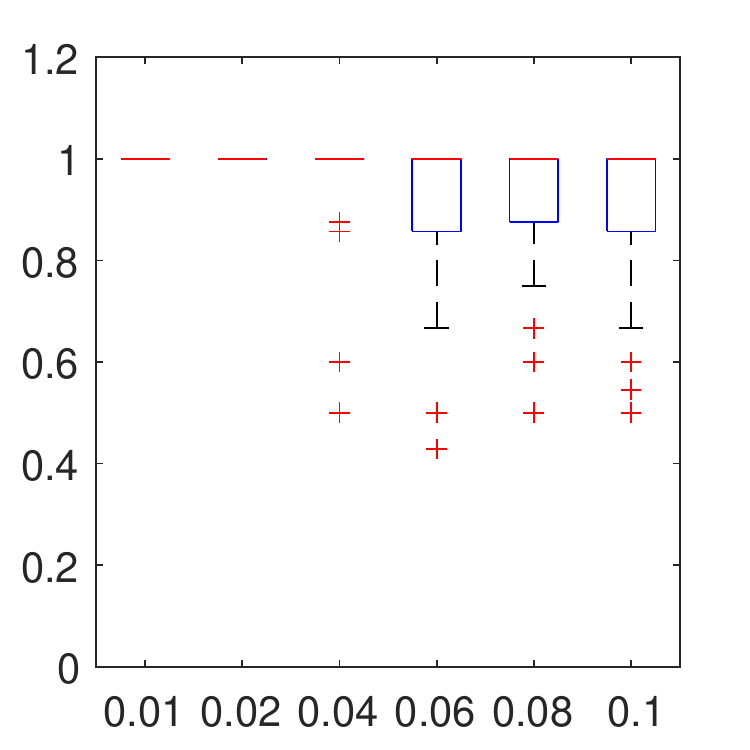} \\    
    & (c1) & (c2) & (c3) & (c4)\\
    SINDy \cite{brunton2016sparse} &
     \includegraphics[width = 0.17\textwidth, align = c]{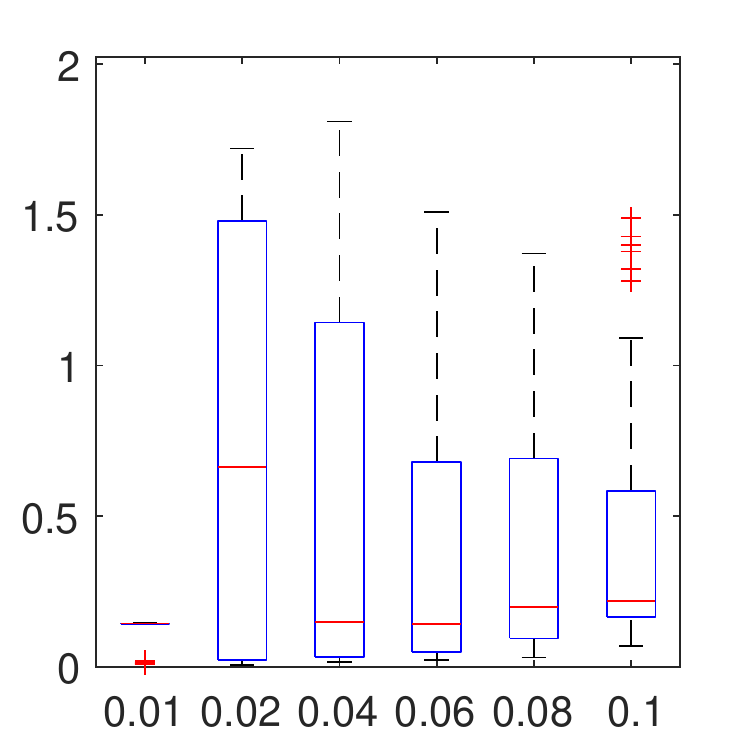} &
     \includegraphics[width = 0.17\textwidth, align = c]{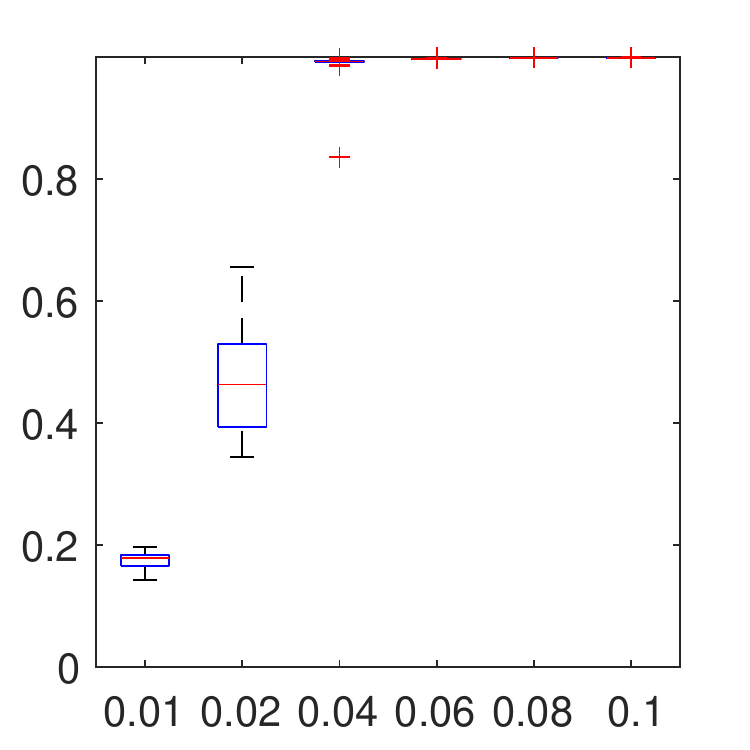} & 
     \includegraphics[width = 0.17\textwidth, align = c]{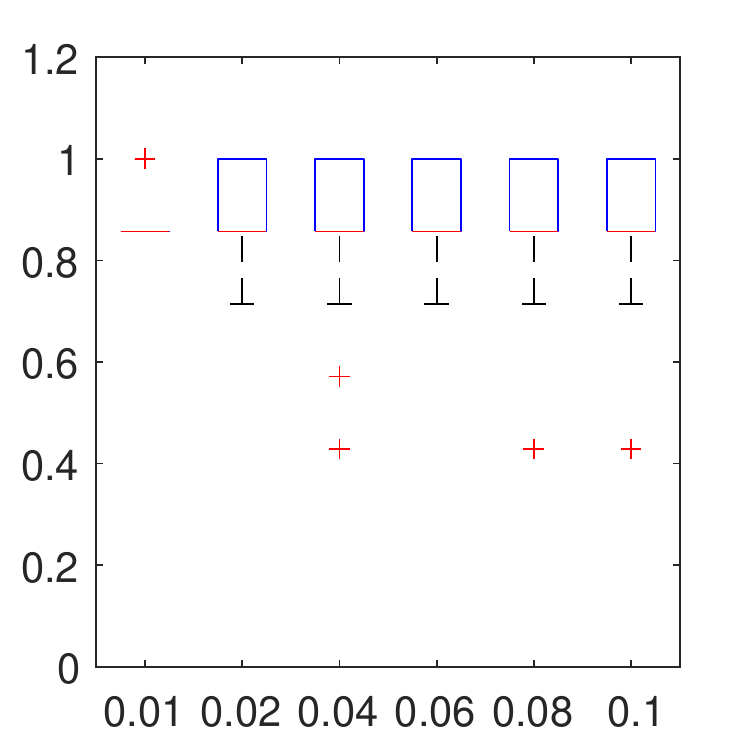} & 
     \includegraphics[width = 0.17\textwidth, align = c]{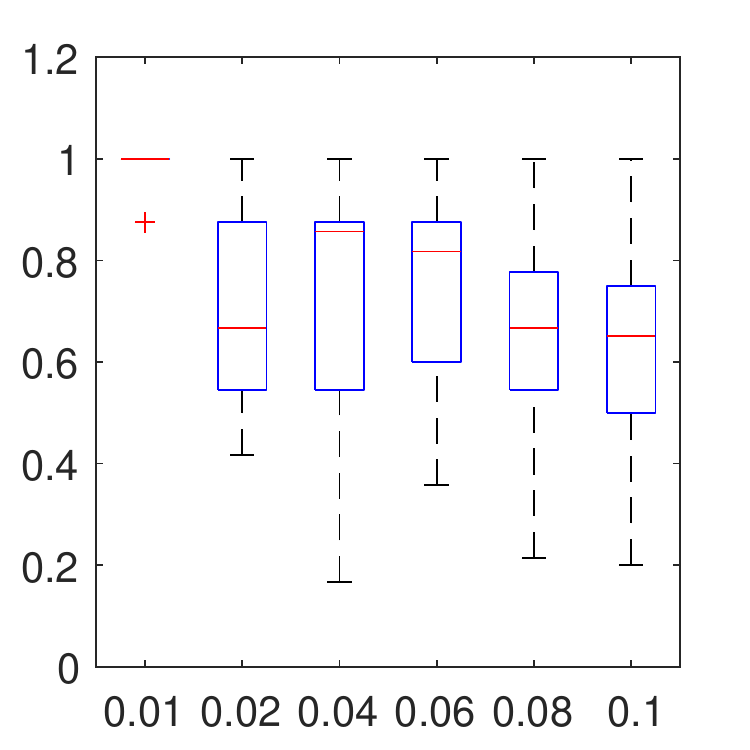} \\
     & (d1) & (d2) & (d3) & (d4) \\
     SC \cite{he2020robust} & \includegraphics[width = 0.17\textwidth, align = c]{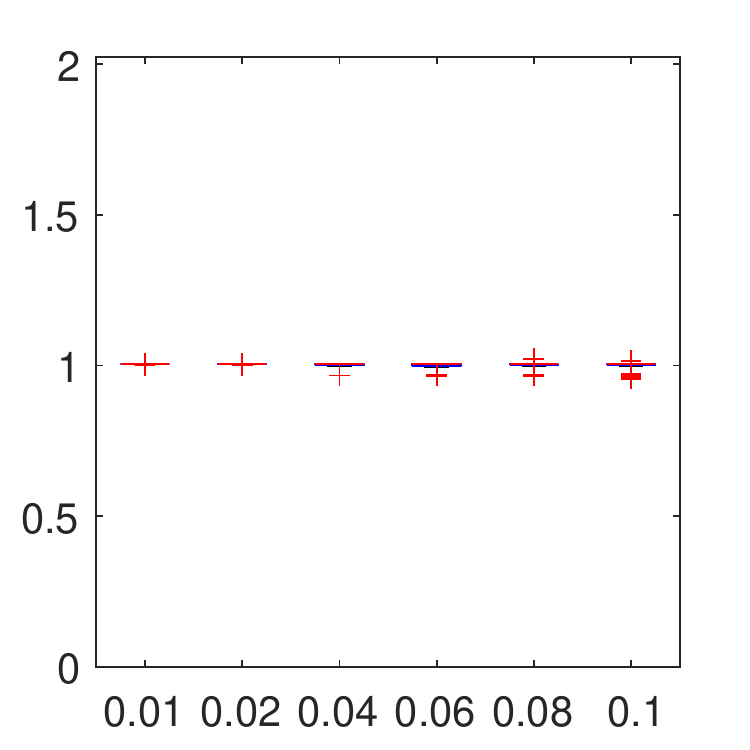} &
     \includegraphics[width = 0.17\textwidth, align = c]{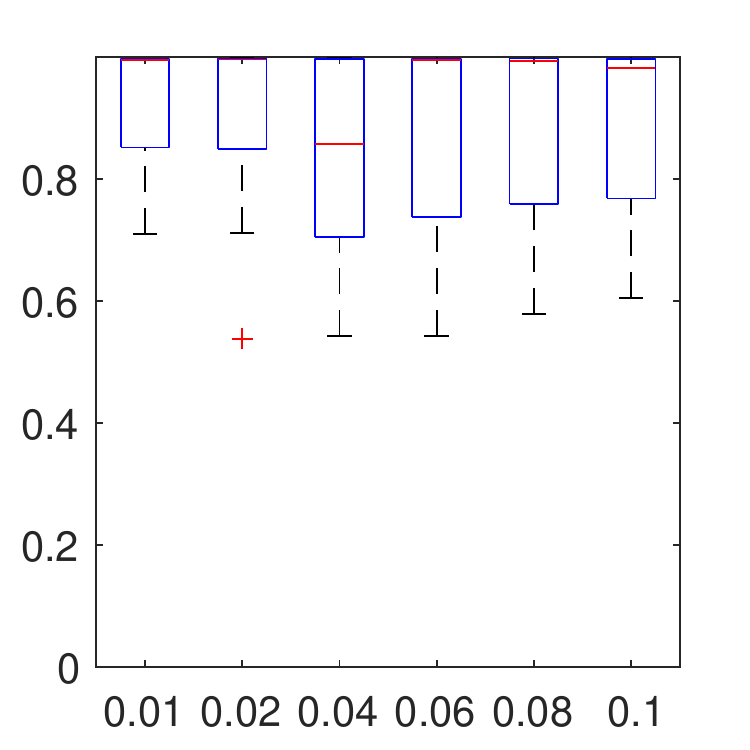} & 
     \includegraphics[width = 0.17\textwidth, align = c]{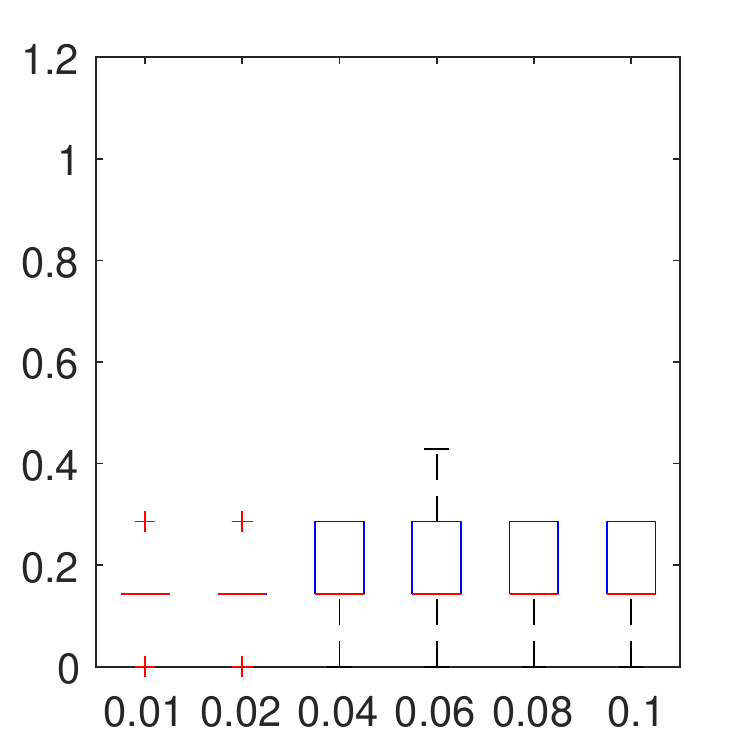} & 
     \includegraphics[width = 0.17\textwidth, align = c]{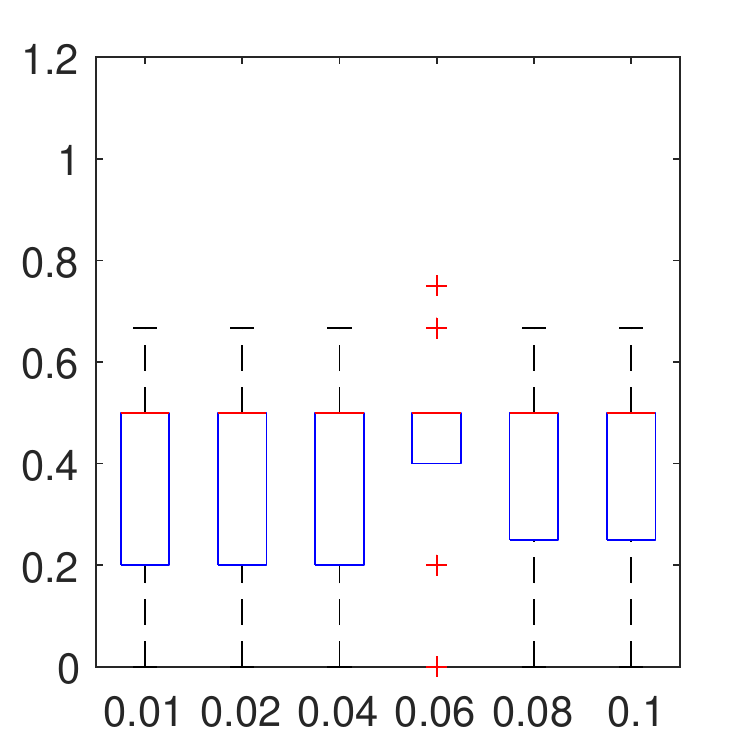} \\
     & (e1) & (e2) & (e3) & (e4) \\
     ST \cite{he2020robust} &
     \includegraphics[width = 0.17\textwidth, align = c]{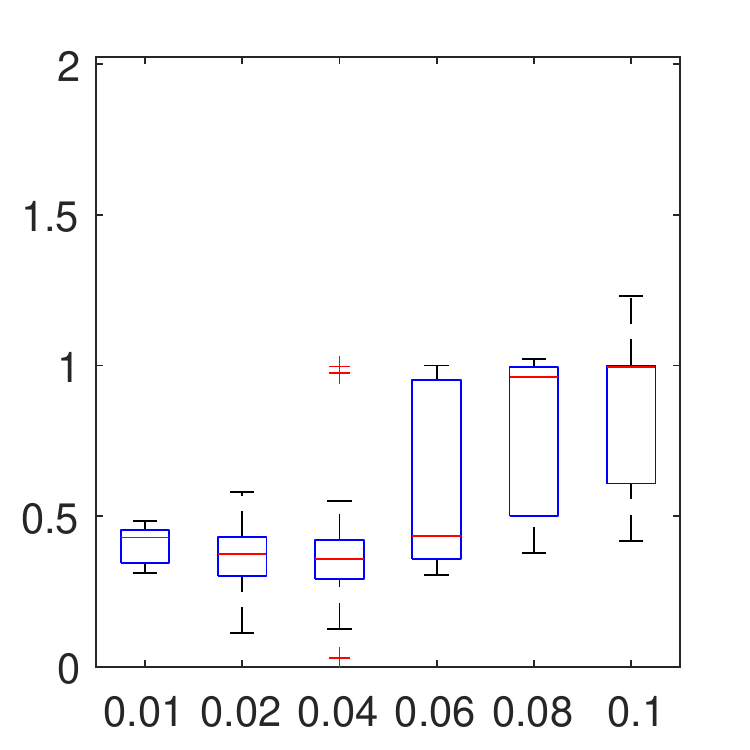} &
     \includegraphics[width = 0.17\textwidth, align = c]{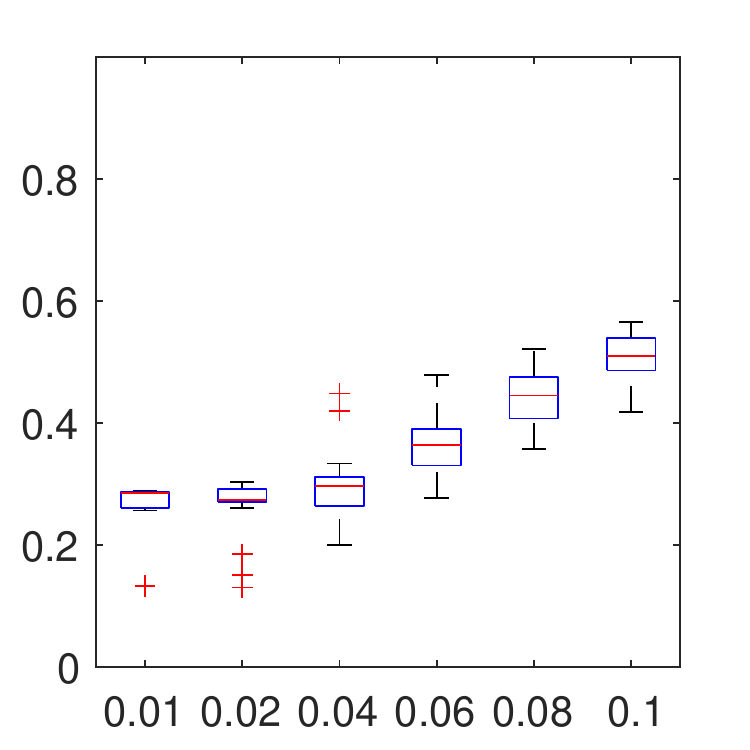} & 
     \includegraphics[width = 0.17\textwidth, align = c]{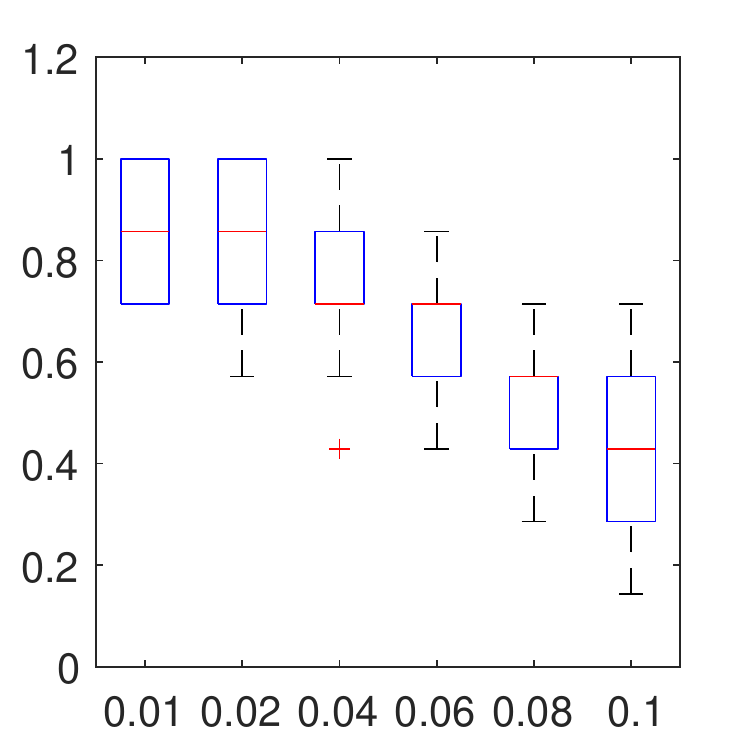} & 
     \includegraphics[width = 0.17\textwidth, align = c]{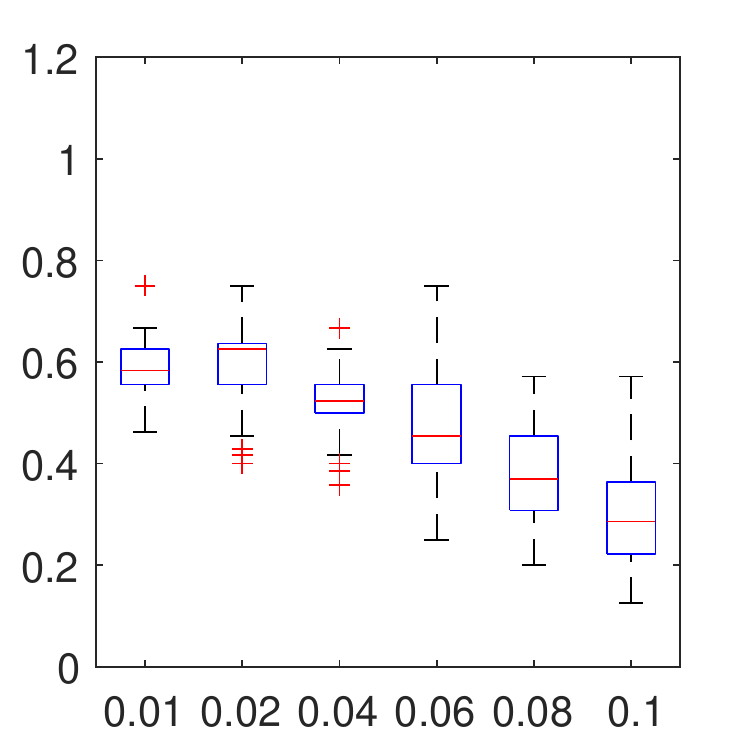} \\     
     \bottomrule
    \end{tabular}
    \caption{The Lorenz equation (\ref{e: ode, lorenz}), statistical comparisons:  WeakIdent (a1)-(a4), WODE \cite{messenger2021weak} (b1)-(b4), SINDy \cite{brunton2016sparse} (c1)-(c4), SC \cite{he2020robust}(d1)-(d4) and ST \cite{he2020robust}(e1)-(e4). 
    The $E_2, E_{\rm res}$ errors , TPR and PPV are shown from 50 experiments for each  $\sigma_{\rm NSR}\in \{ 0.01,0.02,  ,...,0.1\}$  using box-plots.  
    The  $E_2$ error given by WeakIdent  is  lower than others with less variations, and the TPR and PPV by WeakIdent are closer to 1 compared to other methods.}
    \label{fig: ode - box plot - Lorenz system}
\end{figure}

Table \ref{F: an ode example - recovering noisey data - Nonlinear Lorenz} shows the detailed results for the Lorenz system \eqref{e: ode, lorenz}. The data set is the same as the one in  Figure \ref{fig: vis dynamic of odes}(e) with $\sigma_{\rm NSR}=0.1.$.  The noisy data  $\hat{X}, \hat{Y}, \hat{Z}$ are displayed in (a), (b), (c) respecitvely. Table \ref{F: an ode example - recovering noisey data - Nonlinear Lorenz}  shows the recovered equations.  Table \ref{F: an ode example - recovering noisey data - Nonlinear Lorenz}  provides more details associated with  Figure \ref{fig: ode - box plot - Lorenz system} where we present statistical comparisons using 50 experiments for various noise level when $\sigma_{\rm NSR}$ varies from  $0.01$ to $0.1$.  The $E_2$ error by WeakIdent is lower with less variations, and the TPR and PPV values are closer to 1 compared to other methods.

\begin{table}
\begin{center}
\begin{tabular}{ccc}
(a)$\hat{X}$ & (b) $\hat{Y}$ &(c)$\hat{Z}$\\
 \includegraphics[width = 0.25\textwidth]{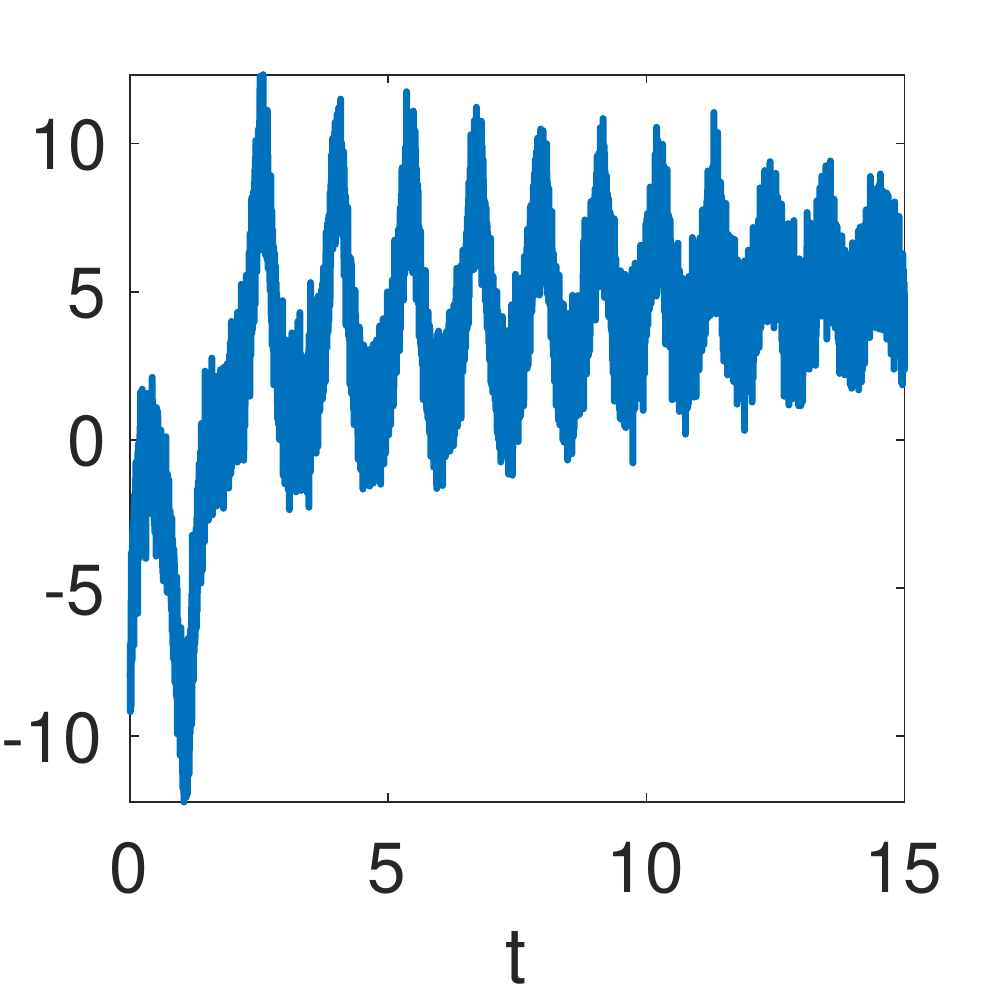}    &
  \includegraphics[width = 0.25\textwidth]{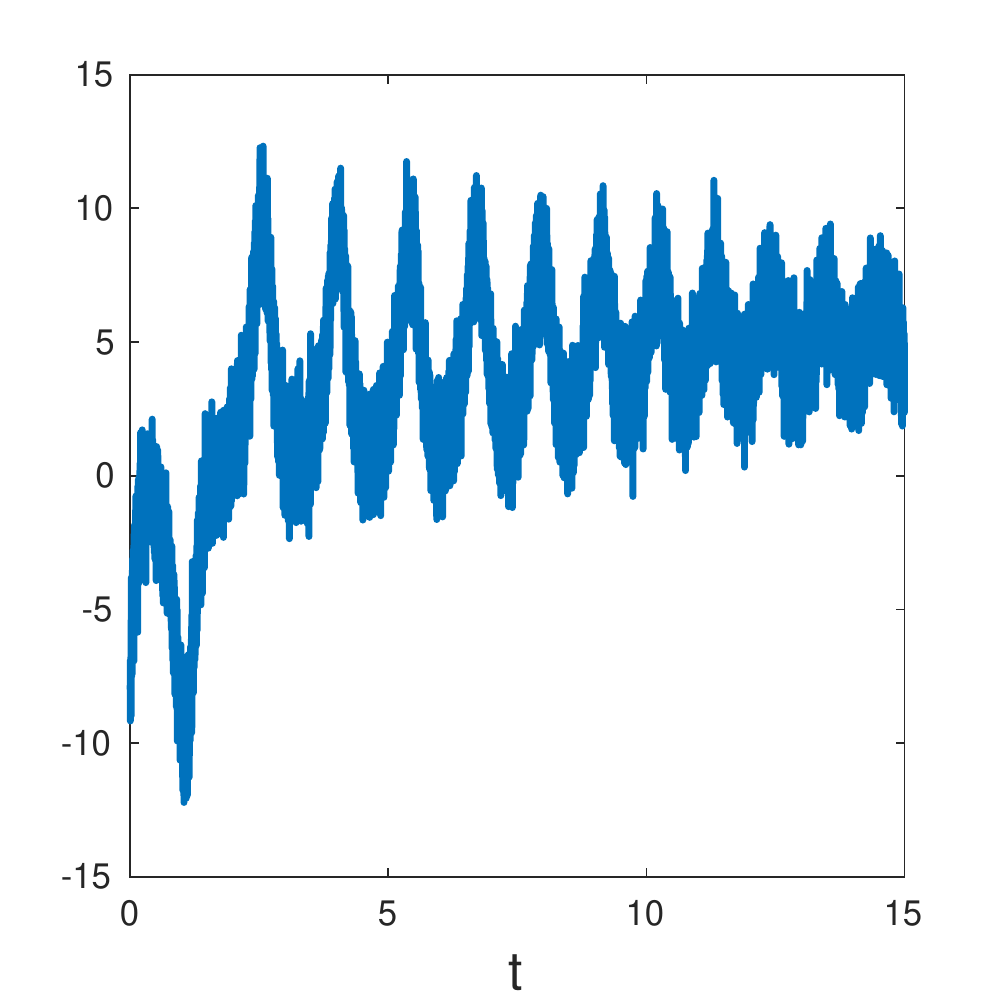}  & 
  \includegraphics[width = 0.25\textwidth]{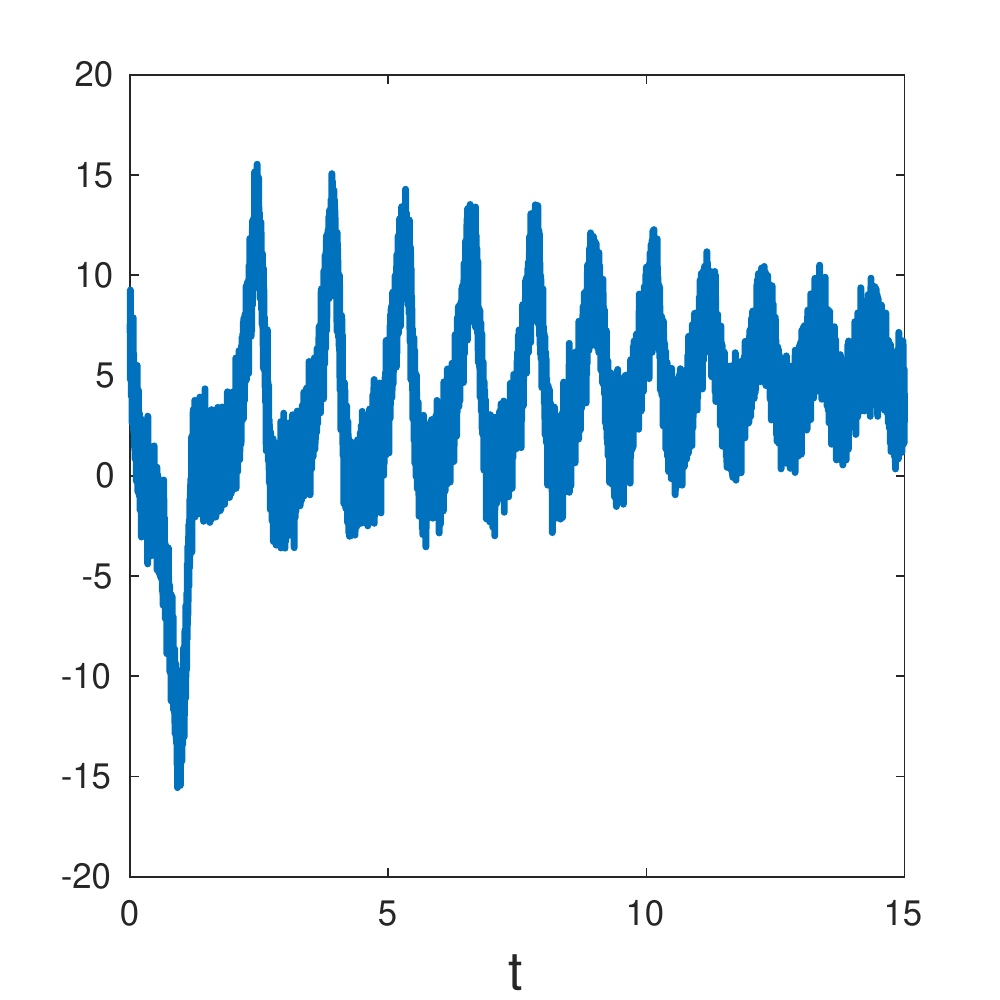} 
\end{tabular}
\begin{tabular}{|l|p{8cm}|r|r|r|}
\multicolumn{5}{c}{(d)}\\
\toprule
Method                         & Equation(s)  & $E_2$  
& TPR & PPV \\ \hline
True equation                           & $\dot{x}  =  +5.00000 y -5.00000 x$ & & &\\
& $\dot{y}  =  -1.00000 y +15.00000 x -1.00000 xz$& & &\\
& $\dot{z}  =  -2.00000 z +1.00000 xy$ & & &\\
\hline
\textbf{WeakIdent}         &              $\dot{x}  =  +4.97854 y -4.97406 x$					&		\textbf{0.011}& 	\textbf{1.00}& 	\textbf{1.00}\\
& $\dot{y}  =  -0.97267 y +14.88230 x -0.99353 xz$			& & &			\\		
& $\dot{z}  =  -1.96169 z +0.98105 xy$						& & &				\\
\hline
WODE\cite{messenger2021weak} & $\dot{x}  =  -4.12201  +2.77136 y +0.39253 y^2 -2.20585 x -0.89285 xy +0.45249 x^2$			&					0.018 &	1	& 0.88\\
& $\dot{y}  =  -11.16919  +3.46066 z -0.27237 z^2 +3.89361 y -0.35785 yz +6.97398 x -0.47972 xz -0.55649 xy +1.16625 x^2$   & & &	\\									
& $\dot{z}  =  -12.63547  +0.71043 y -1.07591 x +1.36546 xy -0.72687 x^2$		& & &		\\			
\hline
SINDy\cite{brunton2016sparse} &
$\dot{x}  =  +4.94736 y -4.91610 x$			&			1.188	 &  0.86 &	0.3\\
& $\dot{y}  =  -0.99564 y +14.88425 x -0.99267 xz$	& & &\\								
& $\dot{z}  =  +0.29463  -2.01828 z +1.00702 xy$								& & &	\\
\hline
SC\cite{he2020robust}      &  $\dot{x}  =  +0.02392 xy^2$						&	1.000	& 0.00&	0.00\\
                              & $\dot{y}  =  -0.00906 xz^2$                                                          & & &\\
                              & $\dot{z}  =  -0.27194 z^2 +0.01103 z^3 +0.05702 xyz$                                                           & & &\\
                              \hline
ST\cite{he2020robust}      &    $\dot{x}  =  +0.03245 x^2y$					&	1.000	& 0.00&	0.00\\
& $\dot{y}  =  -0.00906 xz^2$& & &	\\									
& $\dot{z}  =  +0.02673 y^2z$		& & &								   \\
                               \bottomrule
\end{tabular}

\end{center}
\caption {The Lorenz equation  \eqref{e: ode, lorenz} with $\sigma_{\rm NSR}=0.2$. This experiment shows the comparisons of WeakIdent, WODE, SINDy, SC and ST using the same data set from Figure \ref{fig: vis dynamic of odes}(e). WeakIdent gives rise to the best recovery. }
\label{F: an ode example - recovering noisey data - Nonlinear Lorenz}
\end{table}

\section{Test functions and feature construction}  \label{Asec:testfunction}

We give an outline of the construction of test functions in  \eqref{e: test function}. Specifically, we discuss how to choose the parameters $m_x, m_t, p_x, p_t$ (simply $m$ and $p$ below) according to  \cite{messenger2021weakPDE}: 

\textbf{(1) Frequency consideration:} Given the  data $\{U_i^n\}$, we consider the Fourier transform of data in each dimension. For example, $\mathcal{F}_x(U)$ is the Fourier transform of $U$ is the spatial domain.  We next find a junction point $k_x^*$ by fitting the cumulative sum of the vectorized data $|\mathcal{F}_x(U)|$ by a piecewise linear polynomial with one junction point. The $k_x^*$ minimizes the $L_2$ fitting error. 

\textbf{(2) Fit with a Gaussian distribution:}
The test function $\phi$ is matched to Gaussian for a denoising effect, i.e., 
$\phi_p(x) = C\left(1-\left(\frac{x}{m\Delta x}\right)^2 \right)^{p} \approx \rho_{\sigma}(x)$
where 
$ \rho_{\sigma}(x) = \frac{1}{\sqrt{2\pi} \sigma} e^{ -\frac{1}{2}(\frac{x}{\sigma})^2}$ and $
\sigma = \frac{m \Delta x}{\sqrt{2p + 3}}.
$ 
Here $\phi_p(x)$ matches $\rho_{\sigma}$ up to the third moment such that $ |\hat{\phi}_p(\xi) -\hat{\rho}_{\sigma}(\xi)| \leq \mathcal{O}(|\xi|^4(m\Delta x)^4 p^{-3})$ and $C$ is a constant such that $||\phi_p||_1 = 1$ \cite{messenger2021weakPDE}. Here $\hat{\phi}_p$ and $\hat{\rho}_\sigma$ denotes the Fourier transform of $\phi_p$ and $\rho_\sigma$ respectively. 
To suppress the noise, the high frequency components of data with the mode larger than $k_x^*$ or smaller than $-k_x^*$ are set to be within the 5\% tail of the Gaussian. This gives $ \frac{2\pi}{\mathbb{N}_x\Delta x} k_x^* = \frac{2}{\sigma}$ from  the property of cumulative distribution function of $\hat{\rho}_{\sigma} = \hat{\rho}_{1/\sigma}$, and relating this to $p$ and $m$ gives the first condition: $\frac{2\pi}{\mathbb{N}_x \Delta x} k_x^* = \hat{\tau} \frac{\sqrt{2p + 3}}{m \Delta x}$, where $\hat{\tau}$ is a parameter \cite{messenger2021weakPDE}.

\textbf{(3) Vanishing of $\phi$ on the boundary } To guarantee  the  decay of $\phi$ in each  spatial domain, $p$ and $m$ are set to satisfy the second condition:
$ \phi_p((m-1)\Delta x) \leq 10^{-10}$ and $p>\alpha_x+1$ where $\alpha_x$ is the highest order  derivative in the $x$ direction for all features.  Using the first and the second conditions above, $p$ and $m$ are determined.

\bigskip

\begin{figure}
    \centering
    \includegraphics[width = 0.6 \textwidth]{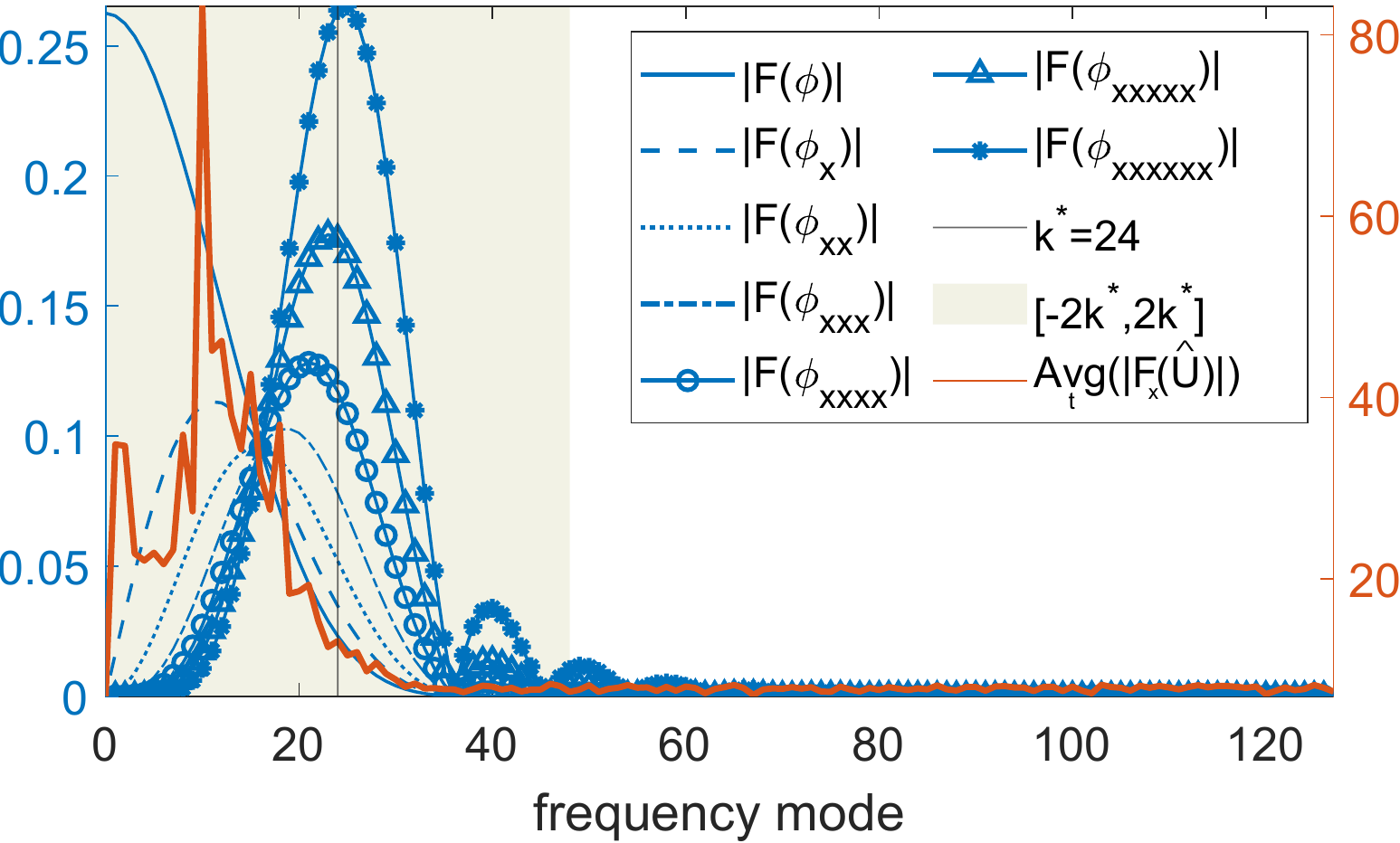}
     \caption{An example of the test function $\phi$ and an illustration in the frequency domain, using data from the KS equation  \eqref{e: pde KS} with $\sigma_{\rm NSR} = 0.5$. 
     This figure shows the test function $\phi(x)$ (solid blue), its derivatives with respect to $x$ (blue lines) in the frequency domain, and the Fourier transform of the given noisy  data averaged over time (red), i.e. $\frac{1}{\mathbb{N}_t}\sum_{n=1}^{\mathbb{N}_t}|{\mathcal{F}_x}(\hat{U}(x,t^n))|$ . 
     The vertical line denotes the location of the junction point $k^*_x$. The colorized region is where the signal dominates and the rest of region is where the noise dominates.}
    \label{fig: visualization of transition point}
\end{figure}

Figure \ref{fig: visualization of transition point}  shows an example of the test function $\phi$ (solid blue), its spatial derivatives up to order 6 (blue lines) and {$\frac{1}{\mathbb{N}_t}\sum_{n=1}^{\mathbb{N}_t}|{\mathcal{F}_x}(\hat{U}(x,t^n))|$}  in  the  frequency domain (red)  for the noisy data of  the  KS equation \eqref{e: pde KS} with $\sigma_{\rm NSR}=0.5$.  
The vertical line denotes the location of  the transition point $k^*_x$ . We colorize the region where the frequency mode is below $k_x^*$ such that the signal dominates in this region.
In this case,  the  junction point is $k^* = 24$,  $p_x = 10$, and $m_x = 17$.
The shape of $\mathcal{F}(\phi)$ (solid blue) and $\frac{1}{\mathbb{N}_t}\sum_{n=1}^{\mathbb{N}_t}|{\mathcal{F}_x}(\hat{U}(x,t^n))|$ (red) demonstrate the denoising effect of using the test function $\phi$. 
The difference among the blue curves  shows that integral forms suppress signals more for the higher order features. 

 \bigskip

\textbf{The computation of the features in \eqref{e: approx of bW}} is done by convolution in each dimension:
\begin{equation*}
    w_{h(x_i,t^n),l} 
    = (-1)^{\alpha_l} \left( U * \frac{\partial^{\alpha_l}\phi }{\partial x^{\alpha_l }}\right)(x_i,t^{n}) 
    =(-1)^{\alpha_l} \sum_{k = n-m_t}^{n+m_t} \sum_{j =i -m_x}^{i+m_x} U_j^k  
    \frac{\partial^{\alpha_l} }{\partial x^{\alpha_l}}\phi(x_j-x_i,t^k-t^n), 
\end{equation*}
and there is a similar form for $b_{h(x_i,t^n)}$
for each $(x_i, t^n)$ in the domain.
FFT is applied to compute the convolution.
For the convolutions in the $x$ direction when $t = t^n$, we use the vector  $\boldsymbol{\phi} = ( \phi(-m_x\Delta x, t^n), ..., \phi(m_x\Delta x, t^n))$ for convolution. Near the boundary, we pad the data by zeros for the computation of convolutions.

\end{document}